\documentclass[12pt,review]{elsarticle}
\usepackage{amsmath,epsfig}
\usepackage[latin1]{inputenc}
\usepackage{subfigure}
\usepackage{epsfig}
\usepackage{multirow}
\usepackage{amssymb,amsfonts,amsmath}
\usepackage{graphicx}
\usepackage{float}
\usepackage{algorithm}
\usepackage{algorithmic}
\usepackage{array}
\usepackage{multirow}
\usepackage{MnSymbol}
\usepackage{lineno}
\usepackage{bm}
\usepackage{color}
\usepackage{hyperref} 
\usepackage[margin=2.5cm]{geometry}
\usepackage{mathrsfs,amsmath}
\usepackage{tipa,mathtools}
\usepackage{wrapfig}
\usepackage{mdframed}
\usepackage{fancybox}

\usepackage[edges]{forest}
\useforestlibrary{linguistics}
\forestapplylibrarydefaults{linguistics}

\newcounter{linecounter}
\newcommand{\linenumbering}{\makebox[2em][r]{\arabic{linecounter}:}}
\renewcommand{\line}[1]{\refstepcounter{linecounter}
\linenumbering}

\def\sgn{{\ensuremath {\rm sgn}}}

\newcommand{\remove}[1]{}

\newcolumntype{L}[1]{>{\raggedright\let\newline\\\arraybackslash\hspace{0pt}}m{#1}}
\newcolumntype{C}[1]{>{\centering\let\newline\\\arraybackslash\hspace{0pt}}m{#1}}
\newcolumntype{R}[1]{>{\raggedleft\let\newline\\\arraybackslash\hspace{0pt}}m{#1}}
\graphicspath{{Figures/}}

\newcommand{\subsubsubsection}[1]{\paragraph[#1]{#1}\mbox{}\vspace{3mm}}

\begin{document}
\begin{frontmatter}

 \title{Persistence in Complex Systems}

\author[sancho]{S. Salcedo-Sanz} \cortext[sss]{Corresponding author: Sancho Salcedo-Sanz. Department of Signal Processing and Communications, Universidad de Alcal\'a. 28871 Alcal\'a de Henares. Madrid. Spain. Ph: +34 91 885 6731, sancho.salcedo@uah.es}
\author[david]{D. Casillas-P\'erez}
\author[javi1,javi2]{J. Del Ser}
\author[carlos]{C. Casanova-Mateo}
\author[sancho]{L. Cuadra}
\author[gus]{M. Piles}
\author[gus]{G. Camps-Valls}

\address[sancho]{Department of Signal Processing and Communications, Universidad de Alcal\'a, Alcal\'a de Henares, Spain}
\address[david]{Department of Signal Processing and Communications, Universidad Rey Juan Carlos, Fuenlabrada, Spain}
\address[javi1]{TECNALIA, Basque Research and Technology Alliance (BRTA), 48160 Derio, Bizkaia, Spain}
\address[javi2]{University of the Basque Country (UPV/EHU), 48013 Bilbao, Bizkaia, Spain}
\address[carlos]{Universidad Polit\'ecnica de Madrid, Madrid, Spain}
\address[gus]{Image Processing Lab (IPL), Universitat de Val\`encia, Val\`encia, Spain}

\begin{abstract}
Persistence is an important characteristic of many complex systems in nature, related to how long the system remains at a certain state before changing to a different one. The study of complex systems' persistence involves different definitions and uses different techniques, depending on whether short-term or long-term persistence is considered. In this paper we discuss the most important definitions, concepts, methods, literature and latest results on persistence in complex systems. Firstly, the most used definitions of persistence in short-term and long-term cases are presented. The most relevant methods to characterize persistence are then discussed in both cases. A complete literature review is also carried out. We also present and discuss some relevant results on persistence, and give empirical evidence of performance in different detailed case studies, for both short-term and long-term persistence. A perspective on the future of persistence concludes the work.
\vspace{0.4cm}
\end{abstract}


%
%

\end{frontmatter}


\clearpage
\tableofcontents
\clearpage
\section{Introduction}

How long will this pandemic, fog event or strong wind last? Will soil retain enough water before it rains again? Is there a relationship between today's stock exchange and yesterday's, or last year's? These, and many other similar scientific questions, drive our daily endeavour to understand complex physical systems. They are all related to the concept of {\em persistence}. Persistence is an important characteristic of many complex systems in nature, in which the variables exhibit some correlation with their past values, or they remain at a certain state for a period of time. In general, it is difficult to give a unique definition of persistence, and its study very often depends on the specific system or process considered. In addition, persistence has been used with different meanings in the literature \cite{Batabyal03}: ``survival'' in ecological systems, ``stability'' in economic systems, or even ``durability'' in chemical and education systems. However, the most useful definition of persistence is related to time persistence of complex systems, i.e. how long a certain variable under study lasts before it changes to a new value or, more specifically, how long does it take for a certain variable to change from one state to another. 

Persistence is closely related to different statistical mechanics properties of system's time series: long-term persistence of different phenomena \cite{Witt13} can be studied by means of correlation-based methods \cite{Bunde02}, Hurst coefficient \cite{Graves17} or Detrended Fluctuation Analysis (DFA) \cite{Peng94,Peng95,Hu01}, which in turn have connections with concepts of time scale invariance \cite{Lesne08} and on-off intermittency \cite{Platt93}, common to a large amount of physical systems or phenomena. Short-term persistence is related to correlation of time series, Markov analysis, AutoRegressive (Integrated) and Moving Average (AR(I)MA) methods for prediction \cite{Box16} and, when exogenous variables are taken into account, it is also related to important properties such as {\em concept drift} in ML prediction systems \cite{concept1,concept2}. On the other hand, the study of persistence has attracted researchers from very different areas in the last years. This concept arises naturally in many complex systems \cite{Comina20}, such as solar-terrestrial physics \cite{Panchev07} or atmospheric and oceanic sciences \cite{Bunde02}, since many related phenomena (rainfall, wind, soil moisture, air and ocean temperature, etc.) are intrinsically persistent. It has recently been also applied to other systems in economy \cite{Grau01,Canarella17}, complex networks \cite{Zou19,Pastor15}, renewable energy \cite{Kosak08,Voyant18}, hydrology \cite{Pelletier97,Yang19}, systems in non-equilibrium thermodynamics \cite{majumdar1999persistence}, optimization \cite{brown1997optimization}, Machine Learning (ML) \cite{concept1}, geophysics \cite{dmowska1999advances,jimenez2005testing}, engineering optimization \cite{brown1997optimization}, or biomedical applications \cite{Depetrillo99}. 

In this paper we review the most important concepts, methods and existing work on persistence in complex systems. We focus both on short-term and long-term temporal persistence, concepts and methods, and on their specific characteristics and  application possibilities in the frame of prediction and analysis problems in complex systems. We pay special attention to how persistence-related concepts and methods can improve our knowledge of a given physical system, and on the computational performance of algorithms and computational methods in specific prediction and analysis problems in complex systems. We also discuss some case studies on the analysis and prediction of persistence-based systems in atmospheric and Earth sciences, renewable energy and machine learning methods, among others. Finally, we give a perspective on the topic, and outline the most promising future lines of research for the broad study of complex systems and nonlinear physical processes.
\begin{figure}[!ht]
\centering 
\includegraphics[width=\columnwidth]{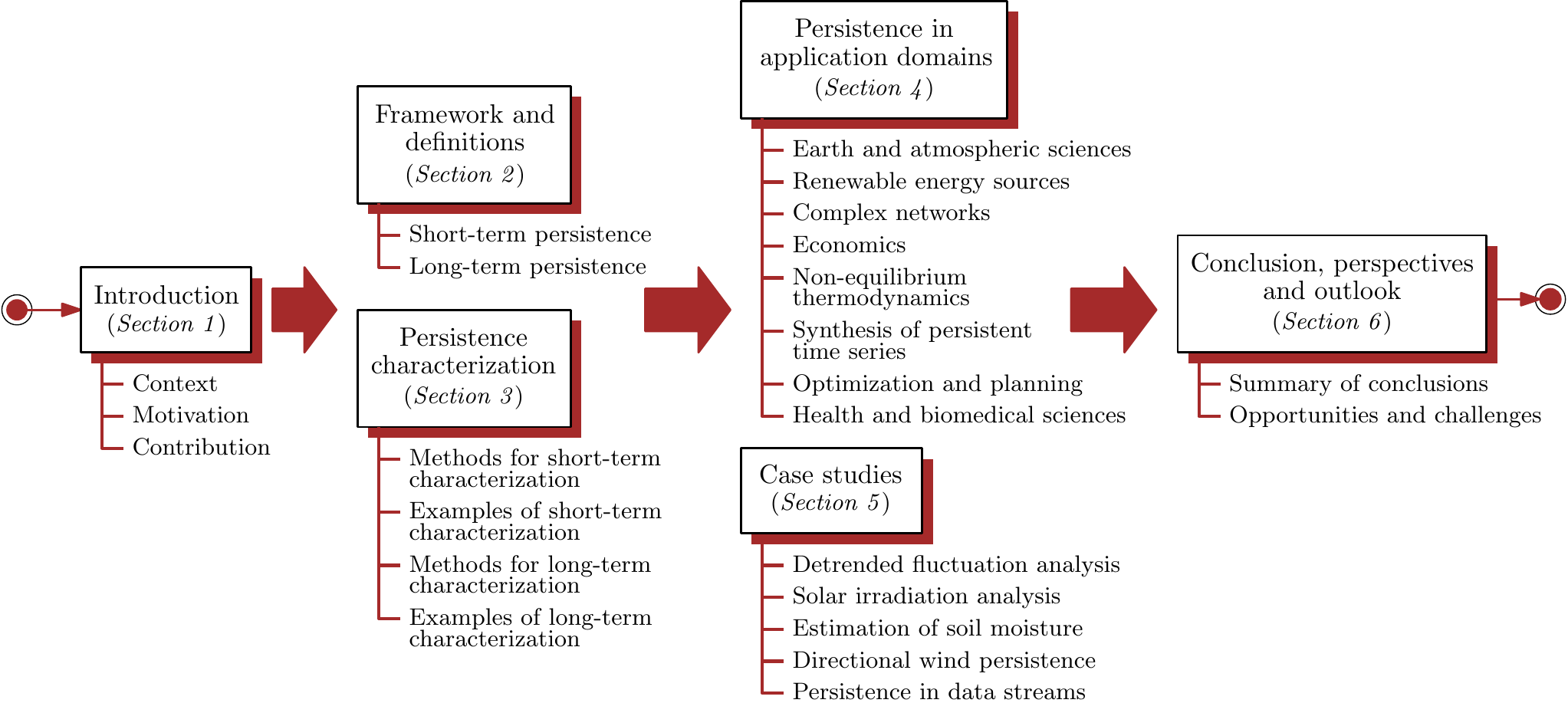}
\caption{Schematic diagram depicting the structure of the paper and a preview of the contents of its sections.}
\label{fig:structure} 
\end{figure}

The remainder of the paper has been structured as shown in Figure \ref{fig:structure}. First, Section \ref{Framework} defines the framework of this study on persistence in complex systems, and formally presents several definitions commonly employed for persistence quantification, both for the short and long-term cases. Section \ref{Persistence_Characterization} discusses different methods for persistence characterization, also for short-term and long-term cases. We support the description of these methods for persistence characterization with some illustrative examples. Section \ref{Literature_Review} presents a complete literature review on persistence in application domains. Key topics such as persistence in Earth and climate sciences, renewable energy, complex networks, economics, non-equilibrium thermodynamics systems, optimization and planning and health and biomedical applications are discussed. In Section \ref{Case_Studies}, we discuss some specific case studies with relevant results on persistence in different complex systems: DFA and time scales, solar radiation prediction based on persistence, e-folding time of soil moisture, directional persistence of wind speed in wind farms, and persistence in ML approaches over data streams. Section \ref{Conclusions} closes the paper with some final conclusions, outlook and future perspectives of this field.

\section{Persistence: framework of the study and definitions}\label{Framework}

Statistical modeling of time series has traditionally considered two types of persistence \cite{Witt13}:
\begin{enumerate}
\item {\em Short-term persistence}, in which time series values are correlated to other values that are in a close temporal neighbourhood with one another.
\item {\em Long-term persistence}, in which all or almost all time series values are correlated with one another, i.e., values are correlated with one another at very long lags in time \cite{Beran94}.
\end{enumerate}
Both the mathematical treatment (methods applied) and characterization of these persistence types are different, and thus deserve a different treatment. 

\subsection{Short-term persistence}\label{short-term-definitions}
Let us consider a complex system (physical, ecological, socio-economical, etc.) that can be characterized by a time-dependent, $t$, variable under study and subject to the influence of external/exogenous drivers $\boldsymbol{\theta}$, that is $x(t,\boldsymbol{\theta})$. 

Formally, the considered variable $x(t,\boldsymbol{\theta})$ is a continuous function defined as follows:
\begin{equation}
    \begin{aligned}
            x(t,\boldsymbol{\theta}): \mathbb{R}\times\mathbb{R}^n &\to \mathbb{R}\\
            (t,\mathbf{\boldsymbol{\theta}}) &\mapsto x(t,\boldsymbol{\theta})
    \end{aligned}
\end{equation}
Generally, we do not count on the whole continuous variable, which includes the evaluation of a dense set of points, but we have access to a sequence of its samples, equally acquired with a period $T_s$, which constitutes the temporal series $x[n,\boldsymbol{\theta}]$ associated with $x(t,\boldsymbol{\theta})$:
\begin{equation}
    \begin{aligned}
            x[n,\boldsymbol{\theta}]: \mathbb{N}\times\mathbb{R}^n &\to \mathbb{R}\\
            (n,\boldsymbol{\theta}) &\mapsto x[n] = \left.x(t,\boldsymbol{\theta})\right|_{t=nT_s}.
    \end{aligned}
\end{equation}
This temporal sequence $x[n,\boldsymbol{\theta}]$ allows making discrete the original continuous domain of the variable $x(t,\boldsymbol{\theta})$. The period $T_s$ should be carefully chosen according to the speed rate of variation of the variable $x(t,\boldsymbol{\theta})$ fulfilling the Nyquist-Shannon theorem, which results in no loss of information~\cite{nyquist1928certain}.
For convenience, in the notation hereafter we will skip the dependencies of the variable $x(t,\boldsymbol{\theta})$ on parameters $\boldsymbol{\theta}$. Consequently, we will simply write $x(t)$ for the continuous variable or $x[n]$ for its associated discrete sequence. 

Let us now define the set of possible {\em states} of a system $\mathcal{S}=\lbrace s_i\,|\, 1\leq i \leq N_s\rbrace\subset\mathbb{N}$ as a subset of integers which represent finite configurations of the system, where it can remain during a period of time. We define the state sequence $s[n]$ (the state variable) associated to the variable $x[n]$ as the following discrete function:
\begin{equation}
    \begin{aligned}
            s[n]: \mathbb{N} &\to \mathcal{S}\\
            n &\mapsto s[n] = \delta(x[n]),
    \end{aligned}
\end{equation}
where the threshold function $\delta(x)$ maps each value of the variable $x[n]$ into the state space $\mathcal{S}$. In general, the threshold function $\delta(x)$ is built as a piece-wise function from a set of ordered thresholds $\lbrace \delta_1\leq \cdots \leq \delta_i \leq \cdots \leq \delta_{N_s-1}\rbrace\subset\mathbb{R}$ as follows:
\begin{equation}\label{eq:thresholdF}
    \begin{aligned}
            \delta(x): \mathbb{R} &\to \mathcal{S}\\
            x &\mapsto \delta(x) = \left\lbrace
            \begin{aligned}
                s_1 & \text{ if } x \leq \delta_1\\
                s_2 & \text{ if } \delta_1 < x \leq \delta_2\\
                \vdots & \text{    } \qquad\vdots\\
                s_{N_s} & \text{ otherwise.}
            \end{aligned}\right.
    \end{aligned}
\end{equation}

Note that this state variable $s[n]$ is more useful than $x[n]$ to define the concept of persistence, since it is related to system's state changes. Also, the residual series associated with $x[n]$, namely, $r[n]=\left|x[n+1]-x[n]\right|$, can be useful to detect state changes, in such a way that if $r[n]<\epsilon  \Rightarrow s[n]=s_i$, but if $r[n]>\epsilon \Rightarrow s[n]=s_j$ ($j \neq i$). In words, we consider both the variable $x[n]$ to set the actual state of the system $s[n]$ (depending on its current value using the threshold function~\eqref{eq:thresholdF}), and the variations larger than $\epsilon$ in $r[n]$  to detect a change in the state of the system.

As an example of system state construction, we consider a case using time series of visibility values. Here $x[n]$ stands for a temporal sequence of visibility values (two complete years, 2018 and 2019, of hourly data) from a measurement station deployed on the A8 highway in Lugo, Spain. The station counts with a visibilimiter device with a maximum value of $2000$m from the viewpoint. The zone under study suffers from low visibility episodes due to persistent orographic fogs. We can define a state variable for this phenomenon, considering a binary variable ($\mathcal{S}=\lbrace 0,1\rbrace$) for the fog, so if we have a fog event, $s[n]=1$, whereas $s[n]=0$ otherwise. Since we consider two states for the variable, the threshold function $\delta(x)$ only uses one single threshold $\delta_1$, which in this case represents a binary threshold function:
\begin{equation}\label{binarizing_function}
s[n]=\left\{
\begin{array}{l l}
1&\text{if}~~x[n] < \delta_1,\\
0&\text{otherwise.}\\
\end{array}
\right.
\end{equation}

\begin{figure}[!ht]
\begin{center}
\subfigure[]{\includegraphics[draft=false, angle=0,width=8cm]{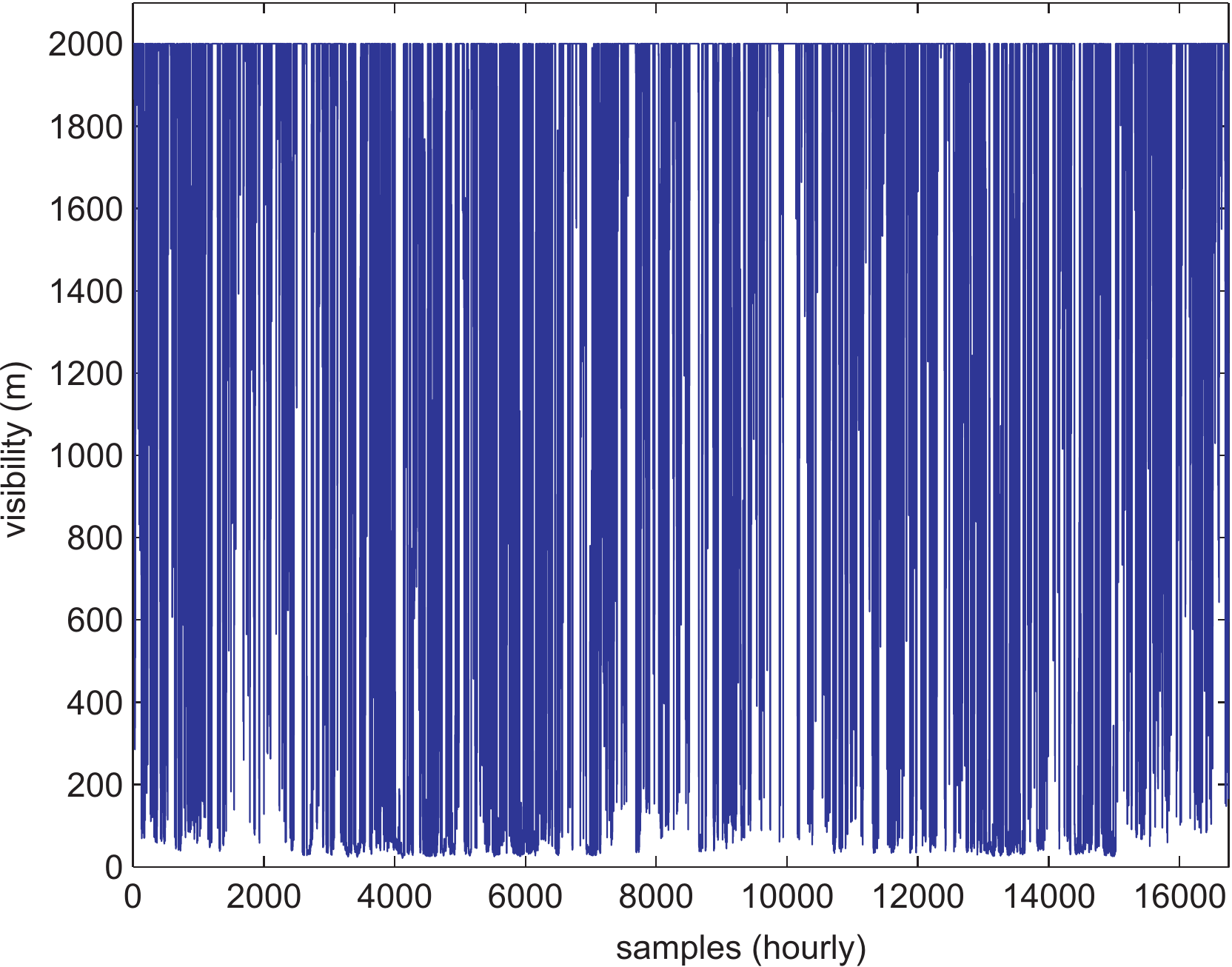}}
\subfigure[]{\includegraphics[draft=false, angle=0,width=8cm]{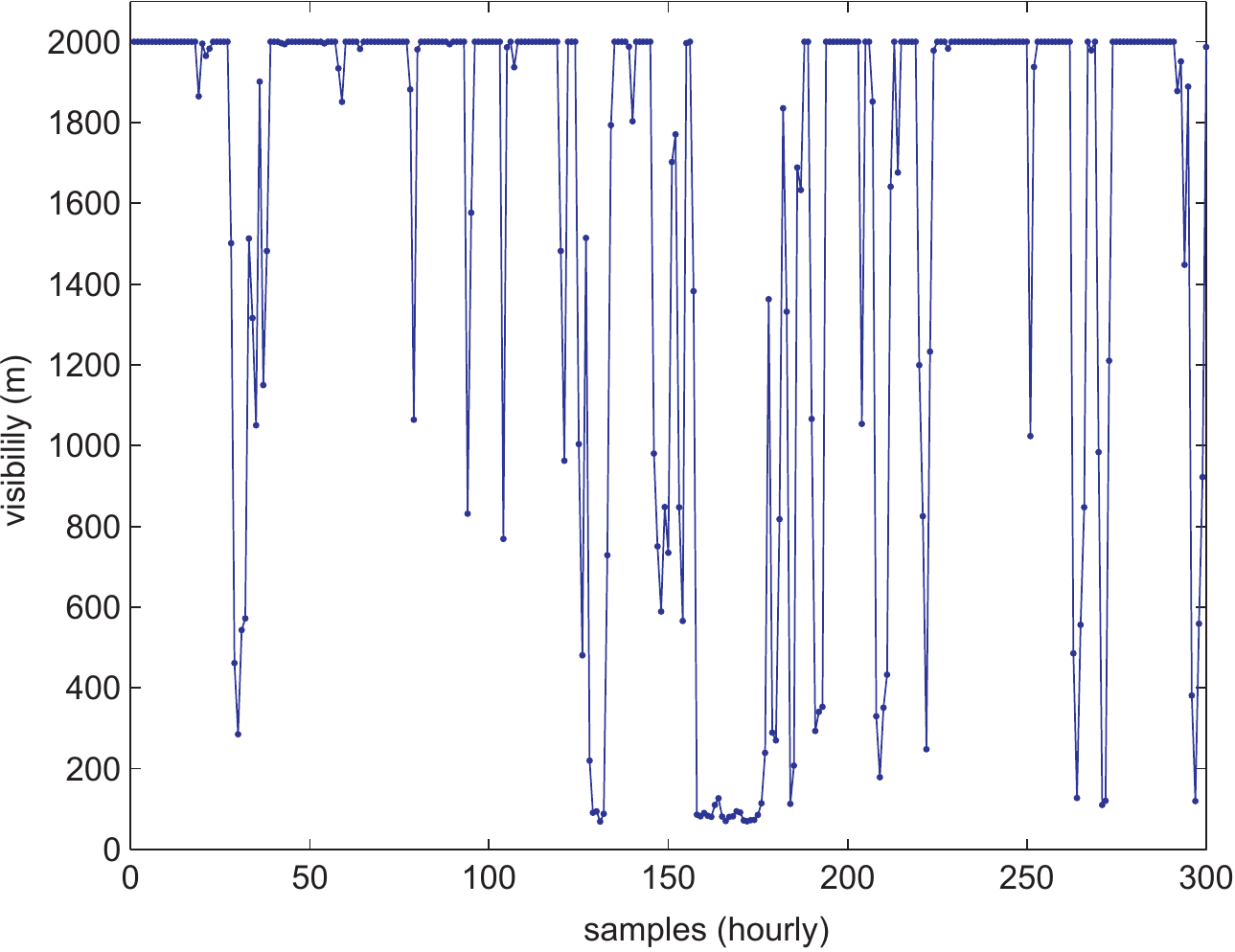}}
\end{center}
\caption{\label{Ejemplo_Sk_visib_x} Visibility time series (hourly resolution) at Lugo station, Galicia, Spain; (a) Complete $x[n]$ time series (two years of hourly measurements); (b) zoom at the first 300 hours of the time series.}
\end{figure}

\begin{figure}[!ht]
\begin{center}
\subfigure[~$\delta_1=1950$]{\includegraphics[draft=false, angle=0,width=7cm]{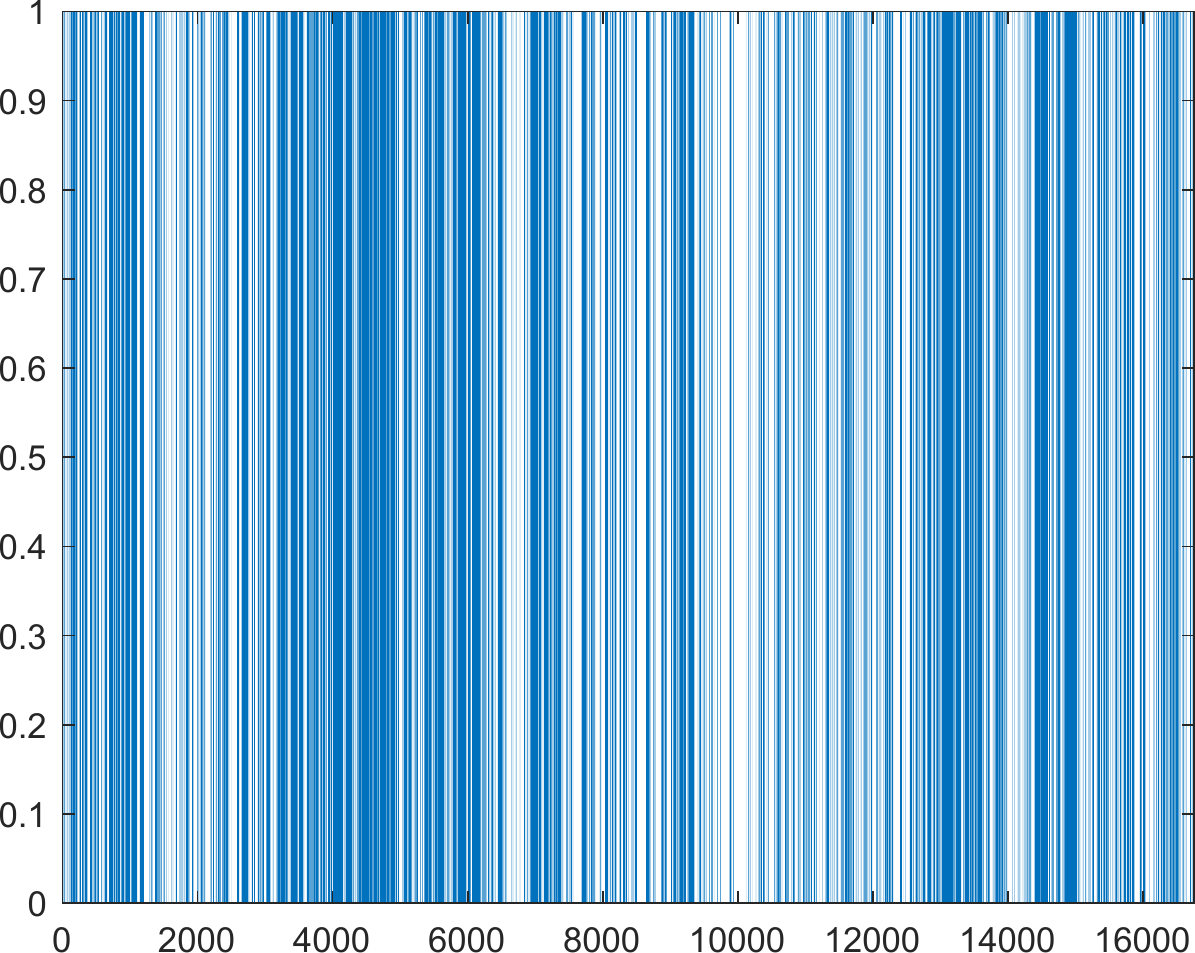}}
\subfigure[~$\delta_1=500$]{\includegraphics[draft=false, angle=0,width=7cm]{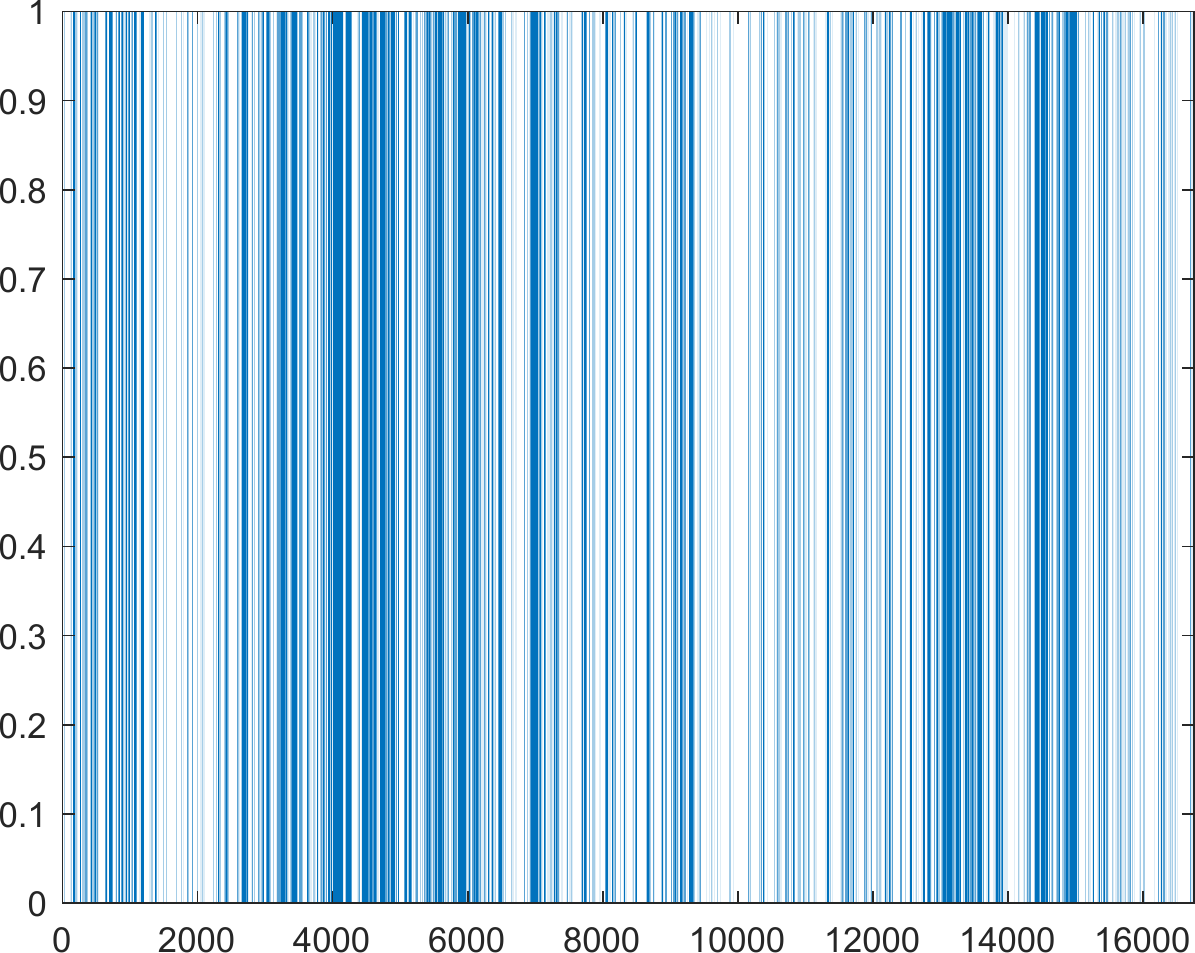}}
\subfigure[~$\delta_1=100$]{\includegraphics[draft=false, angle=0,width=7cm]{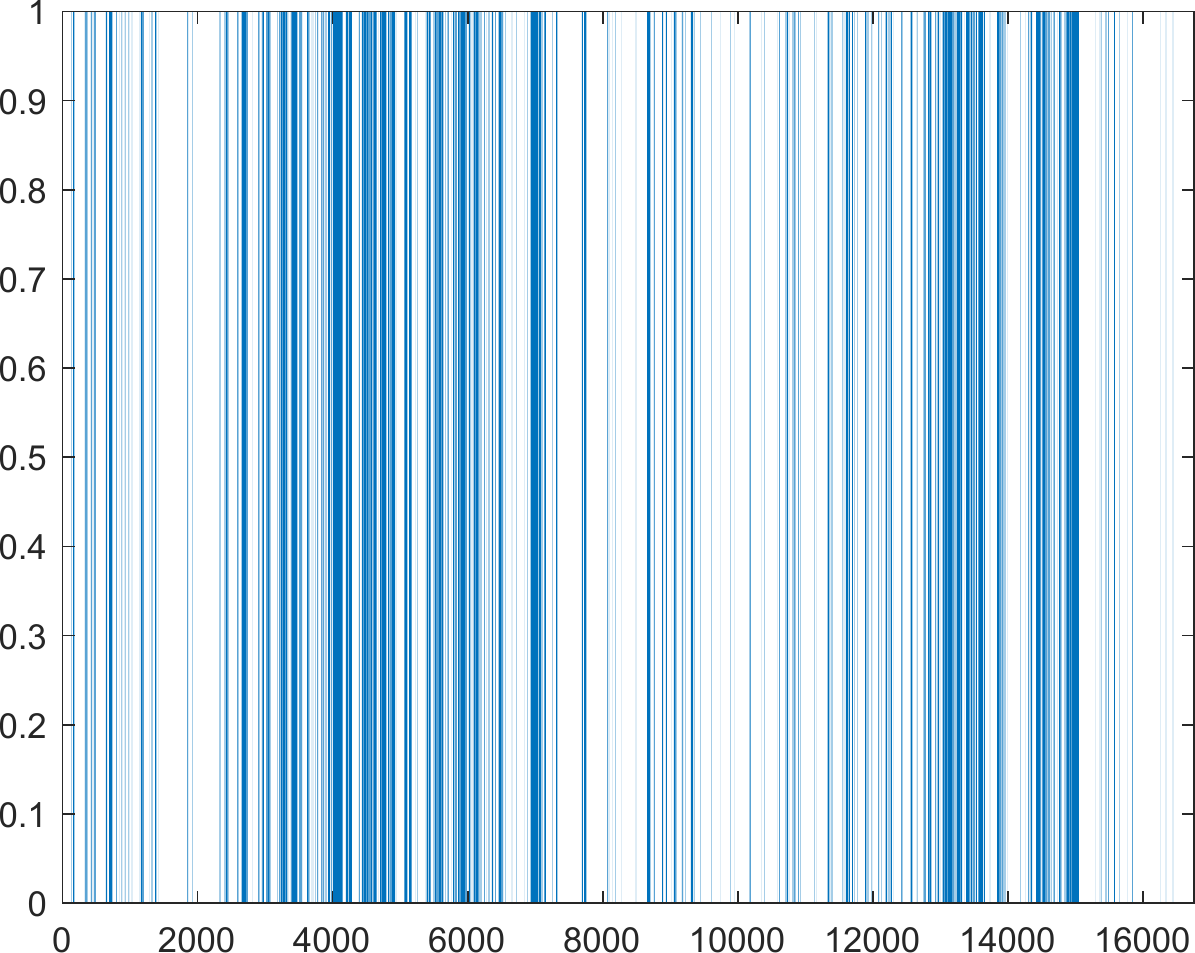}}
\end{center}
\caption{\label{Ejemplo_Sk_visib}
State variables $s[n]$ time series (binary state $s_1=1,s_2=0$) obtained from the original visibility time series $x(t)$ by applying Expression (\ref{binarizing_function}) with different values of the threshold $\delta_1$. We have represented as blue lines the values of the state variable $s[n]=1s$, and as white lines the values $s[n]=0$; (a) $\delta_1=50$; (b) $\delta_1=1500$; (c) $\delta_1=1900$.}
\end{figure}

We can set the threshold $\delta_1$ to different values, obtaining this way different definitions of the system state. Different $\delta$ values mean different definitions of what we consider a low-visibility situation. Figure \ref{Ejemplo_Sk_visib_x} shows the initial time series of visibility values $x[n]$. Figure \ref{Ejemplo_Sk_visib} shows time series of different states variables (binary series) obtained by varying the threshold $\delta_1$. We can consider this way different situations as low-visibility event (depending on $\delta_1$). For example, in Figure \ref{Ejemplo_Sk_visib} (a), $\delta_1=1950$, so we consider as low-visibility event every value of $x[n]$ under $1950$ m. On the other hand, in Figure \ref{Ejemplo_Sk_visib} (c) $\delta_1$ is set to $100$, so we only consider as low-visibility event every value of $x[n]$ under $100$ m. As can be seen in these figures, the state variables time series are completely different from each other. Based on these state variables we can work towards obtaining short-term persistence measures.

Considering the previous definitions, the state sequence $s[n]$ associated to the variable $x(t)$ of a complex system (sampled with a period $T_s$) can be represented as:
\begin{equation}
    s[n]=\left[\ldots, s_i, s_i,\ldots, s_i, s_j, s_j, \ldots ,s_j,s_k,s_k,\ldots\right].
\end{equation}
with $s_i,s_j,s_k\in\mathcal{S}$. In practice, the state sequence $s[n]$ will have a finite length denoted as $N$. We refer as $N_i$ to the number of samples over which the state sequence $s[n]$ remains at state $s_i$.

Intuitively, the persistence of the system (noted as $\mathcal{P}$ in the rest of the manuscript) is better characterized by considering the time sequence of the system's states rather than using the original continuous-variable time series $x(t)$. Before discussing the different possible definitions for the persistence of a system, it is important to note that the majority of works in the literature treats the variable under study $x(t)$ and its corresponding time $x[n]$ or the state sequence $s[n]$ as stochastic processes, under a statistical framework. The persistence of the variables $x(t)$ or $x[n]$ can be determined in terms of usual statistical operators, such as probability vectors or expected time of change among the states. Some works assume important statistical properties, such as the stationary, wide-stationary or Markov property driving to different persistence definitions.

Therefore, there exist several definitions of short-term persistence in the literature. The majority of them estimate the expected time $T_i$ in which the system remains at a given state $s_i$ before changing to another state $s_j$ ($i \neq j$), so that an average for all states $\mathcal{S}$ is provided. A plausible definition, previously used in the literature \cite{Batabyal03}, is to measure the persistence of a system as the expected time that a system remains at a given state:
\begin{equation}
\mathcal{P}^E= \frac{1}{N}\sum_{i=1}^N  \operatorname{E}[T_i],
\end{equation}
where $\operatorname{E}[\cdot]$ is the expectation operator, and $T_i$ stands for the time that the system remains at a given state $s_i$. The larger $\mathcal{P}^E$, the more persistence the system shows. 

Another persistence definition takes into account the {\em burst-length} in the state sequence $s[n]$. Here the persistence is given by the maximum burst length as:
\begin{equation}\label{P_i^b}
    \mathcal{P}_i^b = \frac{\max(T_i)}{N_i},
\end{equation}
where $N_i$ stands for the number of times that the state sequence is in state $s_i$.

An alternative definition considering all the states in the system ${\mathcal S}$ is possible, such as the maximum $\mathcal{P}_i^b$ over the $|{\mathcal S}|$ possible states of the system is considered:
\begin{equation}\label{P^b_1}
\mathcal{P}_{\max}^b = \frac{1}{|\mathcal{S}|} \max \left(\mathcal{P}_i^b\right) = \frac{1}{|\mathcal{S}|} \max \left(\frac{\max(T_i)}{N_i}\right),
\end{equation}
or its mean value:
\begin{equation}\label{P^b_2}
\mathcal{P}_{\mathcal{S}}^b=\frac{1}{|\mathcal{S}|}\sum_{i=1}^{|\mathcal{S}|} \frac{\max(T_i)}{N_i},
\end{equation}
or any other variation of this idea that could be convenient for the problem at hand.

Yet another related definition deals with the probability of the system to remain $N_r$ consecutive times steps in a given state $s_i$ before changing to another state $\mathcal{S}{\setminus}\lbrace s_i \rbrace$, or in other words, the probability of having large values of $N_r$:
\begin{equation}
    \mathcal{P}^P_i(N_r)=P\left(s[0] = s_i,s[1] = s_i, \cdots,s[N_r] = s_i, s[N_r+1] \in \mathcal{S}\setminus\lbrace s_i\rbrace \right).
\end{equation}

Assuming the Markov's property over the state sequence, the persistence can be also defined without even specifying the system's state lifetime $N_r$, just as the probability that the system remains at a certain state (the larger this probability, the larger the expected the persistence of the system):
\begin{equation}
\mathcal{P}^M=P\left(s[n+1] = s_i ~|~ s[n] = s_i\right)
\end{equation}
Recall that there is no standard definition of persistence, so other related definitions different from the described are possible.

As an example of the previous definitions, we consider a Pure Persistent System (PPS), i.e. a system that will remain in the same state $s_i$ indefinitely from the beginning, i.e. $x[n+1]=x[n]~\forall n$. In terms of the residual variable $r[n]$, it is defined as $r[n] < \epsilon~\forall n$. In this particular case (considering $n \rightarrow \infty$), $\mathcal{P}^E \rightarrow \infty$, $\mathcal{P}^b \rightarrow \infty$ for the current state $s_i$, and $0$ for the rest, $\mathcal{P}^b_{max}=\infty$ and $\mathcal{P}^b_{\mathcal{S}}=\infty$, $\mathcal{P}_i^P = 0$ for state $s_i$, and $\mathcal{P}^M = 1$. On the other hand, in systems which are not persistent at all, $\mathcal{P}^E \rightarrow 1$, $\mathcal{P}_i^b \rightarrow 1$, $\mathcal{P}^P \rightarrow 1$ for $n = 1$ (and $\mathcal{P}^P \rightarrow 0$ for $n > 1$) and $\mathcal{P}^M \rightarrow 0$. In systems that are not PPS, but have some degree of persistence, the estimation and characterization of $\mathcal{P}$ is a challenging problem.

\subsection{Long-term persistence}

Long-term persistence of systems are usually defined from the autocorrelation function of the time series $x[n]$ involved. In general, the autocorrelation function for a sequence $x[n]$ is defined as:
\begin{equation}\label{aut_corr}
    r_x\left[k\right]=\operatorname{E}\left[x[n] x[n+k]\right].
\end{equation}
In this context $k$ stands for a given time lag variable. Sometimes, the autocorrelation is standardized by the variance and mean-removed (also called as normalized autocovariance), yielding:
\begin{equation}\label{aut_corr2}
    r_x\left[k\right]=\frac{1}{\sigma_x^2}\operatorname{E}\left[(x[n]-\bar{x}\right)\left(x[n+k]-\bar{x})\right],
\end{equation}
where $\bar{x}$ is the sample mean, and $\sigma_x^2$ is the sample variance.

Following \cite{Witt13,Beran94}, a long-range persistent process exhibits a power-law scaling of the normalized autocorrelation (as given by Equation (\ref{aut_corr2})):
\begin{equation}
    \left|r_x[k]\right| \sim k^{-(1-\beta)},
\end{equation}
which holds for large time lags $k$. The
parameter $\beta\in(-1,1)$ is the strength of long-range persistence, so that $\beta=0$ represents a process that has no long-range persistence between values, $\beta>0$ indicates a long-range persistence, and $\beta<0$ aims at long-range anti-persistence.

For non-stationary times series, i.e. those cases in which $\beta>1$, it is not appropriate to use the autocorrelation function to estimate long-term persistence \cite{Witt13}, since the mean $\bar{x}$ is involved in the autocorrelation calculation. In these cases, an alternative, more efficient way to measure long-range correlations is the semivariogram \cite{Matheron63}, which is given by the following expression:
\begin{equation}
\gamma[k]=\frac{1}{2(N-k)}\sum_{n=1}^{N-k}\left(x[n+k]-x[n]\right)
\end{equation}
For a time series in which $\beta>1$ the semivariogram, $\gamma[k]$, scales with the lag $k$ as:
\begin{equation}
\gamma[k]\sim k^{2\text{\texthvlig}},  
\end{equation}
where \texthvlig~is the Hausdorff exponent. Note that \texthvlig~is related to the persistence strength as $\text{\texthvlig} =\frac{\beta-1}{2}$ \cite{Witt13}.

Long-term persistence can be alternatively defined in terms of the power spectral density $S_x(f)$ of the time series $x(t)$. A process can be defined as long-range persistent if $S_x(f)$ scales asymptotically as a power law for frequencies close to the origin:
\begin{equation}
S_x(f)=\mathscr{F}\left[r[n]\right] \sim f^{-\beta},
\end{equation}
where $\mathscr{F}[\cdot]$ refers to the Fourier transform, and the power-law exponent $\beta$ measures the strength of persistence. This expression is valid for all values of $\beta\in\mathbb{R}^+$.

\section{Persistence characterization: methods and examples}\label{Persistence_Characterization}

In this section, we describe different methods for persistence characterization, both in short-term and long-term cases. We also show several examples of persistence characterization when statistical properties of the time series $x[n]$ or $s[n]$ are known.

\subsection{Methods for short-term persistence characterization}

\subsubsection{Markov chain models for persistence characterization}\label{Markov_chains}

\begin{figure}[t]
 \centering
  \subfigure[]{
   \label{Sys1}
    \includegraphics[width=0.28\textwidth]{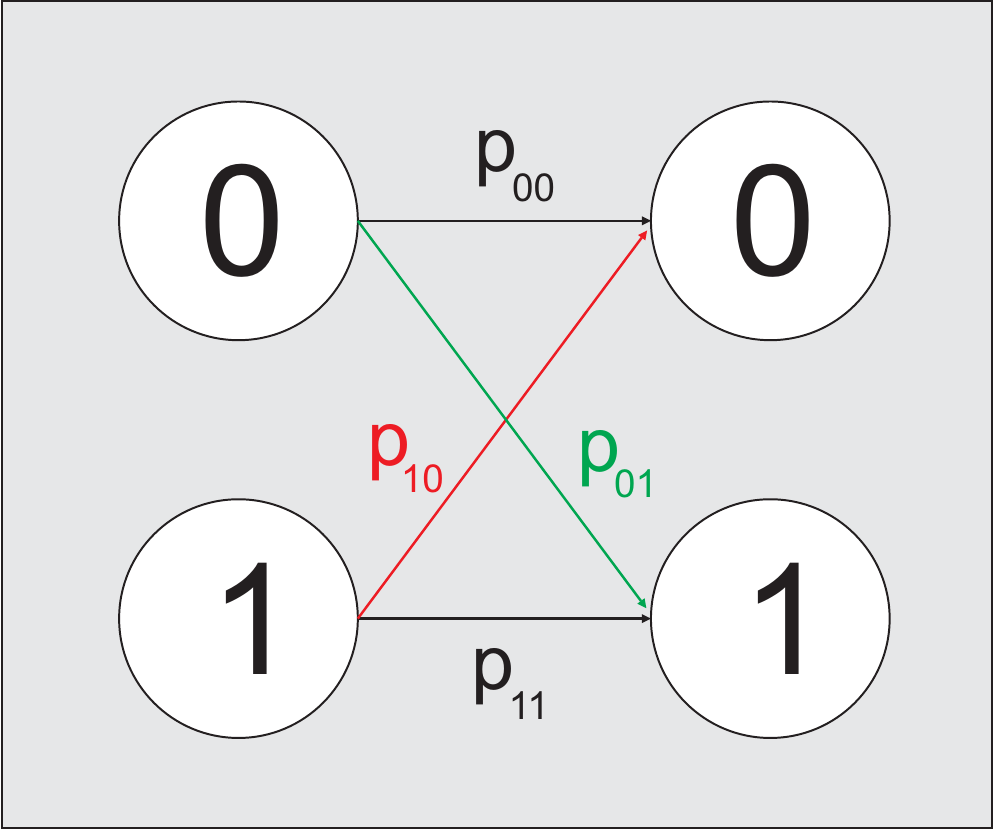}}
  \subfigure[]{
   \label{SysN}
    \includegraphics[width=0.48\textwidth]{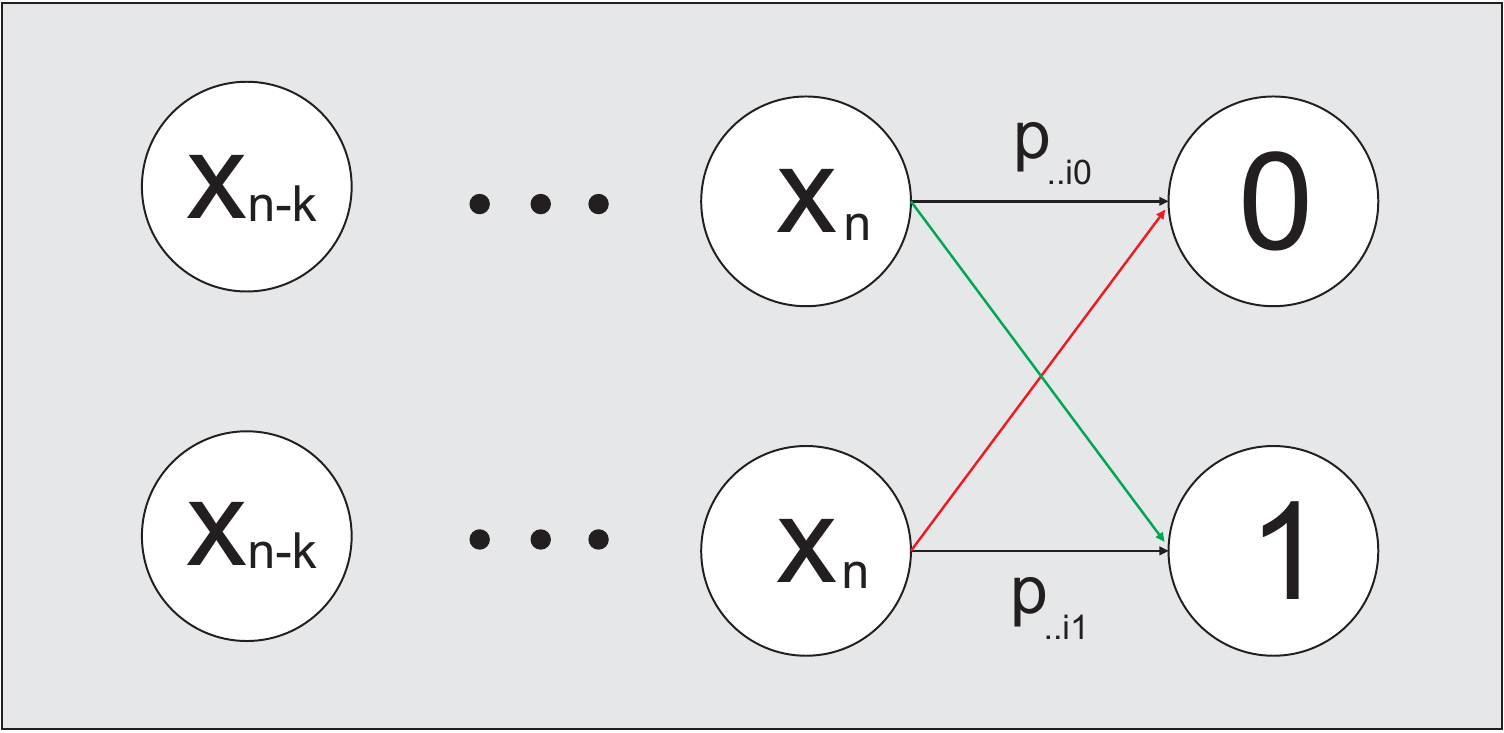}}
 \caption{Diagrams representing a two-state (a) First-order; (b) $k$-order Markov process. }
 \label{DiagProb}
\end{figure}

Markov Chain Models (MCM) are statistical tools useful to analyze persistence in discrete binary time series. Let us consider a system described by means of a binary state variable $s[n]$ which can take values $1$ (positive event), or $0$ (negative event). Assuming that the occurrence of an event in the system at time $n$ depends only on the state at time $n-1$ (first-order Markov chain), the transition probabilities of positive and negative events can be divided into the following four cases:
\begin{equation}
\begin{array}{l}
{p_{00}=\operatorname{P}\left( s[n]=0  \left.\right| s[n-1]=0\right)} \\
{p_{01}=\operatorname{P}\left( s[n]=1  \left.\right| s[n-1]=0\right)} \\
{p_{10}=\operatorname{P}\left( s[n]=0  \left.\right| s[n-1]=1\right)} \\
{p_{11}=\operatorname{P}\left( s[n]=1  \left.\right| s[n-1]=1\right)}
\end{array}.
\end{equation}
An explanatory diagram is represented in Figure \ref{DiagProb} for a first-order and $k$-order MCM.

The estimation of these transition probabilities can be computed through the conditional relative frequencies as:
\begin{equation}
    \begin{array}{l}
        \hat{p}_{00}=\frac{N_{00}}{N_{0}} \hspace{1cm}
        \hat{p}_{01}=\frac{N_{01}}{N_{0}} \\
        \hat{p}_{10}=\frac{N_{10}}{N_{1}} \hspace{1cm}
        \hat{p}_{11}=\frac{N_{11}}{N_{1}}
    \end{array},
\end{equation}
where $N_{ij}$ represents the number of transitions from state $s_i$ to state $s_j$, and $N_i$ is the number of states $s_i$ followed by any other data point, i.e., $N_i = N_{i0} + N_{i1}$. The subscripts refer to the state, $i,j\in\lbrace0,1\rbrace$. Note that the Na\"ive persistence operator is a special case of first-order Markov chain, whose formula $x[n+1]=x[n]$ forces the state preservation at any time, and it can also be described with the following transition probability matrix:
\begin{equation}
\mathbf{P} =
\begin{pmatrix}
1 & 0\\
0 & 1
\end{pmatrix}
\label{Naive}
\end{equation}
Note that $\mathcal{P}^M$ can be easily calculated in this case by means of the transition probabilities $p_{00}$ and $p_{11}$, as:
\begin{equation}
\mathcal{P}^M=\frac{p_{00}+p_{11}}{2},
\end{equation}
and hence the Na\"ive persistence operator characterized by transition matrix of Equation (\ref{Naive}) has $\mathcal{P}^M=\frac{1+1}{2}=1$.

For a higher-order MCM, the transition probabilities take into account the states at the considered time windows. For example, in a second-order Markov chain, the states at times $n-2$ and $n-1$ are considered to predict the state $n$; or in a third-order Markov chain, the states at times $n-3$, $n-2$ and $n-1$ are taken into account to predict state $n$. The transition probabilities for second, third and higher-order models are  respectively defined as:
\begin{eqnarray}
p_{i_2 i_1 i} = \operatorname{P}\left(s[n]=i ~|~ s[n-1]=i_{1}, s[n-2]=i_{2}\right), \quad i, i_{1}, i_{2} \in \lbrace 0,1\rbrace,\\
p_{i_3 i_2 i_1 i} = \operatorname{P}\left(s[n]=i ~|~ s[n-1]=i_1, s[n-2]=i_2, s[n-3]=i_3\right), \quad i, i_1, i_2, i_3 \in \lbrace 0,1 \rbrace,\\
p_{i_1 \cdots i_k i}=\operatorname{P}\left(s[n]=i ~|~ s[n-1] = i_1,\cdots,s[n-k] = i_N\right), \quad i,i_1,\cdots,i_k\in\lbrace0,1\rbrace.
\end{eqnarray}

\subsubsection{Auto-regressive models: ARMA and ARIMA models}

Auto-Regressive-Moving-Average (ARMA) models \cite{Box16} are often used for modeling and forecasting time series, $x[n,\boldsymbol{\theta}]$, in particular when the dependency on the external variables $\boldsymbol{\theta}$ is not evident or unknown, and thus only the values of $x[n]$ are observed. In general, ARMA models provide a parsimonious description of a stationary stochastic process, in terms of two polynomials, one for the dependency on past values of the series (auto-regressive (AR)) and the other one for the dependency on past errors (residuals) of the series (moving-average (MA)):
\begin{equation}
x[n]=a_1 x[n-1]+\cdots+a_p x[n-p]+\varepsilon[n]-b_1 \varepsilon[n-1]-\cdots-b_q\varepsilon[n-q],
\end{equation}
where $\lbrace a_i \rbrace_{i=1}^p$ are the auto-regressive parameters (to be estimated), $\lbrace b_i \rbrace_{i=1}^q$ are the moving average coefficients (also to be estimated), and $\varepsilon$ stands for a series of unknown errors (residuals), usually assumed to follow a zero mean Gaussian probability distribution. This model for a time series is known as an ARMA($p$,$q$) model.

In the description of ARMA models it is customary to use the back-shift operator $\operatorname{B}[\cdot]$ to write the models in a more compact form. This operator has the effect of changing the time $n$ to $n-1$: $\operatorname{B}[x[n]]=x[n-1]$, $\operatorname{B}^2[x(t)]=x[n-2]$ and so on. Note that using this operator, the ARMA model can be rewritten as:
\begin{equation}
\left(1-a_1 \operatorname{B}-\cdots-a_p \operatorname{B}^p\right)[x[n]]=\left(1-b_1 \operatorname{B}-\cdots-b_p \operatorname{B}^q\right)[\varepsilon[n]],
\end{equation}
or, if we define $\Phi_p[\cdot]=\left(1-a_1 \operatorname{B}-\cdots-a_p \operatorname{B}^p\right)[\cdot]$ and $\Theta_q[\cdot]=\left(1-b_1 \operatorname{B}-\cdots-b_p \operatorname{B}^q\right)[\cdot]$ as linear operators, we find:
\begin{equation}
    \Phi_p(\operatorname{B})[x[n])]= \Theta_q(\operatorname{B})[\varepsilon[n]].
\end{equation}
The main advantage of using this notation is that we can detect equivalent models, and therefore simplify them. As an example let us consider the following ARMA(2,1) model:
\begin{equation}
x[n]=0.8 x[n-1]-0.12x[n-2]+\varepsilon[n]-0.6\varepsilon[n-1].
\end{equation}
We can rewrite this model by considering $\Phi_2[\cdot]=(1-0.8\operatorname{B}+0.12\operatorname{B}^2)[\cdot]$ and $\Theta_1[\cdot]=(1-0.6B)[\cdot]$, so:
\begin{equation}
\Phi_2[x[n]]=\left(1-0.8\operatorname{B}+0.12\operatorname{B}^2\right)[x[n]]=
\Theta_1[\varepsilon[n]]=\left(1-0.6\operatorname{B}\right)[\varepsilon[n]].
\end{equation}
Note, however, that $\Phi_2[\cdot]$ can be factorized as $\Phi_2[\cdot]=\left(1-0.2\operatorname{B}\right)\left(1-0.6\operatorname{B}\right)$, leading to:
\begin{equation}
\left(1-0.2\operatorname{B}\right)\left(1-0.6\operatorname{B}\right)x[n]=\left(1-0.6\operatorname{B}\right)\varepsilon[n].
\end{equation}
We can cancel the term $\left(1-0.6B\right)$ from both sides of the equation, thereby obtaining a much simpler equivalent model:
\begin{equation}
x[n] =0.2x[n-1]+\varepsilon[n].
\end{equation}
ARIMA models are applied in some cases where data show evidence of non-stationarity. In these cases, ARMA models cannot be directly applied, and an initial differentiation step (corresponding to the ``integrated'' part of the model) is applied to eliminate the non-stationarity. Following the previous notation, the general form of an ARIMA($p$,$d$,$q$) model is given by:
\begin{equation}
\Phi_p(B)[(1-B)^d [x[n]]]= \Theta_q(B)[\varepsilon(t)].
\end{equation}

\subsubsection{Persistence and memory effects in neural networks}\label{Pers_NN}

When nonstationary and nonlinear processes are present, ARMA and ARIMA models described above are no longer useful to characterize time series data, and thus quantification of persistence from there is often compromised. Other options are needed, and the field of neural networks has appeared as an appealing framework to deal with memory and persistence in such more complex scenarios. In this section we will describe the most common neural networks models with important memory properties which allows persistence modeling and characterization. Note that the definition of persistence from neural networks is not so explicit as in the case of the models defined before, since in many cases neural networks act as blackbox models extracting information from data. Thus, we will focus on defining different types of neural networks which can be useful for persistence analysis of time series, without a fully mathematical characterization of the networks, which is far from the scope of this work.

Artificial Neural Networks (ANNs) are algorithmic models inspired in biological neural structures, and how the brain processes information. A seminal inspiration is found in the neural Hebbian theory \cite{hebb1949organization}, which states that the persistent and repetitive stimulation of axon cells changes their metabolism and gives rise to an adaptation of brain neurons during the learning process. This theory revolutionized the field of neuroscience and behavioural psychology, but also changed the way signal processing was done so far. Actually, inspired by the Hebb hypothesis, in \cite{little1975statistical} a statistical theory of short and long-term memory was proposed for modeling the statistical nature of synaptic signal transmission. The field that studies and develops artificial neural networks drinks from these sources and designs algorithms that (try to) mimic brain processes of learning, adaptation, and memory. Since then, many neural networks have been proposed with all kind of both statistical and biological inspiration to account for information transmission and optimal processing, especially when nonstationary, non-Gaussian and nonlinear processes are involved. 

\subsubsubsection{ARMA modeling with neural networks} 

\begin{figure}[t!]
    \centering
    \includegraphics[width=16cm]{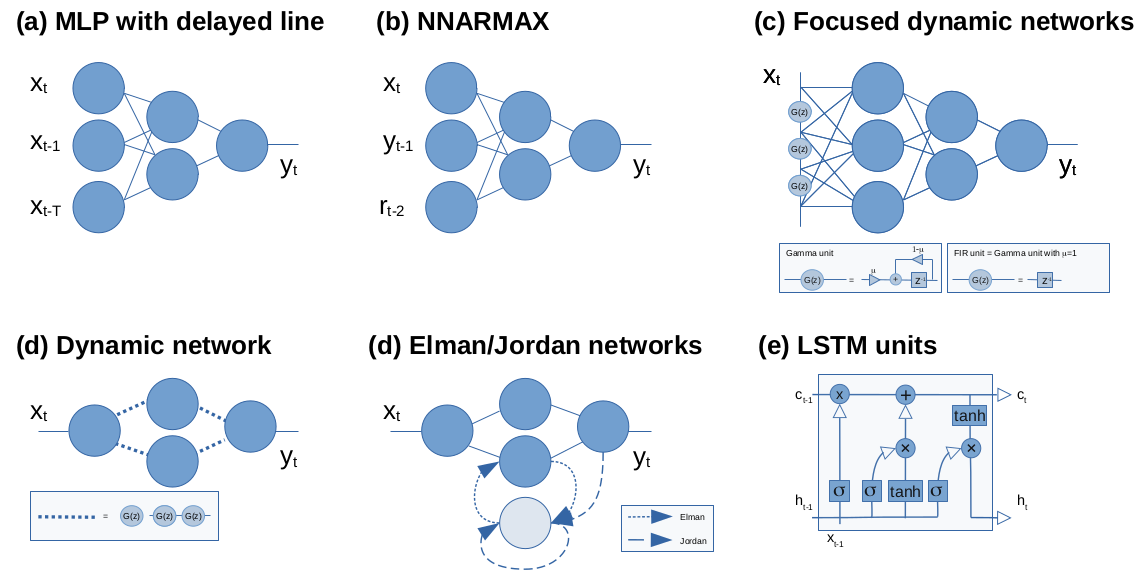}
    \caption{{\bf Taxonomy of dynamic neural networks accounting for memory and recurrent states.} (a) A standard way to consider memory in a neural network is feeding a multilayer perceptron with a tap-delayed line with a prescribed time embedding of the input signal $x_t$; (b) NNARMAX models are inspired in the standard ARMA modeling and essentially feed the network with a time embedding of both the past input values $x_{t-k}$, output values $y_{t-q}$ and residuals/errors $r_t-s$; (C) a focused neural network extends the previous architectures by filtering the input with a particular linear filter (e.g. a gamma filter is a stable instantiation of an IIR filter with restricted memory); (d) dynamic neural networks include such digital filters inside the network itself (replacing the scalar weights by filtered signals) so they can cope with more complex dynamics internally; (d) Elman and Jordan networks include the so-called `contextual neurons' that process and feedback the internal dynamics in the hidden neurons (Elman) or the output of the network (Elman); and (e) long short-term memory (LSTM) units are a complicated recurrent memory composed of a cell state $c_{t-1}\rightarrow c_t$, an input gate takes $h_{t-1}$, and output gate yields $h_t$ and a forget gate that can cope with sophisticated dynamic and nonstationary processes.   }
    \label{fig:nnets}
\end{figure}

The main shortcomings of ARMA models (and hence their treatment of memory, system's dynamics and system's persistence) arise from the non-linearity, uniform sampling, and stationarity assumptions. Hence, many researchers have turned to the use of neural networks, in which few assumptions are typically made and can cope with these problems well. The multilayer perceptron (MLP) is the most commonly used artificial neural network, which is composed of a layered arrangement of artificial neurons where each neuron of a given layer feeds all the neurons of the next layer. This model forms a complex mapping from the $n$-dimensional input to the output $\psi: {\mathbb R}^n \to {\mathbb R}$. Despite its great flexibility, it is a static mapping; one has to pre-define the time embedding beforehand (so-called tapped-delay line) and, after all, there are no internal model's dynamics \cite{bk:weigend94,bk:haykin99}. Figure~\ref{fig:nnets} (a) depicts this primitive form of neural network with memory of the past.

Neural network ARMAX (NNARMAX) modeling is intimately related to the latter approach. In a NNARMAX model \cite{ar:norgaard}, given a pair of input-output discrete time series, an MLP is used to perform a mapping between them, in which past inputs, past outputs and past residuals can be used as inputs. Selecting a model structure is much more difficult in the nonlinear case than in the linear case (classical ARMA modeling). Not only is it necessary to choose a set of regressors but also a network architecture, and several choices are available. For instance, in NNARMAX2 the regression vector is formed by past inputs, past outputs and past residuals, in Neural Network in State-Space Innovation Form (NNSSIF) models the regressor is in the form of state space innovations, while in Neural Network Output Error (NNOE) models the regression vector is formed by past inputs and previous estimates \cite{ar:norgaard}, see Figure~\ref{fig:nnets} (b). The use of NNARMAX models in control applications and nonlinear system identification has expanded in the last decades as its main advantage is the use of a non-linear regressor (usually an MLP) working on a fully tailored ``state'' vector \cite{bk:ljung99}. This makes the model particularly well-suited in problems where one can design the endowed input state vector to accommodate non-stationary system's dynamics by adding for example more ``memory'' in the form of error terms when the prediction error is not substantially decreased. 

\subsubsubsection{Synapses as digital filters} 

The use of MLP as regressor on top of a handcrafted tapped-line (or memory unit) is a naive implementation, and after all the regressor does not represent internal dynamics. In order to account for memory {\em inside} the neural network, and hence endorse the modeling with dynamic capabilities, one has to redesign the concept and operation of individual neurons. A standard approach is to substitute the static synaptic weights for dynamic connections, which are usually linear filters. The FIR neural network models each synapsis as a Finite Input Response (FIR) filter \cite{th:wan93},  Figure~\ref{fig:nnets} (d). There are striking similarities between this model and the MLP. Notationally, scalars are replaced by vectors and multiplications by vector multiplications. These simple analogies carry through when comparing standard backpropagation for static networks with {\em temporal backpropagation} for FIR networks \cite{bk:weigend94}. FIR neural networks are appropriate to work in non-stationary environments or in the presence of non-linear dynamics. This is because {\em time} is treated naturally in the synapsis itself. In fact, they have demonstrated good results in problems with those characteristics, such as speech enhancement \cite{ar:waibel} or time series prediction \cite{ar:wan93_model}, among others.

Extension of the FIR network have also considered including $\gamma$-filters, a class of Infinite Input Response (IIR) filters with restricted stability, which give raise to the gamma network \cite{ar:devries2}, Figure~\ref{fig:nnets} (d). In this structure, the FIR synapsis that uses the standard $z$-Transform delay operator $B=z^{-1}$ is replaced by the gamma operator $G(z)$:
\begin{equation}
G(z)  = \frac{\mu}{z-(1-\mu)},
\end{equation}
where $\mu$ is a real parameter which controls the memory depth of the filter. As pointed out in \cite{Principe93}, gamma filters are theoretically superior to standard FIR filters in terms of number of parameters required to model a given dynamics. The filter is stable if $0<\mu<2$, and $G(z)$ reduces to the usual delay operator for $\mu=1$. This filter also provides an additional advantage: the number of degrees of freedom (order $K$) and the {\em memory depth} remain decoupled \cite{Principe93}. A proposed measurement of the memory depth of a model, which allows quantifying the past information retained, is given by $\frac{K}{\mu}$ and has units of time. Hence, values of $\mu$ lower than the unit increase the memory depth of the system. The gamma structure can be used replacing each scalar weight in an MLP with a gamma filter, or simply used as the first layer of a classical MLP, which yields the so-called {\em focused gamma network}, Figure~\ref{fig:nnets} (c). In general, the gamma network can deal efficiently with complex dynamics and a low number of network parameters.

\subsubsubsection{Recurrent neural networks}

Inspired by the hypothesis of multi-temporal scale processing in the brain \cite{hochreiter1997long}, its reinforced persistence \cite{hebb1949organization} and its plasticity \cite{takeuchi2014synaptic}, the concept of {\em recurrency} has been explored to design new ANNs. The general idea of recurrent ANNs is to construct loops in the connections between neurons or layers of the network, with the aim to increase adaptation and flexibility of the networks to cope with nonstationarities, recurrency and more complex dynamics. The Elman's recurrent network is a simple recurrent model with feedback connections around the hidden layer. In this architecture, in addition to the input, hidden and output units, there are also context units, which are only used to memorize the previous activations of the hidden units \cite{ar:elman88}, Figure~\ref{fig:nnets}(d). The application of recurrent neural networks has traditionally been linked to applications such as speech and language processing. Additionally, Elman networks can result in efficient models to both detect or generate time-varying patterns \cite{ar:elman88}.

The inclusion of feedback loops inside a neural architecture has given rise to an entire family of recurrent neural networks beyond Elman's and Jordan's architectures, cf. Fig.~\ref{fig:nnets}(d). These recurrent connections allow the network developing a set of internal dynamics over time that enable the nets learning to exploit the temporal context of individual data samples. Since recurrent neural network architectures contain feedback loops, gradient-based parameter optimization based on back-propagation through the particular layers is neither applicable nor sufficient right away because their activations and outputs depend on data from different time steps. A solution reformulates the recurrent net as a feed-forward net with each layer representing another time step. Then, applying back-propagation on this ``dual architecture'' gives rise to the so-called \emph{back-propagation through time}  (BPTT) algorithm~\cite{Rumelhart86:LRB,Werbos90:BTW}.
Recurrent neural network models are able to capture the dynamics of time series driven by complex latent dynamical systems. The active steering of gradient flows while training increases the representative power of these models which can yield more stable predictions over longer periods. 

Conventional recurrent nets, however, present instabilities and problems during the training phase because back-propagated gradients tend to fade over time, which produces difficulties with learning long-term dependencies\footnote{This is an important flaw of recurrent networks due to the fact that with the increasing of time steps, the gradients may vanish due to the cumulative multiplication of decimal numbers in the activation function. This results in a virtual collapse or explosion of the weight updates and thus to stopping actual learning.}. Long short-term memory (LSTM) networks are a special recurrent hidden unit that was proposed to deal with the {\em vanishing gradient problem} in recurrent networks and learn long-term dependencies \cite{hochreiter1997long}, Figure \ref{fig:nnets}(e). LSTM networks incorporate a series of steps to decide which information is going to be stored (memorized), and which deleted (forgotten). Thus, the network has a certain internal memory. 

In a non-causal setting, memory provides information on both past and future states, and it is convenient to learn from the whole time series. LSTM units can be combined to obtain a bi-directional long short-term memory (BiLSTM) network. The BiLSTM networks are formed by two LSTM units per time step, and take into account not only past temporal dependencies but also information of future time states. Hence, BiLSTM networks learn from the complete time series at each time step thus having a global view of the sequences \cite{schuster1997bidirectional}.

\subsubsection{Advanced memory-based neural networks}

Despite their excellent capabilities to deal with complex time series and dynamic processes, recurrent nets also have a memory problem: they cannot explore and remember long-range dependencies of, for instance, {\em words} in a natural language processing problem or {\em temperatures} in distant times or regions in a geoscientific/climate problem. This is because recurrent networks put too much emphasis on instances being close to one another, so a limited {\em context} (nearby instances in space or time) is considered. The concept of {\em attention} can address this problem. Other approaches have considered including explicit memory scales and mechanisms to augment memory capacities. In the following we summarize recent developments, and emphasize its, again, neuroscientific inspiration. 

\subsubsubsection{Transformers and reformers}

Even when space and time are considered altogether, the vicinity where (recurrent) neural models look at is quite limited since convolutions and memory units are, by construction, FIR filters.
On the one hand, recurrent and LSTM networks cannot remember longer sentences (sequences), mainly due to the vanishing/exploding gradient problem. Several architectures currently try to tackle the myopic problem of neural networks. 
A promising architecture combines standard convolutional networks with human attention mechanisms in neurosciences: the so-called {\em attention deep learning networks}~\cite{vaswani2017attention} can look at the information in distant regions driving processes and activations in a local window. Attention networks have excelled in classification problems where salient patterns of a target are fired by activity in distant regions~\cite{wang2018scene}. 
The idea has been also developed in the field of natural language processing with networks called {\em transformers}~\cite{jaderberg2015spatial}, which efficiently tackle the problem of sequence transduction with self-attention and focused recurrences. Yet, transformers have been limited to look at limited ranges too, but recently the so-called {\em reformers} \cite{kitaev2020reformer}\footnote{https://ai.googleblog.com/2020/01/reformer-efficient-transformer.html} promise to break the long-range dependency problem.

\subsubsubsection{Deep persistent memory network}

The problem of long-range interactions is even stronger when deep networks are used, and actually it is seldom addressed in very deep neural network models. This makes that early layers do not actually impact on the learning of the subsequent ones, so at a certain extent learning complex dynamic processes is actually compromised. The MemNet network presented in \cite{tai2017memnet}
is motivated by the concept of {\em persistency}, which introduces a memory block similar to LSTM (consisting of a recursive unit and a gate unit) to explicitly exploit persistent memory adaptively. The recursive unit learns multi-level representations of the current state under different scales (receptive fields). The representations and the outputs from the previous memory blocks are concatenated and sent to the gate unit, which adaptively controls how much of the previous states should be reserved, and decides how much of the current state should be stored. The MemNet architecture actually exploits the  gating mechanism to tackle the long-term dependency problem, and resolves some of the problems in convolutional architectures. Even though its use was exemplified for image processing, the network could be eventually applied to other time-dependent data problems too.

\subsubsubsection{Augmented memory networks}

Neural networks with a memory capacity provide a promising approach to meta-learning in deep networks. However, information must be stored in memory in a representation that should be both stable (accessible when needed),  addressable (accessibility of specific pieces of information), and limited in memory demands (long-range constraints). These desirable properties are not typically present (at least naturally) in standard memory architectures like LSTMs. Actually, accounting for memory in arbitrary deep learning architectures is a great challenge, especially because traditional gradient-based networks require a lot of data to learn and thus models need to relearn parameters continuously. This is not only a computational problem but, more importantly, a problem about adequately incorporating the new information without catastrophic interference with the previously learned representations. Architectures with augmented memory capacities, such as Neural Turing Machines (NTMs) and the memory-augmented neural network (MANN) in \cite{santoro2016meta} can efficiently encode, retrieve and ingest new information. Like in other architectures with biological inspiration, the characteristics of accessibility to memory and flexible adaptation to changing environments is a key aspect of human learning \cite{jankowski2011meta}. Such memory-augmented neural networks rapidly assimilate new data, and leverage it to make accurate predictions after observing a few instances only.

\subsubsection{Short-term persistence in non-stationary data streams}\label{Pers_ML}

From a broader perspective, the concept of short-term persistence also emerges in modeling tasks over non-stationary data streams, disregarding which predictive model is in use (i.e. beyond the persistence specifics of neural networks discussed in the previous section). 

In its simplest formulation, the goal of predictive modeling is to build a model that captures the relationship between a set of input features and a target variable. Following the assumptions and notation of preceding sections, in our case the variable of interest $x(t,\boldsymbol{\theta})$ is the target variable to be predicted, whereas $\boldsymbol{\theta}$ denotes the set of input features to the model. The predictive model $M_{\bm{\vartheta}}$ results from the application of a learning algorithm (the ANNs models described above or any other) over a set of supervised data instances $\{(\boldsymbol{\theta}_n^{tr},x(t,\boldsymbol{\theta}_n^{tr}))\}_{n=1}^{Ntr}$, such that the model, once trained over these data instances, provides a prediction $\widehat{x}(t,\boldsymbol{\theta}^{test})=M_{\bm{\vartheta}}(\boldsymbol{\theta}^{test})$ for any given test instance $\boldsymbol{\theta}^{test}$. It is important to note that the so-called model hyper-parameters $\bm{\vartheta}$ drive the generalization properties of the model itself, i.e. the capability of the model to provide accurate predictions for \emph{any} test instance, disregarding whether it was present in the training set or not. Examples of predictive modeling under this supervised learning framework abound nowadays in the literature related to physics, as extensively shown by recent overviews on this matter \cite{carleo2019machine}. Supervised learning can serve as a modeling framework for many classical data-based paradigms, including time series forecasting or curve fitting, to mention a few.

Bearing this framework in mind, it is often acknowledged that one of the core assumptions of predictive modeling is the stationary nature of the pattern to be learned. In mathematical terms, the notion of persistence in predictive modeling refers to the temporal stability of the \emph{concept}, namely, the joint statistical distribution $P_{\boldsymbol{\Theta},X_t}(\boldsymbol{\theta},x_t)$ between the feature variable vector $\boldsymbol{\Theta}$ and the target variable $x_t\equiv x(t,\boldsymbol{\theta})$. Indeed, the ultimate purpose of the predictive model $M_{\bm{\vartheta}}$ is to characterize this concept based on a set of training instances, a model structure (e.g. a network of interconnected neurons, or a tree structure), and a learning algorithm that adjusts the parameters of the model as per a measure of generalization performance. When conceived in a batch data generation scenario in which no data instances are produced over time, the concept of short-term persistence refers to any distribution shift between training and test sets \cite{datasetshift}. However, it is often the case that data instances are produced and collected over time (data streams), which translates into a time-dependent definition of concept given by $P_{\boldsymbol{\Theta},X_t}^{(t)}(\boldsymbol{\theta},x_t)$. In this case, the stability of the concept is what imposes the need for detecting and adapting to eventual changes of the distribution to be learned ({\em drift}), what is known in the literature as {\em concept drift detection and adaptation} \cite{concept1,concept2}. Although we hereafter focus on predictive modeling for the sake of simplicity in the explanations, it is worth mentioning that concept drift also arises in unsupervised learning tasks (e.g. clustering), in which a drift occurs when $P_{\boldsymbol{\Theta}}^{(t)}(\boldsymbol{\theta})\neq P_{\boldsymbol{\Theta}}^{(t')}(\boldsymbol{\theta})$ for $t'>t$.

Contributions elaborating on the detection and adaptation to concept drifts are manifold, specially in what refers to the occurrence of this phenomena in streaming data. In general, concept drift can be characterized in terms of their magnitude (intensity) and duration, which determines the effectiveness of algorithmic adaptation designed to minimize the impact of the drift on the performance of the learned model. Adaptation strategies can be i) passive, e.g. by inducing a controlled amount of diversity in the learning process that preemptively prepares the model for eventual data distribution changes; or ii) active, correspondingly, by detecting drifts from the evolving data and, when they occur, by activating an adaptation mechanism to make the model forget its captured knowledge inasmuch as required by the characteristics of the drift. As such, sudden and intense concept drifts are in general easier to detect than smooth concept transitions. However, the former is often tackled by forgetting mechanisms, whose effect on the model's knowledge should match the intensity of the drift itself. Despite recent efforts in this direction \cite{charact1}, concept drift characterization still remains a subject of active research. This fact, jointly with the varying performance of drift detectors published to date, make passive concept drift adaptation strategies a mainstream in this research area. A special mention is to be paid to ensemble learning, which makes it possible to derive simple albeit effective methods for the insertion of controlled diversity among its constituent learners \cite{ensemble}.

\subsection{Examples on short-term persistence characterization}

\subsubsection{Temporal evolution of state probabilities}
In some complex systems, very often considered in non-equilibrium thermodynamics for example \cite{bray2013persistence,iyer2015first}, it is possible to calculate the temporal evolution of the state probabilities of the system $s(t)$. Supposing a continuum-time domain, let us consider a system with a set $\mathcal{S}=\lbrace s_i\rbrace_{i=1}^{N_s}$ of possible states, where $N_s=\left\lvert\mathcal{S}\right\rvert$ is the number of different states.
The dynamics of the system $s(t)$ are characterized by these states, and the transition among them ($s_i \rightarrow s_j,\ \forall s_i,s_j \in \mathcal{S}$). The probability of transitioning from state $s_i$ to state $s_j$ in an infinitesimal interval, $dt$, is known to be  $\alpha(s_i,s_j)dt$, where $\alpha(s_i,s_j)$ stands for the transition rate. The equations that govern the dynamical evolution of the probability
that the system is in state $s_i$ at time $t$ ($P(s(t)=s_i)$) can be described in different ways. We will summarize here the most general formalism for this, known as the Master Equation \cite{iyer2015first}.

The probability of being in state $s_i$ at a time $t+dt$, $P(s(t+dt)=s_i)$, is the sum of two different terms:
first, the probability, $P(s(t)=s_i)$, that the system was in the state $s_i$ at time $t$, and remained there during $dt$. Second, the probability that the system was originally in some other state $s_j$ at time $t$, times the probability that the system transitioned from $s_j$ to $s_i$ during $dt$. This yields:
\begin{equation}
    P(s(t+dt)=s_i) = P(s(t)=s_i)\left[1-\sum_{s_i\neq s_j}\alpha(s_j,s_i)dt\right]+\sum_{s_j\neq s_i}P(s(t)=s_j)\alpha(s_j,s_i)dt,
\end{equation}

From this equation we can obtain the Forward Master Equation (FME):
\begin{equation}\label{FME}
\partial_t P(s(t)=s_i) = \sum_{s_j\neq s_i}\alpha(s_j,s_i)P(s(t)=s_j) - \sum_{s_j\neq s_i}\alpha(s_i,s_j)P(s(t)=s_i)
\end{equation}
note that the first term in this equation stands for the probability flux into the state $S$, whereas the second term stands for the probability flux out of $S$. The FME can be compactly written as:
\begin{equation}
\partial_t {\bf P}(t)=\mathcal{M}_f {\bf P}(t)
\end{equation}
where ${\bf P}(t)$ is a vector with the components $P(s(t)=s_i)$, and $\mathcal{M}_f$ is a linear operator with the following form:
\begin{equation}
\left(\mathcal{M}_f\right)_{s_is_j}=\alpha(s_i,s_j)-\delta_{s_is_j}\sum_{S_k}\alpha(s_k,s_i).
\end{equation}
where $\delta_{s_is_j'}$ stands for the Kronecker $\delta$.

Note that the FME describes the evolution of the system from an initial state to a state at a later time. The FME can be also obtained as a linear combination of the conditional probabilities $P(s(t)=s_i|x(t_0)=s_0)$, which are solutions to the Master Equation for the special initial conditions, $P(s(t_0)=s_i) = \delta_{s_i,s_0}$ . We can rewrite the FME in terms of these conditional probabilities, in the following way:
\begin{equation}\label{BME}
\begin{aligned}
P(s(t)=&s_i|s(t_0)=s_0)=\left(1-\sum_{s_j} \alpha(s_j,s_0)dt\right) P(s(t)=s_i,t|s(t_0+dt)=s_0)+\\
&+\sum_{s_j}\left(\alpha(s_j,s_i)dt\right) P(x(t)=s_i|x(t_0+dt)=s_j)
\end{aligned}
\end{equation}
Note that this equation is formed by two terms: the first one stands for the probability that the system transitioned to a different state $s_j$ during the time interval $dt$, with the probability
$\alpha(s_j,s_0)dt$, and then evolved to a state $s_i$ by time $t$, with the probability $P(s(t)=s_i|s(t_0+dt)=S_j)$. The second term is the probability that the that the system was still in state $s_0$ at time $t_0+dt$, with the probability $1-\sum_{s_j}\alpha(s_j,s_0)dt$, and then by time $t$ it evolved to the state $s_i$, with the probability $P(s(t)=s_i|s(t_0+dt)=s_j)$.

Following \cite{iyer2015first}, the lack of explicit dependence of the transition rates on $t$ implies that the conditional probability $P(x(t)=s_i|x(t_0)=s_0)$ is a function of $t-t_0$. Thus, Equation (\ref{BME}) can be rewritten as:
\begin{equation}\label{BME2}
\begin{aligned}
P(x(t)=&s_i|x(t_0)=s_0))=\left(1-\sum_{s_j} \alpha(s_j,s_0)dt\right)P(s(t-dt)=s_i|s(t_0)=&S_0)+\\
&+\sum_{s_j} \left(\alpha(s_j,s_0)dt\right) P(s(t-dt)=s_i|s(t_0)=s_j)
\end{aligned}
\end{equation}
Taking $dt \rightarrow 0$ in this equation we obtain:
\begin{equation}\label{BME_Final}
\begin{aligned}
\partial_t P(s(t)=s_i|s(t_0)=s_0) = &\sum_{s_j}\alpha(s_j,s_0)P(s(t)=s_i|s(t_0)=s_j) \\
& - P(s(t)=s_i|s(t_0)=s_0)\sum_{s_j}\alpha(s_j,s_0) \\ \equiv & \sum_{s_j} (\mathcal{M}_b)_{s_0,s_j}P(s(t)=s_i|s(t_0)=s_j).
\end{aligned}
\end{equation}
This equation is known as the Backward Master Equation (BME). The operator $(\mathcal{M}_b)_{s_i,s_j}=\alpha(s_j,s_i)-\delta_{s_i,s_j}\sum_{s_j} \alpha(s_j,s_i)$ is known as the backward master operator. The BME can be also written in compacted form as:
\begin{equation}
\partial {\bf P}^T(t)=\mathcal{M}_b {\bf P}^T(t),
\end{equation}
where ${\bf P}^T(t)$ is a vector whose $i$-th component is $P(x(t)=s_i|s_0)$, i.e.:
\begin{equation}
 {\bf P}^T(t)=(\cdots, P(s(t)=s_i|s_1), P(s(t)=s_i|s_2), \ldots, P(s(t)=s_i|s_l), \cdots).
\end{equation}

The FME and BME are of direct application to different problems related to persistence, such as the {\em First Passage Process} \cite{bray2013persistence,iyer2015first,aurzada2015persistence}, and they can be also applied to direct calculation of $\mathcal{P}^P$ and $\mathcal{P}^M$ in some problems of interest in non-equilibrium problems. In these cases, it is known that $\mathcal{P}^P \sim t^{-\theta}$ in many cases \cite{majumdar1999persistence}.  
More details on further formalisms to define persistence in non-equilibrium problems can be found in \cite{iyer2015first}. A discussion on persistence studies for non-equilibrium thermodynamics systems is carried out in Section \ref{PnonEq}.

\subsubsection{Persistent system with Gaussian noise}

Consider a very simple persistent system in which the residuals are Gaussian (one dimensional Brownian walker \cite{majumdar1999persistence,sire2000analytical}, in this case we consider a discrete problem with state change characterized by a threshold $\epsilon$). Thus, let us suppose a persistent system in which $x[n+1]=x[n]+N_n,\ N_n\sim N(\mu,\sigma)$, i.e., the residual $r[n]$ is a Gaussian process with mean $\mu$ and variance $\sigma$: $r[n]=N_n\sim N(\mu,\sigma)~\forall n$. In this particular case of persistent system, given a time step $T\in \mathbb{N}$:
\begin{equation}
x[n+T]=x[n]+\sum_{i=1}^T N_i,\quad N_i\sim N(\mu,\sigma)\ \forall i,
\end{equation}
for $T+1$:
\begin{equation}
    x[n+T+1]=x[n]+\sum_{i=1}^{T+1} N_i,\quad N_i\sim N(\mu,\sigma)\ \forall i,
\end{equation}
and the residual:
\begin{equation}
    r[T]=\left|x[n+T+1]-x[n+T]\right|=\left|N_{(T+1)}\right|,\quad N_{(T+1)}~\sim N(\mu,\sigma),
\end{equation}
which we compares with the threshold $\epsilon$ for changing to another state $r[t]>\epsilon$.

Then, the probability of changing to another state  can be computed as:
\begin{equation}
P\left(r[T]>\epsilon \right)=\int_{\epsilon}^{\infty} \frac{1}{\sqrt{2 \pi} \sigma} e^{-\frac{(r-\mu)^2}{2 \sigma^2}} dr+\int_{-\infty}^{-\epsilon} \frac{1}{\sqrt{2 \pi} \sigma} e^{-\frac{(r-\mu)^2}{2 \sigma^2}} dr=2 Q\left(\frac{\epsilon-\mu}{\sigma} \right)
\end{equation}
where $Q(\cdot)$ stands for the $Q$-function (recall that
$Q(x)=\frac{1}{2}\operatorname{erfc}(\frac{x}{\sqrt{x}})$). It is easy to see that therefore:
\begin{equation}
\mathcal{P}^P(T)=2Q\left(\frac{\epsilon-\mu}{\sigma} \right)\left(1-2Q\left(\frac{\epsilon-\mu}{\sigma} \right)\right)^{T-1}
\end{equation}

\begin{figure}[!ht]
\begin{center}
\subfigure[~$\mathcal{P}^P$ ($\epsilon=1$)]{\includegraphics[draft=false, angle=0,width=8cm]{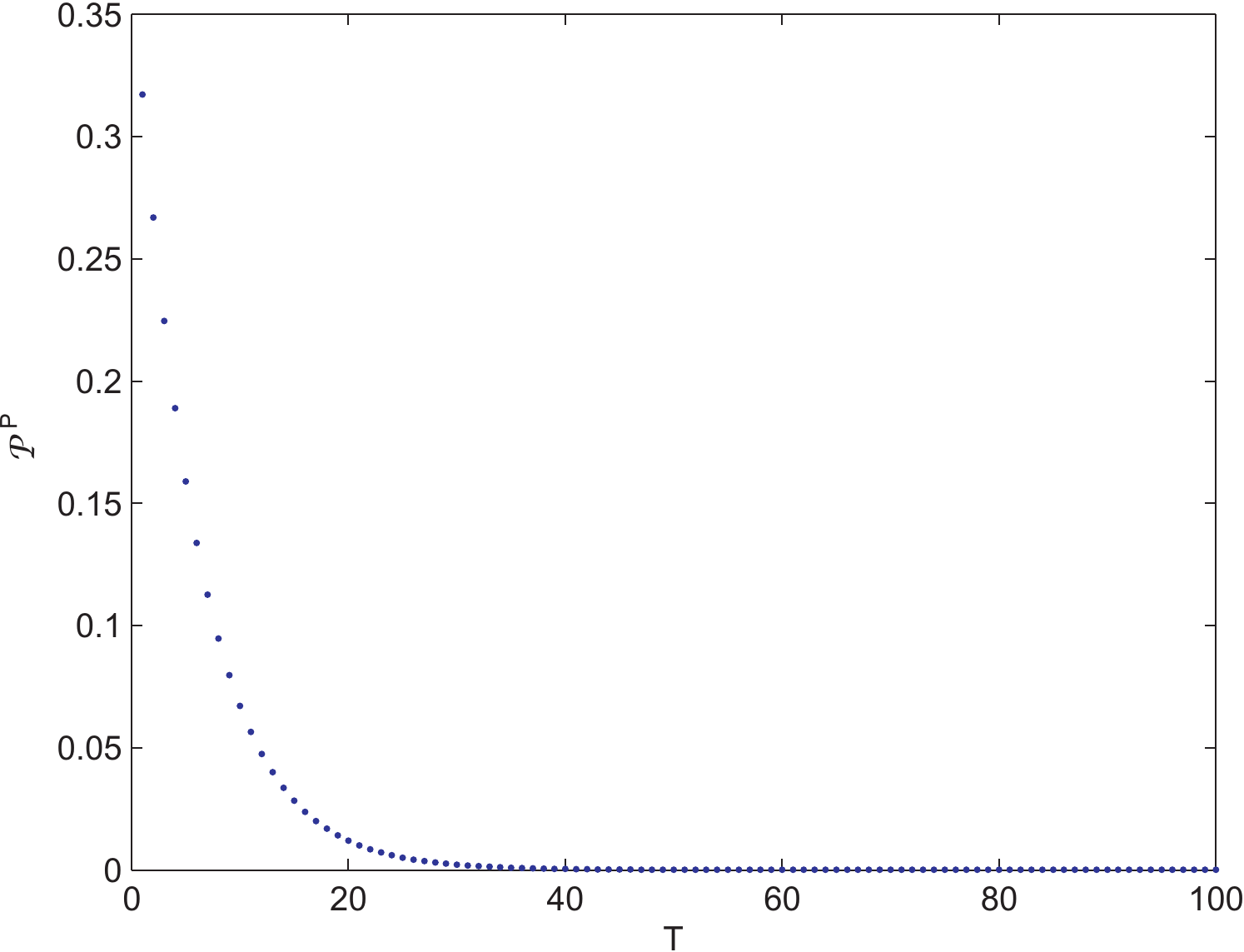}}
\subfigure[~$\mathcal{P}^P$ ($\epsilon=3$)]{\includegraphics[draft=false, angle=0,width=8.2cm]{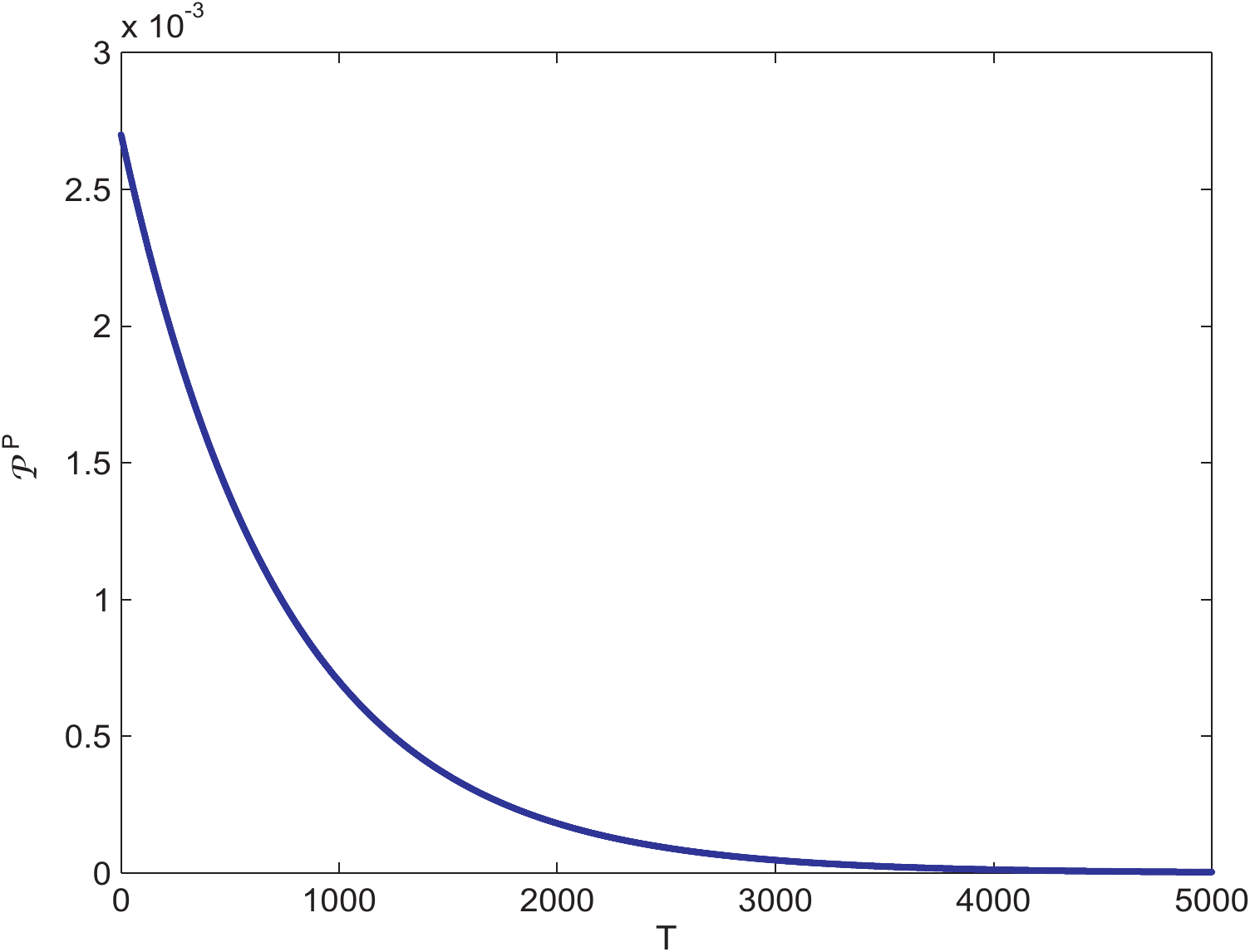}}
\subfigure[~$\mathcal{P}^P$ ($\epsilon=4$)]{\includegraphics[draft=false, angle=0,width=8.2cm]{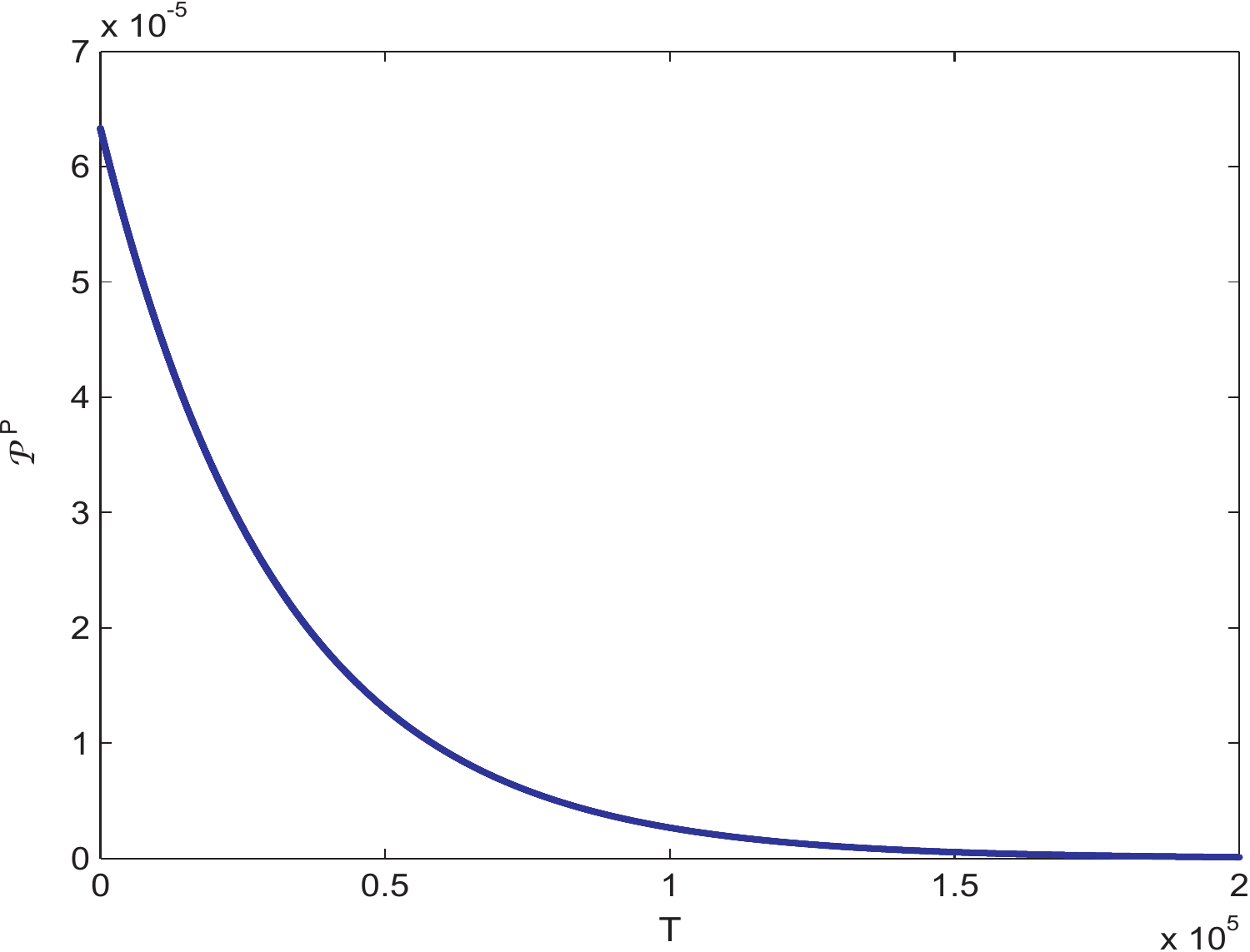}}
\end{center}
\caption{
Persistence estimation ($\mathcal{P}^P$) for a system $x[n+1]=x[n]+N_n,\ N_n \sim N(0,1)$ as a function of $T$; (a) $\epsilon=1$; (b) $\epsilon=3$; (a) $\epsilon=4$.}
\label{Persistencia_PP} 
\end{figure}

Figure \ref{Persistencia_PP} shows an example of $\mathcal{P}^P(T)$ for the system $x[n+1]=x[n]+N_n,\ N_n\sim N(0,1)$, fixing thresholds at $\epsilon=1$, $\epsilon=3$ and $\epsilon=4$. As can be seen, the system with $\epsilon=1$ is not highly persistent, showing negligible probabilities of reaching values of $T>30$. When the threshold is larger, the system is of course more persistent, reaching much larger values of $T$. Note that, in general, it is not a PPS but has a high persistence degree depending on the threshold $\epsilon$. Note that the system becomes a PPS for large values of the threshold $\epsilon$.

This can be better analyzed by calculating the value of $\mathcal{P}^M$ as:
\begin{equation}
    \mathcal{P}^M=\left(1-2 Q\left(\frac{\epsilon-\mu}{\sigma} \right)\right)
\end{equation}

\begin{figure}[!ht]
    \begin{center}
    \subfigure[~$\mathcal{P}^E$ ($\epsilon=1$)]{\includegraphics[draft=false, angle=0,width=8cm]{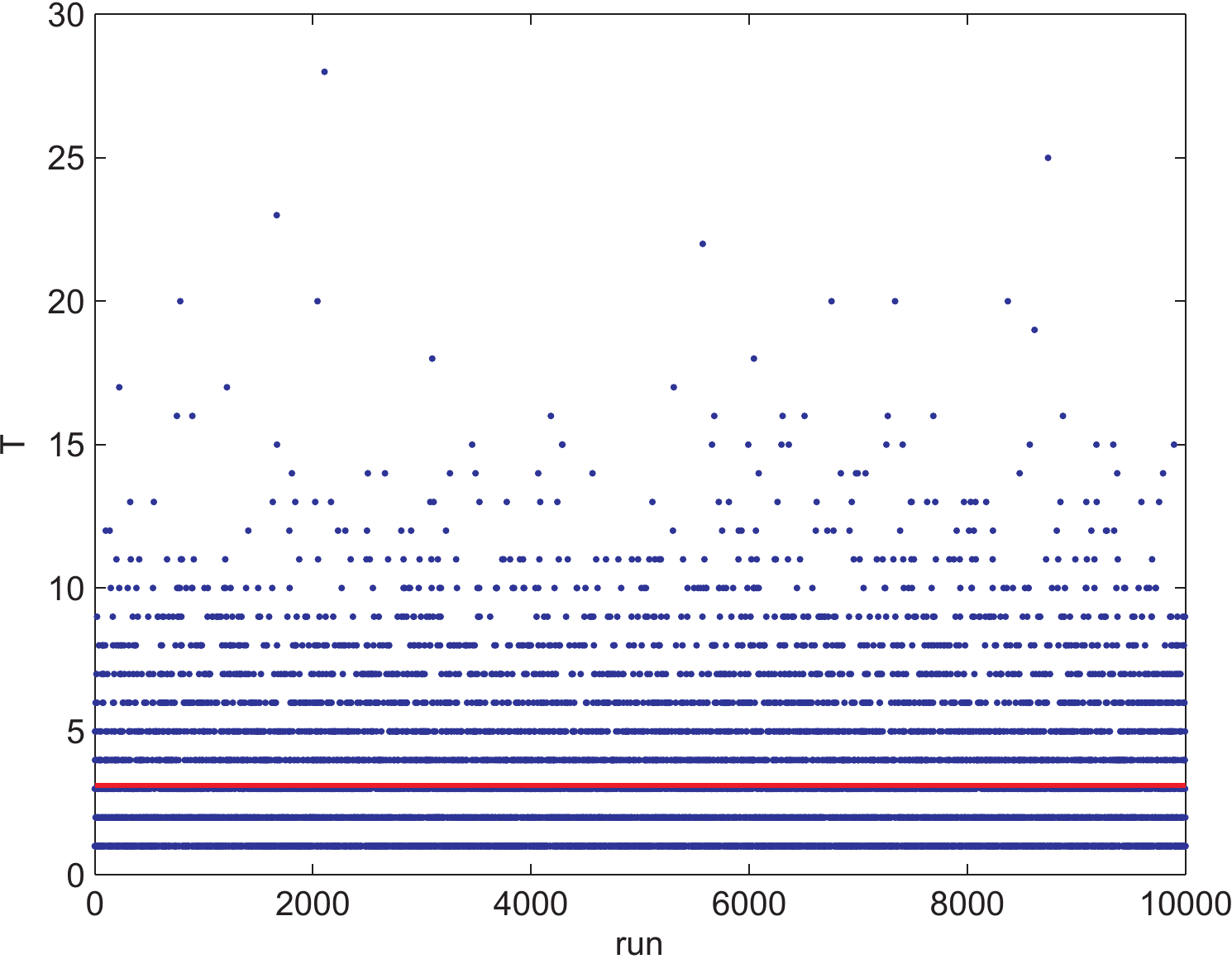}}
    \subfigure[~$\mathcal{P}^E$ ($\epsilon=3$)]{\includegraphics[draft=false, angle=0,width=8cm]{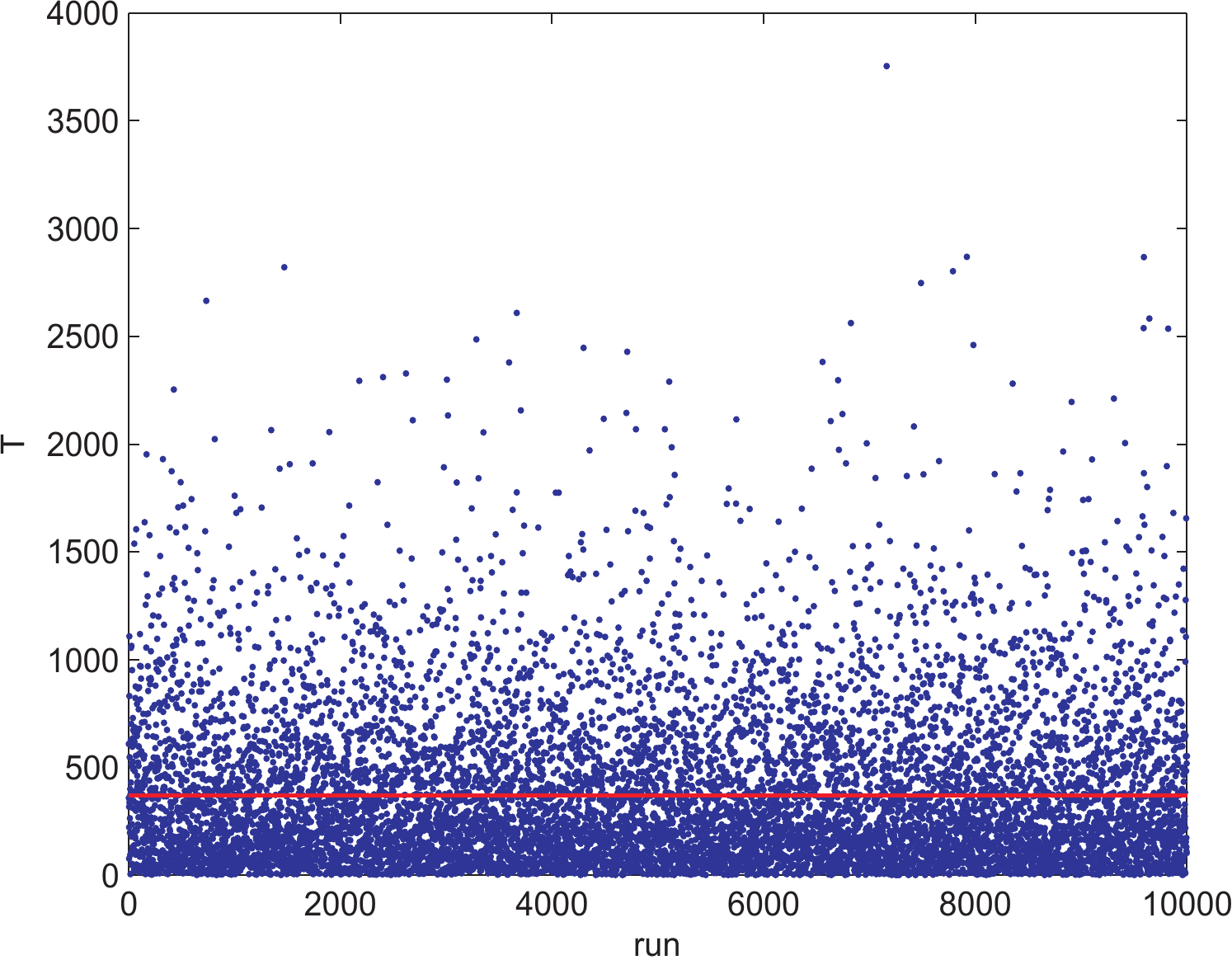}}
    \subfigure[~$\mathcal{P}^E$ ($\epsilon=4$)]{\includegraphics[draft=false, angle=0,width=8.2cm]{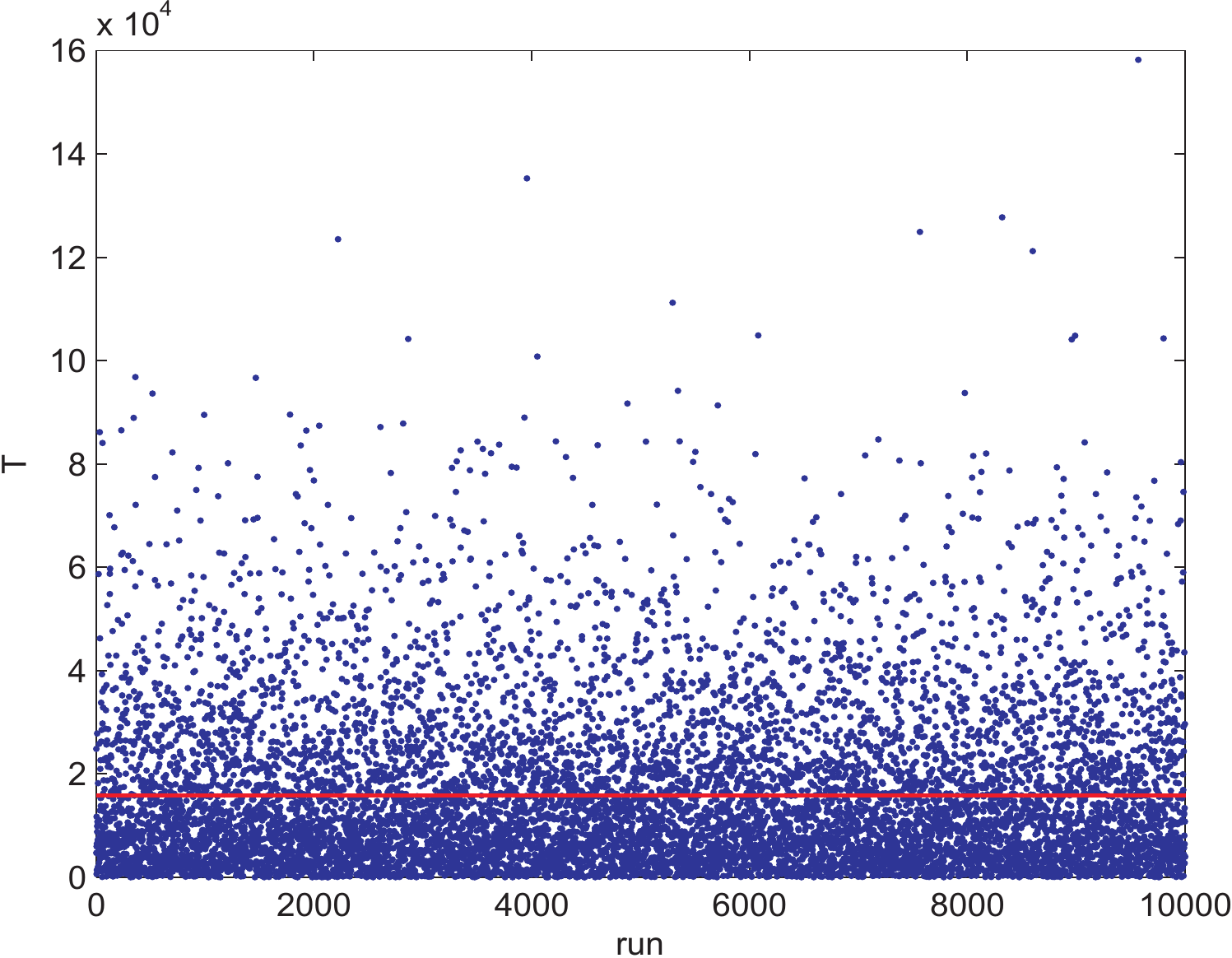}}
    \end{center}
    \caption{
    Persistence estimation ($\mathcal{P}^E=\operatorname{E}[T]$) for a system $x[n+1]=x[n]+N_n,\ N_n\sim N(0,1)$; (a) $\epsilon=1$; (b) $\epsilon=3$; (a) $\epsilon=4$. The estimation of $\mathcal{P}^E$ is depicted as a red line in the plots.}
    \label{Persistencia_E(T)}
\end{figure}

Note that, for the specific case $x[n+1]=x[n]+N_n,\ N_n\sim N(0,1)$, $\mathcal{P}^M=0.6828$ for $\epsilon=1$, $\mathcal{P}^M=0.9973$ for $\epsilon=3$ and $\mathcal{P}^M=0.999936$ for $\epsilon=4$ (system with a high degree of persistence).
The estimation of $\mathcal{P}^E$ can be easily obtained by Monte-Carlo simulation. For example, Figure \ref{Persistencia_E(T)} shows the estimation of $\mathcal{P}^E=\operatorname{E}[t]$ when adding Gaussian Noise $N(0,1)$ for $\epsilon=3$ and $\epsilon=4$, in a Monte-Carlo simulation of $10^4$ epochs. In both cases the estimation of $\mathcal{P}^E$ is depicted as a red line in the plots.

\subsubsection{System with accumulative shocks and damages}

In \cite{Batabyal03} an analysis based on shocks, specific for ecological and economic systems, is carried out. This analysis is extensible to other types of systems. In this analysis, it is assumed that the system changes its state due to the accumulative effect of a number of shocks in the system. In fact, let us assume that the system will remain in a given state until a number $k$ of shocks affect the system. We consider that these shocks follow a Poisson arrival process with rate $\lambda$. In our case, the system will change to a new state when the $k$-th shock hits this system. If we assume that the times between successive shocks are exponentially distributed with mean $\frac{1}{\lambda}$, and from \cite{Batabyal03} we know that the sum of these exponentially distributed times have a Gamma distribution with parameters $k$ and $\lambda$, we can write the density function of the system's lifetime $T$ at a given state $s_i$ as:
\begin{equation}
f_k(t)=\frac{\lambda^k t^{(k-1)} e^{-\lambda t}}{(k-1)!},
\end{equation}
which allows 
the calculation of $\mathcal{P}^E= \operatorname{E}[T]=\frac{k}{\lambda}$.

Following \cite{Batabyal03}, we could elaborate further on this idea, by considering that every shock $i$ arriving at the system produces a {\em damage} $D_i$, measured as the deviation of $r(t)$ from 0. Assuming that the damages from shocks are exponentially distributed random variables, distributed as a density function $g_D(d)=\eta e^{-\eta d}$, with $d>0$, the total damage of $N(t)$ shocks is $X(t)=\sum_{i=1}^{N(t)} D_i$, and the system will change from state $s_i$ to $s_j$ when the total damage arrives the threshold $\epsilon$. The calculation of $\mathcal{P}^E$ in this case is as follows:
\begin{equation}
\mathcal{P}^E= \operatorname{E}[T]=\int_0^\infty  \sum_{n=0}^\infty \frac{(\lambda t)^n e^{-\lambda t} G^n(\epsilon)}{n!} dt
\end{equation}
where $G^n(\epsilon)=P(X(t) \leq \epsilon)=P(D_1+D_2+ \cdots + D_n \leq \epsilon)$. This expression yields:
\begin{equation}\label{refG}
\mathcal{P}^E=\operatorname{E}[T]=\frac{1}{\lambda} \sum_{n=0}^\infty G^n(\epsilon)
\end{equation}
If we further consider that the sum of the exponentially distributed damage random variables has a Gamma distribution, we obtain:
\begin{equation}
\sum_{n=0}^\infty G^n(\epsilon)=\sum_{n=0}^\infty \sum_{k=n}^\infty \frac{(\eta\epsilon)^k e^{(-\eta \epsilon)}}{k!}=1+\eta \epsilon,
\end{equation}
and then, using this expression into Equation (\ref{refG}), we obtain the final expression for the persistence in the case of shocks and their corresponding inflicted damages in the system:
\begin{equation}
\mathcal{P}^E=\operatorname{E}[T]=\frac{1+\eta \epsilon}{\lambda}.
\end{equation}

\subsubsection{Persistence in ARMA processes and the e-folding time} \label{ssec:efolding}

The AR(1) is a simple analysis tool which has been previously used to model persistence in fields such as climate science \cite{Mudelsee10}. The AR(1) model can be seen as a generalization of the na\"ive persistence model, and also of the previously analyzed persistent system with random Gaussian noise.  Let us consider the following AR(1) system:
\begin{equation}\label{AR(1)PGS}
x[n+1]=a_1 x[n] +  \sigma \varepsilon[n]
\end{equation}
where $0<a_1<1$, and $\varepsilon$ is a random noise coming from a zero mean, unit variance Gaussian distribution. This model tends to the persistent system with random Gaussian noise when $a_1  \rightarrow 1$.

A related concept to persistence in AR(1) processes is the so-called {\em e-folding time}. The concept of e-folding time may be used in the analysis of systems kinetics and to study the memory persistence in linear time-invariant systems. Consider a single process A, which decays linearly as:
$$\dfrac{dA}{dt} = -k A,$$
where $k$ has the units of $[t^{-1}]$. This is a first-order linear ordinary differential equation (ODE), whose solution is trivial:

$$\dfrac{dA}{A} = -k dt \to
\int_{A_i}^{A_f} \dfrac{dA}{A}= -k \int_{t=0}^t dt \to
\boxed{\dfrac{A_f}{A_i}=\exp(-kt)},
$$
which tells us that the ratio between the final and initial states of the system follows an exponential function, of which $e$ is the base. So, if one sets the constant to the $k=1/t$, then
$$\dfrac{A_f}{A_i} = \exp(-1) = 0.37,$$
which means that the decay is assumed to be of first order. Many natural processes can be assumed to follow such description, but one has to estimate the $k$ constant from observations, which is typically done via an autoregressive fitting.

\begin{figure}[!ht]
\begin{center}
\includegraphics[draft=false, angle=0,width=9cm]{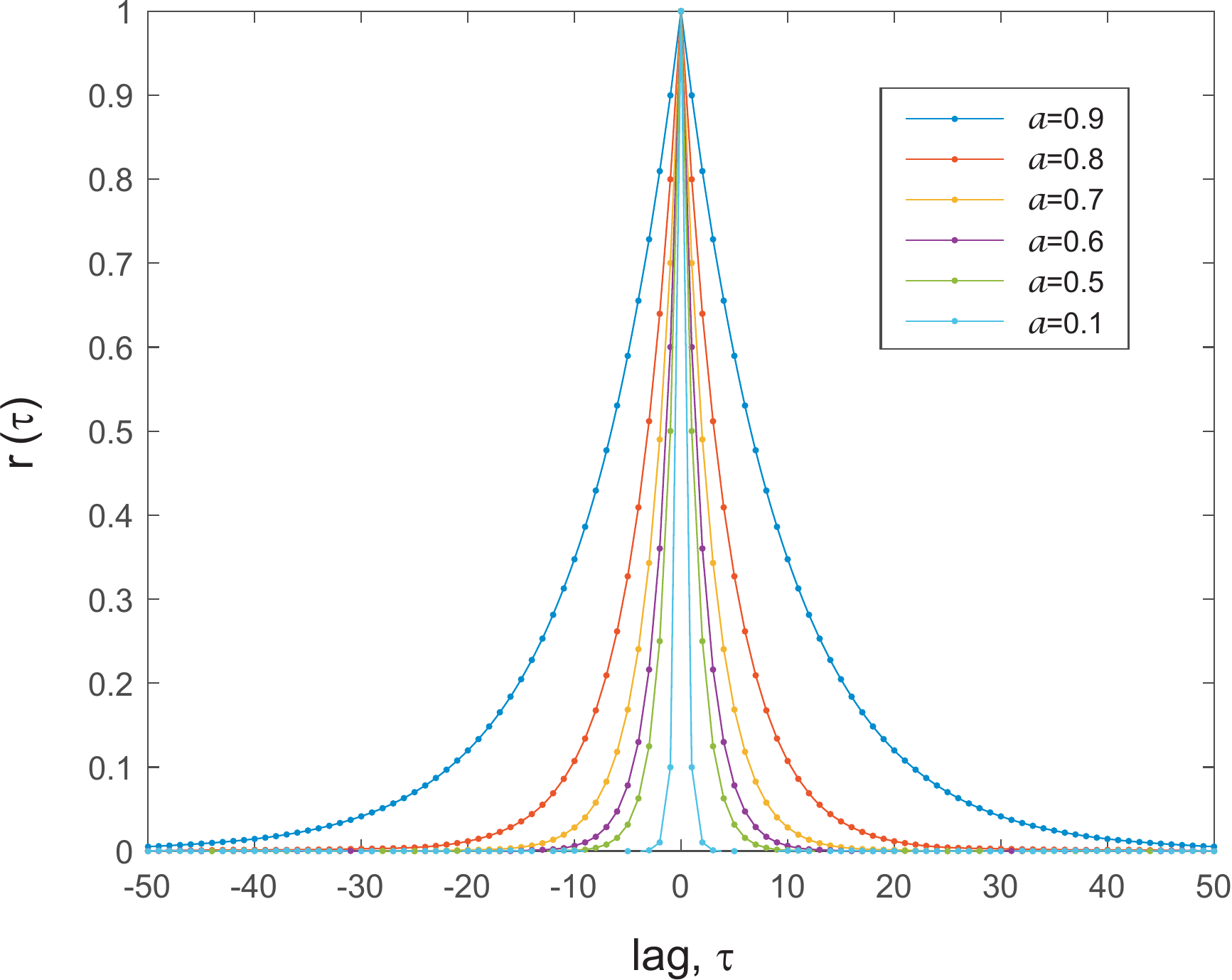}
\end{center}
\caption{Autocorrelation function of the system $x[n] = ax[n-1] + \sqrt{1-a^2}\varepsilon[n],$ for different values of parameter $a$ (that is for different persistence models).}
\label{AR(1)xcorr} 
\end{figure}

Let us assume now a red noise contaminated signal:
\begin{equation}\label{AR(1)Naive}
x(t) = ax(t-\tau) + \sqrt{1-a^2}\varepsilon(t)
\end{equation}
where $0<a<1$ is the degree of memory present from previous states, $\varepsilon$ is a random noise coming from a zero mean, unit variance Gaussian distribution, $\tau$ is the time interval between data points, and let us assume that the signal is standardized first (zero mean, unit std). Note that this is an AR(1) process with $a_1=a$, which tends to the na\"ive persistent system when $a  \rightarrow 1$. Moreover,  parameter $a$ can be estimated from the discrete time series $x[n]$ associated to continuous time series $x(t)$ as follows \cite{Mudelsee10}:
\begin{equation}
\hat{a}=\frac{\sum_{t=2}^N x[n] x[n-1] }{\sum_{t=2}^n x^2[n]}
\end{equation}
It easy follows that the autocorrelation function of the red noise is:
\begin{equation}
    r_{\varepsilon\varepsilon'}(\tau) = \exp\left(-\frac{\tau}{T_e}\right),\quad \text{where } \epsilon'=\epsilon(t-\tau),\ \forall t
\end{equation}
where $T_e$ is called the $e$-folding decay time.\\
Using the autocorrelation function of a wide-sense stationary (WSS) process $x(t)$ which is defined as follows:
\begin{equation}
    r_x(\tau) = \operatorname{E}[x(t)x(t-\tau)],
\end{equation}
directly drives to:
\begin{equation}
\begin{split}
r_x(t,t-\tau) = \operatorname{E}[(ax(t-\tau) + \sqrt{1-a^2}\varepsilon(t))x(t-\tau)] =\\
=  \operatorname{E}[ax(t-\tau)x(t-\tau) + \sqrt{1-a^2}\varepsilon(t)x(t-\tau)] = a
\end{split}
\end{equation}
hence the autocorrelation of $x(t)$ at $\tau$ is $a$.
Following by induction, one can demonstrate that for a red noise AR(1) time series, the autocorrelation at a lag of $n$ time steps is equal to the autocorrelation at one lag, raised to the power $n$. Since a function that has this property is the exponential function, i.e. $({e^{x}})^n = e^{nx}$, one concludes that:
\begin{equation}
r_x(n\tau) = \exp\left(-n\frac{\tau}{T_e}\right).
\end{equation}
We can simulate the autocorrelation of a model similar to that of Equation (\ref{AR(1)Naive}) by means of a Monte-Carlo simulation, obtaining $r_x(\tau)$ for different values of $a$. Figure \ref{AR(1)xcorr} shows the autocorrelation function obtained for different values of the parameter $a$ (different persistence models), with a Monte-Carlo simulation of $10^4$ epochs (we show the mean $r_x(\tau)$ obtained). As can be seen, the smaller the parameter $a$, the faster is the memory decreasing in the system. When $a \rightarrow 0$, $r_x(\tau)=\delta(t)$ (independent samples with no persistence at all in the system).

\begin{figure}[!ht]
\begin{center}
\subfigure[~$\mathcal{P}^E$ ($a=0.9$)]{\includegraphics[draft=false, angle=0,width=8cm]{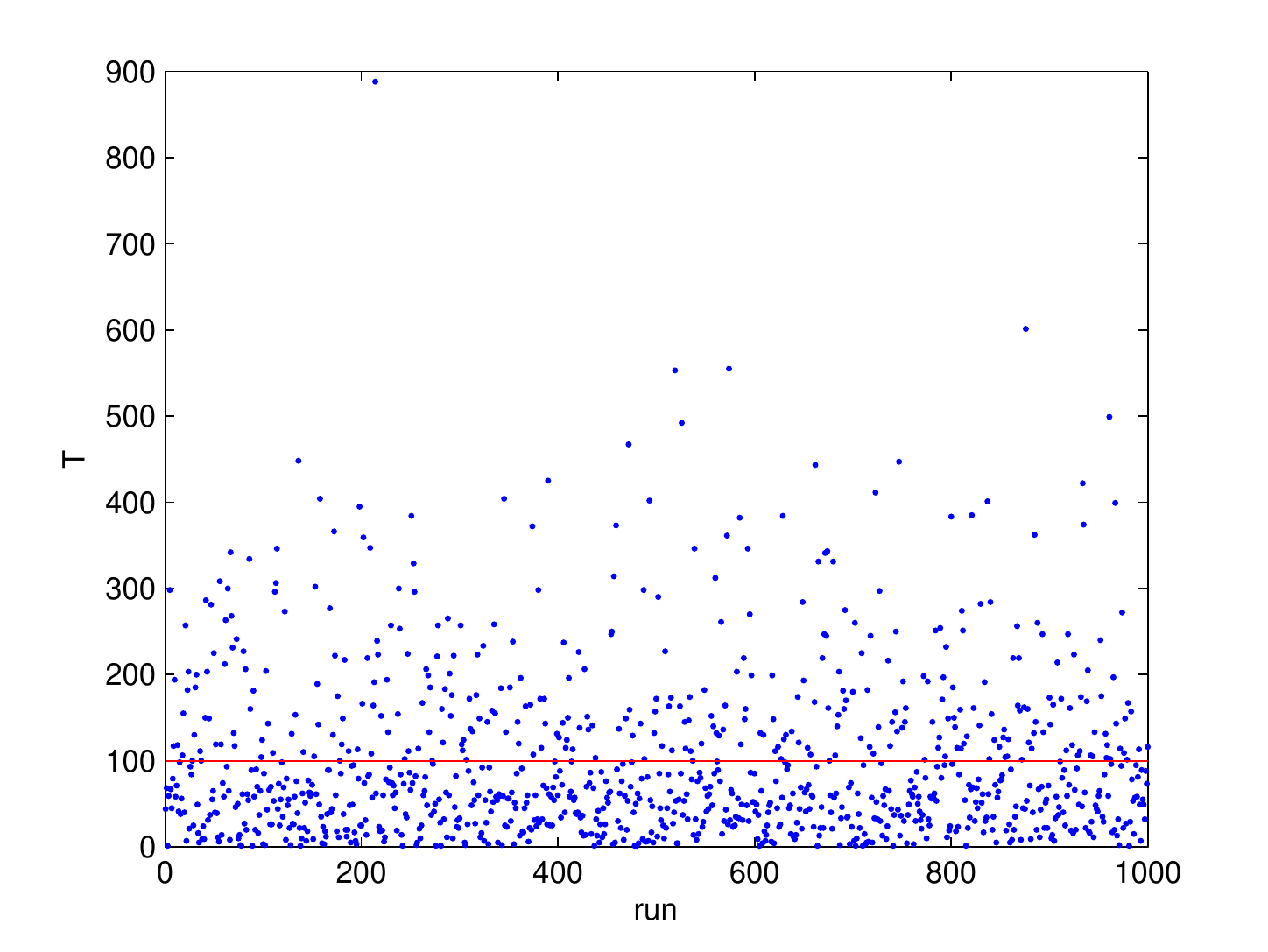}}
\subfigure[~$\mathcal{P}^E$ ($a=0.7$)]{\includegraphics[draft=false, angle=0,width=8cm]{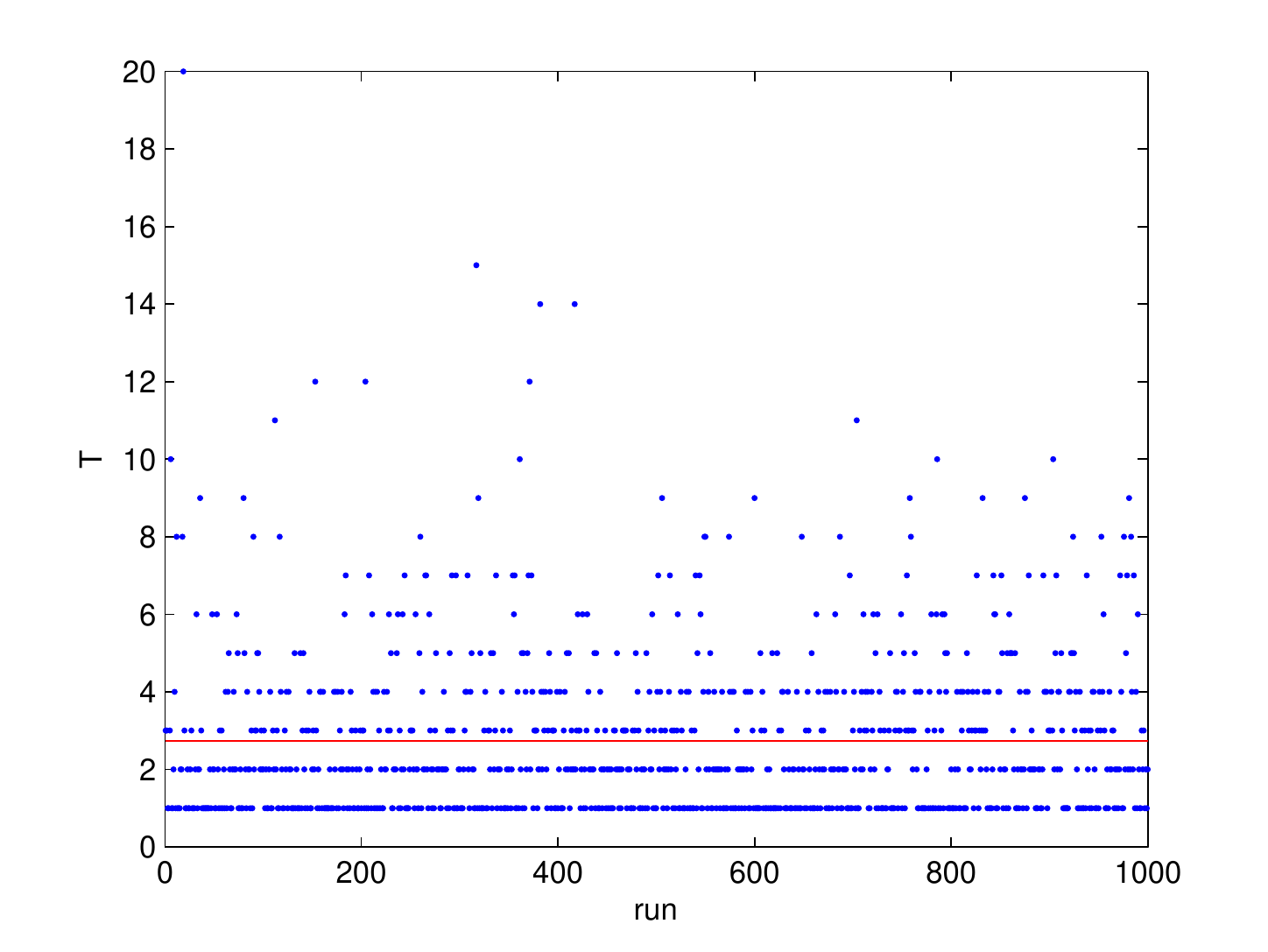}}
\subfigure[~$\mathcal{P}^E$ ($a=0.5$)]{\includegraphics[draft=false, angle=0,width=8cm]{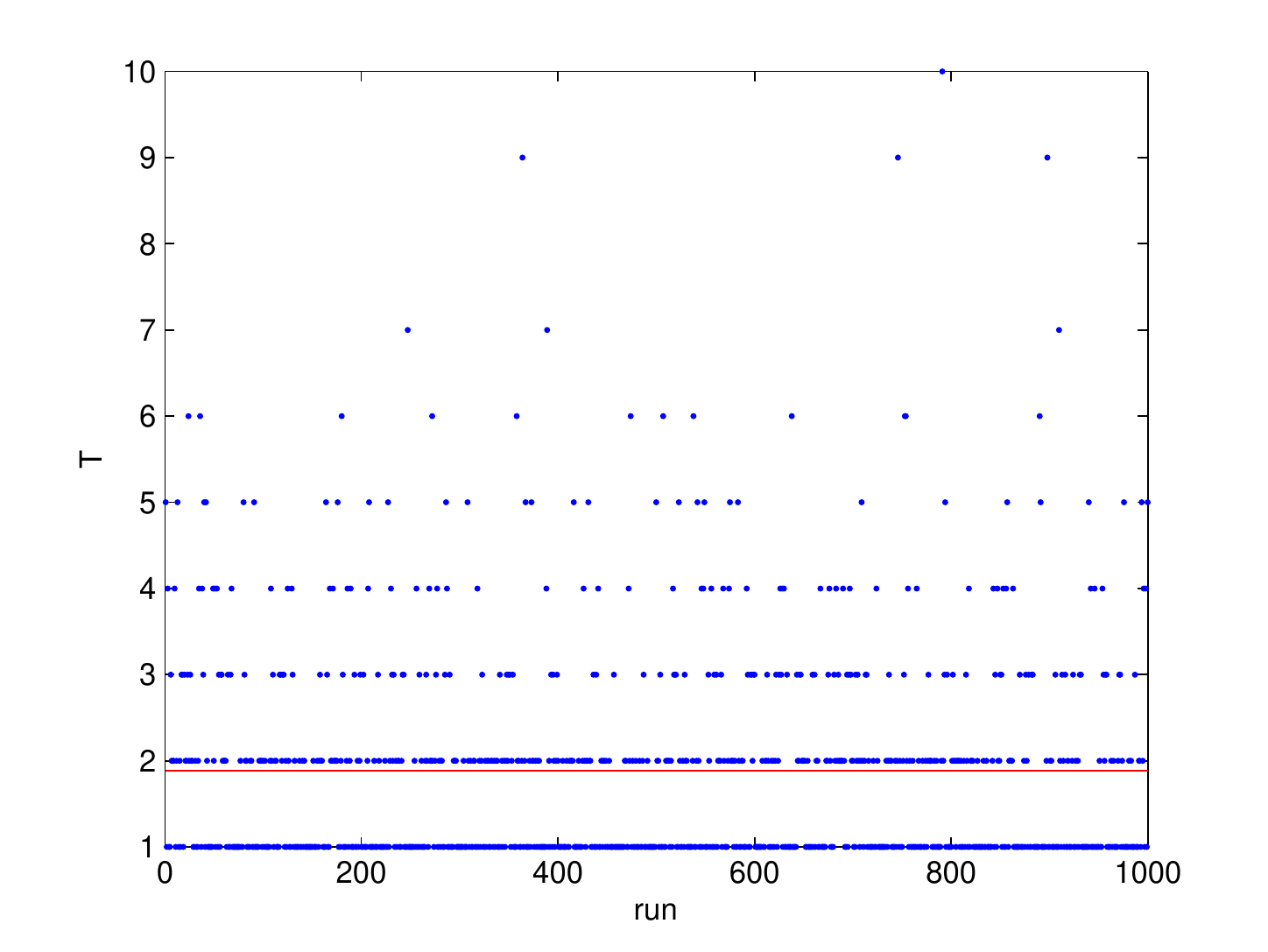}}
\subfigure[~$\mathcal{P}^E$ ($a=0.1$)]{\includegraphics[draft=false, angle=0,width=8cm]{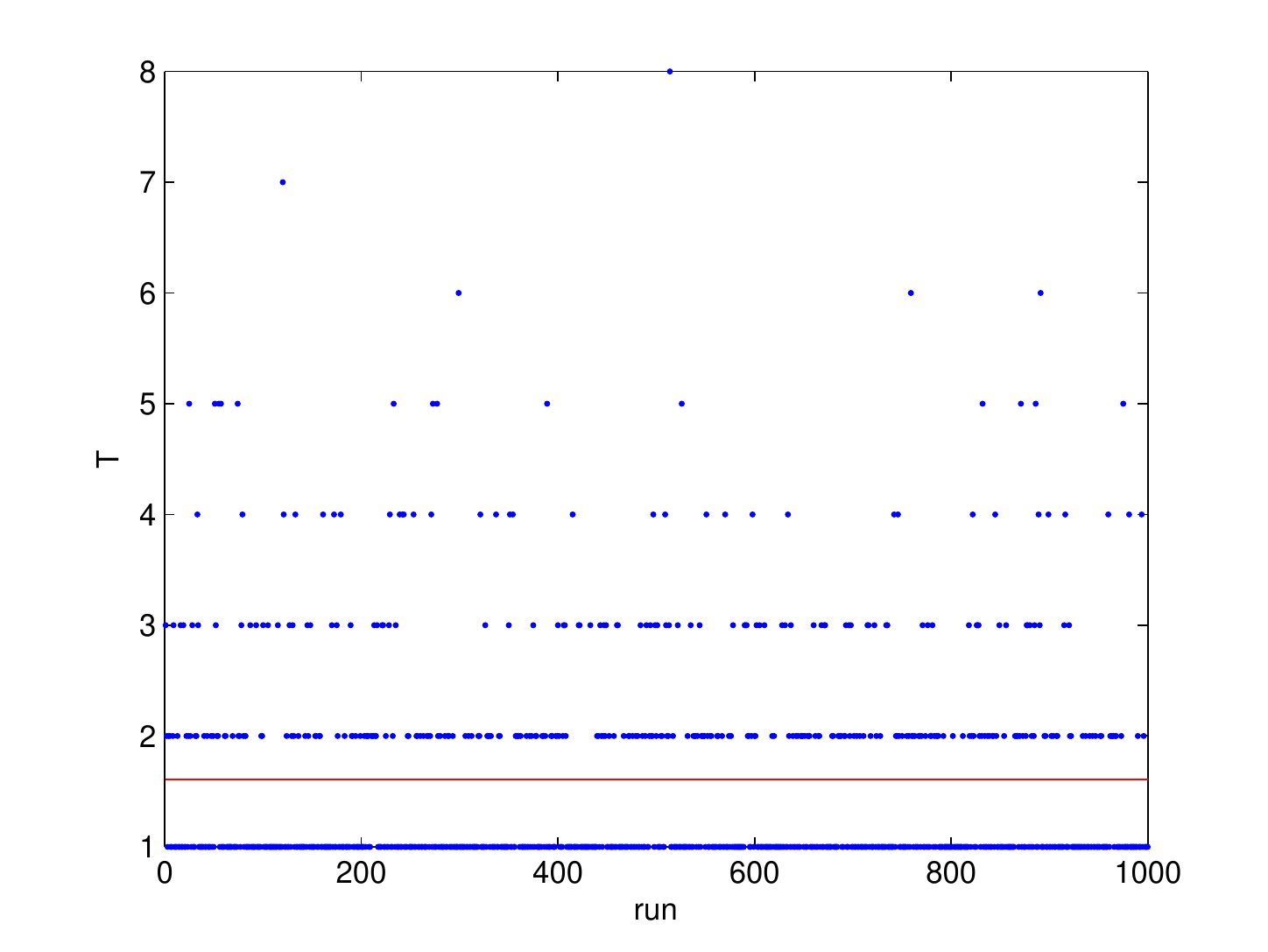}}
\end{center}
\caption{\label{Persistencia_AR(1)E(T)} Persistence estimation ($\mathcal{P}^E=\operatorname{E}[T]$) for the AR(1) system given by Equation (\ref{AR(1)Naive}), $\epsilon=1$, and different values of the parameter $a$.}
\end{figure}

The persistence of the AR(1) system given by Equation (\ref{AR(1)Naive}) can be also analyzed in terms of the measure $\mathcal{P}^E=\operatorname{E}[T]$ by Monte-Carlo simulation. Let us consider a case in which $\epsilon=1$, and $10^4$ epochs of the Monte-Carlo simulation, for different values of $a$: $0.9$, $0.7$, $0.5$ and $0.1$. Figure \ref{Persistencia_AR(1)E(T)} shows the value of $\mathcal{P}^E$ in each case. As can be seen, the system shows a higher degree of persistence when $a \rightarrow 1$, as previously discussed.
In a similar way, we can examine the persistence of an ARMA(1,1) system, given by the following equation:
\begin{equation}\label{ARMA(1,1)Naive}
    x[n]=a_1 x[n-1]-b_1 \varepsilon[n-1]+\varepsilon[n]
\end{equation}
Figure \ref{Persistencia_ARMA(1,1)E(T)} shows the persistence estimation of the ARMA(1,1) system in this case, for different values of $a_1$ and $a_1$. This experiment shows that a higher persistence degree is obtained for the ARMA(1,1) system when $a_1 \rightarrow 1$ and $b_1 \rightarrow 0$, as expected.
\begin{figure}[!ht]
\begin{center}
\subfigure[~$\mathcal{P}^E$ ($a_1=0.9,b_1=0.9$)]{\includegraphics[draft=false, angle=0,width=8cm]{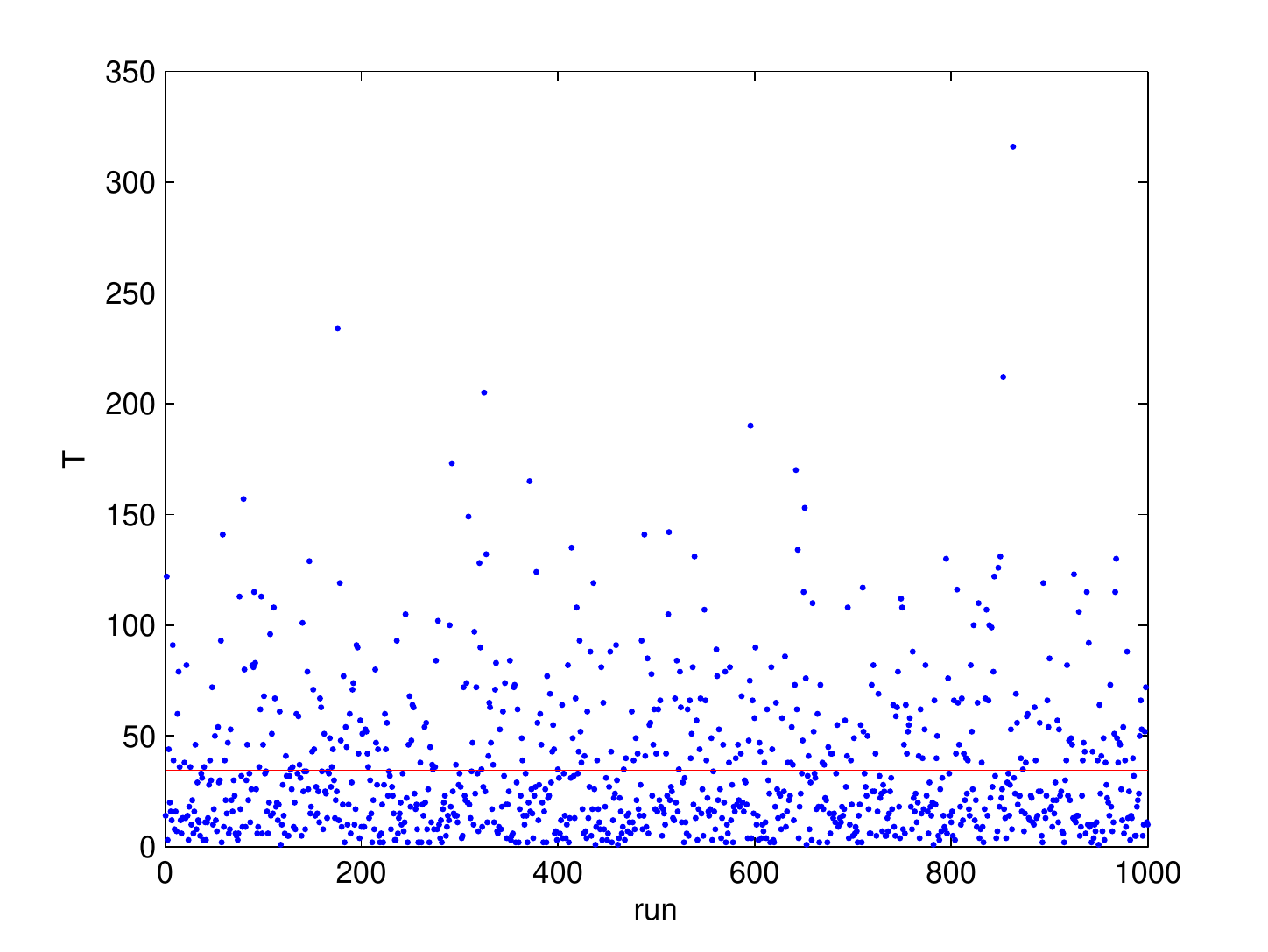}}
\subfigure[~$\mathcal{P}^E$ ($a_1=0.9,b_1=0.5$)]{\includegraphics[draft=false, angle=0,width=8cm]{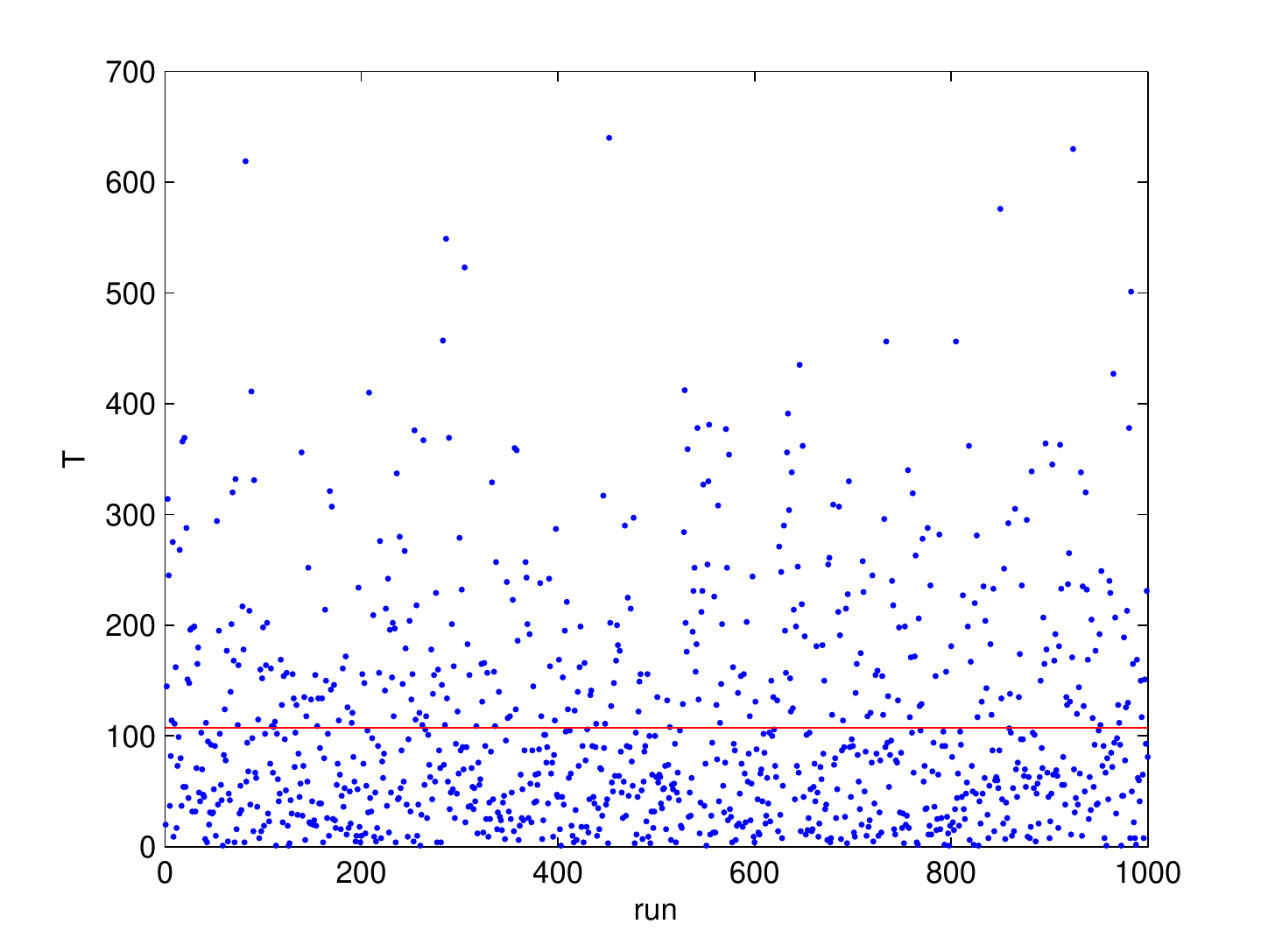}}
\subfigure[~$\mathcal{P}^E$ ($a_1=0.9,b_1=0.1$)]{\includegraphics[draft=false, angle=0,width=8cm]{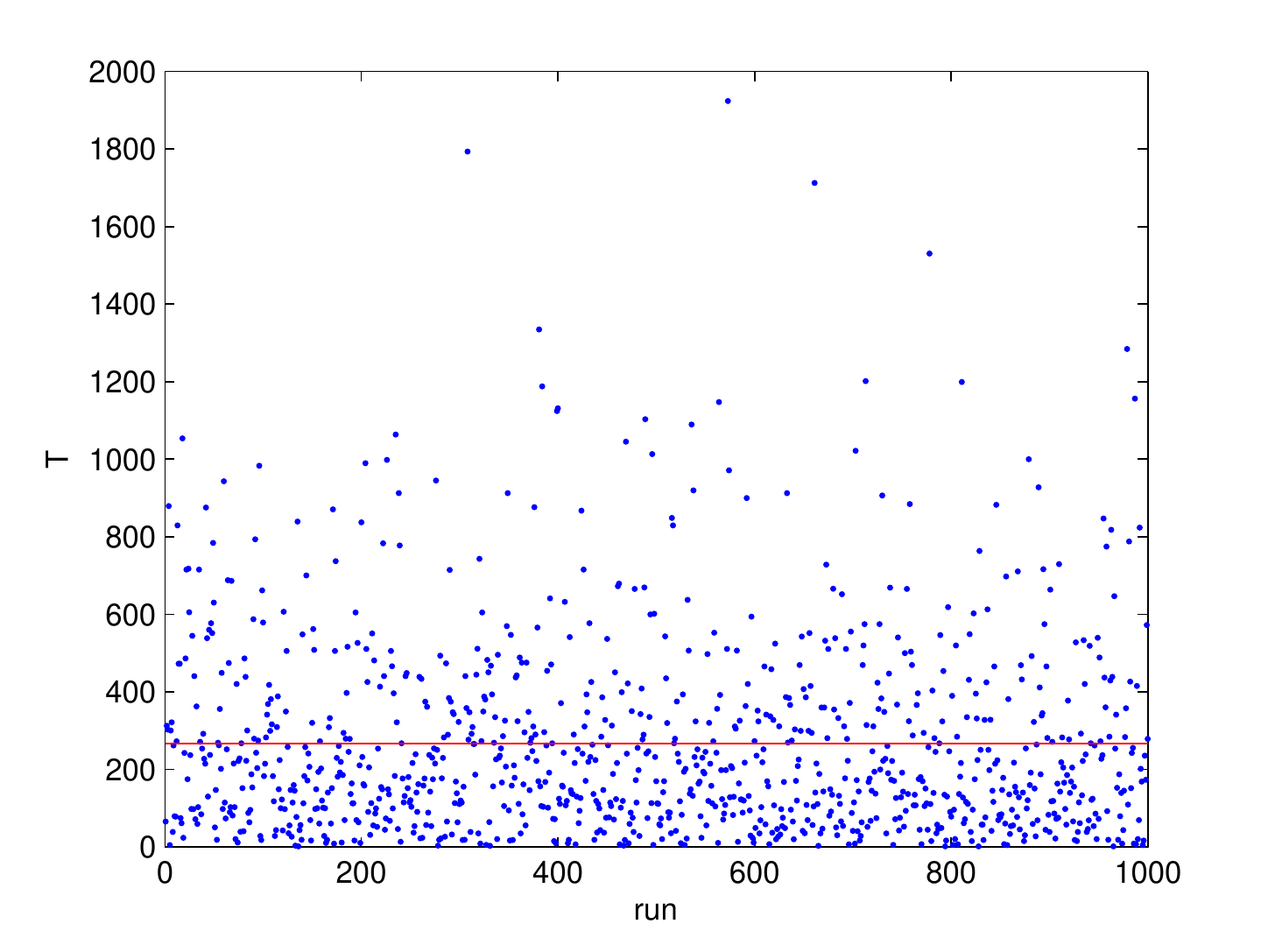}}
\subfigure[~$\mathcal{P}^E$ ($a_1=0.5,b_1=0.9$)]{\includegraphics[draft=false, angle=0,width=8cm]{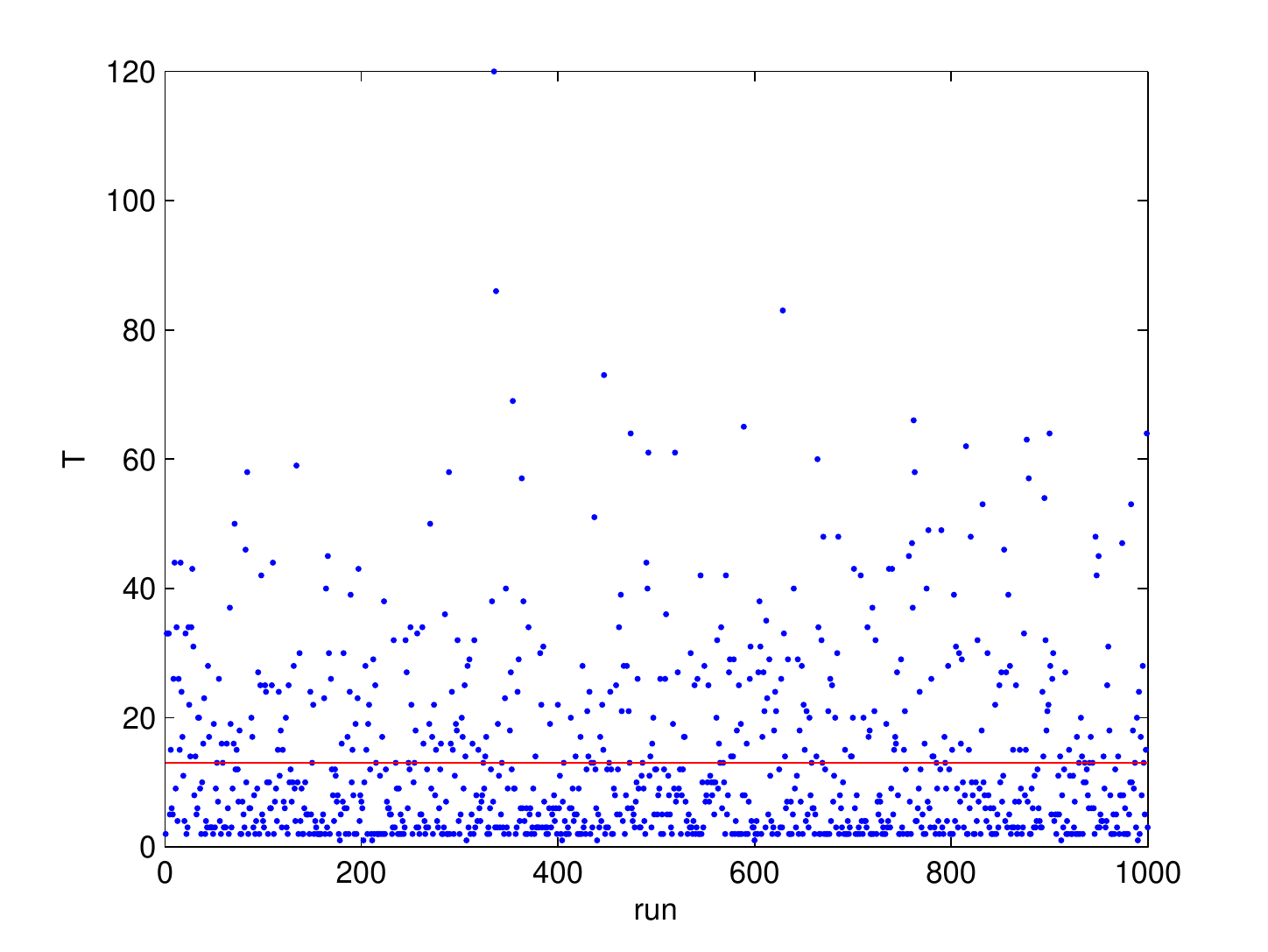}}
\end{center}
\caption{\label{Persistencia_ARMA(1,1)E(T)} Persistence estimation ($\mathcal{P}^E=\operatorname{E}[T]$) for the ARMA(1,1) system given by Equation (\ref{ARMA(1,1)Naive}), $\epsilon=1$, and different values of the parameters $a_1$ and $b_1$.}
\end{figure}

\subsubsection{Memory in ANNs: Monitoring patients with dynamic neural networks}
ANNs are widely used in pharmacy and in daily clinical practice to monitor patients status and health. Despite progress with newer agents, cyclosporine A (CyA) is still the cornerstone of immuno-suppression in renal transplant recipients and its use continues to expand globally \cite{ar:belitsky00}. CyA is generally considered a critical dose drug. Its narrow therapeutical range is an important issue in the clinical management of transplant patients, whereas underdosing may result in graft loss and overdosing causes kidney damage, increases opportunistic infections, systolic and diastolic pressure, and cholesterol. An intensive strategy of therapeutic drug monitoring is necessary during the early post-transplantation period \cite{ar:levy02}. Several studies have been done on the CyA blood concentration prediction \cite{ar:brier95cya,CampsIEE02}. Limitations such as non-uniform sampling during routine clinical data collection, the presence of non-stationary pharmacokinetic processes and the high variability in the CyA blood concentration time series, led to the use of ANNs for this problem. Recurrency and memory are important characteristics of the generating mechanism that should be taken into account. For this, both time-embedding MLP, dynamic FIR networks, and (Elman's) recurrent neural networks were used in \cite{Camps-Valls2003442}.

Thirty-two renal recipients treated in the Dr. Peset University  Hospital in the city of Val{\`e}ncia (Spain) were included in this study  \cite{Camps-Valls2003442}. The patients received a standard immuno-suppressive regimen with a microemulsion lipidic formulation of CyA (Sandimmun Neoral${^{\circledR}}$), mycophenolate mofetil (2 g/d) and prednisone (0.5-1 mg/kg/d). The initial oral dose of CyA (5 mg/kg b.i.d) was reduced according to the measured CyA blood concentration and the desired target range (150-300 ng/mL) \cite{ar:oellerich95}. The patients were randomly assigned to two groups: twenty-two patients were used for training the models (364 samples) and the other 10 patients constituted the validation set (217 samples). Patients were closely monitored  
and a set of variables were collected: 
CyA blood concentration (C [ng/mL]), daily dosage of CyA (DD [mg/Kg/d]), creatinine levels (CR [mg/dL]), anthropometric factors (age, AG [yr], gender, GE, and total body weight, WE [Kg]) and the post-transplantation days, PTD [d].

ARMA modeling showed poor results due to the high variability in the time series, the presence of non-stationary pharmacokinetic processes in the early post-transplantation days, and the direct impact of the non-uniform sampling and nonlinear processes. Results did not improve when more AR terms were considered. The Elman's network was slightly more accurate (-0.019 ng/mL) than the MLP (-0.279 ng/mL) and the FIR net (0.176 ng/mL). Precision for the FIR network was slightly better as can be observed from the RMSE measurements in Table \ref{tab:cc}.

\begin{table}[t!]
\begin{center}
\small
\caption{Models comparison for blood concentration prediction in the validation set. Architectures of the networks are indicated in the form {\sf Input $\times$ Hidden $\times$ Output} nodes. The root-mean-square error (RMSE), the mean error (ME) and the correlation coefficient (r) are given. 95\% CI are given in brackets and were calculated using bootstrap methods.}
\vspace{0.1cm}
\begin{tabular}{lccc@{\vspace{0.1cm}}}
\hline 
\hline
{\bf Performance} & {\bf r} & {\bf ME (ng/mL)} & {\bf RMSE (ng/mL)} \\
& & {\bf ($\pm$95\% CI)} &  {\bf ($\pm$95\% CI)$^\dag$}\\
\hline
\hline
{\bf AR(3)} & 0.325 & 9.878 & 90.96 \\
 & & (-8.907,10.849) & (79.65,102.27)\\
\hline
{\bf MLP (14$\times$7$\times$1)} & 0.728 & -0.279 & 58.69 \\
 & & (-8.106, 7.547) & (52.43,64.91)\\
\hline
{\bf ELMAN} & 0.762 & -0.019 & 55.19 \\
{\bf (14$\times$3$\times$1, CU=3)} & & (-7.379,7.342) & (49.43,61.62)\\
\hline
{\bf FIR (7$\times$4$\times$4$\times$1, } & 0.785 & 0.176 & 52.80 \\
{\bf Taps per layer: 1:1:1) } & & (-6.865,7.217) & (47.79,58.27)\\
\hline
\hline
\end{tabular}
\label{tab:cc}
\normalsize
\end{center}
\end{table}

\begin{figure}[t!]
	\centerline{\includegraphics[width=6cm,trim=0 0 400 0,clip]{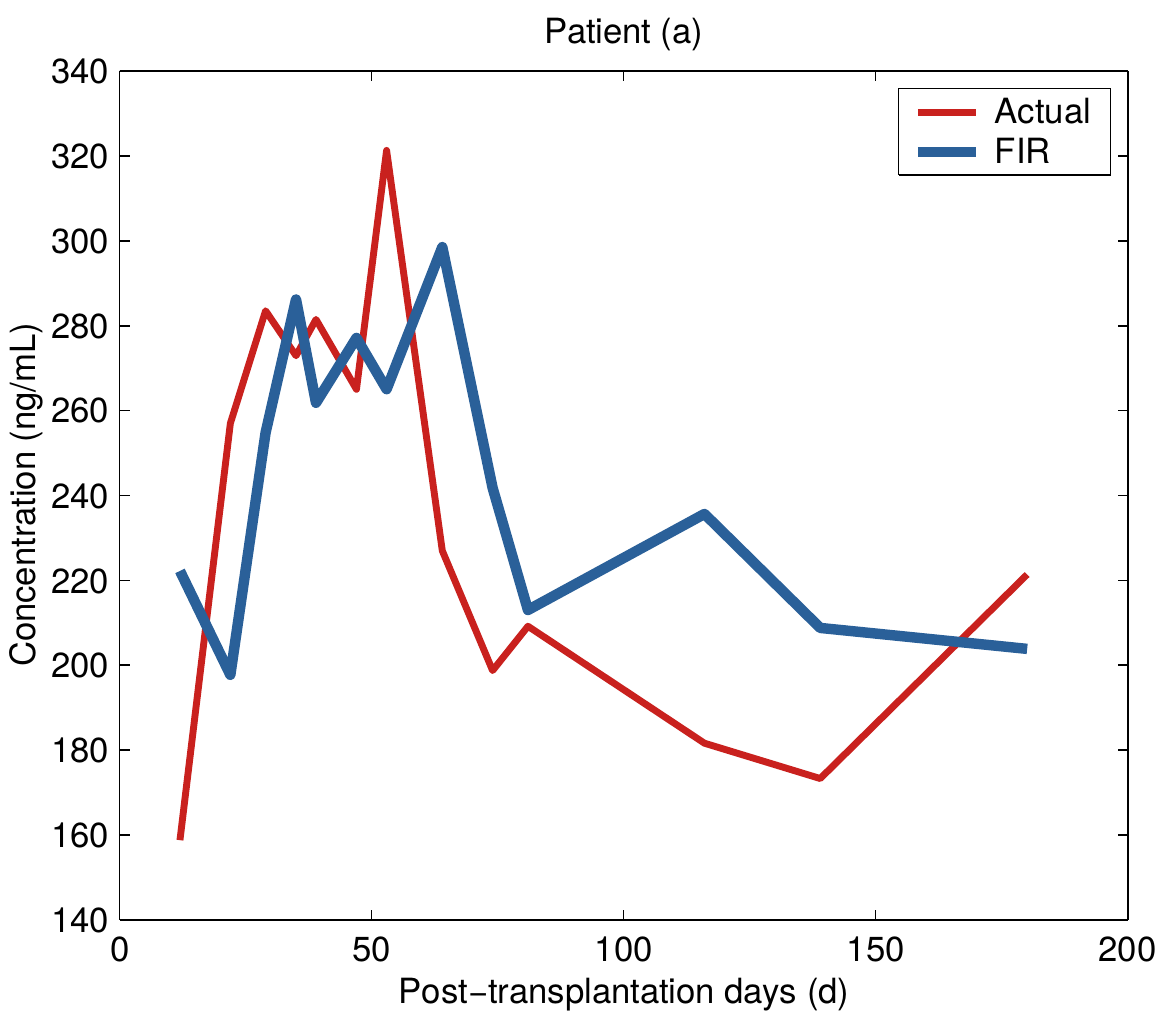}
	\includegraphics[width=6cm,trim=0 0 400 0 ,clip]{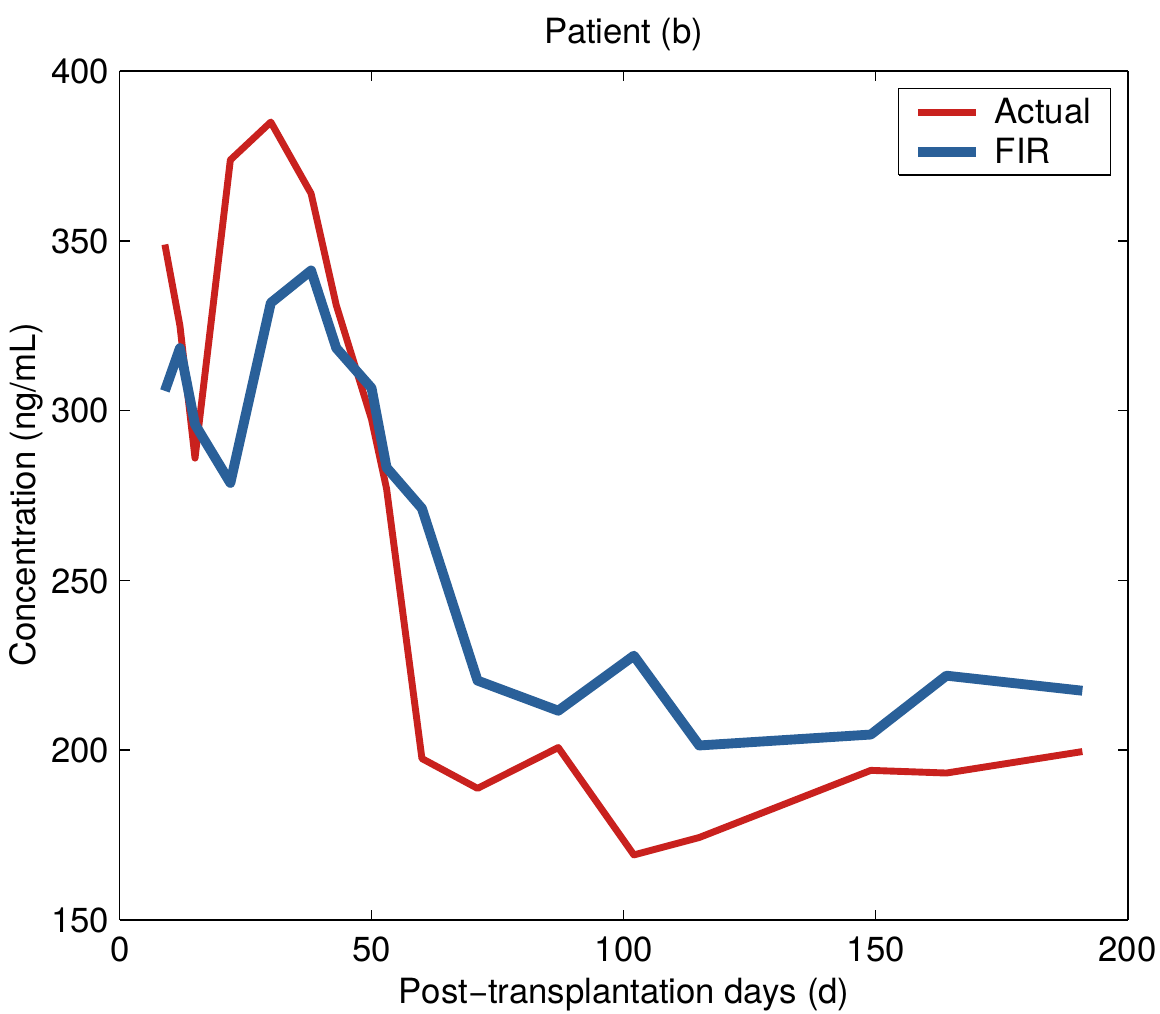}
	\includegraphics[width=6cm,trim=0 0 400 0,clip]{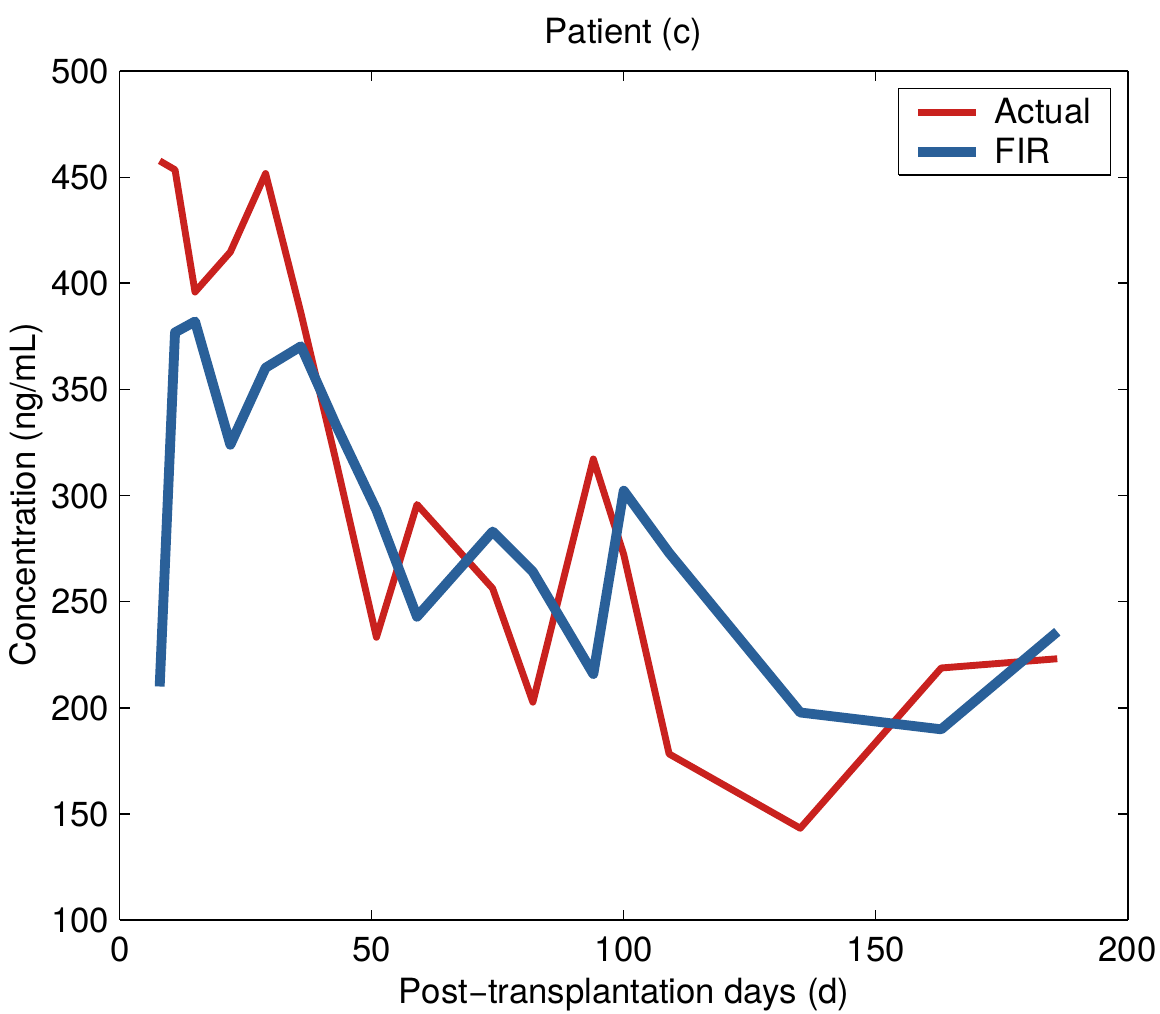}}
	\caption{Plots of observed and predicted CyA concentration (ng/mL) versus postoperative day in individual patients showing (a) good, (b) acceptable and (c) unsatisfactory predictive performance. Only the Elman network for dosage prediction and the FIR network for concentration prediction are shown.\label{fig:cya}}
\end{figure}

The models captured abrupt changes in the time series of blood concentration in the patients, and all of the models performed similarly in the early post-transplantation period. Nevertheless, 21\% of the patients had very poor predictions (RMSE$>$65 ng/mL) which can be due to errors in drug dosage administration, to the inter- and intrasubject variability in the drug absorption process, in the recording of blood sampling times or abrupt changes in each patient's clinical condition. As an example of these situations, we show three patients with good (Figure \ref{fig:cya}a), acceptable (Figure \ref{fig:cya}b) and poor (Figure \ref{fig:cya}c) predictive performances. FIR network predictions captured the complex and heterogeneous memory processes better, and yielded improved results. 

\subsection{Methods for long-term persistence characterization}

\subsubsection{Hurst Rescaled Range (R/S) analysis}\label{H_estim}
The first approach to characterize long-term persistence is due to Hurst \cite{Hurst51}, who carried out pioneering studies on the record of floods and droughts in the Nile river. Hurst introduced the first method for long-term persistence characterization, known as Rescaled Range (R/S) analysis. R/S analysis can be described as follows:
\begin{itemize}
\item[1)] R/S analysis first considers the integrated time series from original $x[n]$:
\begin{equation}\label{Eq:Y_j}
    y[n]=\sum_{i=1}^n x[i].
\end{equation}
The profile $y[n]$ is divided into $N_s = \left\lfloor \frac{N}{s} \right\rfloor$ non-overlapping segments $Y^s=\lbrace y_k^s[n] ~|~ 1 \leq k \leq N_s \rbrace$ of equal length $s$. 
\item[2)] A range $R_{m,s}$ is used to describe the dispersion of these values, looking at the maximum and minimum $Y_j$ values within each segment $m$ of length $s$. It is defined as:
\begin{equation}
R_{m,s}=\max\left[y_{m}^s[n]\right]-\min\left[y_{m}^s[n]\right]
\end{equation}
\item[3)] For each segment $m$ of length $s$, the variance of the original $x[n]$ values in that segment is computed as:
\begin{equation}
S_{m,s}=\sigma_x\left[x_{m}^s\right]
\end{equation}
where this expression indicates taking the standard deviation over each segment of length $s$ in the time series $x[n]$.
\item[4)] Mean values of the range $R_{m,s}$ and the standard deviation $S_{m,s}$ for segments of length $s$ are determined:
\begin{equation}
R_s=\bar{R}_{m,s}=\frac{1}{N_s}\sum_{m=1}^{N_s} R_{m,s}
\end{equation}
and
\begin{equation}
S_s=\bar{S}_{m,s}=\frac{1}{N_s}\sum_{m=1}^{N_s} S_{m,s}
\end{equation}
The ratio $\left(\frac{R_s}{S_s}\right)$, exhibits a power-law scaling as a function of segment length $s$, with a power-law exponent called the Hurst exponent ($H$), i.e.
\begin{equation}
\left(\frac{R_s}{S_s}\right) \sim \left(\frac{s}{2}\right)^H
\end{equation}
\end{itemize}
Note that the Hurst exponent $H$ is related to the strength of long-range persistence $\beta$ as $\beta = 2H-1$ \cite{Malamud99}.

\subsubsection{Detrended Fluctuation Analysis (DFA)}\label{DFA}

Since the DFA algorithm was proposed in two seminal works \cite{Peng94,Peng95}, it has been frequently used to analyze long-term persistence of time series. In the last years, DFA has been much more used than classical R/S methodology, as proved to be highly effective in the analysis of long-term persistence of time series \cite{Kantelhardt01,Hu01}, though some criticisms have also been raised \cite{bryce2012revisiting}. The DFA algorithm consists of three main steps \cite{Hu01}:
\begin{itemize}
\item[1)] 
The periodic annual cycle of the time series is first removed, following the procedure explained in detail in \cite{Yang19}. The process consists on standardizing the input time series $x[n]$ of length $N$ as follows:
\begin{equation}
    \hat x[n] =\frac{x[n] - \mu_x}{\sigma_x},
\end{equation}
where $x[n]$ stands for the original time series, $\mu_x$ represents the mean value of the time series and $\sigma_x$ its standard deviation.
\item[2)] Then, the time series profile $y[n]$ (integrated time series) is computed as follows:
\begin{equation}\label{Eq:Y_j2}
    y[n]=\sum_{i=1}^n \hat x[i].
\end{equation}
The profile $y[n]$ is divided into $N_s = \left\lfloor \frac{N}{s} \right\rfloor$ non-overlapping segments $Y^s=\lbrace y_k^s[n] ~|~ 1 \leq k \leq N_s \rbrace$ of equal length $s$. 
For each segment $y_k^s[n]$, we calculate the local least squares fit $Z^s_k$ of the time series, which measures its local trend.  Note that linear, quadratic, cubic, or higher order polynomials can be used in the fitting procedure, which sets the {\em order} $m$ of the DFA. Usually an order $m=2$ (DFA-2) is used. As a result, we obtain a piece-wise function $\tilde z^s[n]$ compounding each fitting:
\begin{equation}
\tilde z^s[n]=\begin{bmatrix}
Z_1^s & \cdots & Z_k^s & \cdots &  Z_{N_s}^s
\end{bmatrix},
\end{equation}
where the superscript $s$ refers to the time window length used to the perform the fitting of each segment.
\item[3)]
We then obtain the so-called \emph{fluctuation} as the root-mean-square error from this piece-wise function $\tilde z^s[n]$ and the profile $y[n]$, varying the time window length $s$:
\begin{equation}\label{DFA_error}
F(s) = \sqrt{\frac{1}{N}\sum_{k=1}^{N} (\tilde z^s[k]-y[k])^2}.
\end{equation}
At the time scale range where the scaling holds, $F(s)$ increases with the time window $s$ following a power law $F(s) \propto s^{\alpha}$. Thus, the fluctuation $F(s)$ versus the time scale $s$ is depicted as a straight line in a log-log plot. The slope of the fitted linear regression line is the scaling exponent $\alpha$, also called correlation exponent. The scaling exponent $\alpha$ in the DFA method is related to the strength of long-term persistence by $\beta = 2\alpha-1$, so $\alpha=H$ in this context.

Note that when the coefficient $\alpha = 0.5$, the time series is uncorrelated, which means that there is no long-term persistence in the time series. For larger values of $\alpha$ ($0.5 < \alpha \leq 1$), the time series is positively long-term correlated, which also means the long-term persistence exists across the corresponding scale range. When $0 < \alpha \leq 0.5$ the process is anti-persistent. For $\alpha > 1$, the persistence becomes so extreme, that the time series exhibits non-stationary behavior. 
\end{itemize}

\subsubsection{Multi-fractal DFA version}

In \cite{Kantelhardt02} a Multi-Fractal DFA (MF-DFA) approach is introduced as a generalization of the standard DFA described above. The reason behind this generalized version is that many records do not exhibit a simple mono-fractal scaling behavior, but a much more complicated behaviour, with several crossover points (characteristic times) with different scaling exponents for different parts of the series. In these cases a multi-fractal analysis must be applied to better describe the long-term persistence in the time series.

\begin{itemize}
\item[1)]
The periodic annual cycle of the time series is first removed:
\begin{equation}
    \hat x[n] =\frac{x[n] - \mu_x}{\sigma_x},
\end{equation}
where $x[n]$ stands for the original time series, $\mu_x$ represents the mean value of the time series and $\sigma_x$ its standard deviation.
\item[2)]
The time series profile $y[n]$ (integrated time series) is computed:
\begin{equation}
    y[n]=\sum_{i=1}^n \hat x[i].
\end{equation}
The profile $y[n]$  is divided into $N_s = \left\lfloor \frac{N}{s} \right\rfloor$ non-overlapping segments $Y^s=\lbrace y_k^s[n] ~|~ 1 \leq k \leq N_s \rbrace$ of equal length $s$. 
As in the standard DFA, for each segment $y_k^s[n]$, one calculates the local least squares straight-line $Z_k^s$ and we obtain a piece-wise function $\tilde z^s[n]$ compounding each fitting:
\begin{equation}
\tilde z^s[n]=\begin{bmatrix}
Z_1^s & \cdots & Z_k^s & \cdots &  Z_{N_s}^s
\end{bmatrix},
\end{equation}
Note that up until now, the procedure is exactly the same as the standard DFA algorithm. The multi-fractal version of the DFA is defined in the next step of the algorithm.
\item[3)] We then obtain now the \emph{generalized fluctuation} as follows:
\begin{equation}\label{DFA_error_MF}
F_q(s) = \left(\frac{1}{N}\sum_{k=1}^{N} (\tilde z^s[k]-y[k])^{q}\right)^{1/q}.
\end{equation}
Note that when $q=2$, the standard DFA approach is obtained. It is possible to determine now the scaling behavior of the fluctuation functions by analyzing log-log plots $F_q(s)$ versus $s$ for each value of $q$. If the series $x_i$ are long-range power-law correlated, $F_q(s)$
increases, for large values of $s$, as a power-law:
\begin{equation}
F_q(s) \sim s^{h(q)}
\end{equation}
where $h(q)$ is a generalized Hurst exponent. Note that for mono-fractal time series $h(q)$ is independent of $q$, and only if small and large fluctuations scale differently, there will be a significant dependence of $h(q)$ on $q$, which can be used to characterize time series with multi-fractal properties. More details on specific aspects of the MF-DFA, such as the cases $q=0$ and $q<0$ and how they can be treated, can be found in \cite{Kantelhardt02}.

\end{itemize}

\subsubsection{Estimating the persistence strength ($\beta$) with the wavelet transform}\label{Wavelet}

An alternative way of estimating the persistence strength $\beta$ is to consider the wavelets transform of a given time series $x[n]$, instead of the power spectral analysis of the series. We present here the estimation of $\beta$ considering discrete wavelet transform, as shown in \cite{Witt13}. The procedure is as follows:

\begin{itemize}
\item[1)] 
Consider a discrete time series $x[n]$ (usually $N$ is taken as a power of $2$ for convenience).
\item[2)]
Choose a discrete valued mother wavelet function $\Psi[n]$. There are different possibilities for this, such as the classical ones proposed in \cite{grossmann1984decomposition,daubechies1988orthonormal}.
\item[3)]
Determine the wavelet basis functions
\begin{equation}
    \Psi_{kl}[n]=2^{-k/2} \Psi\left(2^{-k}(n-l) \right),~1\leq n \leq N
\end{equation}
where $2^k$ is the scale of the wavelet transformation, $k$ is the level $(1 \leq k \leq K)$, with $K=\log_2(N)$, and $l$ is the number of the wavelet coefficient $(1 \leq l \leq L(k))$, with $L(k)=2^{K-k+1}$ is the maximum number of wavelet coefficients per level
\item[4)]
Perform a wavelet transform on the discrete time series $x_t$, as:
\begin{equation}
x[n]=\sum_{k=1}^k\sum_{l=1}^{L(k)} w_{kl} \Psi_{kl}[n],~1\leq n \leq N
\end{equation}
where $w_{kl}$ are the wavelet coefficients, given by the following expression:
\begin{equation}
    w_{kl}=\sum_{n=1}^N x_t \Psi_{kl}[n]=2^{-k/2} \sum_{n=1}^N x[n] \Psi\left[2^{-k/2} (n-l) \right]
\end{equation}
\item[5)]
Calculate the variance of the wavelet coefficients for each scale:
\begin{equation}
    S^2_W[k]=\sigma^2\left(\{w_{k1},w_{k2}, \ldots, w_{kL}\} \right),~ 1\leq k \leq K
\end{equation}
\item[5)]
The plot of $\log_2\left(S^2_W[k]\right)$ versus the wavelet scale $k$ is approximately a line, which slope approximates $\beta$.
\end{itemize}

\subsection{Examples on long-term persistence characterization}

\subsubsection{DFA performance in persistent systems with random noise}
\begin{figure}[!ht]
\begin{center}
\subfigure[]{\includegraphics[draft=false, angle=0,width=8cm]{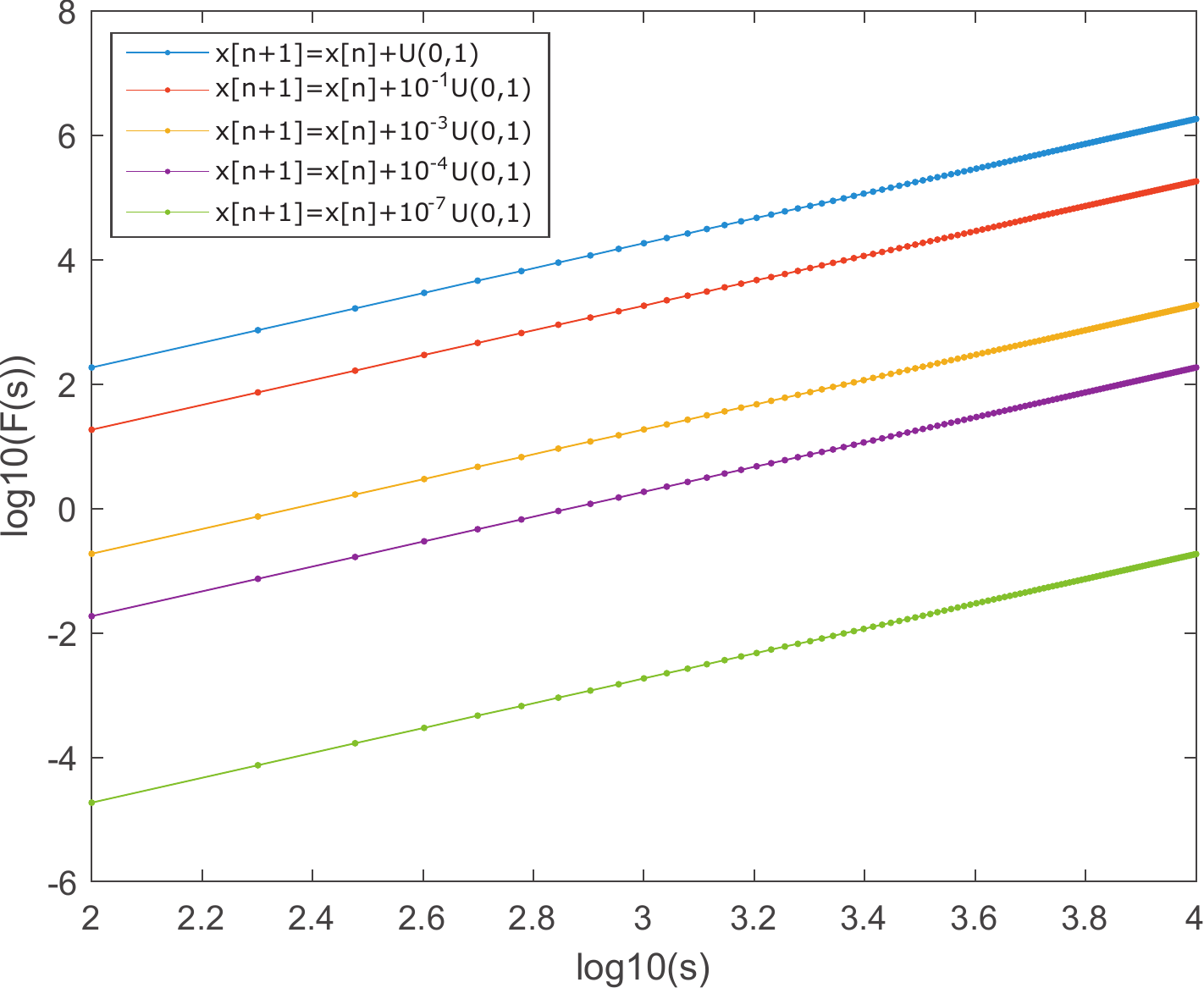}}
\subfigure[]{\includegraphics[draft=false, angle=0,width=8cm]{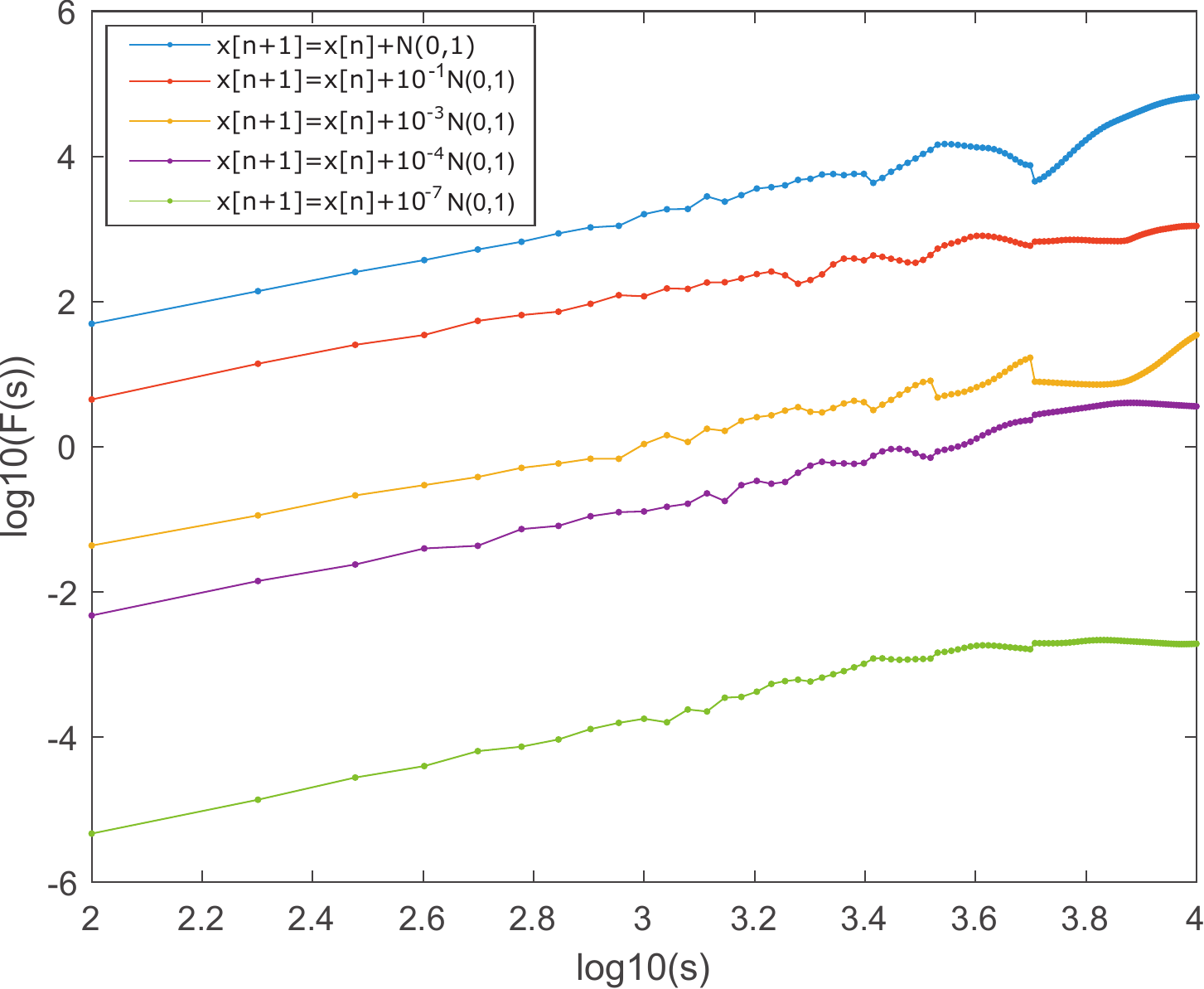}}
\end{center}
\caption{\label{fig:DFA} 
DFA algorithm applied to the persistent systems $x[n+1]=U_n,\ U_n\sim U(0,1)$ and $x[n+1]=N_n,\ N_n\sim N(0,1)$, leading to $\alpha \approx 2$; (a) $x[n+1]=x[n]+U_n,\ U_n\sim U(0,1)$; (b) $x[n+1]=x[n]+N_n,\ N_n\sim N(0,1)$.}
\end{figure}

\begin{figure}[!ht]
\begin{center}
\subfigure[]{\includegraphics[draft=false, angle=0,width=8cm]{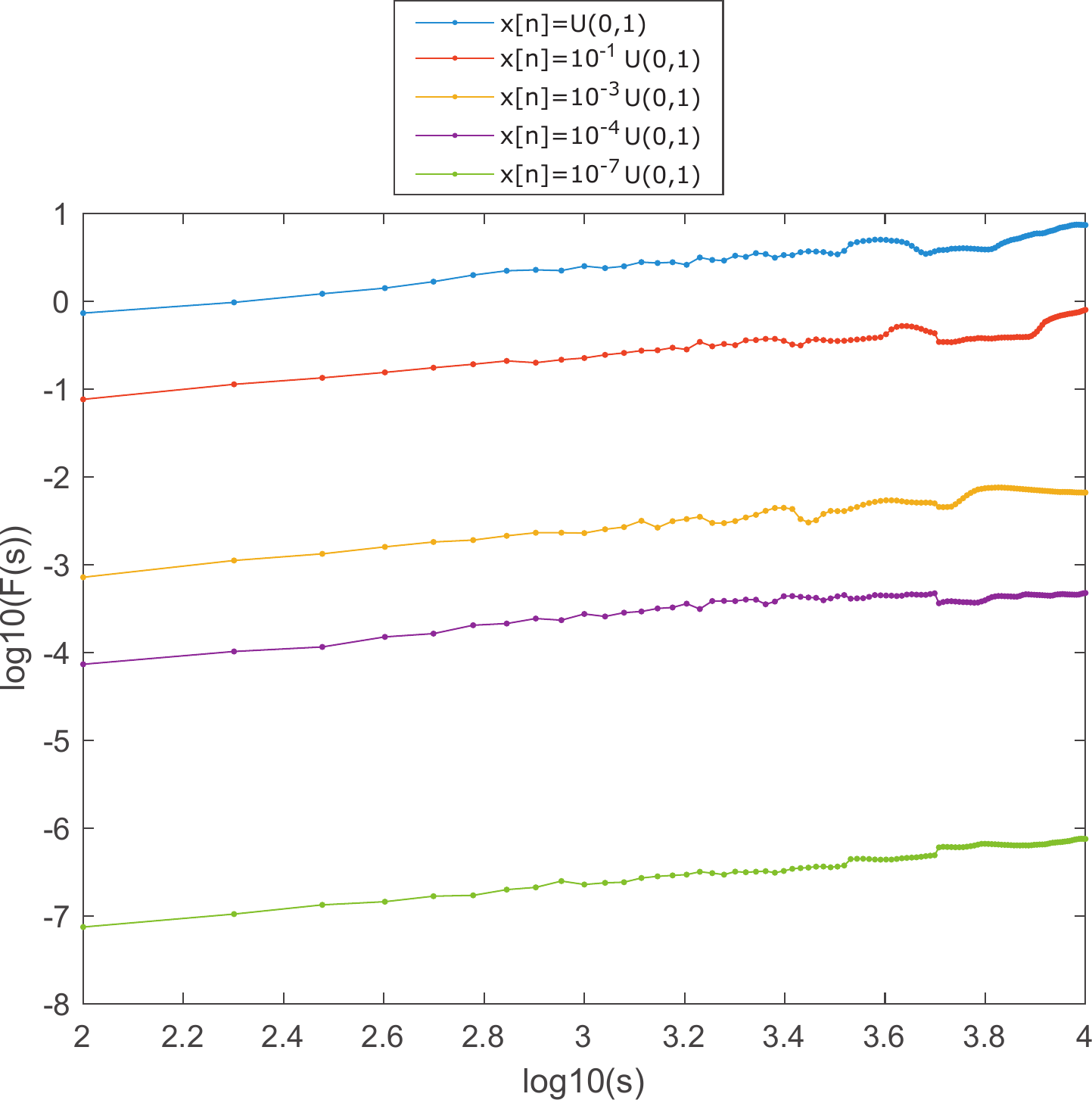}}
\subfigure[]{\includegraphics[draft=false, angle=0,width=8cm]{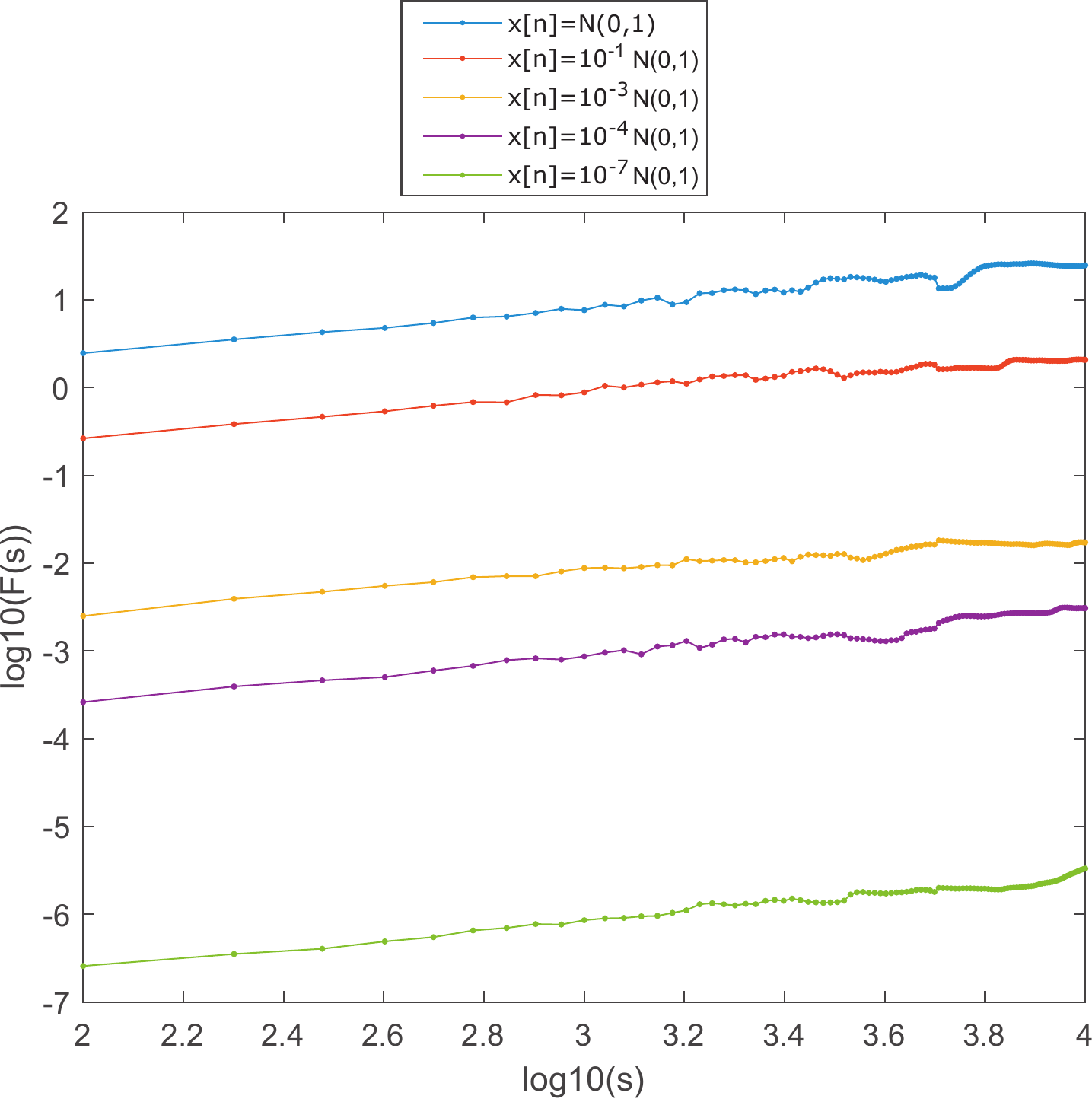}}
\end{center}
\caption{\label{fig:DFA2}
DFA algorithm applied to the uncorrelated systems $x[n+1]=x[n]+U_n,\ U_n\sim U(0,1)$ and $x[n+1]=x[n]+N_n,\ N_n\sim N(0,1)$, leading to $\alpha \approx 0.5$; (a) $x[n]=U_n,\ U_n\sim U(0,1)$; (b) $x[n]=N_n,\ N_n\sim N(0,1)$.}
\end{figure}

We can easily illustrate the performance of the DFA algorithm with some straightforward computational experiments in persistent systems under the presence of random noise. Let us consider two different persistent systems, first $x[n+1]=x[n]+U_n,\ U_n\sim U(0,1)$, where $U(0,1)$ stands for a sample from a Uniform distribution between $0$ and $1$, and second $x[n+1]=x[n]+N_n,\ N_n\sim N(0,1)$, where $N(0,1)$ stands for a sample from a Gaussian distribution with mean $0$ and variance $1$. If we apply the DFA algorithm to both systems ($10^4$ samples each), we obtain the results shown in Figure \ref{fig:DFA}. In both cases the DFA approach indicates a long term correlation of the signals (graphs $F(s)$ versus $s$) with $\alpha \approx 2$. In the case of the uniform noise, 
the DFA obtains a clear correlation even for large values of $s$, but in the case of the Gaussian noise the DFA algorithm with $s>1000$ presents some problems and instability. If we consider uncorrelated systes defined as the corresponding noise sources $x[n]=U_n,\ U_n\sim U(0,1)$ and $x[n]=N_n,\ N_n\sim N(0,1)$ (Figure \ref{fig:DFA2}), the DFA shows that there is not long-term persistence in any case, with linear graphs $F(s)$ versus $s$ with $\alpha\approx 0.5$.

\subsubsection{Long-term persistence and fractal dimension}
In principle, long-term persistence, characterized by $H$, and the fractal dimension of a time series $D$ \cite{falconer2003fractal,li2010fractal}, are independent. Fractal dimension is a local property of time series, whereas long-term memory measures persistence, which can be assumed to be a global characteristic of the system under study \cite{mandelbrot1983fractal}. However, many works have considered a close relationship between $H$ and $D$ \cite{bassingthwaighte1991fractal}. A possible reason for this is that, in self-affine processes \cite{malamud1999self}, local properties are projected into global ones, reaching the well-known relationship between $H$ and $D$,
\begin{equation}\label{H+D}
H+D=n+1,
\end{equation}
where $n$ is the dimension of the considered space. For one-dimensional time series, Equation \eqref{H+D} results in $D=2-H$ \cite{mandelbrot1983fractal}, which has been applied in different works in order to estimate $D$ from long-term persistence estimation methods such as the R/S or the DFA. There are some other works which has pointed out the strict differences in local and global characteristic of $D$ and $H$, such as \cite{gneiting2004stochastic}, where stochastic models are proposed in which $H$ and $D$ are clearly separated, and can be combined arbitrarily. In spite of this there are different works which have exploited the relationship between $D$ and $H$ in the analysis of fractal dimension for long-term persistence characterization of fractal time series. In \cite{north1994bias} a discussion on a bias introduced by the R/S method in the calculation of the fractal dimension is discussed. 

The fractal dimension and Hurst coefficient of wind speed time series and their application in wind speed prediction have been discussed in \cite{liang2015analysis}. In turn the fractal dimension and long-term persistence of wind speed series for variability studies have been discussed in \cite{cadenas2019wind}. In \cite{Breslin99} different methods for calculating the fractal dimension are presented and applied to monthly precipitation time series in Australia, including the R/S method. In \cite{correa2017long} the long-term persistence and fractal dimension of Southern Oscillation Index (SOI) were analyzed by means of the R/S method to calculate the Hurst exponent. The results showed that the considered SOI time series exhibit a chaotic behaviour, with a Hurst index $H=0.56$, which indicates a small long-term persistence of the time series. The relationship between chaotic time series, strange attractors and the Hurst coefficient has also been studied in the literature. In general, non-linear dynamical systems with chaotic behaviour present strange attractors in the phase space. The characterization of strange attractors is usually carried out by calculating the {\em correlation dimension} of the attractor, related to the fractal dimension, usually by applying the Grassberger-Procaccia method \cite{grassberger1983characterization,grassberger1983measuring}. In \cite{de1998r} the R/S method is applied to characterize chaotic signals through their Hurst exponent $H$, as an alternative to other measures to characterize strange attractors. The relationship of $H$ with some other measures to characterize strange attractors are further discussed.

\begin{figure}
    \centering
      \subfigure[]{
        \includegraphics[width=0.45\textwidth]{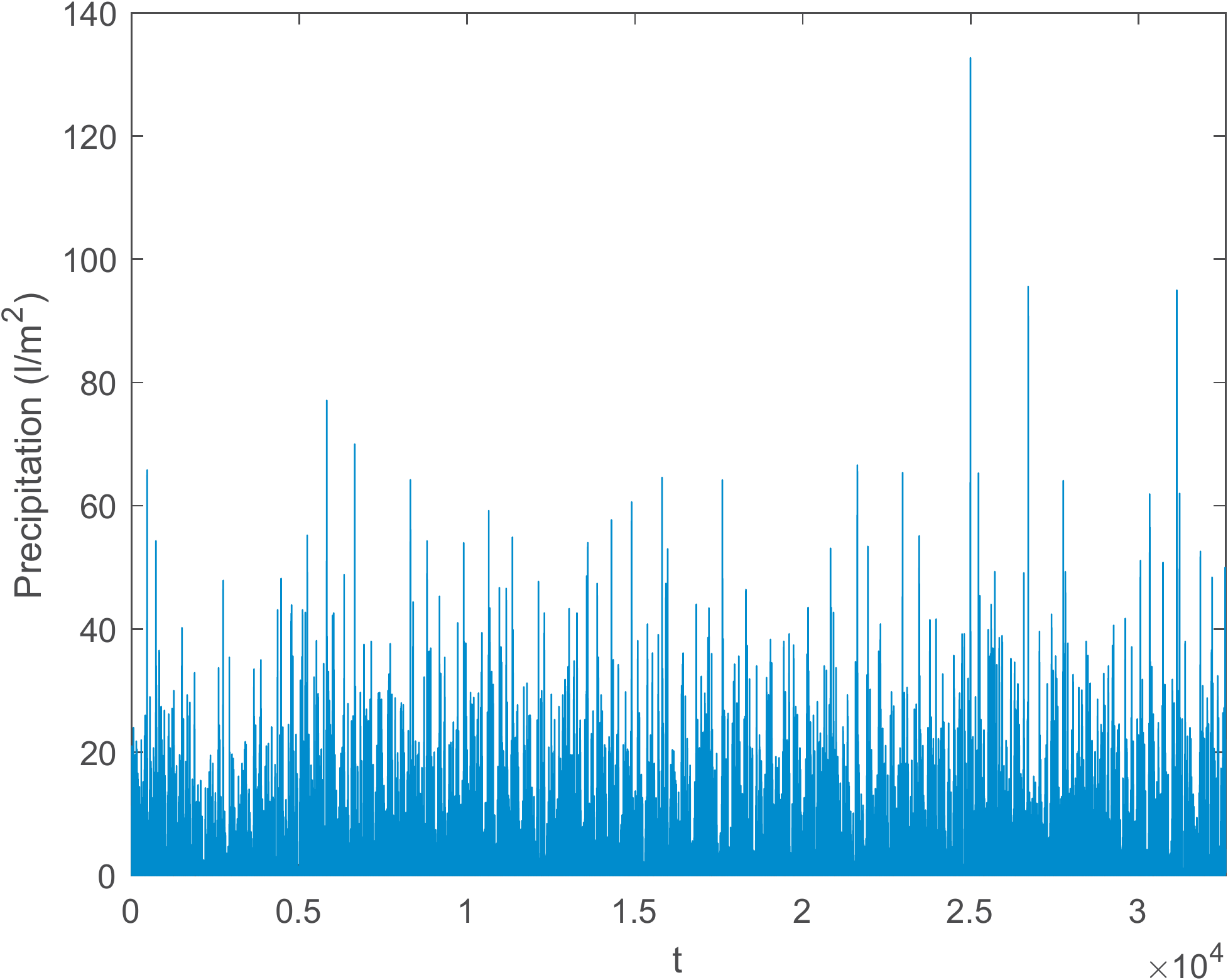}}
      \subfigure[]{
        \includegraphics[width=0.45\textwidth]{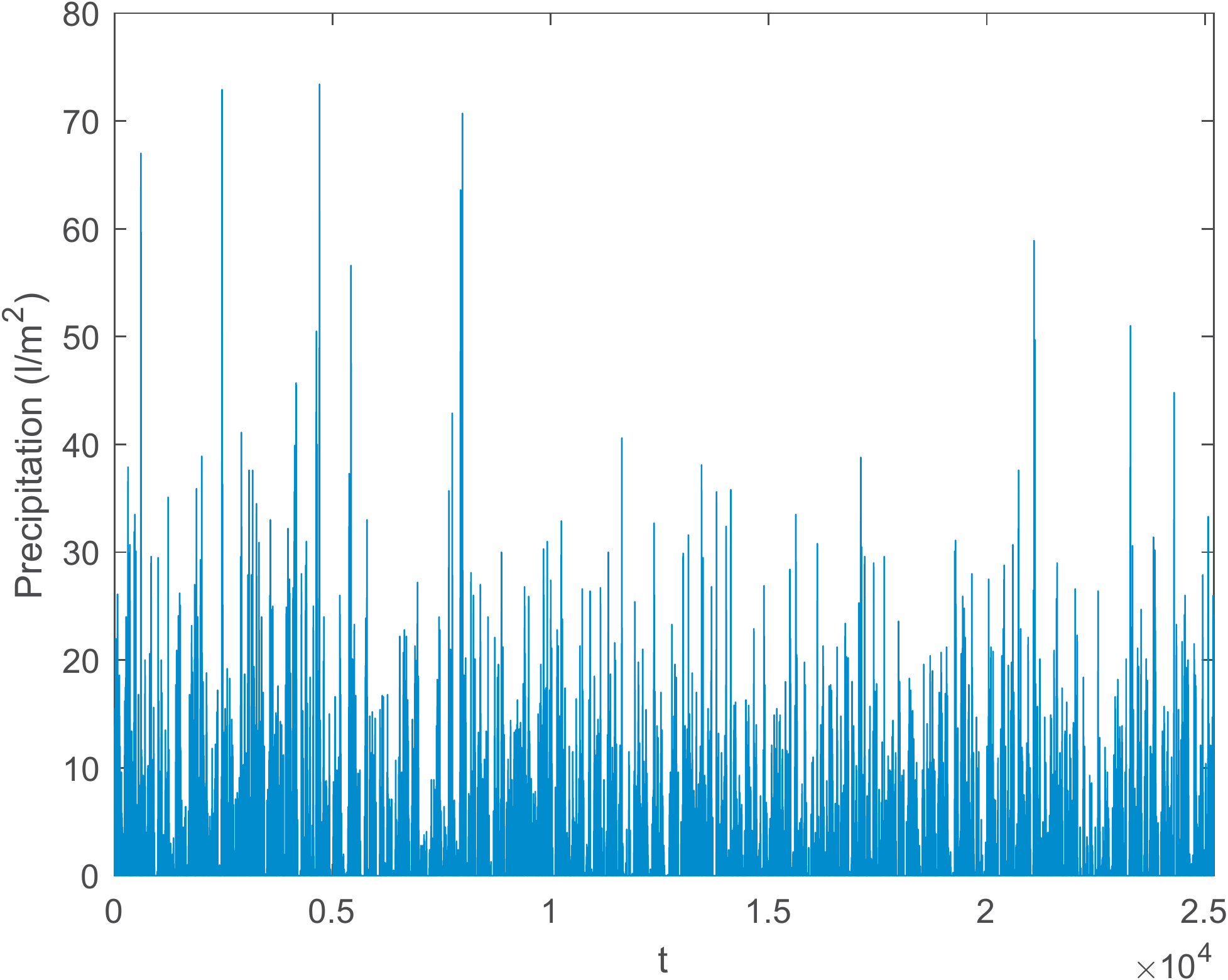}}
        \subfigure[]{
        \includegraphics[width=0.45\textwidth]{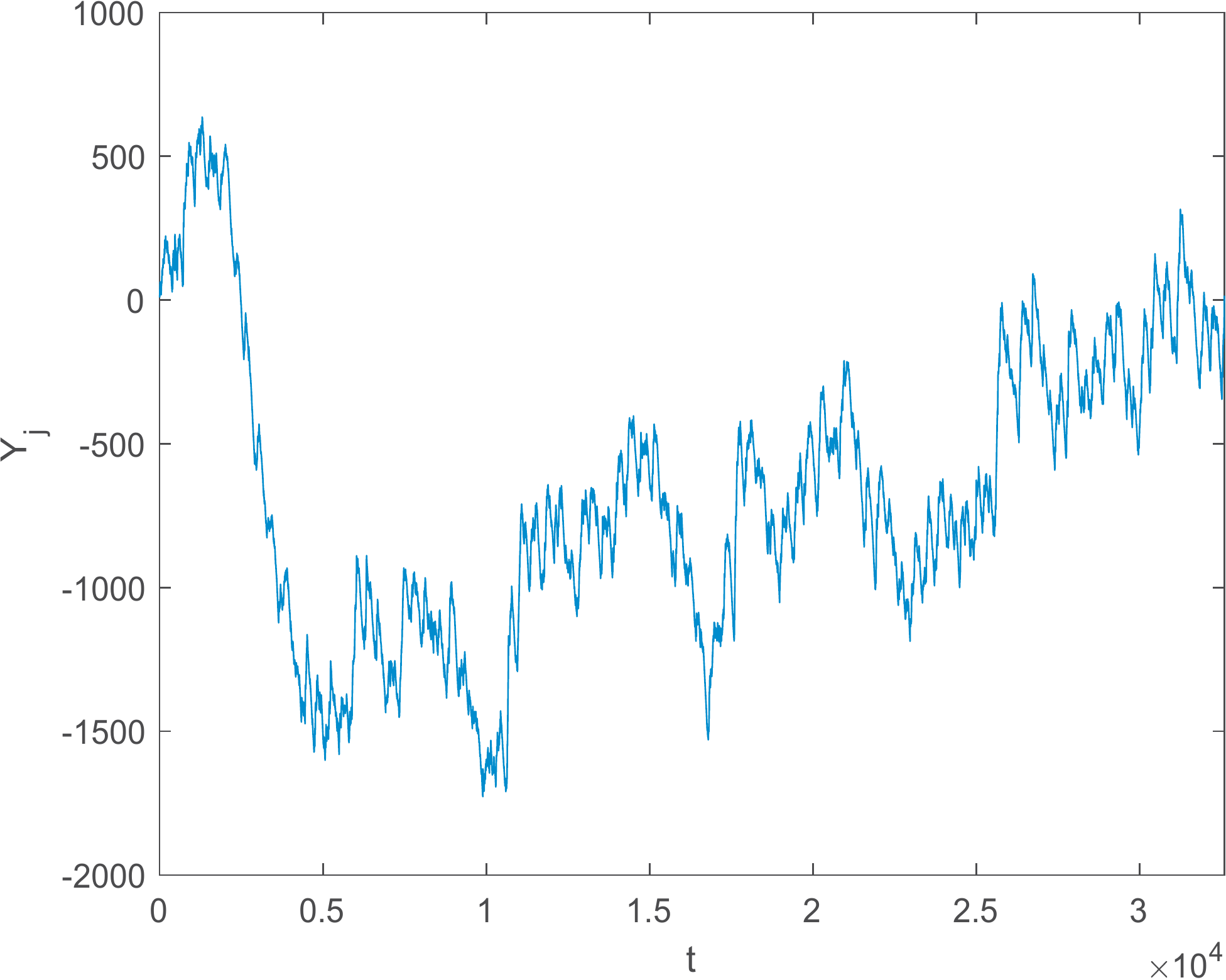}}
      \subfigure[]{
       \includegraphics[width=0.45\textwidth]{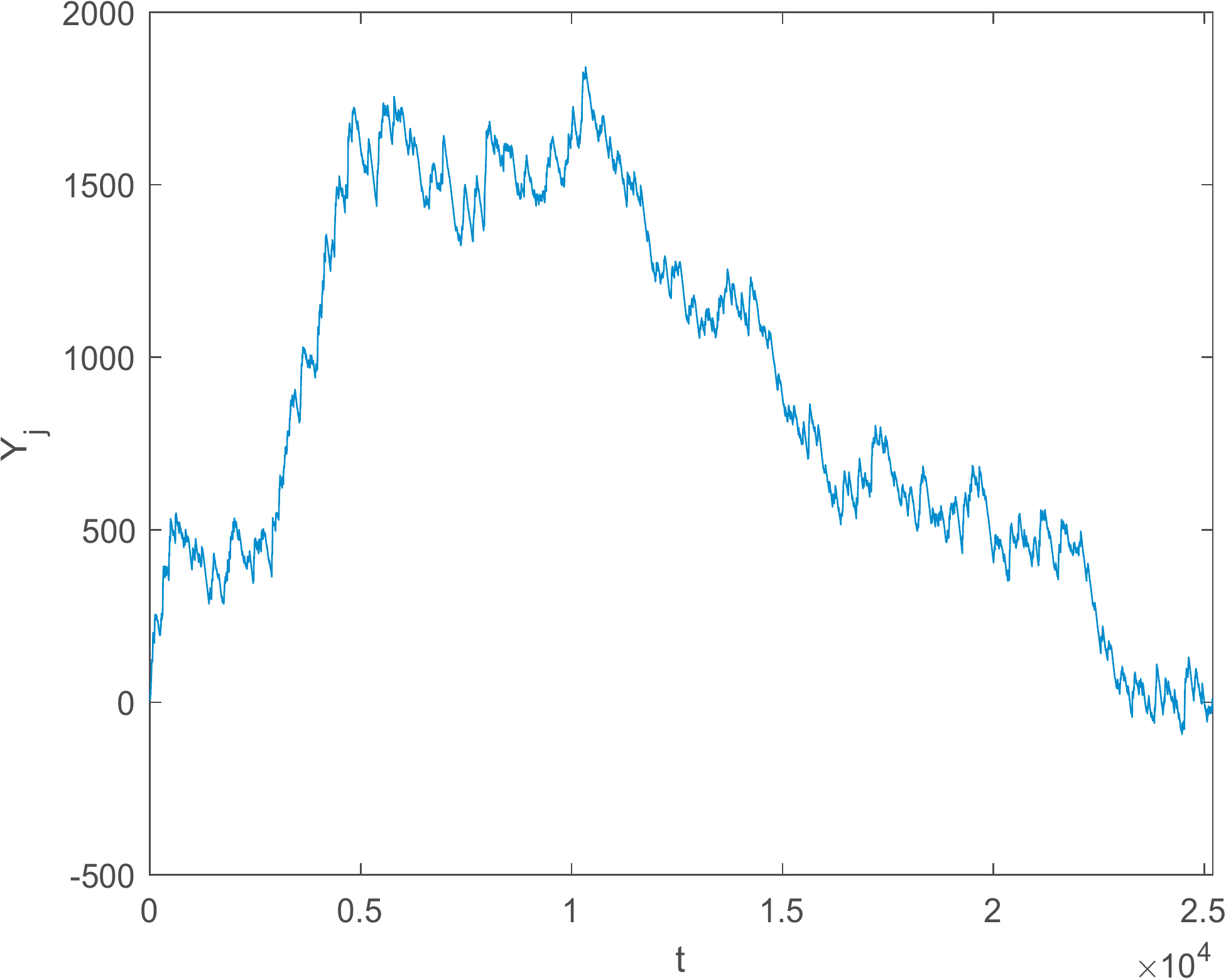}}
    \subfigure[]{
        \includegraphics[width=0.45\textwidth]{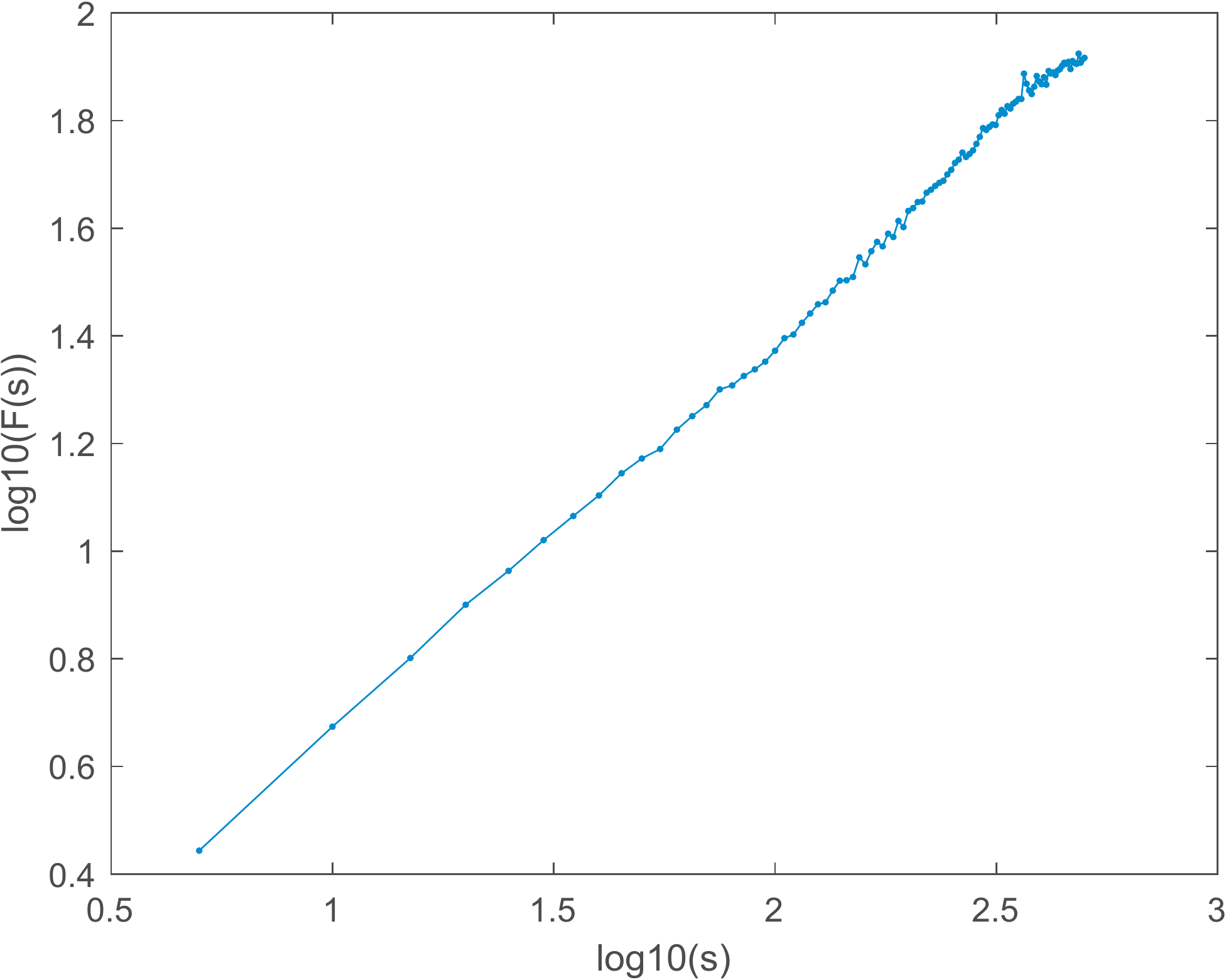}}
      \subfigure[]{
        \includegraphics[width=0.45\textwidth]{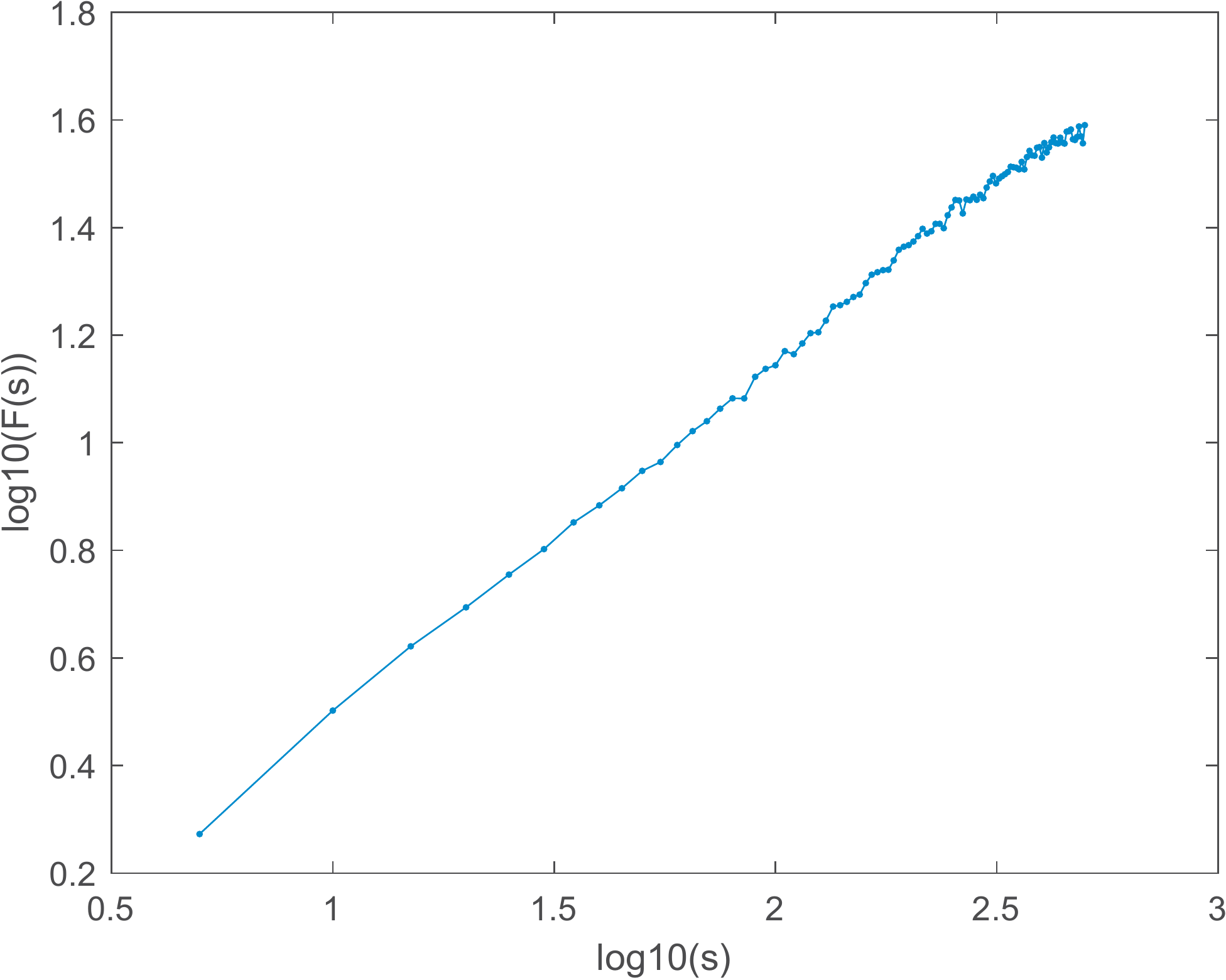}}
    \caption{Precipitation time series considered (original and integrated series ($Y_j$)), and DFA results (log-log plots $F(s)$ vs. $s$); (a) Precipitation time series from A Coru\~na; (b) Precipitation time series from Madrid-Barajas; (c) Integrated time series from A Coru\~na; (d) Integrated time series from Madrid-Barajas;  (e) DFA result for A Coru\~na; (f) DFA result for Madrid-Barajas.}
    \label{fig:x-y_precipitacion}
\end{figure}

Here we present an example of the fractal dimension calculation ($D$) of precipitation time series, following the approach in \cite{Breslin99}. In this case, instead of using the R/S algorithm, the DFA algorithm is considered, which directly provides the Hurst coefficient $H$ by means of the slope of the log-log plots $F(s)$ vs. $s$ graph ($\alpha$). Specifically, we consider two measurement stations in Spain: A Coru\~na, at North-West and Madrid-Barajas, in the center of the country. Daily precipitation data are considered, and long time series are available: A Coru\~na from 01/10/1930 to 31/12/2019 and Madrid-Barajas from 01/01/1951 to 31/12/2019. Figures \ref{fig:x-y_precipitacion} (a)-(d) show the precipitation time series considered, original and integrated series (Equation \eqref{Eq:Y_j2}). Figures \ref{fig:x-y_precipitacion} (e) and (f) show the DFA results (log-log plots $F(s)$ vs. $s$) for the two measuring stations considered. These DFA results suggest a long-term correlation of the time series in both cases. Table \ref{tab:H_precipitation_D} shows the $H$ obtained in the two cases considered from the $\alpha$ coefficient of the DFA, and the corresponding fractal dimension $D$ for them. Note that the considered precipitation time series have a fractal dimension between 1.2 and 1.4, which matches the results obtained in \cite{Breslin99} for Australia.

\begin{table}[]
    \centering
    \begin{tabular}{|c|c|c|}
        \hline
        Station & $H$ (from $\alpha$ in the DFA) & $D$ \\ \hline
         A Coru\~na & $0.77$ & $1.23$ \\ \hline
         Madrid-Barajas & $0.66$ & $1.34$ \\ \hline
    \end{tabular}
    \caption{Hurst coefficient ($H$) and fractal dimension ($D$) of daily precipitation time series in Spain with the DFA method.}
    \label{tab:H_precipitation_D}
\end{table}

\section{Persistence in application domains}\label{Literature_Review}

\begin{figure}[h!]
		\resizebox*{!}{\dimexpr\textheight-5\baselineskip\relax}{\begin{forest}
		        upper style/.style = {draw,top color=white, bottom color=black!20},
                lower style/.style = {draw, thin, align=left},
                where level<=1{%
                    upper style
                }
                {
                    lower style,
                },
                for tree={
				l sep=15em, s sep=1em,
				child anchor=west,
				parent anchor=east,
				grow'=0,
				line width=0.75mm,
				minimum size=1cm,
				anchor=west,
				draw,
			    }
				[{\scalebox{2}{\fontsize{32pt}{0pt}\selectfont Persistence in}}\\{\scalebox{2}{\fontsize{32pt}{0pt}\selectfont application domains}}, for tree={fill=white,top color=white,bottom color=white}
				[{\scalebox{2}{\fontsize{25pt}{0pt}\selectfont Earth and}}\\{\scalebox{2}{\fontsize{32pt}{0pt}\selectfont atmospheric sciences}}, bottom color=blue!20
				[{{\huge Precipitation}} \vspace{3mm}\\
    				{\huge \cite{Pelletier97,Yang19,Katz1977,Mimikou84,Jimoh96,Cazacioc2005,Matalas03,Lennartsson08,guilbert2015characterization,Paschalis14,kumar2013evaluation,Matsoukas00,thyer2000modeling,Markonis2016scale,Dey18b}}\\
    				{\huge \cite{Dey18b,vogel1998regional,sagarika2014evaluating,iliopoulou2019revealing,Wu19,jovanovic2018long,martinez2021modified,zhang2019modified,sarker2021detrended,chakraborty2021exploring}}\\
    			]
				[{{\huge Temperature}} \vspace{3mm}\\
				    {\huge \cite{Bunde02,kumar2013evaluation,bloomfield1992trends,Zekai90,Raha20,Bunde01,rybski2006long, Capparelli11,Weatherhead10,vecchio2010amplitude,Govindan01,blender2003long,Zhu2010demonstration,Mann2011long}}\\
				    {\huge \cite{Vyushin12,moon2018drought,Monetti03,Fraedrich03,Gan07,Breaker19,Zhang15_physa,Luo15,deng2018impact,wang2021unnatural,li2021widespread,lopez2021effect,sarvan2021classification}}\\
				]
				[{{\huge Soil moisture}} \vspace{3mm}\\
				    {\huge \cite{Delworth88,Delworth89,Manabe90,Delworth93,Liu98,Xu01,lorenz2010persistence,Salcedo20,Shen18,ghannam2016persistence,nicolai2016long,sanz2021generalized,tatli2020long,millan2021hurst}}\\
				]
				[{{\huge Climate processes}} \vspace{3mm}\\
    				{\huge \cite{Graves17,MacDonald92,Feldstein00,Tsonis99,keeley2009does,rybski2008long,sreedevi2022spatiotemporal,adarsh2021multifractal,sanchez2021diversity}}
				]
				[{{\huge Fog events}} \vspace{3mm}\\
    				{\huge \cite{Cornejo20,Perez18,Salcedo21}}\\
				]
				[{{\huge Sea and reservoir level}} \vspace{3mm}\\
    				{\huge \cite{Barbosa06,dangendorf2014evidence,agarwal2012trends,wu2018multifractality,Castillo20,stratimirovic2021changes,san2021long}}\\
				]
				[{{\huge Atmospheric pollution}} \vspace{3mm}\\
    				{\huge \cite{Chelani16,Chelani09,vyushin2010statistical,fioletov2003seasonal,Varotsos05,Varotsos06,Varotsos06ACP,kiss2007long,sun2017impact,Liu15Pollution,Lu12}}\\
				]
				[{{\huge Geophysics and seismology}} \vspace{3mm}\\
                    {\huge \cite{Witt13,dmowska1999advances,jimenez2005testing,Chamoli07,lee1999persistence,aggarwal2015multifractal,telesca2016multifractal,fan2017multiscale,kataoka2021detrended,flores2015multifractal,telesca2004fluctuation,varotsos2009detrended,varotsos2011scale}}\\
				]
				]
				[{\scalebox{2}{\fontsize{25pt}{0pt}\selectfont Energy resources}}, bottom color=red!20
				[{{\huge Solar}} \vspace{3mm}\\
					{\huge \cite{Voyant18,Antonanzas16,Pedro12,huertas2019using,Yang19b,Driemel18,Yang19b,Driemel18,Lipperheide15,Harrouni09,Dos15,liu2021use,yelchuri2021short}}\\
    			]
				[{{\huge Wind}} \vspace{3mm}\\
    				{\huge
    			\cite{Kosak08,Shirvaikar72,Poje92,Gadian04,Jiang18,de2021long,santos2019analysis}}\\
				]
		    	]
				[{\scalebox{2}{\fontsize{25pt}{0pt}\selectfont Complex Networks}}, bottom color=green!20
				[{{\huge \emph{On} CN}} \vspace{3mm}\\
				{\huge Links:  \cite{papadopoulos2019link}, }
				{\huge Epidemics:} {\huge\cite{pastor2001epidemic,bohme2013emergence}, }
    			{\huge Message, memes:} {\huge \cite{cui2014message,weng2012competition}}\\
				]
				[{{\huge \emph{Of}  CN}} \vspace{3mm}\\
    			{\huge Links:  \cite{nicosia2013graph,barucca2018disentangling,li2021percolation}, }
    			{\huge Nodes:  \cite{valdano2015predicting}, \vspace{1mm}}
    			{\huge Hubs: \cite{banerjee2020persistence}, \vspace{1mm}}
    			{\huge Communities: \cite{li2018persistent},  \vspace{1mm}}\\
    		    {\huge Person's social signature: \cite{saramaki2014persistence}, \vspace{1mm}}
    			{\huge Patterns: \cite{morse2016persistent,zhang2014dynamic},}\\
    		    {\huge Processes Rumors: \cite{vazquez2013spreading}, \vspace{1mm}}
    		    {\huge Epidemics: \cite{sun2015contrasting}, \vspace{1mm}}
    			{\huge Emotions: \cite{garas2012emotional} \vspace{1mm}}\\
				]
				]
				[{\scalebox{2}{\fontsize{25pt}{0pt}\selectfont Economics and}}\\{\scalebox{2}{\fontsize{25pt}{0pt}\selectfont Market Analysis}}, bottom color=orange!20
				[{{\huge Inflation}} \vspace{3mm}\\
				    {\huge \cite{mcgahan1999persistence,fuhrer1995inflation,Sbordone07,Gaglianone18,pivetta2007persistence,Sbordone07,Meenagh09,Gaglianone18,koenker2001quantile,Tule19,pivetta2007persistence}}
				]
				[{{\huge Expected Returns}} \vspace{3mm}\\
    				{\huge \cite{Priestley01,Dichev09,Frankel09,Dichev09,Wu12,Pla19}}\\
				]
				[{{\huge Exchange Rates and Crypto-currencies}} \vspace{3mm}\\
				    {\huge \cite{Curran19,Caporale18,al2018efficiency,takaishi2021time,costa2019long,vaz2021price,luis2019drivers,quintino2020efficiency,alvarez2018long,stosic2019multifractal,stosic2019multifractal,david2021fractional}}\\
				]
				[{{\huge Stock Market}} \vspace{3mm}\\
    				{\huge \cite{Grau01,Cajueiro05,constantin2005volatility,Oh08,bollerslev1986generalized,baillie1996fractionally,Bentes14,Granero08,Lu13,milocs2020multifractal,stovsic2015multifractal,yin2013modified}}\\				
				]
				[{{\huge Market Analysis}} \vspace{3mm}\\
    				{\huge \cite{mensi2021does,cerqueti2021long,cai2019exploring,delbianco2016multifractal,gu2010multifractal,ftiti2021oil,zhang2021cross,david2020measuring,yin2021market,feng2021multifractal,raza2021multifractal,nejad2021multifractal,gorjao2021change}}\\
    				{\huge \cite{han2021complexity,ali2021modeling,fan2015multifractal}}
				]
				] 
				[{\scalebox{2}{\fontsize{25pt}{0pt}\selectfont Non-equilibrium}}\\{\scalebox{2}{\fontsize{25pt}{0pt}\selectfont Thermodynamics}}, bottom color=cyan!20
				[{{\huge Non-equilibrium systems}} \vspace{3mm}\\
    				{\huge \cite{majumdar1999persistence,bray2013persistence,iyer2015first,aurzada2015persistence,Majumdar01,bray2002theory,villain1984nonequilibrium,Derrida1994non,stauffer1994ising,sire1995coarsening,derrida1995exponents,majumdar1998persistence,tam2002cluster}}
				]
				]
				[{\scalebox{2}{\fontsize{25pt}{0pt}\selectfont Synthetic time series}}, bottom color=magenta!20
				[{{\huge Fourier filtering}} \vspace{3mm}\\
    				{\huge \cite{peng1991directed,prakash1992structural,makse1996method,eichner2006extreme,halley2009using}}\\
				]
				[{{\huge Hydrology}} \vspace{3mm}\\
                    {\huge \cite{Efstratiadis14,boughton2003continuous,boughton2003continuous,Ilich14,Castalia2,Castalia,koutsoyiannis2000generalized}}
				]
				[{{\huge Wind}} \vspace{3mm}\\
                    {\huge \cite{tsekouras2014stochastic}}\\
				]
				]
				[{\scalebox{2}{\fontsize{25pt}{0pt}\selectfont Optimization and planning}}, bottom color=yellow!20
				[{{\huge Robust Optimization}} \vspace{3mm}\\
                    {\huge \cite{brown1997optimization,ben2009robust,gabrel2014recent,kalsi2001comprehensive,allen2006robust,morrison2010new,petit2019enriching}}\\
				]
				[{{\huge Planning}} \vspace{3mm}\\
                    {\huge \cite{brown1996scheduling,brown1997optimizing,Borthen2019}}\\
                ]
				]
				[{\scalebox{2}{\fontsize{25pt}{0pt}\selectfont Health and biomedical sciences}}, bottom color=gray!20
				[{{\huge Biomedicine}} \vspace{3mm}\\
                    {\huge \cite{Depetrillo99,de1999long,de1999persistence,thurner2003scaling,rahmani2018dynamical}}\\
				]
				[{{\huge Sport}} \vspace{3mm}\\
                    {\huge \cite{Gomez20}}\\
				]
				] 
				] 
		\end{forest}}   
		\centering
		\caption{Application domains of persistence discussed in this paper.}
		\label{fig:taxonomy}
	\end{figure}
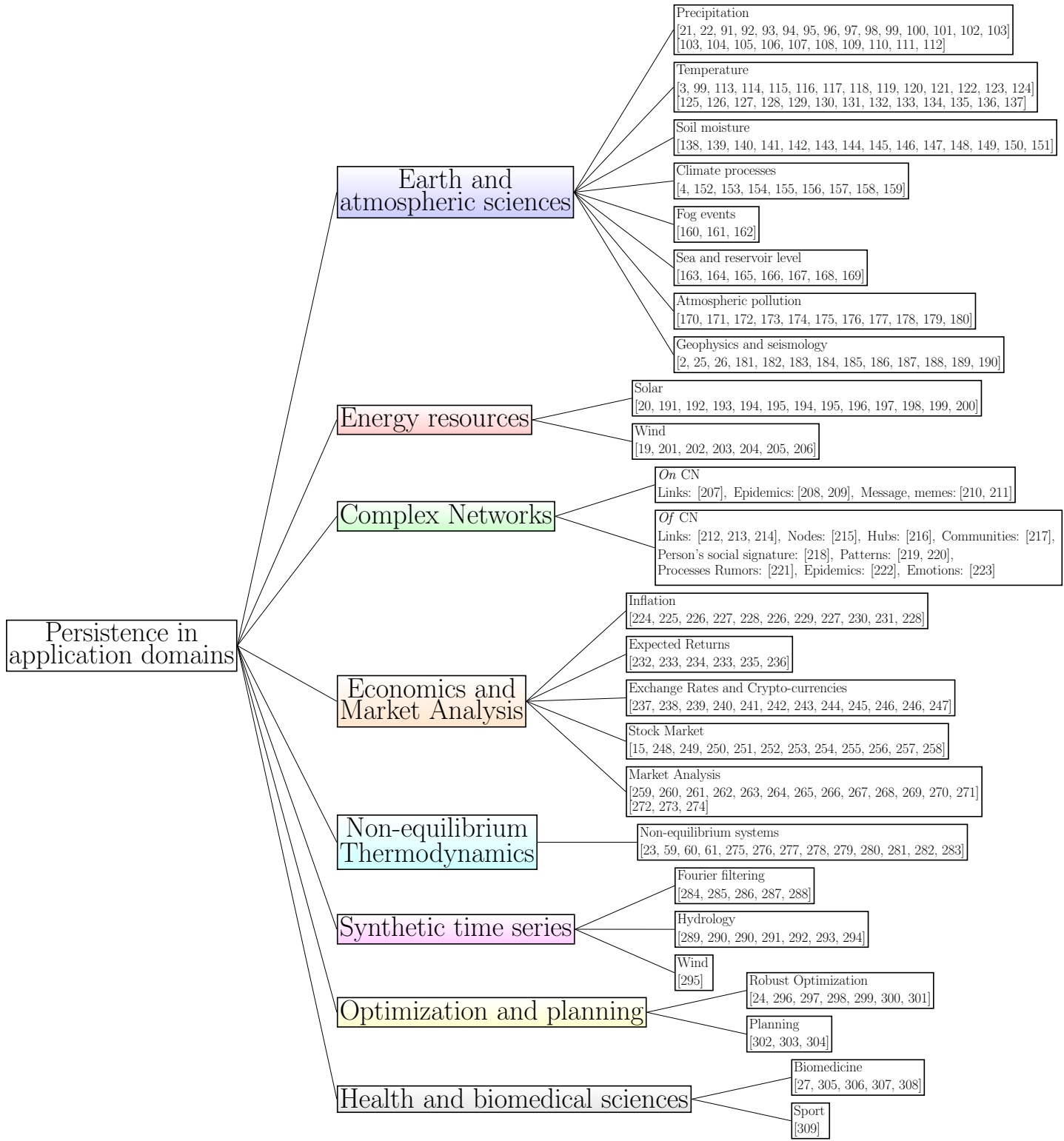

Persistence and related concepts are ubiquitous in a wide range of application domains. We here review their use in paradigmatic research areas where persistence analysis has been considered.
Persistence studies of Earth, atmospheric and climatic systems has been traditionally carried out, and it is therefore the research area where persistence has been studied the most. Closely connected to atmospheric processes, the persistence of renewable energy resources has also been important in the last years. 
We then revise different persistence studies carried out in Complex Networks, a multidisciplinary field to study complex systems that consist of many elements interacting with each other. We also review persistence works in economics, a field with a long tradition on studies in persistence. 
Non-equilibrium thermodynamics systems is another research field where persistence studies have been important in the last years. Many of the previous research areas have relied on synthetic time series; the generating processes to account for different persistence properties are also revised later. Persistence applications in optimization and planning is following taken into account, an interesting and novel way of looking to persistence of systems in engineering. We end up this section by discussing some persistence applications in health and biomedical applications. A visual summary of the application domains of persistence discussed in this paper can be seen in Figure \ref{fig:taxonomy}.

\subsection{Persistence studies in Earth and atmospheric sciences}

Persistence of weather and climate phenomena has been well studied using different techniques since the 1950 \cite{Graves17}. In recent years, different works have tackled the persistence analysis of weather (short-term) and climate (long-term) regimes. The typical time scale for weather phenomena depends a lot on the specific event studied, but may vary in short-term studies from hours in the local and meso-scale scales (fog and mist formation, storms, etc.) to days, until approximately one week (general weather regimes in synoptic-scale phenomena). Long-term persistence involves synoptic regimes (up to time scales of several weeks), usually associated with {\em blocking} phenomena (very stable high pressure over a region, which produces persistent weather on that region). Even longer persistence phenomena exist over many months or seasons, years, and even over several decades \cite{MacDonald92}. These are usually associated with external forcing and climate-related events such as El Ni\~no Southern Oscilation (ENSO phenomenon) and other decadal and multidecadal variations of the climate system, usually characterized as climate indices. Persistence studies in geophysics have been less frequent in the literature, but they also have importance in fields such as geomagnetism or seismology. In this section, we review relevant works dealing with persistence, both at short-term and long-term scales, in Earth and atmospheric science, including meteorology, climate and also geophysics.

\subsubsection{Short-term persistence in atmospheric processes}\label{Persistence_atmospheric_Processes}

There is an important number of works which study short-term persistence in atmospheric processes. The majority of these works apply Markov chain models to model phenomena such as precipitation occurrence, drought and heatwave persistence, soil moisture and streamflow, pressure anomalies analysis or fog occurrence.

Markov chain models were one of the first persistence-based approaches to model precipitation occurrence with probabilistic approaches. Early works applied Markov chains for this purpose. For example \cite{Katz1977}, where a Markov chain model was used to model daily precipitation occurrence in Pennsylvania, USA. Other early work on Markov Chain analysis of precipitation is \cite{Mimikou84}, where daily precipitation occurrences were modelled with a Markov chain in Greece. The order of the chain depends on the season, the station's meteorological conditions and the geographical location of the measuring station. In \cite{Jimoh96} Markov models are applied to model daily precipitation occurrence in Nigeria. This work mainly discusses different statistical criteria to set the Markov chain model order, such as the Akaike or the Bayesian Information Criteria (AIC and BIC, respectively). A similar work is carried out for Romania in  \cite{Cazacioc2005}, also using Markov Chain models. In \cite{Matalas03} simulated hydrologic sequences (as realization of Markov process, ARMA(1,1) or fractional Gaussian noise) are used to study the trend in this type of series, via estimation of the regression coefficient. In \cite{Lennartsson08} a Markov chain model is applied to model the probability of occurrence of precipitation at a weather station in Sweden. This approach is combined with a Gaussian process to model the amount of precipitation collected when precipitation at the station occurs. In \cite{guilbert2015characterization} the short-term persistence of precipitation is studied by using Markov chain models transition probabilities. This work showed an increasing persistence in daily precipitation in the northeastern United States in the last years, which suggests that global circulation changes are affecting regional precipitation patterns. In \cite{Paschalis14} the performance of composite stochastic models for precipitation occurrence is discussed. Different models are considered to model precipitation over several time scales, from daily (with Markov models) to minutes, using other models such as Multiplicative Random Cascade models. This approach can satisfactorily reproduce precipitation properties across a wide range of temporal scales.

Global temperature, droughts and heatwaves analysis have also been tackled with persistence-based methods. In \cite{bloomfield1992trends} ARMA models are applied to study the global temperature variation of the last century. In \cite{Zekai90} it is shown that a second-order Markov chain is an efficient model for probability drought generation mechanism, which was tested over three representative annual flow series and from different parts of the world. In \cite{Raha20} heatwave duration is studied by means of direct analysis of real data in the USA, in which calculation of several persistence measurements (see Section \ref{short-term-definitions}) is carried out with probabilistic inference models.

In \cite{Delworth88}, the authors used a first-order Markov model together with simulations produced with the Geophysical Fluid Dynamics Laboratory (GFDL) General Circulation Model (GCM), to study the timescale of decay of soil moisture disturbances. A significant persistence in soil moisture disturbances, especially at high latitudes, and during the winter season, was found in this study. Further studies on soil moisture persistence, decay and its effect in climate and atmospheric variability have been studied with similar  \cite{Delworth89,Manabe90,Delworth93}, and more advanced land-atmosphere analytical model \cite{Liu98}. Different timescales for soil moisture (persistence characteristics) are obtained in this latter work, namely seasonal, monthly, weekly, and daily, characterized by different hydrological and thermal processes.
In close relation with soil moisture, streamflow modeling has been also tackled with Markov chains, such as in \cite{Xu01}, where a stochastic model for the description of daily streamflow at multiple sites of Wuper river (Germany) is presented. It is an extension of a Markov chain process commonly used to model daily precipitation series.

In \cite{Feldstein00} the timescales of prominent pressure low-frequency anomalies (the North Atlantic oscillation (NAO), the Pacific-North American (PNA), the West Pacific (WP) oscillation and the ENSO oscilation) are studied. It is shown that several of the prominent low-frequency anomalies are well described as a Markov process (AR1 process). In fact, the temporal evolution of the NAO, PNA, and WP anomalies can be interpreted as being a stochastic (Markov) process with an e-folding timescale between 7.4 and 9.5 days. The time series corresponding to the ENSO spatial pattern did not match that of a Markov process, and thus a well-defined timescale could not be specified.

Finally, in \cite{Cornejo20} Markov Chain models with different order are applied to the probability of fog occurrence at Valladolid airport, Spain. The Markov model are then hybridized with different machine learning regressors in order to form a Mixture of Expert (MoE) algorithms to improve the prediction of occurrence in the zone. This prediction approach had been previously discussed for low-visibility records \cite{Perez18}.

\subsubsection{Long-term persistence in atmospheric and climatic processes}

Long-term persistence analysis of atmospheric and climatic processes mainly involve the analysis of power-laws by means of the application of the DFA algorithm or correlation-based approaches. Long-term persistence studies have been applied to different time series, mainly air temperature, sea surface temperature anomalies, atmospheric  pollutants and hydrology-based processes.

One of the first works discussing long-term persistence in atmospheric circulation was \cite{Tsonis99}. In that work, 500-hPa measurements were used to show the existence of scale invariance in the variability of extra-tropical atmospheric circulation anomalies, over the
whole range of timescales resolved by the available data (from weeks to decades). In \cite{keeley2009does} the persistence of the North Atlantic Oscillation (NAO) index on intra-seasonal
time-scales (summer and winter) is studied by means of the  autocorrelation of the NAO time series. The estimation of the e-folding time for the NAO obtained in this work is approximately of 5 days. In \cite{Bunde01,Bunde02} it is shown that there is a universal power-law in correlation of temperature variations (long-term persistence in temperature variations), found in many different atmospheric temperature records, and also for ocean temperature records \cite{Bunde02}. In these works, it is shown that temperature fluctuation follows a universal correlation power-law with exponent 0.7 for continental measuring stations and 0.4 for ocean stations. Also dealing with temperature records, in \cite{rybski2006long} a DFA is applied to evaluate six different Northern Hemisphere atmospheric temperature datasets, finding
that all are governed by long-term persistence. In \cite{Capparelli11} the spatial distribution of the DFA exponent is used to quantify the degree of memory of surface temperature data in the USA. This study shows that the long-term persistence is related to the presence of several climatic regions, sensitive to climatic phenomena such as the ENSO. In \cite{Weatherhead10} the temperature persistence variation in the Canadian Arctic is analyzed, by using one day autocorrelation measurements from air temperature time series at Baker Lake, Nunavut, Canada. In \cite{vecchio2010amplitude} a DFA approach is applied to study the long-term persistence of temperature anomaly from two historical temperature records lasting about 250 years measured at Milan and Prague stations. A long-term persistence in these temperature records for time scales between 3 and 10 years are detected in this work. In \cite{deng2018impact} the impact of inter-annual variability of annual cycle on long-term persistence of surface air temperature is studied by means of a DFA approach. Long-length air temperature series over 150 years are considered in Italy, Germany, Czech Republic, Sweden, Austrian and Croatia. 
In \cite{wang2021unnatural} a DFA approach is applied to global land surface air temperature in the period (1951-2015) to separate unnatural trends in the series. The work shows discriminated evidence that observed air surface temperature change exceeds its natural internal variability, pointing to a unnatural trend warming process, specially strong in some specific areas of the world. In \cite{lopez2021effect} the performance of DFA on air temperature data at Laguna del Maule (Chile) is analyzed in series with missing data, with and without applying imputation techniques. The results showed that the DFA technique is highly  robust when it is applied to a short time series with missing data and without data imputation. In \cite{sarvan2021classification} DFA, wavelet and Hurst analysis methods are applied to characterized the long-term persistence of temperature anomalies records in the UK, from long-length gridded data in the UK.
In \cite{li2021widespread} the changes in the persistence of temperature events under climate change are studied. The persistence of the series are evaluated by means of 10-day autocorrelation of daily mean temperatures. The results obtained indicate that climate change simulations are marked by important changes in surface temperature persistence.

The persistence of Sea Surface Temperature Anomalies (SSTA) has been analyzed in different works. In \cite{Gan07} a DFA analysis is applied to study the long-term persistence of SSTA in the South China Sea. The analysis is carried out on several time series from different measuring stations, obtaining scaling coefficients over 0.8 in all them. More recently, in \cite{Breaker19} a DFA analysis of SSTA in the west coast of California was presented. In \cite{Zhang15_physa} an assymmetric DFA algorithms is applied to study sea surface temperature anomaly. Global distributions of the scaling exponents are shown in order to quantify the degree of memory of the system, on time scales in the range from 10 months to 20 years. In \cite{Luo15} a DFA analysis showed a clear two-ranges pattern in SSTA, with a characteristic time (crossover point) around 52 months, which seems to be modulated by the ENSO oscillation. This behaviour is shown to be similar in most ocean basins of the Earth.

Some other works have studied the persistence associated with variables from climate models. For example in \cite{Govindan01} the trends and temporal correlations in the monthly mean temperature data of Prague and Melbourne from four state-of-the-art general circulation models is studied. For this, DFA and correlation-based approaches are used. In \cite{blender2003long} the analysis of long-term correlation in surface temperature from the NCEP re-analyses simulations is carried out. An analysis based on spectral power density of the time series is carried out, showing values of the persistence strength $\beta=-1$ over the oceans, $\beta=0$ in the inner continents, and $\beta=0.3$ in the coastal areas. In \cite{rybski2008long} the appearance of long-term persistence in simulated temperature records is analyzed. Specifically the 
the global coupled general circulation model ECHO-G was considered to simulate the temperature records, with a historical simulation of the last 1000 years. The DFA method is applied to analyze the long-term persistence of this simulated record. In \cite{Zhu2010demonstration} the existence of long-term correlation in the surface temperature field was shown. The work applies a DFA algorithm over grid points of a millennium control simulation, from a state-of-the-art global circulation model. The results obtained shown a clear long-term persistence in surface temperature records in high-latitude oceans, which is most pronounced in the North Atlantic and North Pacific with a DFA scaling parameter close to 1. In \cite{Mann2011long} the long-term persistence of surface temperatures is analyzed by using the Hurst coefficient $H$ as a persistence measurement. This analysis is carried out in synthetic temperature time series from a simple model, and compared with observational records. In \cite{lorenz2010persistence} the role of soil moisture in heat wave persistence is investigated. The CLM RCM regional model (a climate version of the COSMO-model \cite{Salcedo20}) is used to perform different simulations for this study. The results obtained show that simulations in which soil moisture is fixed to a constant value or prescribed seasonal cycle, present a lower intrinsic heat wave persistence than simulations in which interactive soil moisture is present. In \cite{Vyushin12} the performance of persistence-based models AR1 and power-law to model climate variability is studied. Specifically, both persistence-based techniques are fitted to observed surface air temperature at different spatial points. The parameters estimated from observations are found in general to be well captured in the Coupled Model Inter-comparison Project 3 (CMIP3) simulations. It is shown that AR1 representation captures the tendency for the variability of some climatic time series  annual or multi-annual variability and the power law representation is appropriate for very long-term instrumental i.e. decadal or multi-decadal timescales. In \cite{kumar2013evaluation} the Hurst exponent is used in order to study the long-term persistence of the twentieth-century temperature and precipitation from 19 climate models of the Coupled Model Inter-comparison Project. The results show that the studied models capture the long-term persistence in temperature
reasonably well, observing an areal coverage of observed long-term persistence of 78\%. The areal coverage of observed long-term persistence in precipitation is smaller, about 60\%. In \cite{moon2018drought} the persistence of drought events and how good they are simulated by Global Circulation Models is studied. A methodology based on dry-to-dry transition probabilities at monthly and annual scales is considered. In \cite{Monetti03} the scaling effect of SSTA is studied for the Atlantic and Pacific Oceans. The study is carried out by applying a DFA algorithm, which reveals two different long-term persistence regimes with a characteristic time of 10 months. In \cite{Fraedrich03} a DFA approach is applied in SSTA from the coupled atmosphere-ocean model ECHAM4/HOPE, at global level (coarse grained resolution). The results show a clear long-term persistence of this variable all around the globe, with scaling factors between 0.6 and 0.7. Finally, in \cite{sanchez2021diversity} global climate model simulations of the last millennium, and paleoclimate archives have been used to study the variability of the persistence in El Ni\~no events, showing a wide range of variability in the frequency of such events. 

Persistence in sea level records has been also tackled in the literature. In \cite{Barbosa06} long-term persistence in sea level records from the North Atlantic ocean is studied by means of the spectral density function scaling exponent $\beta$. In this case parameter $\beta$ is calculated from the slope of the wavelet spectrum, as shown in Section \ref{Wavelet}. In \cite{dangendorf2014evidence} a study on long-term memory in sea level was carried out by applying a DFA approach. A set of 60 centennial tide gauge records together with ocean reanalysis were used. The DFA approach showed that sea level records
exhibit long-term correlations on time scales up to several decades.
In close connection with sea level studies, the persistence of the arctic sea ice is studied in \cite{agarwal2012trends}. A multi-fractal temporally weighted DFA approach is proposed. The results obtained show a clear long-term persistence of the time series, resulting in a  DFA log-log plot with different slopes (different scaling factors $\alpha$), at time scales ranging from days to several years. 

Other long-term persistence analysis of atmospheric processes involve the study of air pollution indices and pollutants concentration \cite{Chelani16}. In \cite{Chelani09} the persistence of hourly ground level ozone concentration is studied. The work deals with observed hourly tropospheric ozone measurements during 2006 at two sites in
Delhi, applying a rescaled range analysis (R/S method) based on Hurst exponent estimation (see Section \ref{H_estim}). Persistence is observed up to 120 hours ($~5$ days) which corresponds to weekend and week day cycles, pointing out the significance of traffic activities in controlling the correlations in the ground ozone time series. In \cite{vyushin2010statistical} long-term persistence of total ozone anomalies are studied by calculating the Hurst coefficient of several total ozone times series (monthly means) from satellite data (version 8 of TOMS, SBUV and OMI). In \cite{fioletov2003seasonal} temporal autocorrelations of monthly mean total ozone anomalies are analyzed, revealing that anomalies established in the wintertime are persistent until the end of the following autumn. In \cite{Varotsos05} a DFA algorithm is applied hourly air-pollution time-series in Athens (Greece), specifically times series of ozone, nitrogen oxides, and atmospheric particles, from five air-pollution monitoring stations during the period 1987--2003. The work also applies DFA to analyze PM2.5 particle fluctuations in a 6-months data set collected at the University of Maryland in East Baltimore, USA. Persistent power-law correlations were found for ozone concentrations with lag times ranging from 1 week to 5 years. Also, the fluctuations of nitrogen oxide concentrations exhibit a similar behavior. Finally, persistent power-law correlations from about 4 h to 9 months were found in the case of atmospheric particles. In \cite{Varotsos06} a DFA approach is applied to a time series of zonal mean daily Aerosol Index values derived from satellite observations during 1979-2003 period. The results obtained showed that the aerosol index globally obey  persistent long-range power-law correlations for time scales longer than 4 days and shorter than 2 years. In \cite{Varotsos06ACP} a DFA approach is applied to series of Total Ozone Content (TOC) in the atmosphere. The results showed a clear long-term persistence pattern, with two ranges and characteristic time around 28 months. In \cite{kiss2007long} the analysis of global TOC persistence in the atmosphere is carried out by applying a DFA approach to time series from Total Ozone Mapping Spectrometers, on board of Nimbus-7 and EarthProbe NASA satellites. The results obtained showed a clear long-term persistence of TOC with consistent values of scaling coefficient $\alpha$ around 0.7 all around the globe. Recently, a study in \cite{sun2017impact} has investigated how ozone distribution is determined by the persistence of different meteorological conditions. In \cite{Liu15Pollution}, the long-term persistence of air pollutants is studied, focused on a zone with intense air pollution such as Shanghai, China. Analysis with DFA show a temporal scaling following a power-law, with a clear two-range pattern for this phenomenon. The application of DFA to different meteorological variables' volatility is shown in \cite{Lu12}. Specifically, the volatility (magnitude of the increments between successive elements) of daily mean relative humidity, daily mean temperature, daily temperature range, daily maximum and minimum temperature is analyzed. Universal scaling behaviors are found for all the variables considered, whose scaling exponents take similar distributions, mean values and standard deviations. 

In close connection with atmospheric processes, problems related to hydrologic analysis and their long-term persistence have also been discussed in the literature. In \cite{Pelletier97} an analysis of long-term persistence in hydrologic-related time series such as tree-rings, precipitation, rivers' discharge and temperature, considering hundreds of stations around the world, is carried out. The long-term persistence analysis is carried out by means of power-spectrum of each time series, and proving that they follow a power-law in all cases analyzed. In \cite{Matsoukas00} a DFA analysis is used to compare the long-term persistence of precipitation and runoff at different measurement stations of Kansas, USA. This is one of the first works which detects a two-ranges structure in the DFA log-log plot, and associates it with the binary time series structure. In \cite{thyer2000modeling} a hidden Markov model approach is proposed to evaluate long-term precipitation persistence in Australia, showing that this approach is able to closely simulate the influence of the global climatic mechanisms in precipitation. The work in \cite{Markonis2016scale} discusses, in terms of the Hurst coefficient, the differences between persistence in precipitation records from proxies (palaeoclimatic reconstructions) against a period of instrumental records. The main conclusion raised in the study is that instrumental records are too short to be used for inferring long-term properties of the rainfall variability, however, proxy data usually lead to the underestimation of the persistence structure regardless of proxy type, or to overestimation due to low resolution of the time series. In \cite{Yang19} a DFA analysis of rainfall time series over USA is carried out. The analysis showed a clear two-ranges structure in the log-log DFA plot, with characteristic time around 160hours. In \cite{martinez2021modified} a modified MF-DFA analysis with Polynomial and Trigonometric functions is applied to study the long-term persistence of daily precipitation time series in meteorological stations of Tabasco, M\'exico. A similar analysis is carried out in \cite{zhang2019modified}, where the polynomial fitting local trend in traditional MF-DFA has been replaced by an ensemble empirical mode decomposition algorithm to evaluate the local trends in the MF-DFA. This modified approach is applied to evaluate the long-term persistence of precipitation time series in Dongting Lake Basin, China. Also, Indian Monsoon precipitation persistence has been analyzed very recently with DFA \cite{sarker2021detrended} and MF-DFA algorithms \cite{chakraborty2021exploring}.

In \cite{Dey18b} the spatio-temporal persistence of rainfall and streamflow and their joint behaviour is examined by applying a DFA algorithm, for different watersheds in USA and Cauvery river basin in India. In \cite{vogel1998regional} an analysis of the long-term persistence of streamflow (by means of the estimation of the Hurst coefficient) is carried out in the USA, using data from over 1500 measuring gauging stations. In \cite{sagarika2014evaluating} a similar study on the long-term persistence of streamflow time series in USA is carried out by means of applying Mann-Kendall tests to data of 240 streamflow measuring stations over the USA. In \cite{kantelhardt2006long} a DFA and Multi-Fractal DFA algorithms are applied to study the long-term persistence of river runoff records. The study is focused on discussing the multi-fractal temporal scaling properties of precipitation and river discharge records on large timescales. In \cite{iliopoulou2019revealing} the links between clustering of rainfall extremes and long-term persistence (characterized by the Hurst coefficient) is studied. In \cite{Wu19} a DFA and MF-DFA algorithms have been applied to analyze the long-term persistence of river runoff fluctuations, using data of 12 mayor rivers in China. A two-ranges structure is again found in the DFA log-log plot, showing a strong long-term persistence, more accused for river runoff than for precipitation in all rivers considered. In \cite{jovanovic2018long} the Hurst coefficient is used as a measurement of urbanisation, in turn associated with a reduction in the long-term correlation within a streamflow series. Specifically, R/S and MF-DFA are applied, and the work showed that the Hurst exponent $H$ is a useful metric to assess the impact of catchment imperviousness on streamflow regimes. In \cite{wu2018multifractality} a MF-DFA algorithm has been applied to analyze the persistence of streamflow and sediment in the Pearl river delta, Southern China. The results indicate that there was a significant multifractal structure present in the fluctuations of streamflow and sediment in the time series records measured in the Pearl river delta.

Regarding soil moisture, there have been a number of works dealing with persistence analysis in soil moisture time series. In \cite{Shen18}, a DFA is applied to analyze soil moisture over the zone of Babao River Basin (China). The results obtained showed a pattern with three ranges of behaviour, from strong long-term correlation to anti-persistence, with characteristic times varying between 0-19 hours, 19-160 hours and over 160 hours approximately, depending on the depth considered. In \cite{ghannam2016persistence} the differences, underlying dynamics, and relative importance of persistence timescales in root-zone soil moisture are studied. Persistence characteristic below or above some thresholds (characteristic time) are explored in this work for soil moisture. These persistence scales are more indicative of the wet and dry states of soil moisture and are some of the principal measure of land-atmosphere coupling strength.
In \cite{nicolai2016long} a study on the predictability of soil moisture based on persistence is carried out. Linear regression models taken into account persistence only or persistence plus climate indices are considered on the long-term predictability of soil moisture dynamics, with prediction time-horizon from 1 to 5 months. In \cite{san2021long} the soil water content persistence in tillage land with different cover crops (vineyards) are analyzed by means of a MF-DFA method in Mediterranean Spain. The results show a multifractal behaviour, compatible with complex patterns of long-range correlations for the time series of soil water content. In close connection with soil moisture analysis, there are works which have focused on the analysis of long-term persistence properties of vegetation index time series (NDVI), such as \cite{sanz2021generalized}, where a MF-DFA algorithm has been applied to time series of NDVI in arid zones of Southern Spain. 

The long-term persistence of droughts has also been studied with DFA and related approaches. In \cite{tatli2020long} the time-domain characteristics of drought persistence over Turkey is analyzed by applying a DFA algorithm to a time series of the Palmer drought severity index. In \cite{millan2021hurst} a R/S analysis has been applied to obtain Hurst scaling regimes in time series of Standardized Precipitation-Evapotranspiration Index, in order to analyze drought events in data from two states of Brazil. The R/S analysis has also been applied in \cite{sreedevi2022spatiotemporal} over Standard Precipitation Index data to evaluate the long-term persistence of drought in India, by obtaining the Hurst coefficient of the associated time series. Another recent work \cite{adarsh2021multifractal} also analyze drought events in India (Western India in this case), by applying a MF-DFA algorithm to Standard Precipitation Index records.

Finally, in some recent works the analysis of persistence of water level in dammed rivers has been analyzed. In \cite{Castillo20} a DFA analysis of water level at Belesar dam, Galicia, Spain is carried out. The study shows a clear two-ranges pattern with 1 year characteristic time. A persistence-based prediction approach based on ARMA methods and on typical year construction is also proposed in that work. In \cite{stratimirovic2021changes} DFA and wavelets approaches have been applied to investigate the persistence of Danube river level in the vicinity of dams at regrio of Djerdap/Iron Gates (Romania and Serbia). The work reports a partial loss of annual cycles in the upstream stations' data due to the operation of the dammed water reservoir at the zone.

\subsubsection{Persistence in geophysics}

Long-term persistence of geophysical time series have been discussed in the literature. In \cite{dmowska1999advances} several long-term persistence methods are presented for geophysical time series. Applications in climate variability, sedimentation and Earth magnetic field variability are discussed. Similarly, \cite{Chamoli07} discusses several methods to measure the long-term persistence of geophysical series, such as wavelet transform, semivariogram, or the R/S method. In this case these methods are applied to geophysical data of the Bay of Bengal. The gravity, magnetic and bathymetry data indicate the self-affine nature with $H=0.8$, 0.8 and 0.9, respectively, indicating long-term persistence in all cases. In \cite{telesca2004fluctuation} a MF-DFA method is applied to analyze geoelectrical data at measuring station of Giuliano, Italia, a zone of seismic activity. The results show a generalized Hurst coefficient around 0.5, with a slight $q$-dependence. In \cite{varotsos2009detrended} the DFA method is applied to geomagnetic and geoelectric series preceding ruptures. In \cite{Witt13} the results of the application of different long-term persistence methods to geophysical records is reported. First, palaeotemperature time series based on GISP2 bi-decadal oxygen isotopes data for the last $10,000$ years was considered. Methods such as R/S, semivariogram, DFA and power spectral density showed a strong power-law dependence on the segment lengths, lags, and frequencies, which in turn points out to a strong long-term persistence. The second time series considered was the daily discharge of Elkhorn River in Nebraska, USA for the period 1929-2001. Methods such as R/S and DFA showed a two-ranges pattern of persistence with characteristic time of 1 year. Finally, the geomagnetic auroral electrojet index data, sampled per minute for 01 February 1978 and differenced is considered. The different methods to characterize the long-term persistence showed a power-law scaling, with a strong long-range persistence in this time series.

Persistence in seismology has been also tackled in the literature. In \cite{lee1999persistence} a theoretical study on short-term persistence of earthquakes occurrence in a quasistatic model of two parallel heterogeneous faults is carried out. In \cite{jimenez2005testing} the long-term persistence of a seismic catalog in the Iberian Peninsula is investigated by means of the R/S method. A strong long-term correlation is found in the region. The study also shows that small earthquakes are very important to the stress transfers. Some works have applied the MF-DFA approach to discuss the persistence of magnitude series of earthquakes (shallow and deep) at different locations such as western India \cite{aggarwal2015multifractal} Pannonia region (Hungary) \cite{telesca2016multifractal}, California \cite{varotsos2011scale,fan2017multiscale}, Mexico  \cite{flores2015multifractal} or Japan \cite{kataoka2021detrended}.

\subsection{Persistence of renewable energy resources}\label{Persistence_Ren_Enery}

Among renewable energy resources, wind and solar energy are currently the most important ones, and those that are already playing a key role in the energetic mix of developed countries, and pushing the developing of the rest \cite{Ren21}. Maybe the most important issue with renewable sources, specially wind and solar resources, is their intrinsic intermittency, which prevents their further penetration in the energetic mix. This is why prediction systems are key in order to ensure a correct management of the whole electric system when renewable resources are taken intro account \cite{Costa08R,Yang18P}. Wind and solar prediction systems are mainly used at short-term, but there are also a great interest in long-term prediction of these resources. In all the cases, the study of wind and solar persistence is key in order to improve short-term and long-term prediction of these renewable energy resources \cite{Antonanzas16}.

Regarding short-term persistence analysis, it has received much more attention in solar energy systems, due to the solar resource characteristics. Note that solar radiation received at a given point of the Earth's surface has an intrinsic persistence which can be calculated only by astronomic parameters (latitude and longitude) of the point of interest and exact hour and day of the day). This receives the name of {\em clear sky} model, and it represents the amount of solar radiation that would be received in the point if there was no atmosphere on Earth. Of course, atmosphere, clouds and aerosols would affect the amount of solar radiation received, and prediction models must incorporate these points, usually taken into account the intrinsic persistence of the resource, given by the clear sky model. Thus, note that there are efficient short-term prediction methods mainly based on {\em smart persistence} techniques, which include modifications on the clear sky models. For instance, in \cite{Antonanzas16} a review of the most important smart persistence techniques are given. In \cite{Pedro12} different prediction methods without exogenous atmospheric variables are proposed. Persistence-based techniques, ARMA models and neural networks are some of the methods proposed and tested in solar radiation records at California, USA. In \cite{Voyant18} two versions of na\"ive predictors based on persistence are proposed. One of them is based on an additive and the other one a multiplicative scheme from the clear sky and the previous values of solar radiation. These schemes have been tested on solar radiation measurements in France and Greece. 
In \cite{FLIESS2018519} a solar prediction method exclusively based on an additive decomposition of a time series into its mean, or trend, and quick fluctuations around it is proposed. A clear sky model (solar persistence) is also used to improve the prediction in this approach, without including any exogenous prediction variable. In \cite{huertas2019using} a prediction approach based on smart solar radiation persistence is combined with a machine learning method (Random Forest) and past production measurements in a problem of solar power production prediction. Experiments on the performance of this approach in three years of data from six solar PV modules at Faro (Portugal) have been carried out. The work in \cite{Yang19b} also highlights this idea of short-term persistence based approaches to solar radiation prediction, by proposing a combination of persistence-based techniques with other based on climatology measurements, for solar radiation prediction without considering exogenous atmospheric variables. Good results are reported on solar radiation prediction at 32 stations of the Baseline Surface Radiation Network \cite{Driemel18}. Very recently, some works have dealt with physics-based smart persistence for solar radiation prediction, such as \cite{liu2021use}, where physics-aware persistence models have been used to obtain accurate short-term prediction of Global Horizontal Irradiance, Direct Normal Irradiance and Difuse Horizontal Irradiance in a measuring station at the USA. In \cite{yelchuri2021short} smart and physics-based smart persistence models have been successfully applied to a problem of very short-term (intra-hour) solar radiation prediction at 15 Solar radiation resource assessment stations located at different parts of India.  Other works include clouds persistence to improve the capabilities of a nowcasting prediction system, as in \cite{Lipperheide15}, where the clouds persistence are taken into account and included in the prediction system. Results are obtained in a solar power plant at Henderson, NV, USA. 

Studies involving long-term persistence have dealt with solar energy but also with wind energy. In \cite{Harrouni09} a new method to measure long-term persistence based on fractal dimension calculation. The method is applied to evaluate the solar radiation time series persistence in Panam\'a. In \cite{Dos15} the long-term persistence of both solar radiation and wind is studied at the Fernando de Noronha Island, Brazil, by applying correlation-based techniques. 
Long-term persistence of wind time series has been evaluated in different studies. There is an early work which deals with wind speed direction persistence \cite{Shirvaikar72}, where the directional persistence of the wind in three nuclear power station sites in India is discussed. More recently works deal with persistence studies of wind speed in Croatia \cite{Poje92}, or the work in \cite{Kosak08} where several ways and methods for evaluating wind speed persistence are presented, including autocorrelation and probability measures. Results are reported for wind persistence at different stations in Turkey. In \cite{Gadian04} the directional persistence of low wind speed observations is studied, in this case with possible application on micro-wind energy but also on air pollution dispersion and natural ventilation of buildings. Results on real data for Scotland are reported. In \cite{Jiang18} the long-term persistence of mean wind over China is studied by applying a DFA approach. Recently, two works have applied the DFA technique to discuss the long-term persistence of wind speed in different regions of Brazil \cite{de2021long,santos2019analysis}, showing correlations with power laws in two different scales (subdiffusive and persistent flows).

\subsection{Persistence in complex networks}

Complex networks (CN) \cite{newman2018networks} have emerged as a multidisciplinary field for studying complex systems that consist of many elements interacting with each other \cite{barabasi2016network}. The field of CN covers a lot of diverse systems, ranging from artificial systems (such as the world-wide web \cite{albert1999diameter}, power grids \cite{cuadra2015critical,cuadra2017optimizing}, or  social networks \cite{strogatz2001exploring}) to natural systems (e.g. metabolic \cite{jeong2000large}, ecological \cite{montoya2006ecological} or climate  \cite{donges2009complex,fan2021statistical} networks). 
More examples can be found in \cite{newman2018networks,barabasi2016network,boccaletti2006complex} and the references therein. 

Most of these systems exhibiting such a different nature have in common that they can be described in terms of  graphs \cite{albert2002statistical,strogatz2001exploring}. A graph is a set of entities called ``nodes'' (or vertices) that are connected to each other by means of ``links'' (edges) \cite{chung1997spectral}. A node represents a system element (for instance, a user of a social network), which is linked with others via a relationship (e.g. who follows who in the social network) or by the exchange of information (e.g. in a communication network \cite{nicosia2013graph}) or energy (e.g. in a power grid \cite{cuadra2015critical}). This could be argued as a reductionist approach  \cite{bohme2013emergence}. However, the methodology helps capture the essential features and eases its mathematical analysis, making it possible to predict and/or explain complex emergent phenomena, which go beyond the individual behavior of their constituent elements \cite{zou2019complex}. \emph{Persistence} is just one of these phenomena. Understanding persistence within the CN framework requires to previously introduce some essential concepts.

The first concept arises from the very interaction between nodes. When two nodes are directly connected by a link then they are said to be ``adjacent'' or neighboring. The adjacency matrix $\mathbf{A}$ encodes the topology of a network, that is, whether or not there is a link ($a_{ij}=1$ or $a_{ij}=0$) between any two pairs of nodes $i$ and $j$. Sometimes, this binary encoding is not enough and it is necessary to quantify the importance of any link (for example the strength of a tie between two users in a social network, or the flow of electricity between two nodes in a power grid \cite{cuadra2015critical}) by assigning a ``weight'' to each link. In this case, the matrix that encodes the connections is called {\em weighted adjacency matrix}, $\mathbf{W}$ \cite{newman2018networks}. 

A second important concept is the number of links that each node $i$ has, which is called its {\em degree}, $k_i$. The degree distribution of a network captures the probability $P(k)$ that a randomly chosen node exhibits degree $k$. $P(k)$ and its mean value $\langle k \rangle$ (mean degree) are very useful since they quantify to what extent nodes are heterogeneous with respect to their connectivity. In fact, many real-world networks exhibit broad, heterogeneous degree distributions. In a degree-heterogeneous network, the probability to find a node with $k >\langle k \rangle$  decreases slower than exponentially, leading to the existence of a non  negligible number of nodes with very high degrees. A key feature of such degree distributions is the so-called scale-free behavior \cite{barabasi1999emergence,barabasi2003scale}, characterized by a degree distribution $P(k)\sim k^{-\gamma}$. This means that most of the nodes have very few links, while only a few nodes have a large percentage of all links. These most connected nodes are called ``hubs''.

Another important CN feature is the transitivity or clustering property \cite{serrano2006clustering}, which quantifies the probability that two neighboring or adjacent nodes of a given node $i$ are connected. This concept is clear in social networks: the fact that usually ``the friend of a friend is a friend'' leads to high clustering coefficient.

The last key feature we introduce here is the average shortest path length, $\ell$, which quantifies the extent to which a node is accessible from any other node \cite{barabasi2016network}. The length of the shortest path between two nodes $i$ and $j$ in a network is the minimum number of links for going from node $i$ to $j$. Its average is computed over all possible pairs of nodes. When $\ell$ is small when compared to the network size (number of nodes, $N$), the small-world property arises \cite{watts1998collective}. Intuitively this means that any pair of nodes are relatively ``close''. Mathematically, this means that the average shortest path scales logarithmically with the network size \cite{barabasi2016network}: $\ell \sim \ln(N)$.

The complex network concept can be generalized to that of multiplex network, in which  nodes take part in several layers of networks \cite{bianconi2013statistical}, connected by relationships of different kinds \cite{battiston2014structural}. The recent work \cite{papadopoulos2019link} shows that links tend to persist from one layer to another in multiplex networks modeling real systems (Internet, Drosophila, C. Elegans, Human Brain), and points out that the process of link formation in one layer has a key influence in the topology of other layers that form the multiplex network.

The aforementioned CN features help capture not only the topology (adjacency relationships) but also the \emph{dynamics} of phenomena that can occur in the system encoded by the network  \cite{boccaletti2006complex,barrat2008dynamical}. In the same way as many complex systems represented by CN can be very diverse in nature, the dynamic phenomena they experience can also be very different. Representative examples in CN literature \cite{newman2018networks,barabasi2016network} are, for instance, the spreading of outages (in power grids), infectious diseases (in biological networks), malicious software (engineered networks), or memes and ideas (social networks). Despite their differences, these  phenomena can be modeled within the CN framework as epidemic processes  \cite{newman2018networks,perra2012activity,liu2014controlling}. Regarding this, exploring persistent phenomena within CN requires differentiating two classes of dynamics: dynamics \emph{on} CN (Section \ref{sec:on}) and dynamics \emph{of} CN (Section \ref{sec:of}). 

\subsubsection{Dynamics on complex networks and persistence}\label{sec:on}

Dynamics \emph{on} a complex network refers to a process that occurs on its underlying (or static) network structure. Illustrative examples are spreading models for the propagation of epidemics \cite{wang2019coevolution}, information, memes \cite{dawkins1989selfish} or behaviors, and gives the theoretical background for predicting and controlling these processes \cite{pastor2015epidemic}. In this regard, finding those paths causing the spreading on the network is crucial for implementing efficient strategies to either hinder dissemination (in the case of diseases), or speed up spreading (in the case of information) \cite{ren2018structure}. In CN with a broad degree distribution 
\cite{barabasi1999emergence,barabasi2003scale,albert2000error}, hubs are the key agents involved in fast spreading process \cite{pastor2001epidemic,albert2000error,cohen2001breakdown}. Additionally, the importance of a node for spreading is usually associated with the betweenness centrality, a measure of how many shortest paths go across such a node, and helps find the most influential node, for instance, in the context of social networks \cite{freeman1978centrality,friedkin1991theoretical}. 

Focusing first on the persistence of certain diseases, a feasible way to approach these processes is to introduce the concept of states of a node. This means that, at a given time instant $t$, any node $i$ is characterized by a variable $\varepsilon_{i}(t)$ that takes values in a discrete set. This idea is better understood if we focus 
on some of the epidemic models \emph{on} CN \cite{wang2017unification}. For illustrative purposes, we focus here on the so-called susceptible-infected-susceptible (SIS) model  \cite{pastor2001epidemic}. In this approach, the node state $S$ means that a node (individual) is susceptible to being infected, while the node state $I$ encodes that a node is already infected. Node states take values in the set $\{S,I\}$. Susceptible nodes may be infected by others via contacts (links) with a rate $\beta$, while infected nodes can be recover with a rate $\gamma$, becoming thus susceptible again. The key question is whether the disease can spread through the system and become \emph{persistent} when a threshold (``epidemic threshold'', $\epsilon_{Th}$) is exceeded \cite{bohme2013emergence}. An epidemic \emph{persist}  if $\frac{\beta}{\gamma} \geq   \epsilon_{Th}$ \cite{pastor2001epidemic,pastor2001endemic}.
In scale-free networks $\epsilon_{Th}$ is the maximum eigenvalue of the adjacency matrix $\mathbf{A}$. This means that $\epsilon_{Th} \to 0 $ for infinite networks size (thermodynamic limit) \cite{bohme2013emergence}. However, in finite size scale-free networks, an epidemic can always persist, though at low prevalence level \cite{pastor2001epidemic,pastor2001endemic}.

The existence of a threshold that separates two different macroscopic phases (absence/ persistence of a disease) is an instance of percolation transitions on complex networks \cite{li2021percolation}. 
An illustrative example of percolation during the growth of networks is the emergence of electron transport in some networks that model systems of disordered quantum dots (QDs) \cite{cuadra2021modeling,cuadra2021approaching}. In this approach, a QD $-$which confines most of the wavefunction inside it$-$ is encoded by a node, while electron hopping between two QDs (nodes) is represented by a link. Electron hopping is related to the overlapping of the electron wavefunctions in both nodes \cite{cuadra2021modeling,cuadra2021approaching}. In the case in which the dot density is very small, these are disconnected because the wavefunctions do not overlap. As the QD density (or equivalently, the average node degree $\langle k \rangle$)  increases, QDs are closer and closer and electron wavefunctions start to overlap, causing the formation of small, disconnected subnetworks  (called components) in which electrons are localized. There is a threshold value $\langle k_{Th} \rangle$ for which one of the subnetworks becomes dominant (giant component) and begins to grow to the detriment of the others. There exist even a value $\langle k \rangle > \langle k_{Th} \rangle $  for which a single giant component emerges and the network becomes connected \cite{cuadra2021modeling}. Thus, the electron wavefunction is spread throughout the network instead of remaining localized in one or several subnetworks. As shown in the recent review \cite{li2021percolation}, there are many other complex systems whose dynamics can be modeled as percolation phenomena on complex networks, involving not only the aforementioned epidemic spreading but also cascading processes or the persistence of network functionality against deliberate attacks or random failures. The interested reader is referred to \cite{li2021percolation} for further details about network-specific percolation models and their applications. Note that percolation is, together with synchronization, one of the phenomena that can become explosive, i.e. first-order-like \cite{boccaletti2016explosive} in CN.  While explosive percolation \cite{achlioptas2009explosive} is related to an abrupt change in the network structure, explosive synchronization \cite{kuramoto2003chemical} corresponds to the sudden emergence of a collective behaviour in the networks dynamics, which has sparked important debates in the scientific community  \cite{boccaletti2016explosive}. Especially interesting in the context of persistence in CNs is the so-called bootstrap percolation \cite{gao2015bootstrap}. Put it simple, it works as follows: (1) nodes are either inactive or active with a probability $p$; (2) inactive nodes become active if they have at least $k$ active neighbors; (3) once activated, nodes persist. Bootstrap percolation on CNs has been used to explain the emergence of social behaviors including trends, political opinions, beliefs, rumours, memes or cultural fads \cite{gao2015bootstrap}.

Just regarding the emergence and persistence of beliefs, fads and rumours, let us now shift our focus to the spreading of information on CN. A susceptible-infected-recovered (SIR) model has been used to study the message persistence on different CN \cite{cui2014message}. In that work, nodes encode the individuals in a population and links represent the propagation of information among them. Each node can be in one of three node states at any time $\{S,I,R\}$. In the susceptible ($S$) or uninformed state, the node has not yet received the message, or knows it but has not yet decided to disseminate it. In the infected ($I$) state, the node receives the information and forwards it to the adjacent ones. In the recovered ($R$) state, the node, which had already transmitted the information before, will no longer transmit it anymore. This SIR approach \cite{cui2014message} for information spreading exhibits the novelty of considering memory effects, social reinforcement and decay effects together. Memory effect refers to the fact that previous activities can affect the current process \cite{weng2012competition}. Social reinforcement \cite{krapivsky2011reinforcement} consists in that an individual usually has a greater tendency to adopt an idea if many others adopt it since they consider the idea is important or interesting \cite{centola2010spread}. Decay effects refer to the fact that the novelty of a message tends to fade over time \cite{wu2007novelty}, since there are many memes competing for our limited attention \cite{weng2012competition}. The proposed model points out that, on the one hand, higher social reinforcement (arising from larger clustering coefficient) leads to better message diffusion, which can reach (infect) a longer number of susceptible nodes, while, on the other hand, the larger the persistence, the easier the idea spreads. The results successfully explain the phenomenon ``A lie, if repeated often enough, will be accepted as truth''.

Another especially relevant paper is \cite{weng2012competition}, which focuses on the emergence and persistence of memes on Twitter. The authors model the problem as a complex network in which nodes represent Twitter users and links encode retweeted posts that carry the meme under study. This work shows that the combination of the social network structure and the competition among memes (for the finite user attention) is a sufficient condition for the emergence and persistence of memes, regardless the intrinsic quality of the meme. These results do not prove that exogenous characteristics (the inherent value of a particularly successful meme) play no role in determining its persistence. The research \emph{does} show that, at a statistical level, it is not necessary for an idea to be particularly relevant to explain its global dynamics on the network. This is a crucial difference when comparing with biological epidemics in which  the inherent properties of viruses and their adaptation to hosts are strongly important in determining the winning competing strains \cite{weng2012competition}.

\subsubsection{Dynamics of complex networks and persistence}\label{sec:of}

The dynamics \emph{of} a complex network refers to the evolution in time of the network structure, which is no longer static, as the approach in Section \ref{sec:on}. Mathematically, this leads to a time-dependent adjacency matrix $\mathbf{A}(t)$ \cite{liao2017ranking}. In fact, most of real-world networks  change over time: links activate and deactivate over time according to intermittent interactions between nodes. Illustrative examples are online social networks (friends are added and removed) \cite{xiang2010modeling}, mobile networks \cite{morse2016persistent,karsai2014time}, human contact networks \cite{barrat2013temporal,rodriguez2017risk,kim2012centrality}, or animal contacts between farms \cite{koher2019contact}. For instance, in the case of e-mail, a link is active only for the time interval in which the information is transmitted between two nodes (users). These and other CN are called Temporal Networks (TN) \cite{holme2012temporal} since the very  existence of the links evolve over time. TN can be analyzed as a time-ordered sequence of network snapshots over the set of nodes (time-ordered sequences of graphs). Within the framework of TNs, the basic unit of interaction is usually called ``contact'' \cite{holme2015modern} instead of link, more used for static CN, whose main concepts have introduced in Section \ref{sec:on}. 

\begin{figure}[!ht]
\vspace{24pt}
\begin{center}
\includegraphics[draft=false, angle=0,width=12cm]{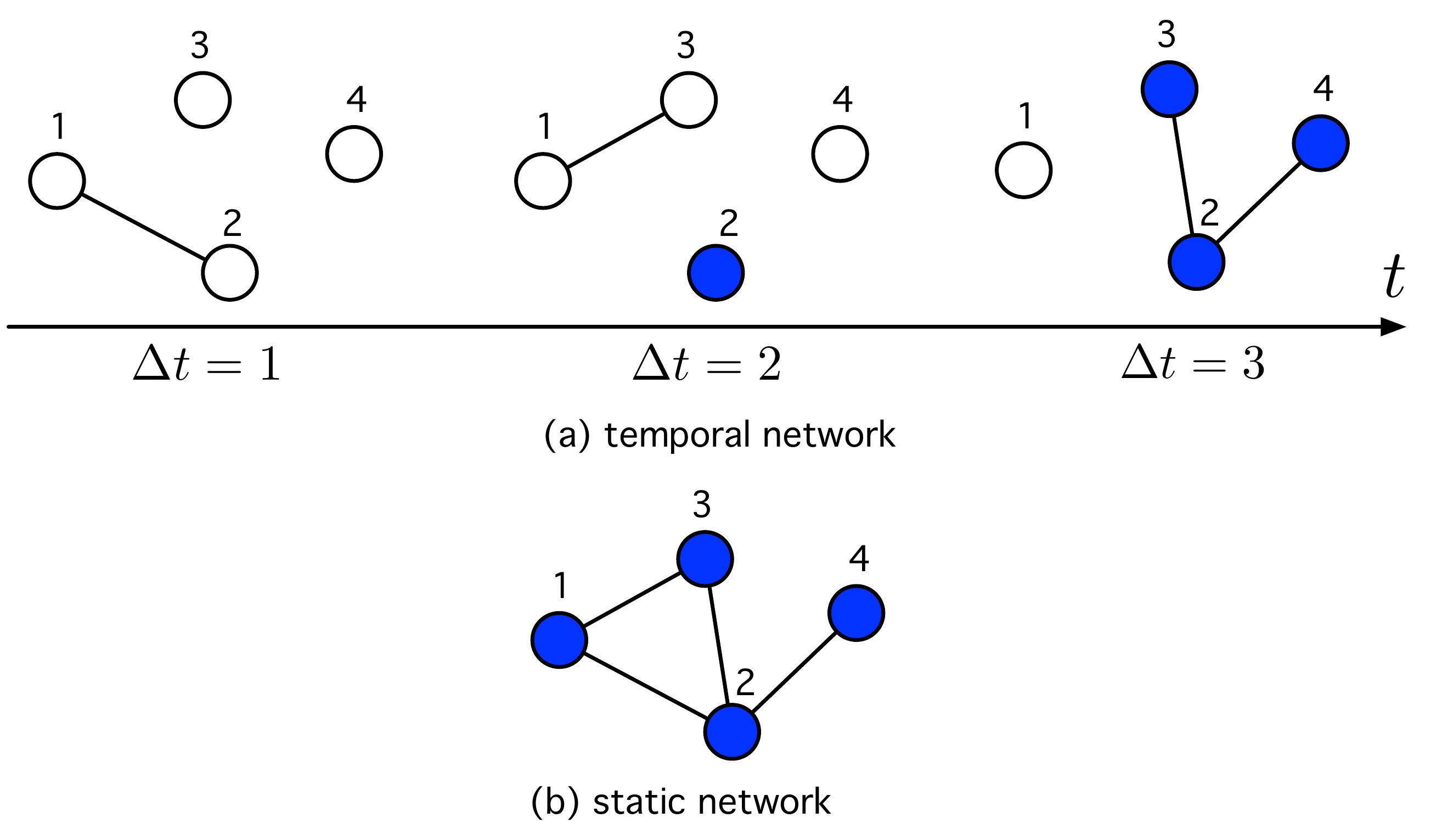}
\end{center}
\caption{\label{TN_Concept} (a) A schematic representation of the impact of the order of contacts for an epidemic process (blue nodes) that is spreading on a network. (b) Corresponding static network in which all four nodes appear to be connected. 
\vspace{24pt}}
\end{figure}

TN are more difficult to study than static CN because the \emph{order} in which contacts are established is as important as their lengths over time. For instance, in the spread of an epidemic, the order of interaction can have a crucial effect on which nodes become infected \cite{colman2018reachability,clauset2012persistence}. Figure \ref{TN_Concept} aims at illustrating this idea for both temporal (a) and static (b) networks. 
In particular, Figure \ref{TN_Concept} (a) shows that those interactions that precede infection time (blue nodes) do not propagate the disease. Note that if nodes 1 and 2 interact in $\Delta t=1$ (before node 2 being infected in $\Delta t=2$) then node 1 is not infected \cite{clauset2012persistence}. That is, in TN, the concept of node adjacency depends on the exact temporal ordering of the links \cite{nicosia2013graph,holme2013temporal}. However, in the corresponding static network in Figure \ref{TN_Concept} (b), all four nodes appear to be connected and infected. Along with order, the time scale in which the system is studied is one of the most important difficulties in TN  \cite{caceres2013temporal}.

Nonetheless, until very recently, the most common approach to study the dynamics of CN has consisted in using time-aggregated representations (static approach), as illustrated in Figure \ref{TN_Concept} (b), in the detriment of their intrinsic temporal nature. In the static approach, the degree of a node $k$ quantifies the total number of connections to other nodes, while the weight $w$ of a link encodes its importance, for instance, the total number of calls between a pair of linked nodes in a telecommunication network \cite{karsai2014time}. However, the static approach may not be accurate enough because in complex, time-varying systems represented by TN, the contacts usually exhibit two features. The first one is that contacts are usually bursty, since they occur in short and sudden episodes within time intervals \cite{min2013burstiness,vazquez2013spreading}. This leads to time-between-contacts distributions with heavy tails. The second important property is that contacts may be \emph{persistent}, that is, if two nodes are linked at a time instant, $t$, there is a non-zero probability that they will still be linked at time $t + \Delta t$. As TN evolve, some subsets of their nodes and links may be more continuously active than others \cite{holme2013temporal}. Persistence in TNs exhibits different facets, depending on the element that persists, be it a link or bigger structures or patterns \cite{holme2013temporal} such as, for example, hubs \cite{banerjee2020persistence} or communities (sets of nodes which are internally densely connected \cite{newman2004finding}).

Focusing first on the persistence of a link, its  probability of persistence can be computed as the temporal-correlation coefficient \cite{nicosia2013graph,tang2010small}. This measure is based, in turn, on the average topological overlap (over all the temporal snapshot, $N_{S}$) of the neighbourhood of node $i$, which has the following expression in a discrete-time domain:
\begin{equation}
    \bar{C_i}[n,n+1]  = \frac{1}{N_{S}-1}   \sum_{k=1}^{N_{S}-1} C_i[k,k+1],
\end{equation}
where $C_i[n,n+1]$ is the topological overlap of the  neighbourhood of node $i$ in two consecutive time steps $n$, $n+1$:
\begin{equation}
C_i[n,n+1] =  \frac{\sum_{j=1}^{N} a_{ij}[n] a_{ij}[n+1] }{  \sqrt{\sum_{j=1}^{N} a_{ij}[n] \sum_{j=1}^{N}a_{ij}[n+1] }     } .
\end{equation}
The average of $\bar{C_i}$, computed over all the $N$ nodes in the network, is the temporal-correlation coefficient \cite{tang2010small}, and it is a measure of the overall average probability for a link to persist. A similar metric for nodes, called loyalty, has been defined in 
 \cite{valdano2015predicting}. The interested reader is referred to \cite{nicosia2013graph} to understand how the metrics used for static complex networks have to be redefined for TN to take into account the effects of time ordering. 
There are other works that study the persistence of structures that are bigger than links or nodes  \cite{barucca2018disentangling,li2018persistent,saramaki2014persistence,morse2016persistent}. For instance, \cite{li2018persistent} focuses on the problem of finding communities that are persistent over time in a temporal network, in which every link is associated with a timestamp. In 
 \cite{saramaki2014persistence} the person's social signature (the pattern of his/her interactions with different friends and family members) has been explored, leading to the conclusion that it tends to persist over time. Similarly, \cite{morse2016persistent} explores persistent patterns of mobile communication among users. \cite{zhang2014dynamic} explores dynamic motifs (statistically significant subgraphs or patterns of a larger graph) in socio-economic networks, and shows that loyalty leads to heavy tails of the duration of business contacts \cite{holme2013temporal}.
 
A key point related to persistence on TNs arises from the fact that the dynamics of many complex social, technological and economic systems, are driven by individual human actions, which follow non-Poisson statistics characterized by bursts of quickly occurring events, separated by long periods of inactivity, as mentioned before. Regarding this, \cite{vazquez2013spreading} explores how human bursty activity patterns has a key influence on the dynamics of spreading processes in social and technological systems, leading to the conclusion that such non-Poisson nature results in prevalence decay times significantly larger than those predicted by models assuming Poisson processes. Thanks to this slow dynamics, the spreading entity (rumor, virus, etc.) can persist in the system. In the case of biological diseases, epidemic thresholds (determining whether or not a disease becomes persistent) on TN have been explored in \cite{leitch2019toward}. In particular, 
\cite{sun2015contrasting} has studied SIS epidemic models on activity driven networks with and without memory (Markovian and non-Markovian, respectively). Activity-driven models \cite{perra2012activity,perra2012random,liu2014controlling} associate each node $i$ with an activity $a_i$, a feature that represents the tendency of a node to create links per unit of time: at any time step $t$, nodes create links with a probability that is proportional to their activity. This  model \cite{perra2012random,liu2014controlling} has been enriched with some improvements, such as memory \cite{laurent2015calls}, aging effects caused by bursty social interactions \cite{moinet2015burstiness} and attractiveness \cite{alessandretti2017random}.
The attractiveness of a node is a measure of its tendency to attract contacts, and models the fact that some nodes may be more popular \cite{ghoshal2006attractiveness,starnini2013modeling}. In the particular case of SIS processes studied in \cite{sun2015contrasting}, memory reduces the epidemic threshold and increases the fraction of nodes in which the disease persists, the virus surviving in tightly connected local clusters that acts as reservoir.

Finally, persistence of emotions in online chatting communities has been explored in 
\cite{garas2012emotional} within the framework of TN. At any time step, each user generates a post expressing a positive, negative, or neutral emotion. These TN have activity patterns and emotional patterns. On the one hand, the activity patterns are described by two parameters. The first one, called inter-activity time or waiting-time, $\tau$, is the time interval between two consecutive posts of the same user in the same channel. It has a power-law distributed $P(\tau) \sim \tau ^{-\alpha}$, already found in different human activities \cite{barabasi2005origin,malmgren2008poissonian,crane2010power}. 
The second activity parameter, called inter-event time, $\omega$, is the time interval between two consecutive posts in the same channel. Note that the inter-event time is independent of any user and must be distinguished from the inter-activity time, which characterizes a single user. The inter-event time $\omega$ has a distribution $P(\omega)$ that is exponential and suggests the existence of long range correlations in the conversation activity. Applying the DFA method, the authors of \cite{garas2012emotional} have found a Hurst exponent value $H\approx 0.6$, showing the presence of long-term \emph{persistence}. On the other hand, regarding emotional patterns, users' emotional expression has been found to exhibit a Hurst exponent $H>0.5$, which means that most of users are quite persistent regarding their positive, negative or neutral emotions. As mentioned before, this persistence can be also seen as a kind of memory in changing the emotional expression, that is, the following post from the same user is more likely to have the emotion.

\subsection{Persistence in Economics and market analysis}

Persistence is a key point in Economics and market analysis, recurrently studied for different economic processes \cite{mcgahan1999persistence}. Evidence that economic shocks or initial differences in economics outcomes dissipate very slowly (if at all) has been given in the study of processes such as inflation, firms' asset prices and returns. The long-term persistence of other time series such as exchange rates, stock market indices, market analysis or crypto-currency analysis, among others, have been also analyzed in the literature.

Inflation persistence is defined as the time that it takes for an inflation shock to dissipate, i.e. how quickly a stationary inflation process reverts to its initial level or to the long-run equilibrium after a shock \cite{fuhrer1995inflation,Sbordone07,Gaglianone18}. The research on inflation persistence tries to shed light into questions such as whether inflation persistence is an inherent characteristic of the economy, or it instead depends on the specific historical sample considered, or what are the factors that effect changes in the degree of persistence over time, among others. In general, inflation persistence analysis uses short-term persistence techniques analysis, such as autorregressive models. Usually inflation persistence is analyzed by means of univariate analyses, based on the size of the largest AR root, the sum of the coefficients, the largest root, or the half-life (the number of periods for which inflation remains above 0.5 for a unit shock) \cite{pivetta2007persistence}. Multivariate analyses have also been applied to model inflation persistence, based on impulse response functions, when modeling inflation as a highly persistent process (US inflation case, for example) \cite{Sbordone07}. In \cite{Meenagh09} AR and ARMA processes are applied to model the inflation persistence in UK since WWII. In \cite{Gaglianone18} the persistence of Brazilian inflation and its main components is studied by using quantile regression techniques \cite{koenker2001quantile}. In \cite{Tule19} inflation persistence in Nigeria is modelled by using the recently developed fractional co-integration VAR model. In \cite{pivetta2007persistence} the inflation persistence in USA is studied using different approaches, including autorregressive models. All the models tested led to the same conclusion: persistence of inflation has been high and approximately constant over time since 1965.

The persistence of firms' returns and assets prices have been analyzed in different works in the literature. In \cite{Priestley01} the analysis of the persistence of expected returns in terms of its effect on assets prices is carried out. An AR(1) model is used to characterize persistence in the expected holding period. In \cite{Dichev09} the persistence of earnings in a corporation is analyzed, trying to find out whether past earnings volatility has an explanatory power for the persistence of current earnings. In \cite{Frankel09} the same problem is further analyzed in terms of its possible effects, pointing out the underlying lack of a causal theory in \cite{Dichev09} findings. In \cite{Wu12} the persistence of firm's growth type is studied and compared to their leverage (investment strategy using borrowed capital). This work demonstrates that firms rationally invest and seek financing in a manner compatible with their growth types, or in other words, growth type can parsimoniously predict and distinct future leverage ratios. In \cite{Pla19} a DFA analysis is applied to a time series of daily absolute firm's assets and equity returns for a sample of non-financial iTraxx Europe companies. The scaling exponent of the DFA is slightly higher for absolute equity returns than for absolute firm's asset returns (0.85 vs. 0.81), but in both cases it confirms the long-term persistence of these time series. These results agree with a previous FIGARCH model also applied to these time series of firm's assets and equity returns.

The persistence and long-term behavior of exchange rates and crypto-currencies has been very important in the last years. In \cite{Curran19} the analysis of the short-term persistence of real exchange rates across 151 countries is carried out. Several autorregressive models such as AR(1), and variations such as General higher-order autoregressive models or exponential smooth transition autoregressive models, are tested. In \cite{Caporale18} R/S analysis is applied to analyze the long-term persistence of the four main crypto-currencies (BitCoin, LiteCoin, Ripple and Dash) over the period 2013-2017. The results obtained indicate that this market exhibits persistence, i.e. there is a positive correlation between its past and future values, and that its degree changes over time. The conclusion of the paper is that this long-term predictability may lead to trend trading strategies to generate abnormal profits in the crypto-currency market. In \cite{costa2019long} a DFA and detrended crosscorrelation analysis are applied to evaluate the persistence properties of four cryto-currencies (Bitcoin, Ethereum, Ripple, and Litecoin). The results obtained in this case showed that Bitcoin and Ripple seem to behave as efficient financial assets, while Ethereum and Litecoin present some evidences of persistence. In \cite{takaishi2021time} a complete analysis of volatility of daily Bitcoin returns and multifractal properties of the Bitcoin market is carried out. Among other efficiency-related measures, the Hurst exponent, multifractal degree, and kurtosis are obtained for the Bitcoin price, return and volatility time series. In \cite{stosic2019multifractal} an analysis of the persistence of price changes and daily volume changes of 50 crypto-currencies by means of a MF-DFA algorithm is carried out. Results show an absence of correlations in price changes, whereas volume changes present anti-persistence in the long-term. In \cite{david2021fractional} the persistence of different crypto-currencies are analyzed, such as Bitcoin, Litecoin, Ripple, Monero, Ethereum and Ethereum Classic. ARIMA, Auto-Regressive Fractionally Integrated Moving Average and DFA algorithms are used to evaluate the short-term and long-term persistence of the considered crypto-currencies. The results suggest that the analized crypto-currencies exhibit a lower predictability when compared to the Bitcoin. The results for the Bitcoin indicate a clear persistent behavior in the long-term.
In \cite{al2018efficiency} the efficiency of Bitcoin market compared to gold, stock and foreign exchange markets is analyzed by using a MF-DFA approach. In \cite{vaz2021price} a MF-DFA and a Multifractal Regime Detecting Method are applied to study the multi-fractality of bitcoin returns in the period 2013-2020. The results obtained suggest that bitcoin is multifractal, with significant long-range correlations which expose the informational efficiency of the cryto-currency. In \cite{quintino2020efficiency} a DFA analysis of Bitcoin prices from the Brazilian market is carried out. The results suggest that the persistence of Bitcoin in slightly higher for the Brazilian Bitcoin than in other markets. In \cite{alvarez2018long} an efficiency analysis of the Bitcoin market is carried out, by applying a DFA approach. the resuls obtained indicate that the Bitcoin market exhibits periods of efficiency, alternating with periods where the price dynamics are driven by anti-persistence. In \cite{luis2019drivers} the demand drivers of Bitcoin are analyzed in the short-term and long-term horizons by applying AR models and GARCH methodology. The results obtained show that Bitcoin behaves as a speculative asset in the short-term, but in the long term, speculation does not seem to influence demand for Bitcoin. 

The persistence of stock markets indices has been also studied and analyzed in different works. In \cite{Grau01} a long-range power-law correlations using DFA are obtained for different stock market indices in USA, UK, Japan, Germany, France and Spain. A two-ranges pattern is obtained in the log-log plot of the DFA applied to the absolute and squared returns (logarithmic difference in the index). A characteristic time (crossover point) of 41 days is obtained in such a way that $\alpha_1\approx 0.6$ and $\alpha_1\approx 0.8$.
In \cite{Cajueiro05} a R/S analysis and a DFA algorithm are applied in order to investigate long-term persistence of the Brazilian stock market (BOVESPA Index). This study evidences statistically significant correlation between specific variables of the Brazilian firms analyzed and a clear long-range dependence phenomenon present in these stocks.
In \cite{constantin2005volatility} an analysis of stock prices persistence is carried out by direct analysis of the correlation function to obtain the Hurst coefficient, and also using a MF-DFA technique. Results in different stock indices are reported. In \cite{Oh08} the long-term persistence of several stock market indices and foreign exchange rates are studied by using a DFA algorithm. Real data from the KOSPI 1-minute market index in the Korean stock market, the KOSDAQ 1-minute market index from the Korean Stock Exchange, and the foreign exchange rates of six currencies to the US dollar (Euro (EUR), the UK (GBP), Japan (JPY), Singapore (SGD), Switzerland (CHF), and Australia (AUD)) are considered. The results obtained showed that the exponents estimated for this time series by the DFA method were in the range $0.56 \leq H \leq 0.68$, what suggests a long-term memory of the series. This work also tries to explain the causes of the long-term persistence by modeling the series with AR(1), Generalized Autoregressive Conditional Heteroscedasticity (GARCH) models \cite{bollerslev1986generalized} (specifically GARCH(1,1)), and a variation known as Fractional Integrated GARCH (FIGARCH) models \cite{baillie1996fractionally} (specifically FIGARCH(1,d,1)). There are other works in the literature which relays on GARCH and FIGARCH models to analyze the long-term persistence in econometric time series, specially in stock market analysis, such as \cite{Bentes14}. In \cite{Granero08} the R/S methodology is used in order to study the long-term persistence of  several international stock market indices: Cac40, FTSE, Nikkei, Nasdaq, Ibex 35 and S\&P500. This study shows that the standard R/S method has serious problems which may lead to even obtain evidence of long memory in random series. Two new geometrical interpretations of the Hurst index are proposed in this work in order to solve these issues with the R/S method. In \cite{Lu13} a MF-DFA is applied to the specific analysis of the Chinese stock market. The results obtained suggest that the multi-fractality is mainly due to long-range correlations in the stock market indices analyzed.
In \cite{yin2013modified} the DFA algorithm and and multiscale detrended cross-correlation analysis have been applied to analyze the persistence of different US and Chinese stock markets indices during the period 1997-2012. The well-known S\&P500, NQCI, HSI, and the Shanghai Composite Index have been considered as representative stock market indices. The results indicate that US and Chinese stock indices differ in terms of their multiscale auto-correlation and long-term persistence structures. In \cite{stovsic2015multifractal} a MF-DFA and related techniques are applied to analyze the persistence of 13 global stock market indices, considering daily price changes and daily volume changes. The results obtained show that time series of price changes are more complex than those of volume changes, and that large fluctuations dominate the long-term persistence of price changes, while small fluctuations dominate the long-term persistence of volume changes.
Finally, In \cite{milocs2020multifractal} a MF-DFA approach is applied to the analysis of seven Central and Eastern European stock markets, with recent financial data. In many cases it was found that stock indices returns exhibit long-range correlations, indicating inefficiency in the analyzed stock markets. 

The analysis of markets and their persistence properties have been recently studied in a number of works. Electricity market has been analyzed with DFA and MF-DFA and related techniques, at different levels. For example, the long-term persistence characteristics of the European Electricity market spot prices have been studied in \cite{gorjao2021change,han2021complexity}. The US Electricity market persistence has been analyzed in \cite{ali2021modeling} and there are also works dealing with national electricity markets such as in the Czech Republic \cite{fan2015multifractal}. Gold market long-term persistence properties has been analyzed using MF-DFA techniques in \cite{nejad2021multifractal}. Crude oil and derived prices and markets long-term persistence analysis has been discussed in the literature, with the usual techniques already discussed in this peper, such as Hurst coefficient estimation, DFA, MF-DFA and related approaches. Studies on general crude oil markets persistence \cite{gu2010multifractal,delbianco2016multifractal,cai2019exploring,cerqueti2021long,mensi2021does,ftiti2021oil} and specific national markets such as Chinese \cite{zhang2021cross} or Brazilian \cite{david2020measuring} can be found in the literature. Long-term persistence evaluation techniques have also been recently applied to the study of other markets such as housing \cite{raza2021multifractal} or agricultural \cite{yin2021market,feng2021multifractal}.

\subsection{Persistence in non-equilibrium thermodynamics systems}\label{PnonEq}

Most thermodynamics systems in nature are not in equilibrium, i.e. they exchange fluxes of matter or energy with their surroundings, and/or undergo chemical reactions \cite{rupprecht2016fresh}. The study of statistical mechanics in non-equilibrium systems include different models such as phase separation process \cite{Derrida1994non}, simple diffusion equation with random initial conditions \cite{majumdar1996global}, several reaction diffusion systems \cite{cardy1995proportion}, fluctuating interfaces \cite{krug1997persistence}, Lotka-Volterra models of population dynamics \cite{frachebourg1996spatial} and granular medium \cite{swift1999survival}, among others. 

In the last years, the study of persistence concepts in non-equilibrium systems has been intense. The precise definition of persistence in non-equilibrium systems is as follows: Let
$\phi({\bf x},t)$ be a non-equilibrium field fluctuating in space and
time according to some dynamics. Persistence is defined as the probability $P_0(t)$ that, at a fixed point in space, the quantity
$\sgn\left[\phi({\bf x},t) - \langle\phi({\bf x},t)\rangle\right]$ does not change up to time $t$ \cite{majumdar1999persistence}. Note that this definition is one of the described in Section \ref{short-term-definitions} for short-term persistence. It is known that in many non-equilibrium thermodynamics systems this probability decays as a power law $P_0(t) \sim t^{-\theta}$ at late times (non-stationary systems), where the persistence exponent $\theta$ is usually nontrivial \cite{majumdar1999persistence}. In stationary systems, persistence usually decays as $P(\tau) \sim \exp\left(-\theta_s \tau \right)$ \cite{Majumdar01}, where $\tau=t_2-t_1$, stands for time differences.

Paradigmatic non-equilibrium systems where persistence topics have been studied are phase ordering kinetics, i.e. the growth of order through domain coarsening when a system is quenched from a homogeneous phase into a broken-symmetry phase \cite{bray2002theory}. In these coarsening systems, the dynamics are usually characterized by a single length scale $L(t) \sim t^{-z}$, which measures the typical size of the domains \cite{villain1984nonequilibrium}. There have been studies of persistence in similar systems, for example in \cite{Derrida1994non,stauffer1994ising} the persistence properties for the coarsening dynamics of ferromagnetic spin models are discussed. Usually, studies on ferromagnetic systems undergoing phase ordering have focused on the study of persistence $\mathcal{P}^P(t)$ of the local magnetization, which is the probability that the local spin at site $r$ has not flipped between time $0$ and time $t$. An algebraic decay of $\mathcal{P}^P(t) \sim t^{-\theta}$ has been reported for these systems, with a non-trivial exponent depending on the specific problems. Other works on coarsening systems dynamics and their persistence are \cite{sire1995coarsening,derrida1995exponents,majumdar1996survival,majumdar1998persistence,tam2002cluster}. Similarly, in \cite{marcos1995self} a study on statistical mechanics of non-equilibrium dew formation is carried out. In that work, it was shown that the fraction of $f_{dry}(t)$ of the surface which was never covered by any droplet decays as a power law $f_{dry}(t) \sim t^{-1}$. Note that $f_{dry}(t)$ is equivalent to a persistence probability $\mathcal{P}^P(t)$ for the coarsening dynamics of ferromagnetic spin models described before. 

In the last years there have been some good reviews related persistence issues in non-equilibrium thermodynamics systems \cite{bray2013persistence,iyer2015first,aurzada2015persistence,majumdar1999persistence}

\subsection{Constructing synthetic time series with persistence properties}

The study of persistence of different physical phenomena heavily depends on having accurate real time series describing them. In some cases, the availability of these real time series is not immediate, real time series can be incomplete or including missing values, or simply their length is shorter than the one required to apply the mathematical techniques to study their behaviour. In these cases, synthetic time series with certain persistence-related characteristics can be very useful to test mathematical techniques or corroborate findings on real time series.

In \cite{peng1991directed} a method for generating synthetic time series with long-term correlation properties is proposed. Specifically, this method produces a sequence of random numbers with power-law correlations by using the Fourier filtering method (FfM). It consists of filtering the Fourier components of an uncorrelated sequence of random numbers with a suitable power-law filter, in order to introduce correlations among the variables. After applying the inverse Fourier transform, the resulting time series will maintain the desired correlations. In \cite{prakash1992structural} the FfM method is used in order to obtain long-term correlated time series for studying percolation phenomenon on them \cite{li2021percolation}. In spite of the FfM has been successfully used for analyzing and describing the behaviour of different complex systems, it has an important issue: the FfM presents a finite cutoff in the range over which the variables are correlated. This produces undesired effects, for example in the one dimension case, only about 0.1\% of the samples will have the desired correlation. In two dimensions this percentage grows up to 1\%. Note that this limitation makes the FfM not suitable for the study of scaling properties in the limit of very large systems. In \cite{makse1996method} a modification of the proposed FfM method is presented, in such a way that the cutoff for correlated variables is eliminated. This modified version of the FfM method has been used then to study how the presence of long-term correlations in time series affects the statistics of extreme events \cite{eichner2006extreme}. 

In \cite{halley2009using} statistical processes with long-term persistence are use to modeling the natural dynamics of Earth's temperature. Specifically, the family of $1/f$-noises is used to assess the possible natural origin for global warming, parameterized using paleoclimate reconstructions.
In \cite{Efstratiadis14} a robust method for time series synthetic  generation is presented, in the context of hydrology applications, e.g involving river inflow, rainfall, floods simulation, etc. The proposal consists of using a three-level multivariate scheme for stochastic simulation of correlated processes. The proposed method preserves the essential statistical characteristics of historical data at three time scales (annual, monthly, daily) and it also maintains the long-term persistence (Hurst coefficient), its periodicity and intermittency. In fact, the simulation of hydrological systems is a constant source of different methods for synthetic time series construction with persistence properties \cite{boughton2003continuous}. In \cite{tsekouras2014stochastic} real and simulated time series of wind speed and sunshine duration are used study the design and management of renewable energy systems. The paper analyzes both the properties of marginal distributions and the dependence properties of these natural processes, including possible long-term persistence by estimating and analyzing the Hurst coefficient of the series. These two papers (\cite{Efstratiadis14} and \cite{tsekouras2014stochastic}) use the software {\em Castalia} \cite{Castalia2,Castalia} to generate the synthetic time series with persistence properties, which implements the framework defined in \cite{koutsoyiannis2000generalized}. In \cite{Ilich14} a new method is presented to generate stationary multi-site hydrological time series. The algorithm is based on standard decomposition models and the Box-Jenkins approach. 

\subsection{Persistence in optimization and planning}

Optimization and planning problems are areas of major interest in Engineering and Computation. Persistence concepts also arise in optimization and planning problems, within the field of {\em robust optimization} \cite{ben2009robust,gabrel2014recent}. Before defining the concept of robust optimization, let us remind the concept of optimization problem. A general definition for an optimization problem as given in \cite{salcedo2016modern} is the following:
Given a system $\mathcal{S}$ (real), and an associated computational model of this system $\mathcal{M}$ (model), the optimization of $\mathcal{M}$ is defined as:
\begin{equation}\label{Opt_def}
\begin{array}{r r}
\operatorname{Opt}~~ z=f({\bf x})&\\
{\bf x} \in \Omega&\\
\text{subject to}&\\
g_i({\bf x}) \geq 0, & i \in \{1, \ldots, m\}\\
h_j({\bf x})= 0, & j \in \{1, \ldots, p\}\\
\end{array}
\end{equation}
where $\Omega$ stands for a given search space, $f({\bf x}): \Omega \rightarrow \mathbb{R}$, and the $\operatorname{Opt}$ sign stands for minimizing or maximizing (usually a cost function must be minimized and a benefit function will be maximized). Functions $g_i({\bf x})$ and $h_j({\bf x})$ stand for the constraints of the optimization problem.

Let ${\bf x}^*$ be the (optimal) solution for this optimization problem given by model $\mathcal{M}$. The key point to understand the concept of robust optimization is to realize that we are optimizing a model of reality, and therefore, as pointed out in \cite{Beyer2007}, in the majority of cases we do not have a detailed knowledge of the error function of the model $\mathcal{M}$ respect to the real system $\mathcal{S}$. We cannot, therefore, be certain that the solution ${\bf x}^*$ for the model optimum can be mapped to the true optimum. Thus, trying to be too precise in the model solution might waste computational and time resources. Other issues about obtaining a too precise optimal solution for the static optimization problem of Equation \eqref{Opt_def} are related to problems for building the true optimum either because of manufacturing uncertainties or because the required precision during the manufacturing, the existence of dynamic environments in the real system, such as environmental parameters fluctuation, materials wear down, life spam of the design, etc.
Following this discussion, note that systems optimized in the classical sense can be very sensitive to small changes, which in fact are very likely to occur in the real system optimized. 

Robust design optimization deals with these issues, by looking for  solutions and performance results which remain relatively unchanged {\em when exposed to uncertain conditions}, not only sensitivities of objective function with respect to design parameters, but also environmental, manufacturing or any other uncertainty that may arise in the process. Note that, in many cases, coping with uncertainties in such optimizations is subject of robust multidisciplinary design \cite{kalsi2001comprehensive,allen2006robust}.  

Persistence was defined as a subarea of robust optimization in \cite{brown1997optimization}, and this concept has reached their maximum meaning in planning problems. The idea behind persistence in optimization and planning is that a given optimization or planning problem may suffer small variations in its definition in a period of time, due to engineering manager decisions or project necessities. Once the first solution is obtained, the following solutions must be {\em persistent}, i.e. not very different from the first one, so the engineering or planning problem is the less possible affected. Mathematically, the first model for the problem $\mathcal{M}$ is now considered as $\mathcal{M}(t)$, so we have a time series 
\begin{equation}
\left[\mathcal{M}(1), \mathcal{M}(2), \mathcal{M}(3), \ldots,\mathcal{M}(T)\right].
\end{equation}
Let $\left[x_1^*, x_2^*, x_3^*, \ldots, x_T^*\right]$ be the optimal solutions for these models. The idea is that the time series of the optimal solution states is persistent for all times $T$ considered, i.e., all the optimal solutions can be considered to be in the same state:
\begin{equation}
    \left[x_1^* \rightarrow s_1^*, x_2^* \rightarrow s_1^*, x_3^* \rightarrow s_1^*, \ldots, x_T^* \rightarrow s_1^*\right]
\end{equation}
or, in words, the solutions for the different problem's models are similar enough so that small changes must be done in the engineering project or planning at hand.

The first work taking into account persistence of solutions in optimization problems was \cite{brown1996scheduling}, dealing with a scheduling of patrol boats for the US coast guard. The schedules must sometimes be revised due to unexpected issues in boats or crew, so the proposed algorithm was prepared to respond these schedule changes while retaining as much of an already-published schedule as practical. The concept of {\em persistence} was first introduced here. Soon later, this concept was fully presented in \cite{brown1997optimization}, which was considered as one of the founding works for optimization persistence. Different optimization and planning problems where the concept of persistence can be useful are presented, and the some issues of considering persistence in optimization problems are also discussed. In \cite{brown1997optimizing} a problem of submarines berthing is tackled, taken into account persistence in the solutions obtained. Submarine berthing is tackled as a planning problem. A real case study in Naval Submarine Base, San Diego, was taken into account, where changes are often necessary after a berthing plan has been approved, due to changed service requests, delays, early arrival of inbound vessels, etc. The algorithm presented in this work includes a {\em persistence incentive}, implemented as a penalty function which rewards persistence of solutions to future changes. In turn, \cite{morrison2010new} presents a new paradigm for robust combinatorial optimization, focused on persistence of solutions. This work is focused on treating persistence of decisions
as evidence of robustness, i.e. how solution characteristics persist across a given set of ranked solutions (a list of optimal and near-optimal solutions), and how this way it is possible to find solutions that are more robust to uncertainties in the underlying model. In \cite{petit2019enriching} it is also dealt with combinatorial optimization. In this work a general
framework for problems of a combinatorial nature is proposed, taken into account the persistence of solutions. Persistence is introduced in this work by generating a set of of multiple
(near-)optimal, diverse solutions, that can be further infused with desirable features to adapt them to changes in the problems' initial condition. Finally, in \cite{Borthen2019} a problem of offshore supply vessel planning for the oil and gas industry is tackled. A meta-heuristic approach \cite{salcedo2016modern} is proposed (genetic algorithm), which tries to obtain cost-efficient and also persistent solutions, i.e. new weekly plans exhibit few changes from the previous plan.

\subsection{Persistence in health and biomedical sciences}

Finally, we review here some applications of persistence concepts in health and biomedical sciences.

The persistence of biomedical time series have been studied in several works. In \cite{Depetrillo99} a method for estimating the Hurst coefficient $H$ from the fractal dimension $D$ for self-affine time series is proposed. The method is tested in the calculation of $H$ of heart rate variability time series in healthy individuals. The results obtained show that the proposed method based on fractal dimension obtains more accurate values of $H$ than other methods such as power
spectral density or discrete wavelet transform. The proposed method has been tested in cases where $H<0.5$ (anti-persistence), due to the characteristics of the heart rate variability time series. No results are reported for cases $H>0.5$. In \cite{de1999long} a R/S method was applied to phase-shift records of random transition between periodic solutions of a biochemical dynamic system. The Hurst coefficient obtained showed a clear long-term persistence, indicating that each phase-shift values depend on past values. The methodology was also applied to cardiac rhythm time series in order to obtain their $H$ and showing their associated long-term persistence. In \cite{de1999persistence} the R/S method is applied to study the persistence properties of metabolic networks, (dynamical systems formed by the activity of several catalytic dissipative structures, interconnected by substrate fluxes and regulatory signals). Hurst exponents $H < 0.5$ (antipersistence) were found in most cases when the R/S method was applied to this systems. In \cite{thurner2003scaling} a DFA algorithm and calculation of the power spectral density are applied to study study the temporal variability of human brain activity in time series of functional magnetic resonance imaging (fMRI) data. In \cite{rahmani2018dynamical} a R/S method is applied to electroencephalogram (EEG) signals to calculate their $H$ coefficient. The study showed that the Hurst coefficient demonstrate a significant differential response between healthy and post-traumatic stress disorder (PTSD) samples (combat-related cases). The authors pointed out that these results may lead to diagnostic applications of EEG for PTSD. Finally, somehow related to biomedical is sport science. Persistence studies here are very recent, and restricted to specific sports such as badminton. In \cite{Gomez20} a study on long rallies in badminton and their effect on subsequent rally is carried out, taking into account technical parameters of the game, and differencing by sex. Long rallies in badminton are defined per rally duration or number of strokes, using thresholds to separate between long rallies (state 1) and short rallies (state 0). The persistence of long rallies in a badminton game can be considered as another measurement of high performance and game difficulty, similar to that of the number of consecutive points won by a player.

\subsection{Other application domains}

Analysis and characterization of persistence, mainly long-term persistence, have been carried out in alternative application domains. Some examples are solar activity time series, where a long-term analysis with a DFA algorithm has been discussed in \cite{das2021hemispheric}, and the persistence of the sunspot cycle with Hurst coefficient has been discussed in \cite{das2021hemispheric}. Earth's longwave radiation flux long-term persistence has been studied with the Hurst exponent and MF-DFA analysis in \cite{stathopoulos2017long,stathopoulos2019long}. Also, there are recent studies on long-term persistence of completely different time series such as mortality data \cite{peters2021statistical}, toll-free calls \cite{gui2021long}, particle concentration on roads \cite{lu2014detrended} or ship flow sequences of containers ports \cite{liu2021long}.

\section{Case studies and practical applications}\label{Case_Studies}

In this section we show and discuss in depth several case studies, where persistence has been used to analyze different real complex systems. We have chosen significant cases, with both short-term and long-term persistence characteristics. Some of the cases presented here are related to Earth Observation \cite{Salcedo20}, and we have also included cases related to renewable energy and a study related to persistence in machine learning.

\subsection{DFA and time scales: characteristic time invariance}
Long-term analysis of atmospheric-related events with a DFA approach requires selecting the time scale of the procedure. Note that the DFA algorithm incorporates the time scale of the analysis in the length of the non-overlapping segments in which the normalized and integrated time series is divided ($s$, see Section \ref{DFA}). Thus, the DFA is able to analyze the time series at different time resolutions or scales just by setting $s$. In the standard DFA procedure, an initial time scale is set, and then it is increased, usually in regular steps, up to a maximum time scale to be analyzed, in such a way that the total number of segments $N_s$ is large enough to obtain correct estimations of error from Equation (\ref{DFA_error}). There is, however, a second possibility for considering time scales, in this case due to the actual resolution of the times series $x_i$. Note that, due to the discrete nature of the data collection, the time series $x_i$ is stored with resolution $a$ (a sample is obtained each $a$ time units). By increasing this resolution, we can explore different time scales. In this case, we consider averaged series to obtain a range of time scales in some atmospheric-related time-series, i.e. series obtained by averaging the original series $x$ over non-overlapping blocks of
size $b$: Let us consider a time series of $N$ values, associated with a given atmospheric phenomenon, $\{x_i\}_{i=1}^N$. These values are obtained at the highest possible resolution of the measurement instrument for this phenomenon $a$. From this initial time series, we construct some other time series at a different time scale $b a$, i.e. $\{x'_i\}_{i=1}^{\lfloor N/b \rfloor}$, by using the  procedure given in Algorithm \ref{alg:TS}.

\begin{algorithm}[!h]
    \caption{Time scale change procedure in time series.}
   \label{alg:TS}
   \begin{algorithmic}[1]
  \REQUIRE Time series $x(t)$ of length $N$, with time scale $a$ (a sample each $a$ time units).
   \ENSURE A time series $x'(t)$, with different time scale $a b$ and length $\lfloor N/b \rfloor$.
 \STATE Set scale change parameter $b$\;
 \STATE \[
\{x'_i\}_{i=1}^{\lfloor N/b \rfloor}=\left[ x'_1=\frac{1}{b} \sum_{k=1}^b x_k,  x'_2=\frac{1}{b} \sum_{k=b+1}^{2b} x_k, x'_3=\frac{1}{b} \sum_{k=2b+1}^{3b} x_k, \ldots \right.
\]\
\[
\left. \ldots, x'_p=\frac{1}{b} \sum_{k=(p-1)b+1}^{pb} x_k, \ldots, x'_{\lfloor N/b \rfloor}=\frac{1}{b} \sum_{k=(N-b+1)}^{N} x_k\right]
\]\
\STATE Repeat this process to obtain time series of visibility at different time scales.\

\end{algorithmic}
\end{algorithm}

Note that this procedure is a {\em time scale change}, in which the average of $b$ previous values of the series is used to obtain the series at a different sampling rate. 

\subsubsection{DFA and characteristic time}

\begin{figure}[!ht]
\vspace{24pt}
\begin{center}
\includegraphics[draft=false, angle=0,width=10cm]{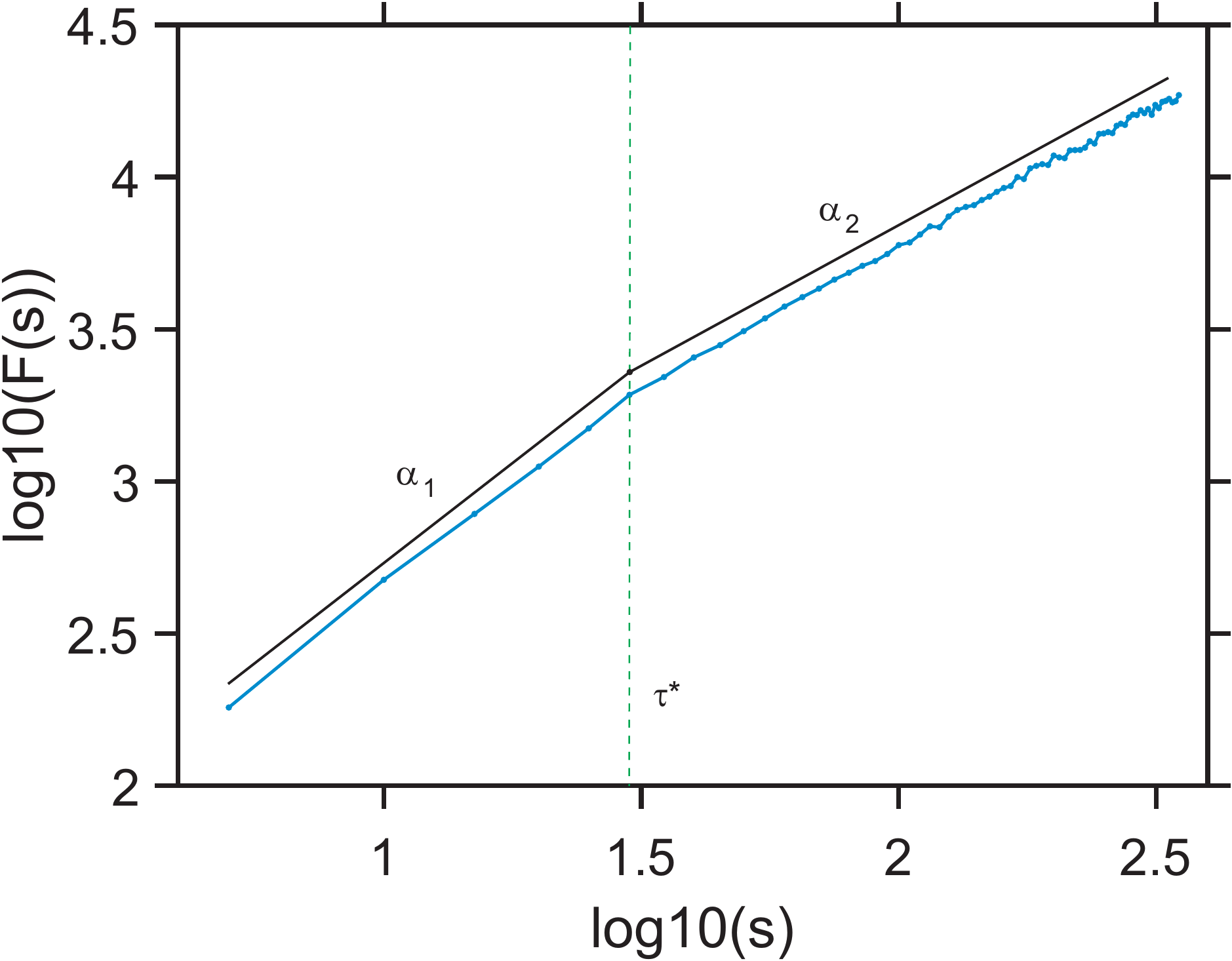}
\end{center}
\caption{\label{DFA_Example} Example of a typical DFA pattern with two temporal ranges with different scaling. The crossover point is a characteristic time of the long-term persistence of the time series ($\tau^*$), and $\alpha_1$ and $\alpha_2$ are the slopes of each segment (DFA scaling exponents or coefficients).
\vspace{24pt}}
\end{figure}

As previously discussed in Section \ref{Persistence_atmospheric_Processes}, works dealing with DFA analysis of meteorological phenomena, such as rainfall and streamflow series or SSTA analysis \cite{Matsoukas00,Luo15,Yang19,Dey18b}, have reported long-term persistence patterns given by two temporal ranges with different scaling, as shown in Figure \ref{DFA_Example}. 

This behavior was theoretically predicted in \cite{Kantelhardt01}, where a study on crossover points in DFA due to different scaling
regimes were carried out by using synthetic time series. In \cite{Gao06} this phenomenon was also analyzed for different synthetic noise time series. Note, however, that crossover points commonly appear in real data, mainly when atmospheric processes are being studied. As previously discussed, in \cite{Matsoukas00} a DFA analysis for rainfall was carried out, which showed that this two-range pattern for rainfall series is associated with the binary structure of the time series (the pattern of alternating dry and wet spells independently of the actual values of rainfall intensity. Recently, in \cite{Rasanen18} it has been shown that rainfall and fog events have a similar structure in terms of on-off intermittency properties. This two ranges pattern with characteristic time has been also shown in other persistence-based studies dealing with DFA such as, orographic fog events \cite{Salcedo21}, SSTA data \cite{Luo15}, rainfall time series \cite{Yang19}, soil moisture \cite{Shen18} or water level in reservoirs \cite{Castillo20}, among others. Note that the crossover point between two temporal ranges in a DFA analysis can be seen as a characteristic time ($\tau^*$) for the meteorological phenomenon under study. This characteristic time is invariant under time scales transformations, providing a measurement which characterizes the scaling of the phenomenon under study, together with the scaling factors (slopes of the DFA log-log plot), before and after $\tau^*$.

\subsubsection{Examples}

\begin{figure}[!ht]
\begin{center}
\subfigure[~original]{\includegraphics[draft=false, angle=0,width=8.2cm]{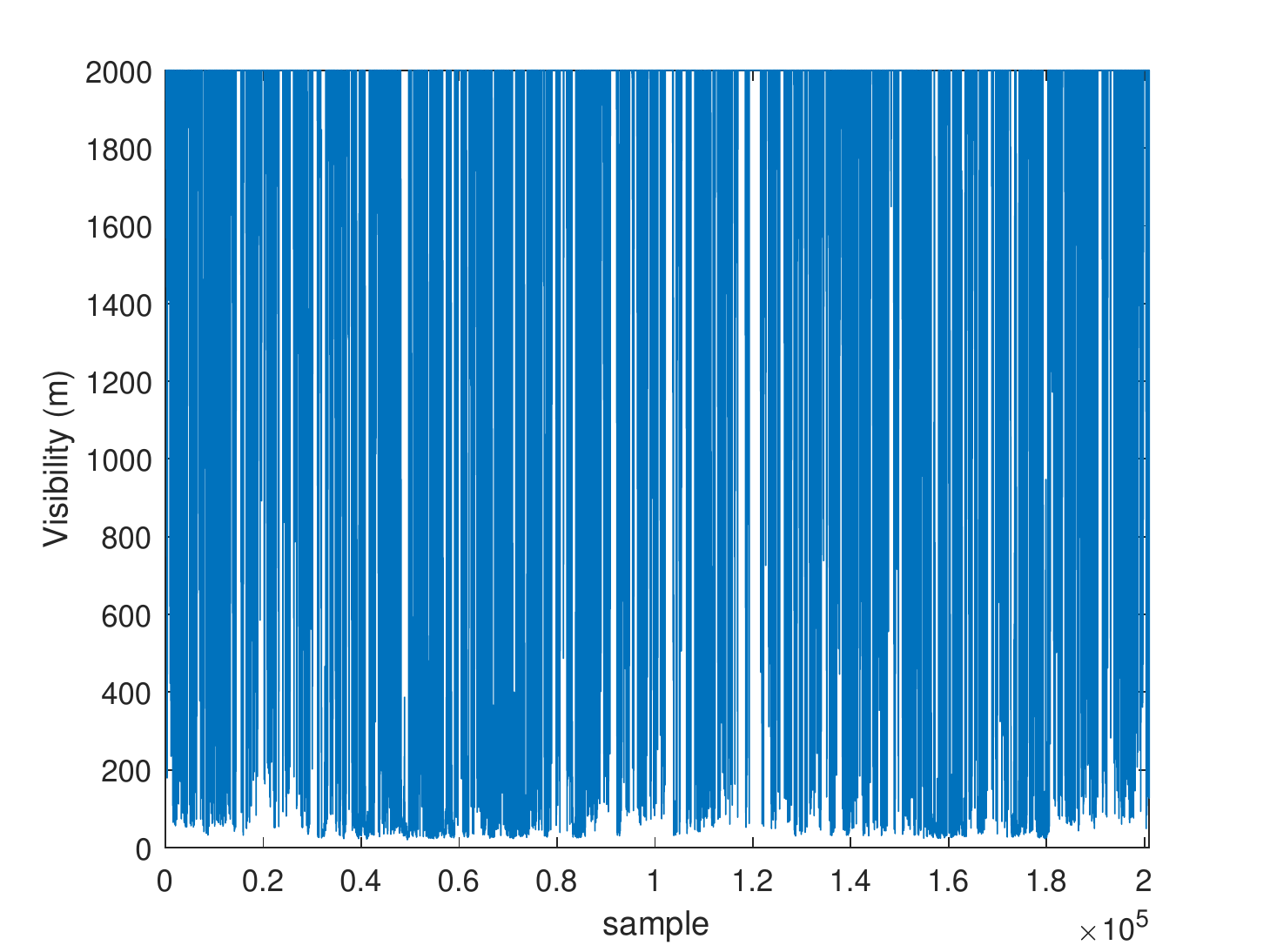}}
\subfigure[~zoomed]{\includegraphics[draft=false, angle=0,width=8.2cm]{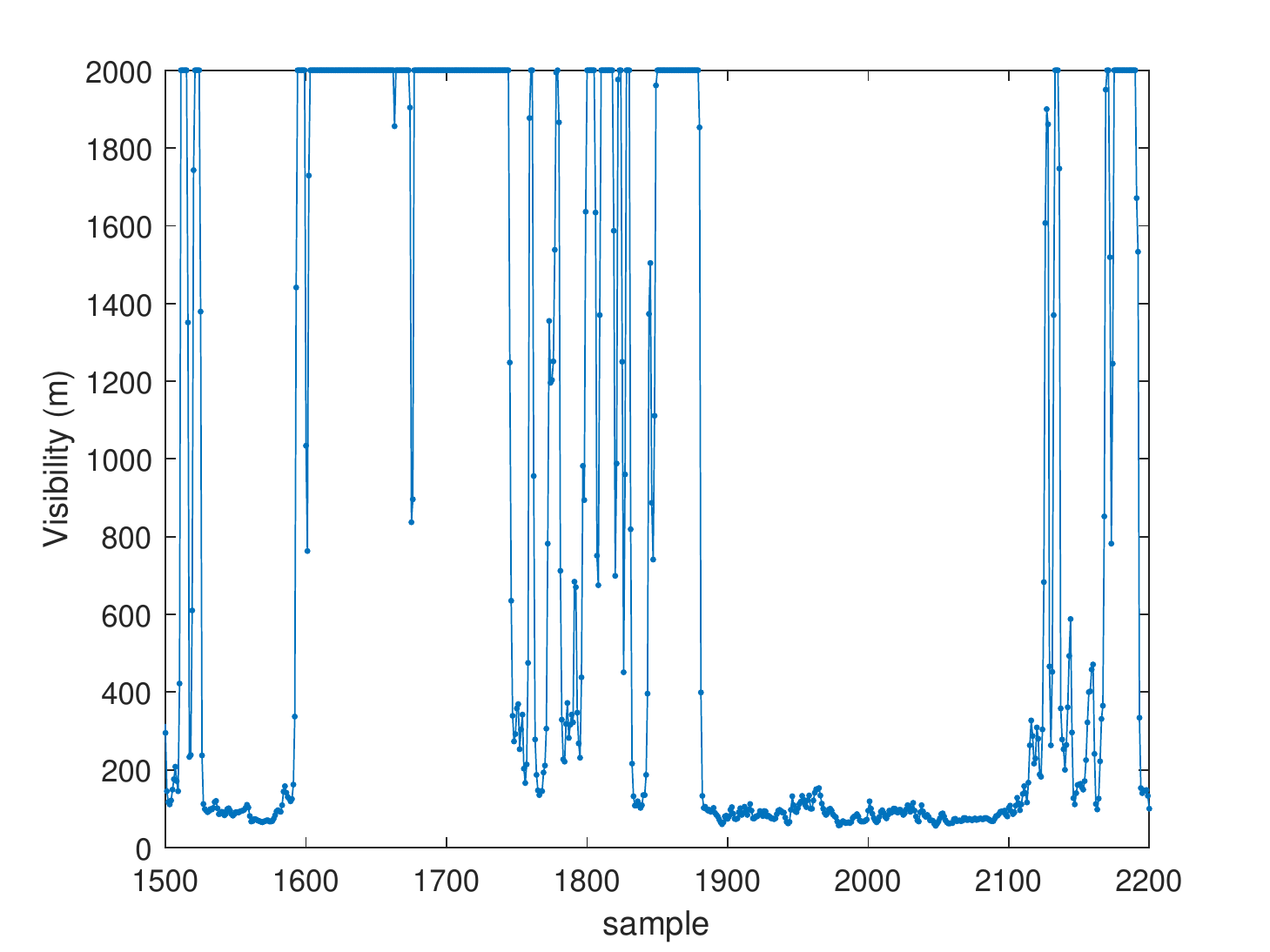}}\\
\end{center}
\caption{\label{A8fog_series} Low-visibility values due to orographic fog at A8 motor-road (Lugo measuring station) with a time resolution of a sample each 5 minutes; (a) Complete time series (23 months); (b) Time series zoom.}
\end{figure}

Let us show some examples of this invariance of characteristic time over time scale changes. We start with an example of DFA application in the analysis of low-visibility events in Northern Spain \cite{Salcedo21}. Specifically, we have collected visibility data from a measuring station at the A8 motor-road, which goes across the Spanish northern coast. This measuring station is located at the province of Lugo, Galicia (43.3841N, 7.3692W), where orographic-type fog events are common, severely affecting the traffic in the motor-road many days every year \cite{cornejo2021statistical}. The measuring station is equipped with a Biral WS-100 visibility sensor. Following the specifications of the visibilimeter, the maximum visibility value is 2000m (full visibility in the zone). A total of 23 months of data with a resolution $a=5$ minutes are available for this study, measured in the station from 1st January 2018 to 30th November 2019. Figure \ref{A8fog_series} shows the low-visibility values at Lugo, with a 5 minutes resolution, and a zoom on the first values of the series.

\begin{figure}[!ht]
\begin{center}
\subfigure[~1 sample/5 min., $s=5$]{\includegraphics[draft=false, angle=0,width=7cm]{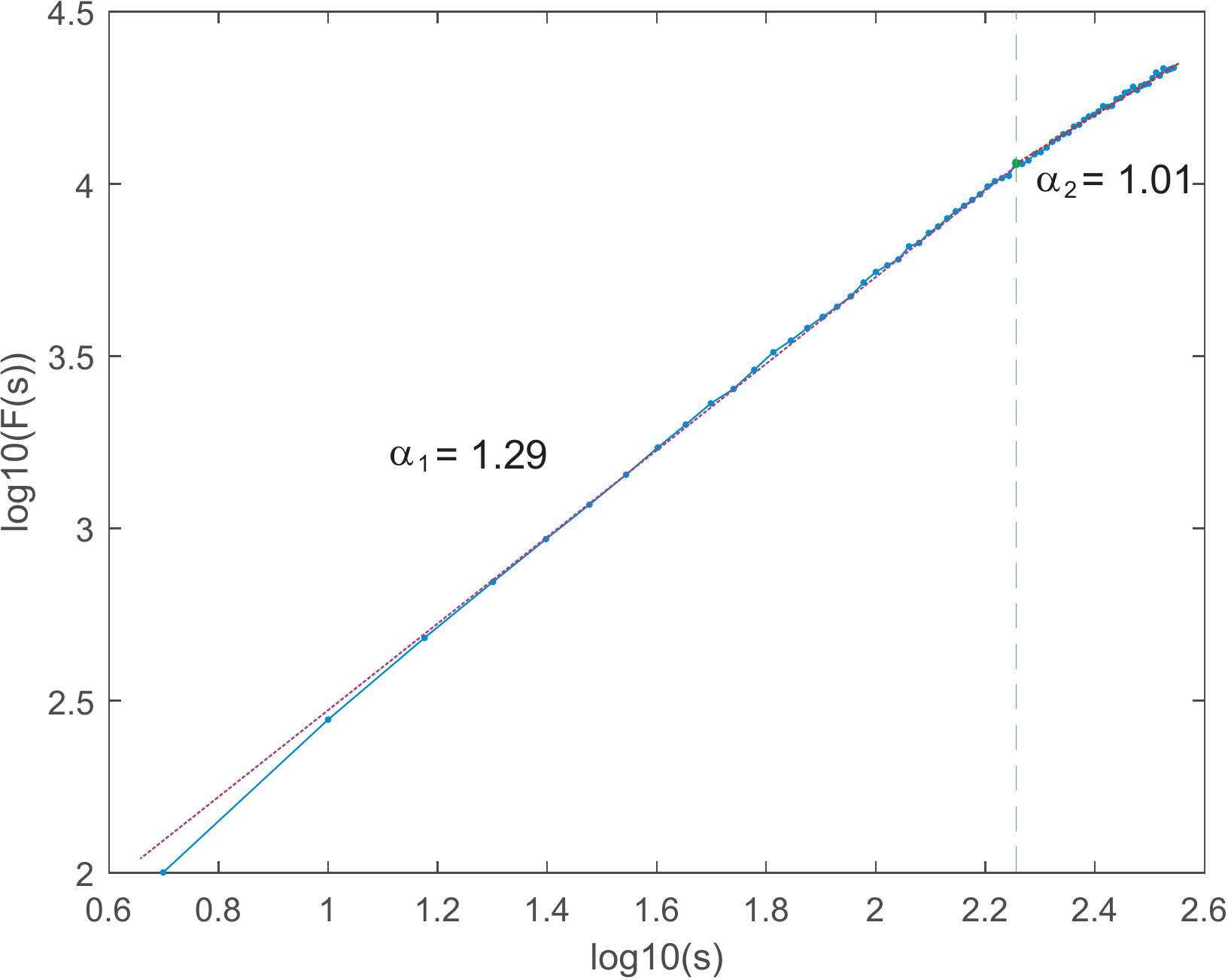}}
\subfigure[~1 sample/15 min., $s=5$]{\includegraphics[draft=false, angle=0,width=7cm]{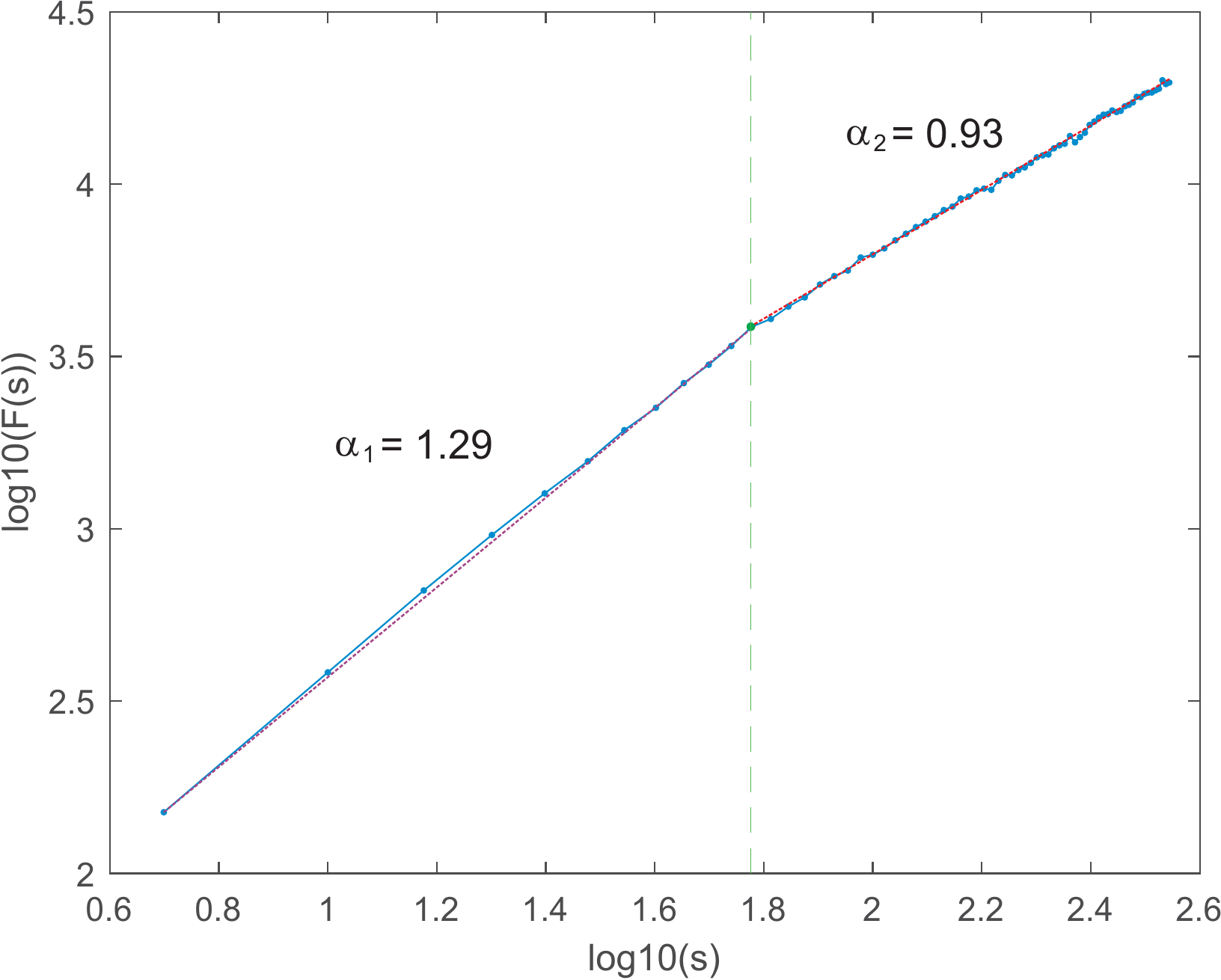}}\\
\subfigure[~1 sample/45 min., $s=5$]{\includegraphics[draft=false, angle=0,width=7cm]{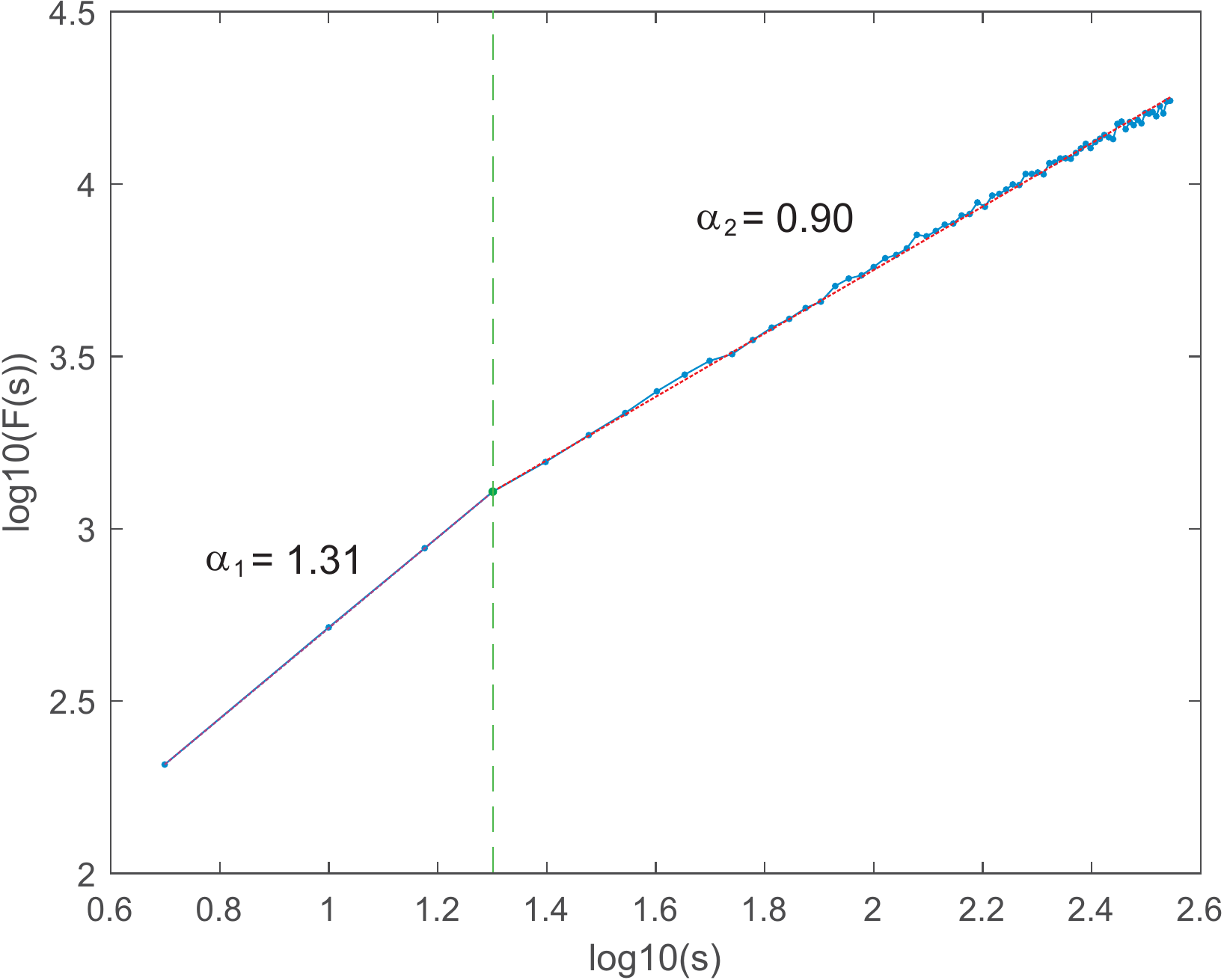}}
\subfigure[~1 sample/60 min., $s=5$]{\includegraphics[draft=false, angle=0,width=7cm]{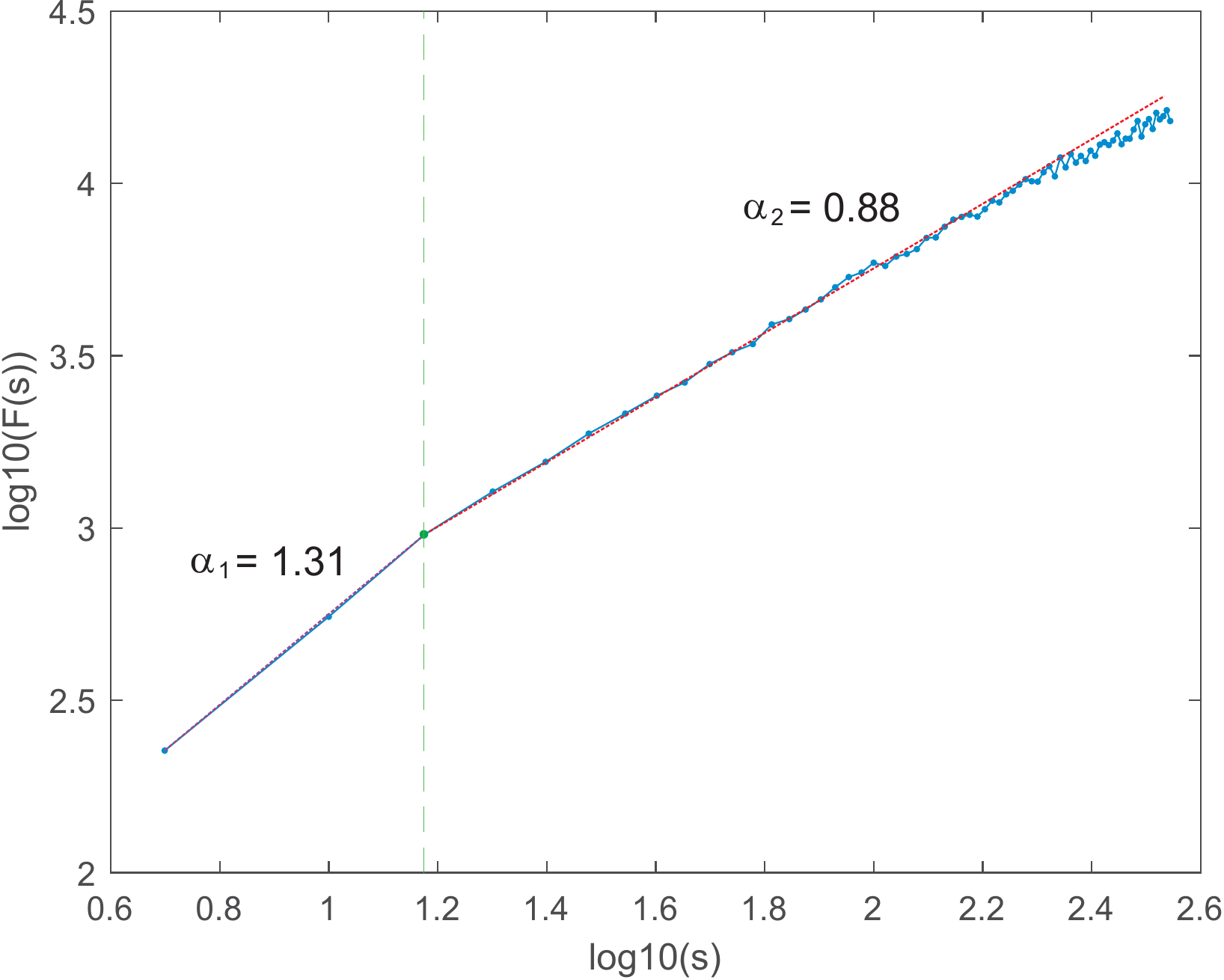}}\\
\end{center}
\caption{\label{A8fog} DFA (Double log plots of $F (s)$ vs $s$) for fog events time series at Lugo measuring station (orographic fog), considering different time resolutions and $s=5$; (a) (a sample each) 5 min.; (b) 15 min.; (c) 45 min.; (d) 60 min.}
\end{figure}

We can apply the time scale changing procedure shown above in order to obtain visibility data at different resolutions ($b$), of one sample each 15, 30, 45, 60 and 75 minutes. Figure \ref{A8fog} shows the results of applying a DFA algorithm to the time series of visibility, with different time resolution and $s=5$. Note that the two-range pattern can be detected in the DFA result in all cases, with a characteristic time $\tau^* \approx 15$ hours independently of the time scale resolution taken into account (invariance of the characteristic time). In all cases the scaling exponents $\alpha_1$ and $\alpha_2$ show a strong long-term persistence of the low-visibility time series, before and after the characteristic time $\tau^*$.

\begin{figure}[!ht]
\begin{center}
\subfigure[~Temperature original]{\includegraphics[draft=false, angle=0,width=8cm]{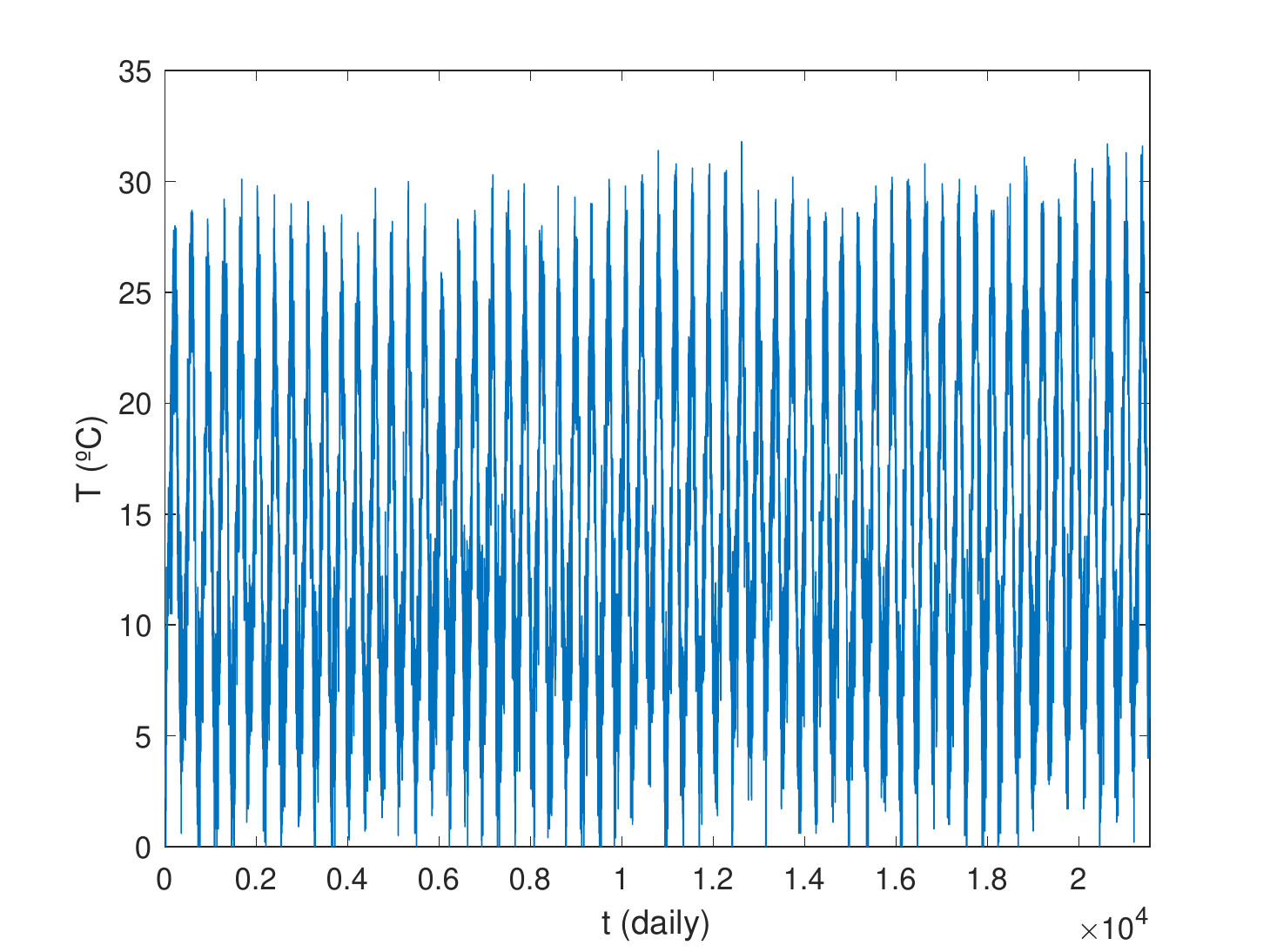}}
\subfigure[~Temperature zoomed]{\includegraphics[draft=false, angle=0,width=8cm]{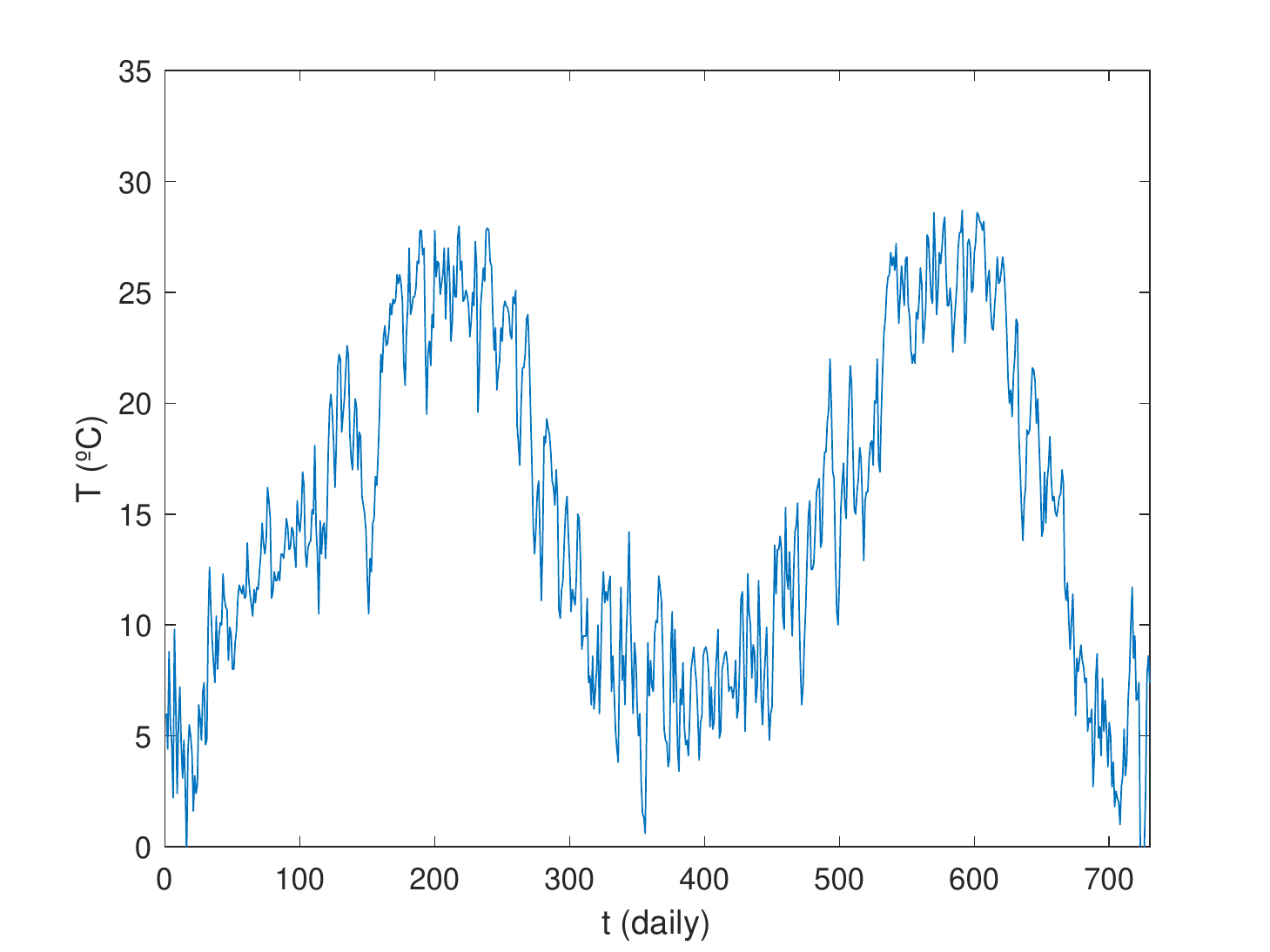}}\\
\subfigure[~Pressure original]{\includegraphics[draft=false, angle=0,width=8cm]{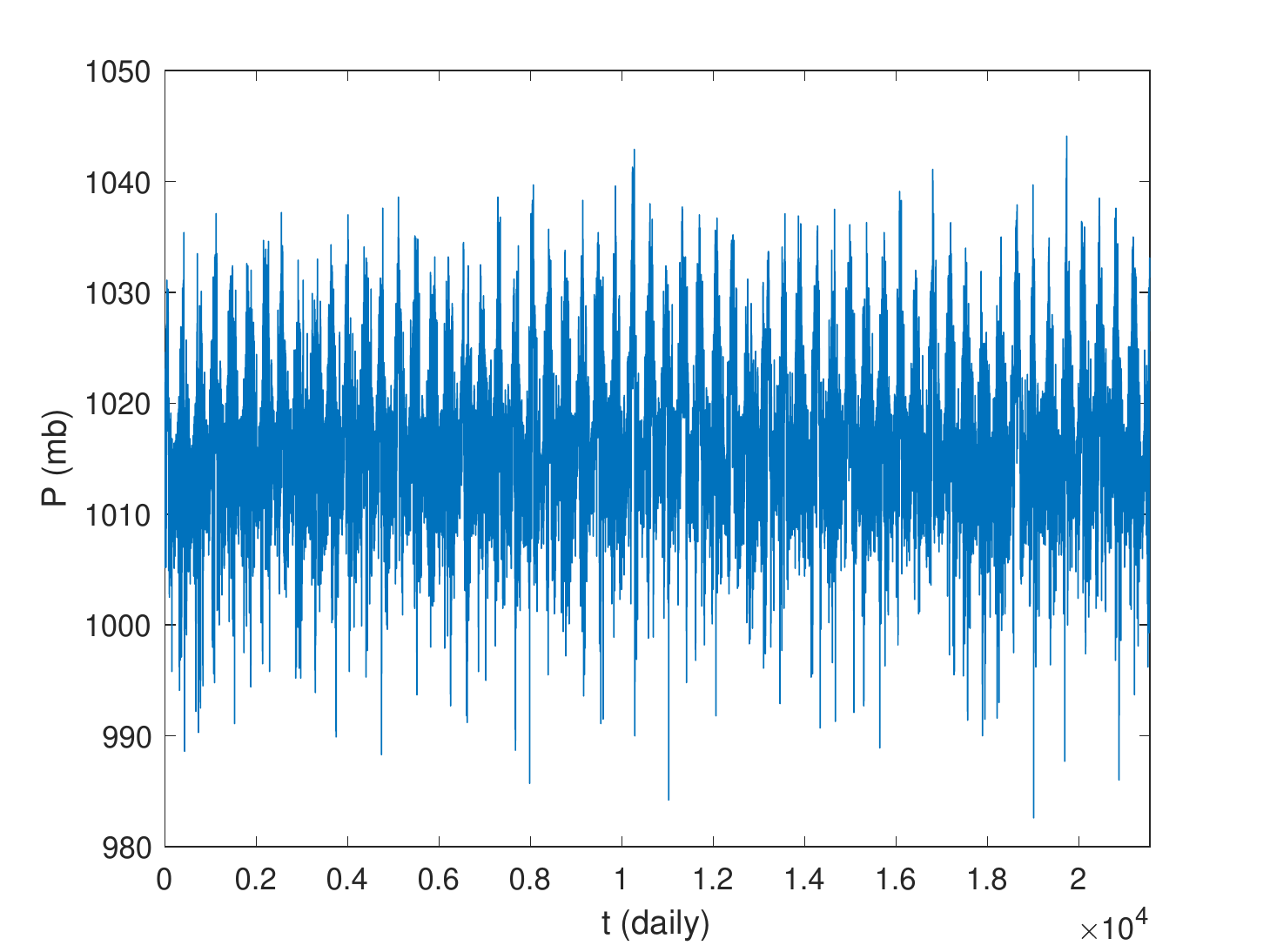}}
\subfigure[~Pressure zoomed]{\includegraphics[draft=false, angle=0,width=8cm]{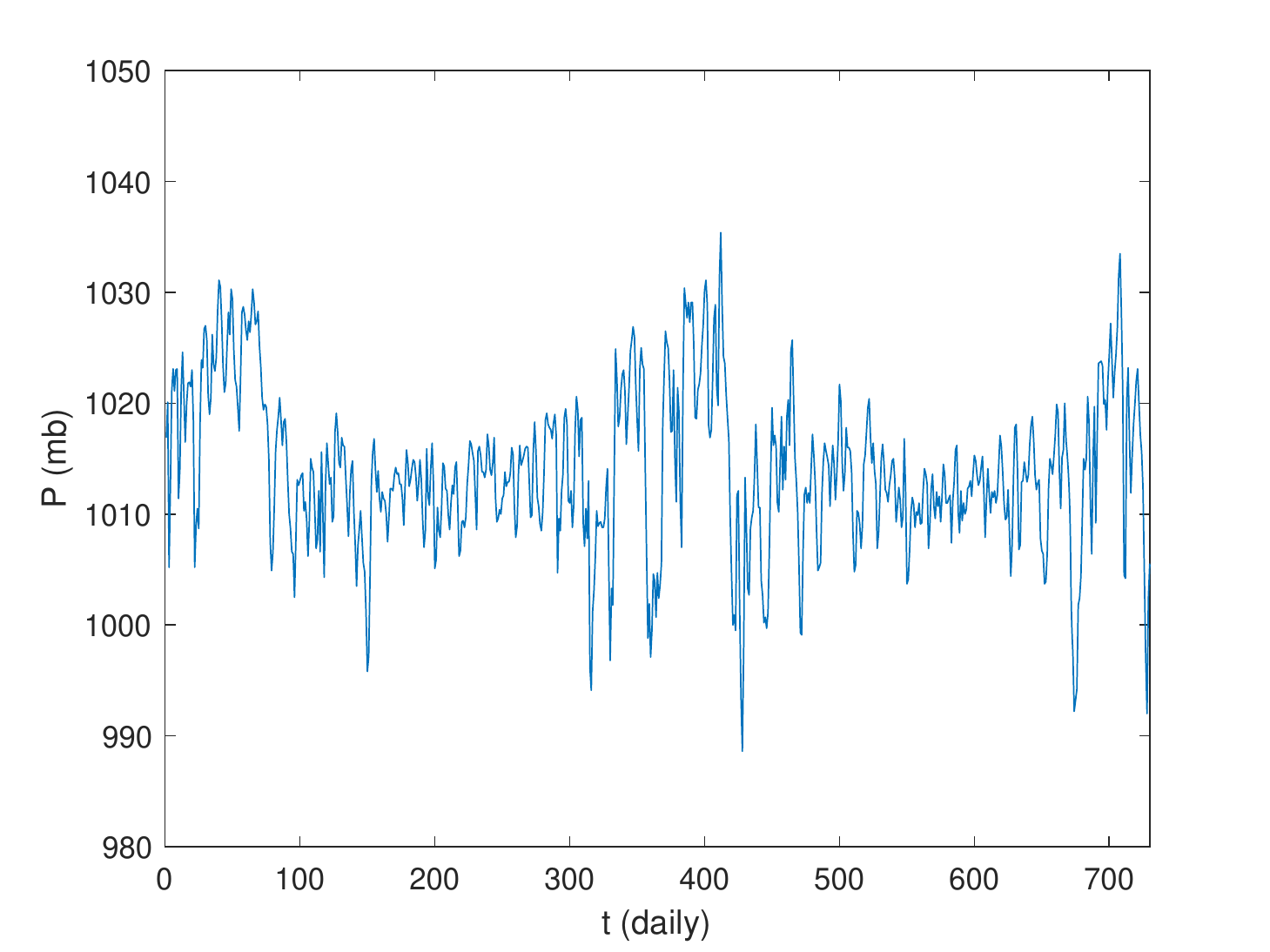}}\\
\end{center}
\caption{\label{PBarajas_series} Daily average temperature and pressure values values at Madrid-Barajas airport in the period 2006-2017); (a) Daily average temperature; (b) Zoom of the first values of the temperature series.; (c) Daily pressure; (d) Zoom of the first values of the pressure series.}
\end{figure}

We have also tested the application of the DFA algorithm to alternative meteorological data, such as temperature, pressure, wind speed and precipitation. The data of temperature and pressure have been obtained from the historic data at Madrid-Barajas airport. Temperature and pressure data consist of 11 years of daily averaged measurements at the airport (2006-2017). Figure \ref{PBarajas_series} shows the time series of daily average temperature and pressure available and a zoom of the first samples in the series.

\begin{figure}[!ht]
\begin{center}
\subfigure[~Temp., sample/day]{\includegraphics[draft=false, angle=0,width=5cm]{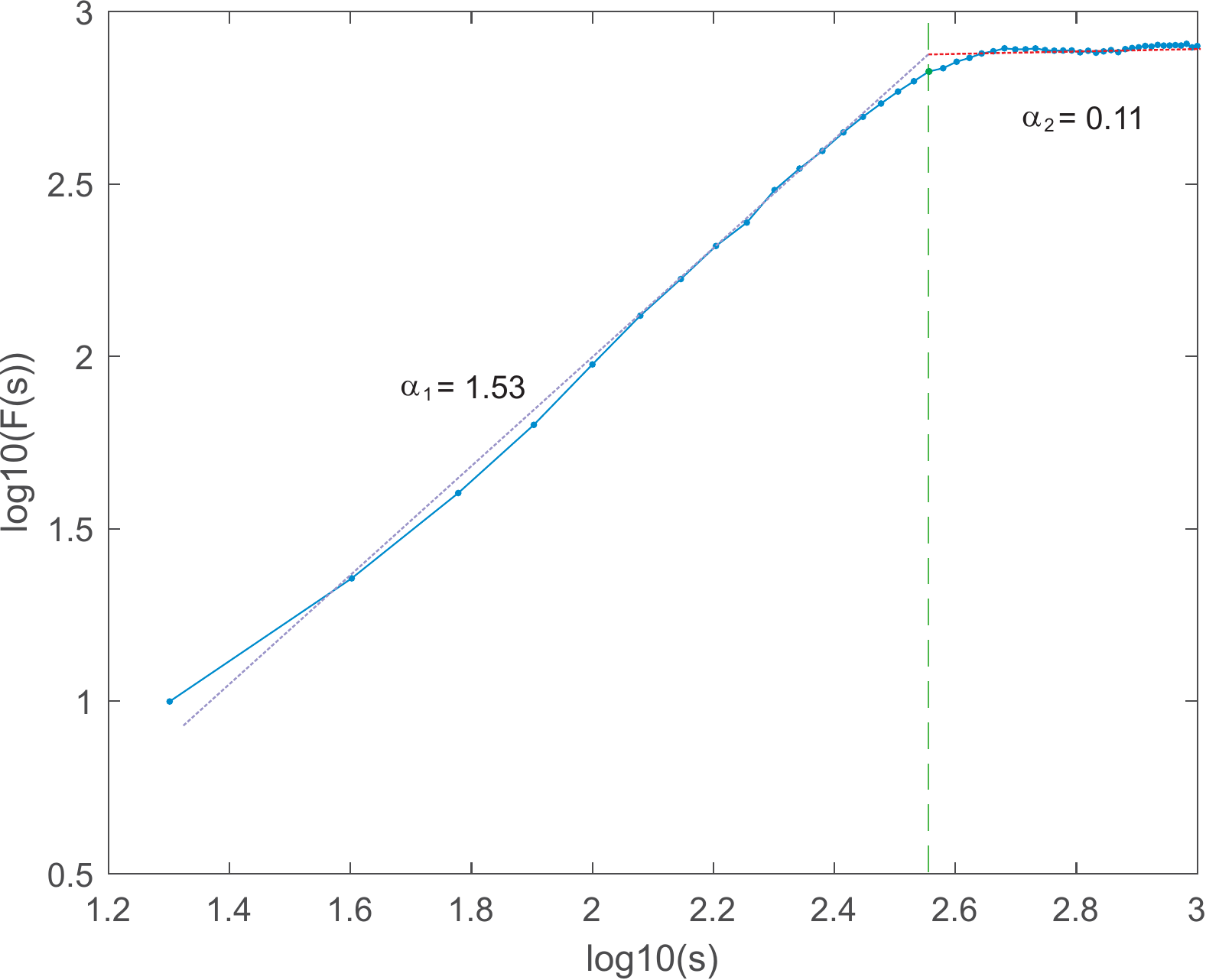}}
\subfigure[~Temp., sample/3 days]{\includegraphics[draft=false, angle=0,width=5cm]{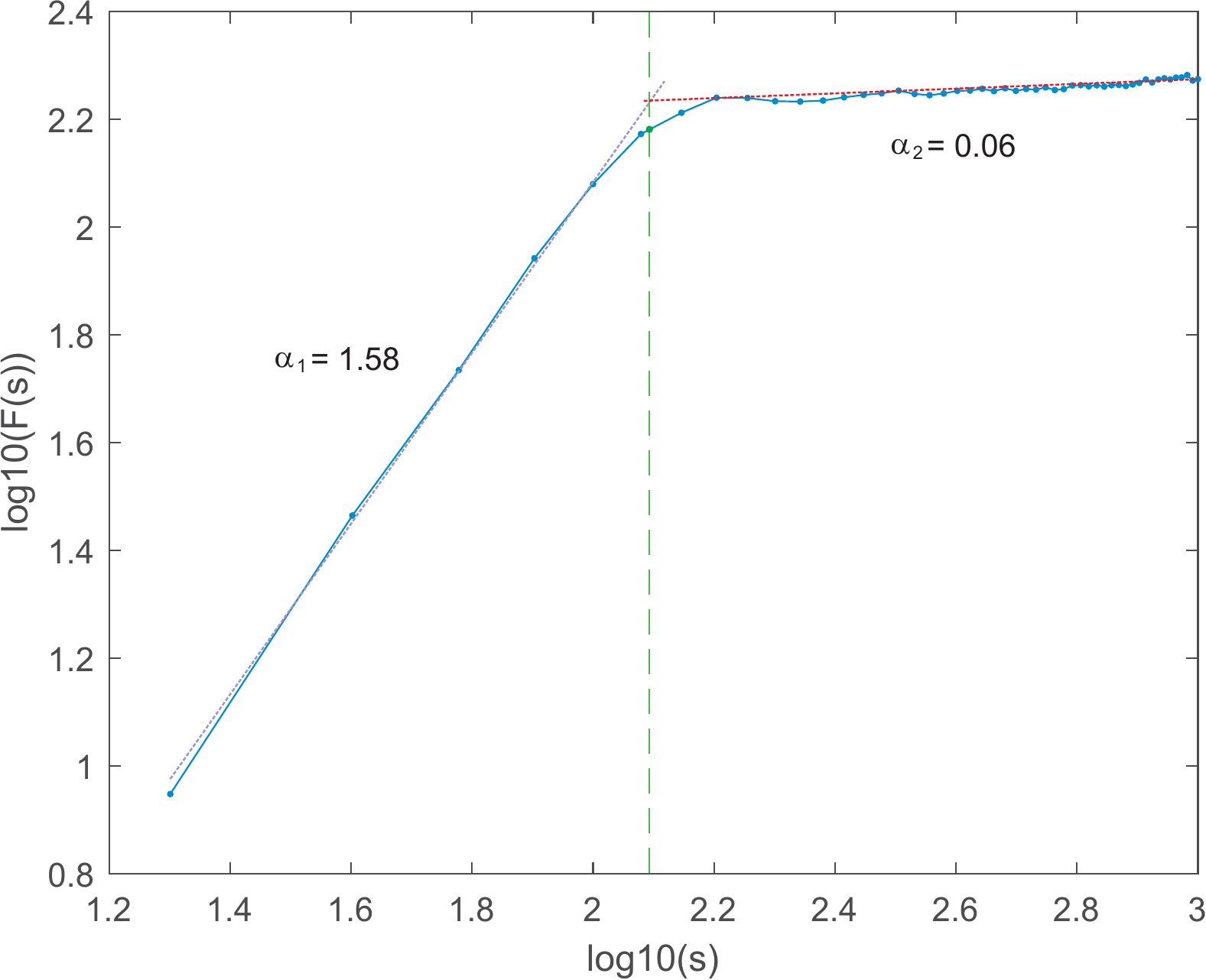}}
\subfigure[~Temp., sample/6 days]{\includegraphics[draft=false, angle=0,width=5cm]{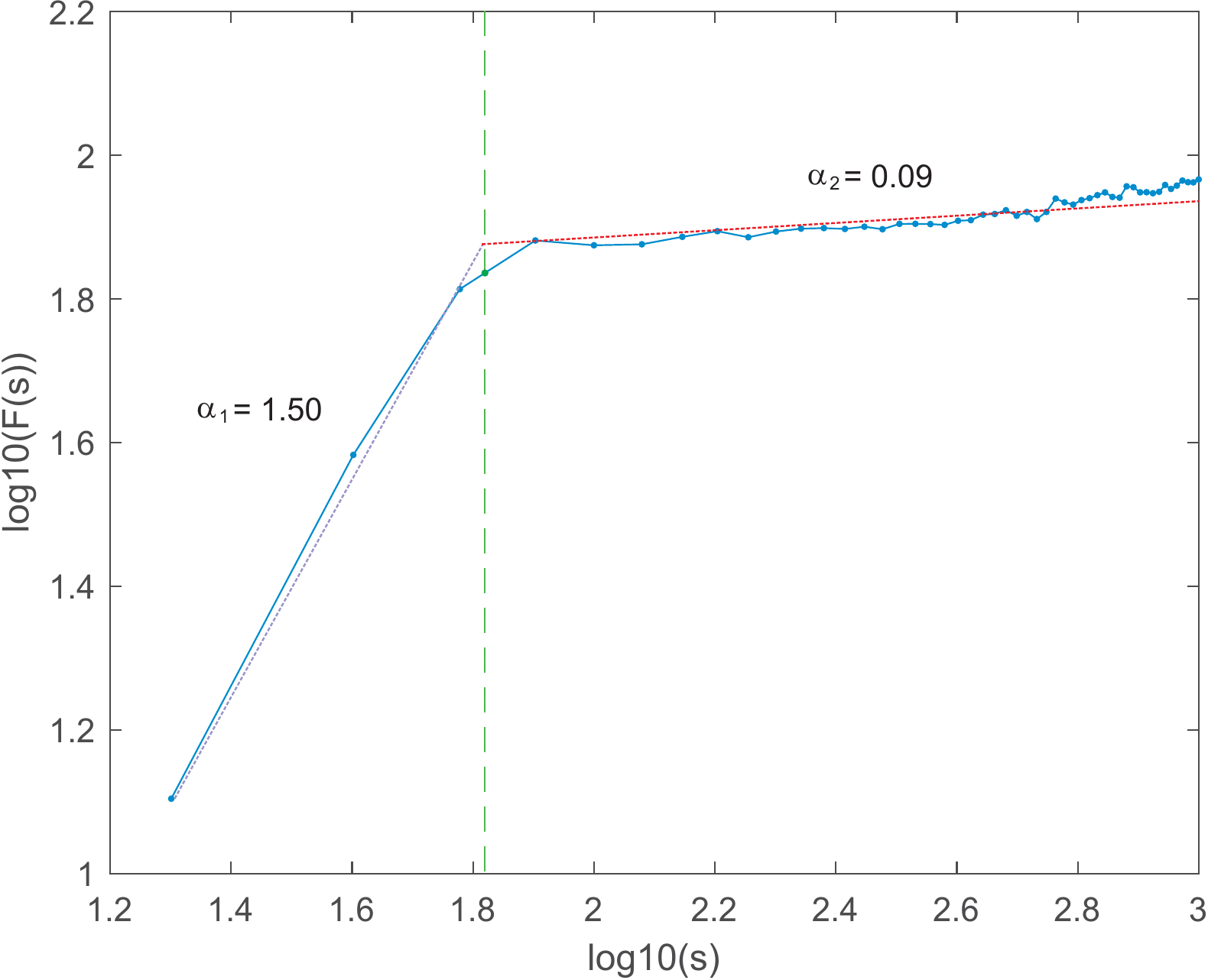}}
\subfigure[~Press., sample/day]{\includegraphics[draft=false, angle=0,width=5cm]{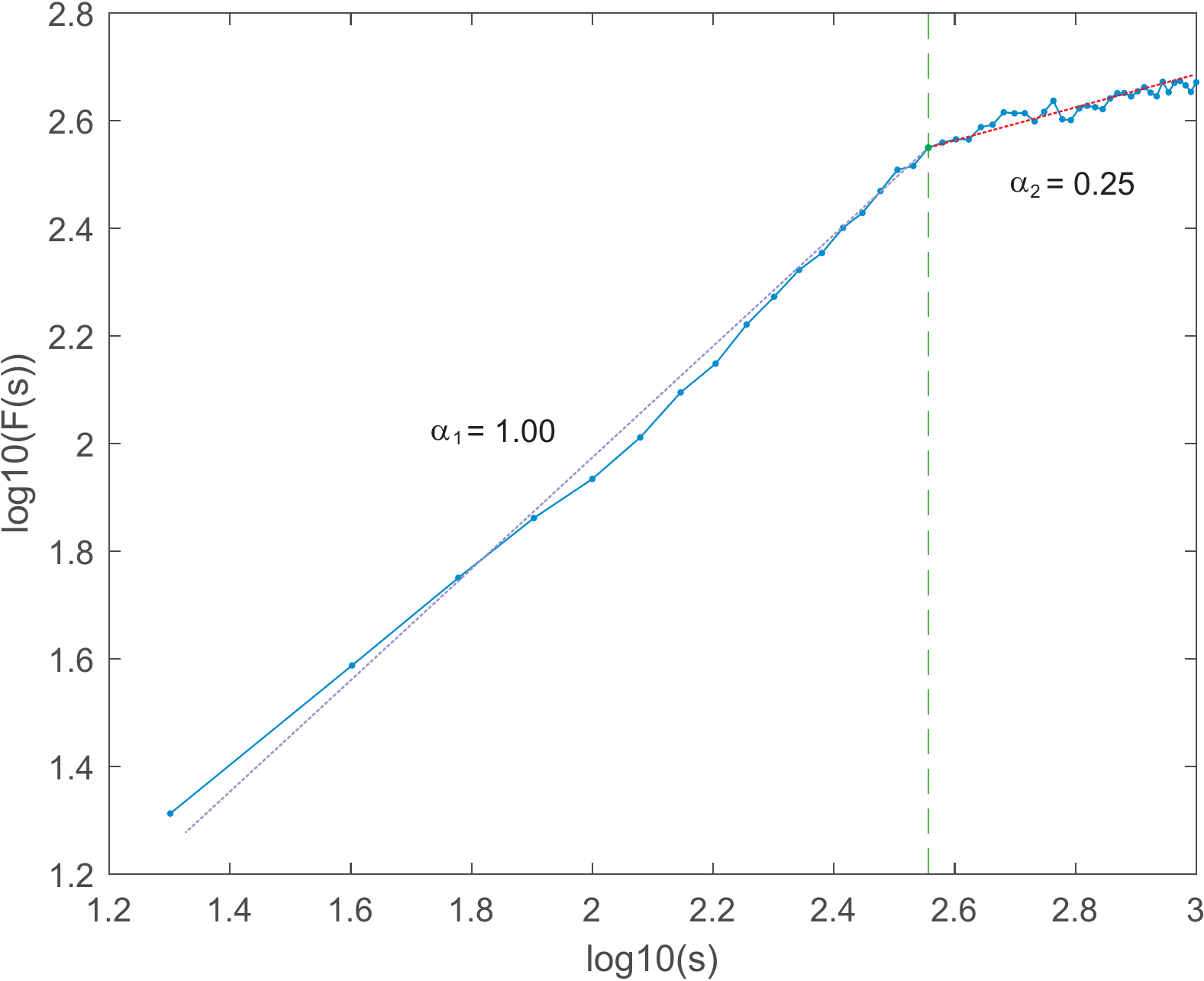}}
\subfigure[~Press., sample/3 days]{\includegraphics[draft=false, angle=0,width=5cm]{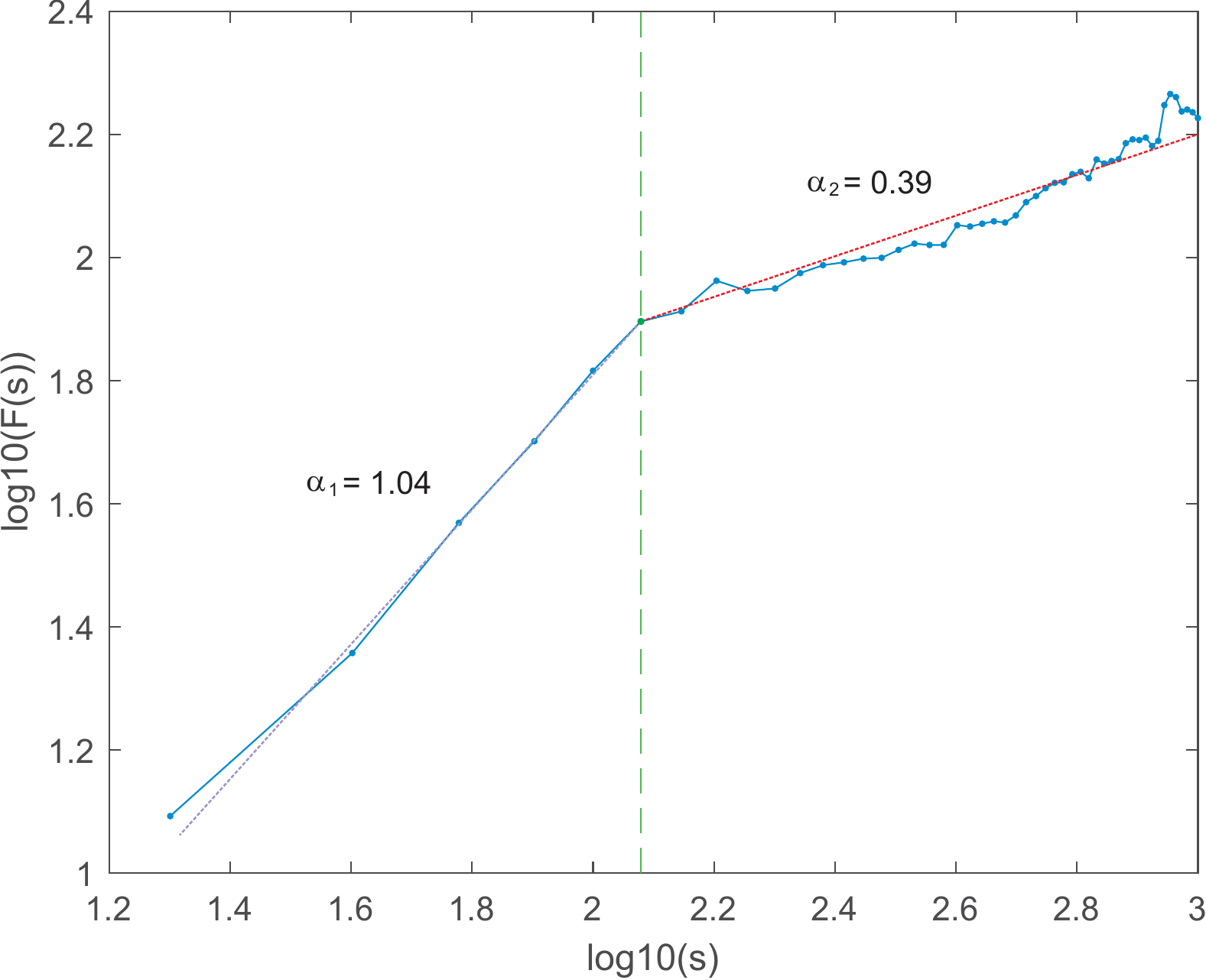}}
\subfigure[~Press., sample/6 days]{\includegraphics[draft=false, angle=0,width=5cm]{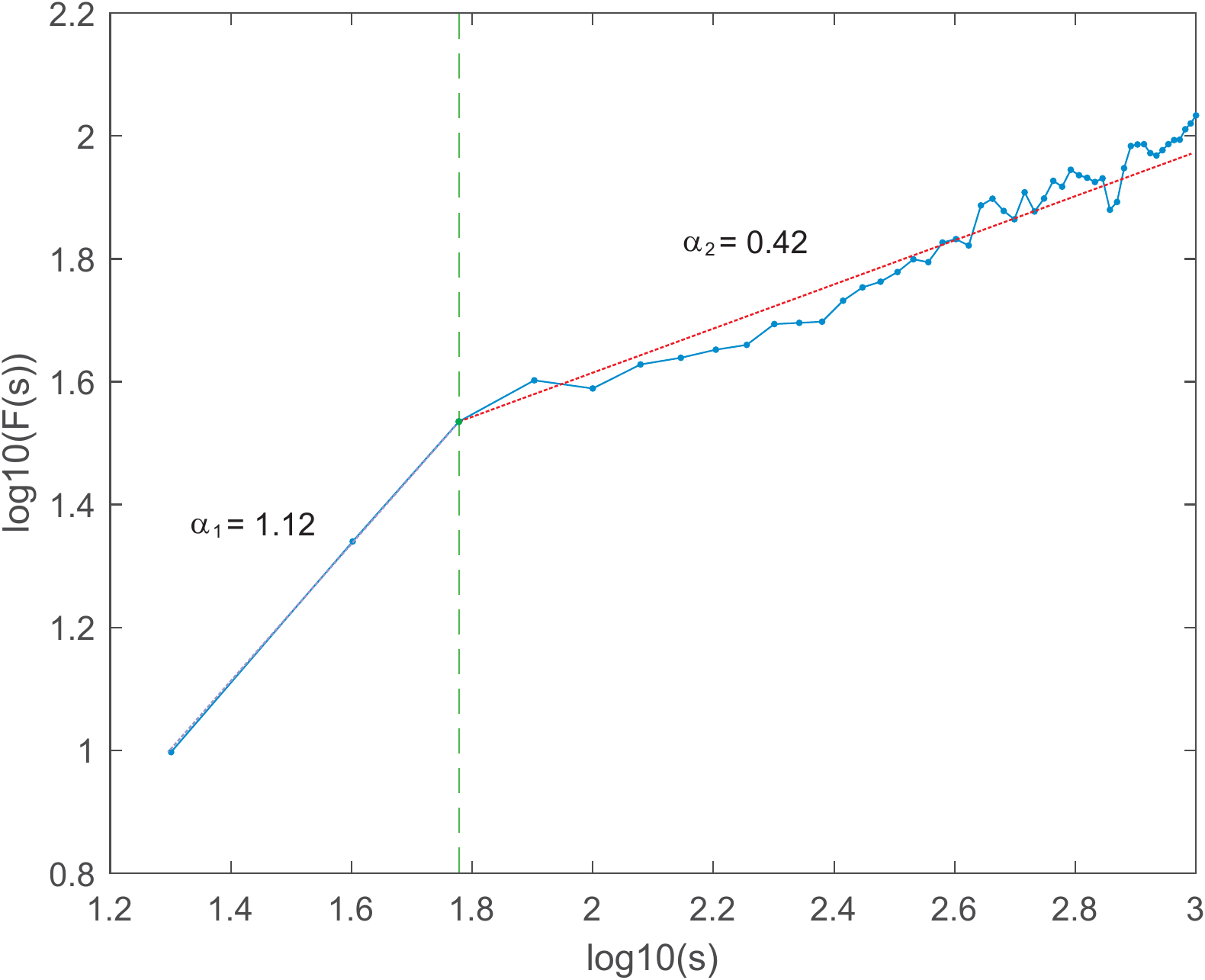}}
\end{center}
\caption{\label{TBarajas} DFA (Double log plots of $F (s)$ vs $s$) for the air temperature time and pressure series at Madrid-Barajas airport,  considering different time resolutions and $s=20$; (a) Temperature sample/day; (b) Temperature sample/3 days; (c) Temperature sample/6 days; (d) Pressure sample/day; (b) Pressure sample/3 days; (c) Pressure sample/6 days.}
\end{figure}

Figure \ref{TBarajas} shows the DFA application (Double log plots of $F (s)$ vs $s$) to the air temperature and pressure time series, considering different time resolutions (daily, 3 days and 6 days). In the case of average temperature we obtain a clear two ranges form, with a clear characteristic time $\tau^*$ of one year. As can be seen, the characteristic time is invariant over time scale changes. In this case, the first part of the curve, before the characteristic time point shows a clear line with $\alpha_1$ around 1.5, which points to a marked persistence of the series up to the characteristic time $\tau^*$. From this point, $\alpha_2 \approx 0$, which indicates an anti-persistent behaviour. In this case, the DFA shows a clear annual persistence pattern of the air temperature time series, as expected. In the case of pressure time series, we again obtain a clear two-ranges form, with a characteristic time $\tau^*$ of one year, and it is possible to see that the characteristic time is invariant over time scale changes. In this case of atmospheric pressure, the first part of the curve, before the characteristic time point shows a line with $\alpha_1$ around 1.0, which again points to a marked persistence of the series up to the characteristic time $\tau^*$. Beyond this characteristic time, $\alpha_2$ takes values from 0.25 to 0.4, which indicates an anti-persistent behaviour, but less marked than for the air temperature time series. The DFA analysis shows a clear annual persistence pattern for the atmospheric pressure, but slightly less pronounced than the one obtained for the air temperature case.

\begin{figure}[!ht]
\begin{center}
\subfigure{\includegraphics[draft=false, angle=0,width=8cm]{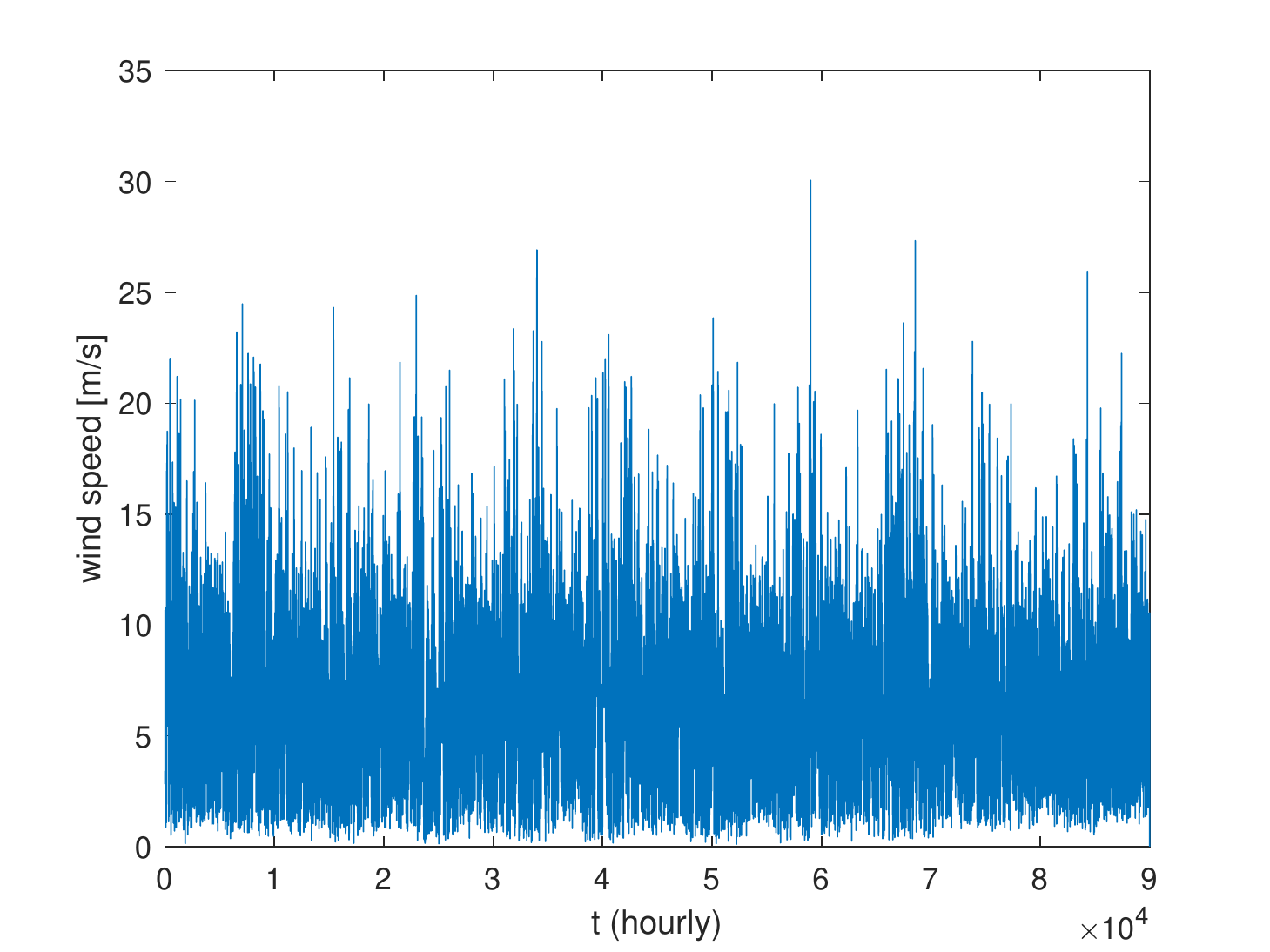}}
\end{center}
\caption{\label{Viento_series} Hourly wind speed time series from a measurement tower in Burgos, Spain.}
\end{figure}

\begin{figure}[!ht]
\begin{center}
\subfigure[~sample/3 hours, $s=5$]{\includegraphics[draft=false, angle=0,width=5cm]{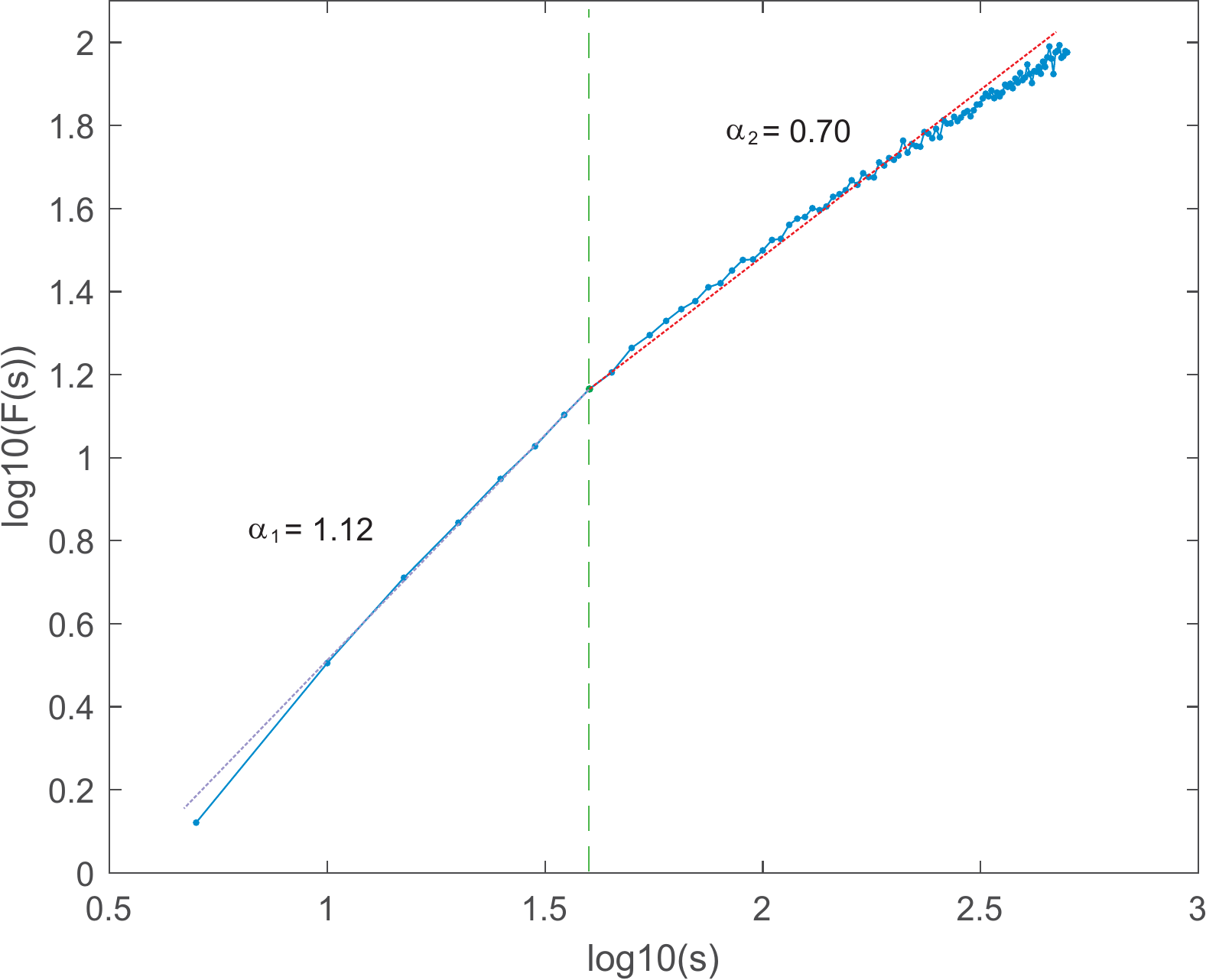}}
\subfigure[~sample/6 hours, $s=5$]{\includegraphics[draft=false, angle=0,width=5cm]{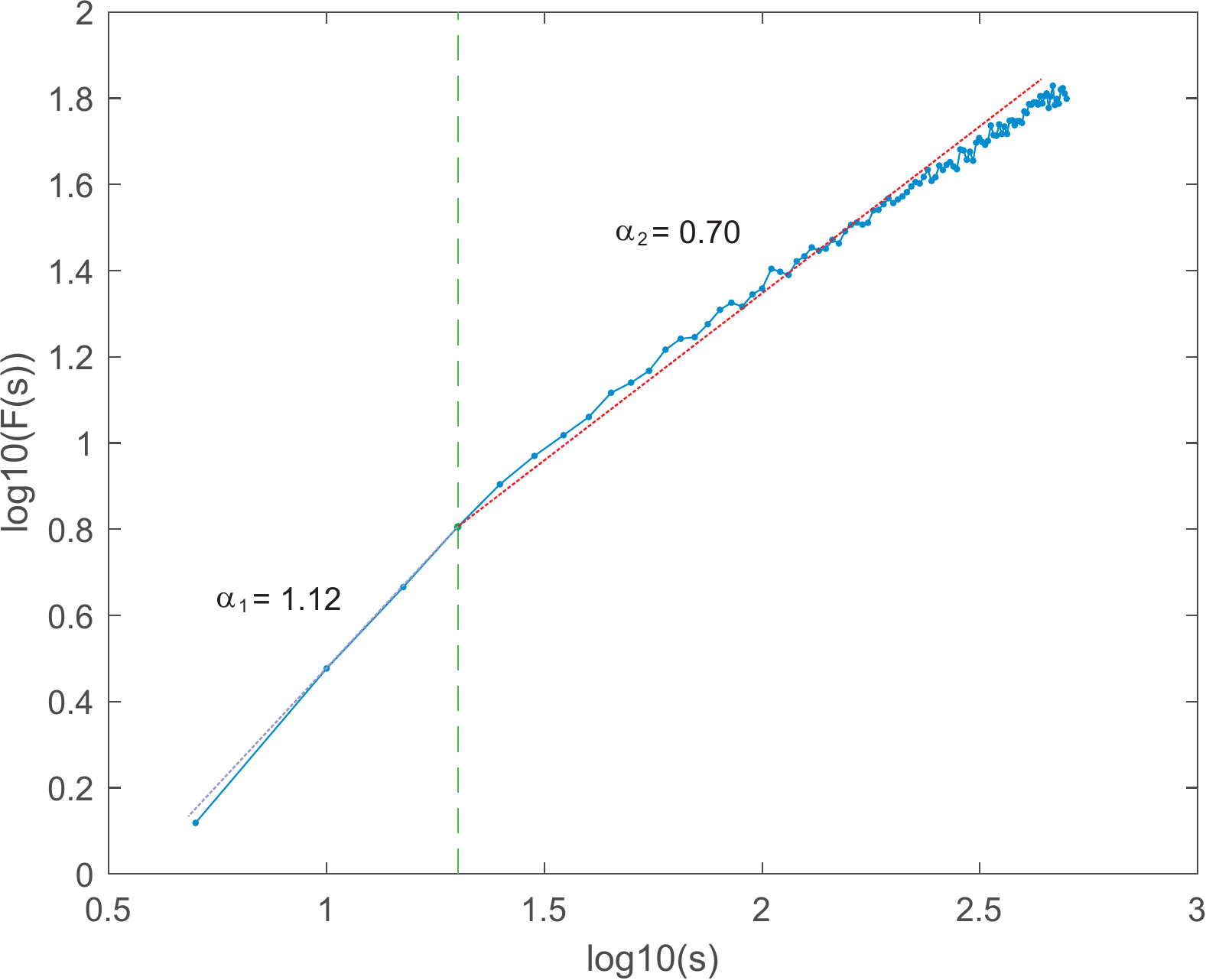}}
\subfigure[~sample/12 hours, $s=5$]{\includegraphics[draft=false, angle=0,width=5cm]{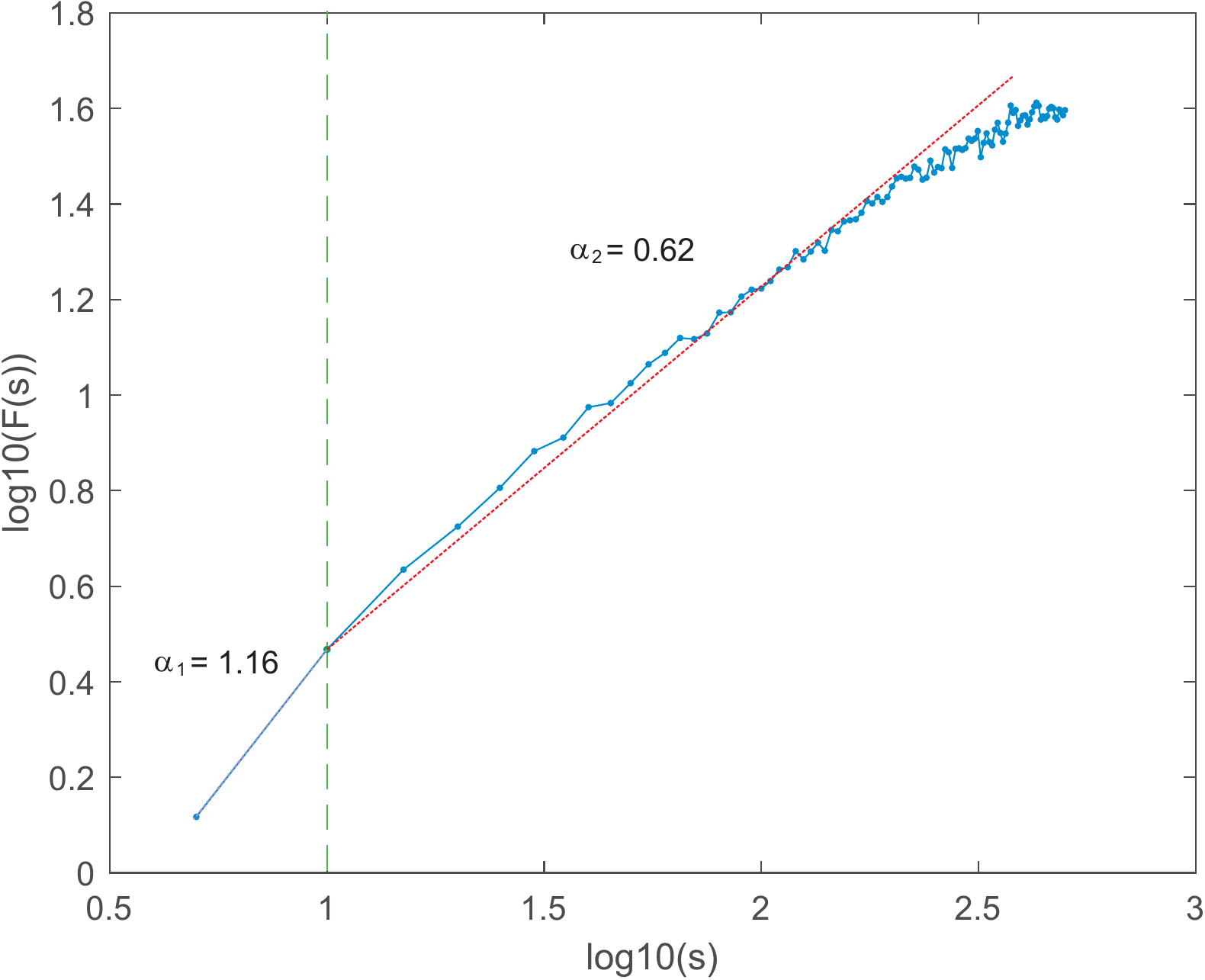}}
\end{center}
\caption{\label{Viento_Mnon} DFA (Double log plots of $F (s)$ vs $s$) for the wind speed time series at Burgos measurement station,  considering different time resolutions and $s=5$; (a) (a sample each) 3 hours; (b) 6 hours; (c) 12 hours.}
\end{figure}

We will further show the two-ranges DFA structure in time series of wind speed and also in precipitation data. We first consider wind speed data in a measurement tower located in Burgos, Spain (42.48N,3.80W). We have 10 years of hourly data, measured in the tower from 1st November 2002 to 29th October 2012. Figure \ref{Viento_series} shows the wind speed time series obtained. Figure \ref{Viento_Mnon} shows the DFA application (Double log plots of $F (s)$ vs $s$) to this time series. As can be seen, the DFA plot shows a two-ranges graph, with a characteristic time of 120 hours, invariant over time scale changes. In this case $\alpha_1 \approx 1.1$ and $\alpha_2 \approx 0.7$ which indicates a strong long-term persistence of this time series.

\begin{figure}[!ht]
\begin{center}
\subfigure{\includegraphics[draft=false, angle=0,width=8cm]{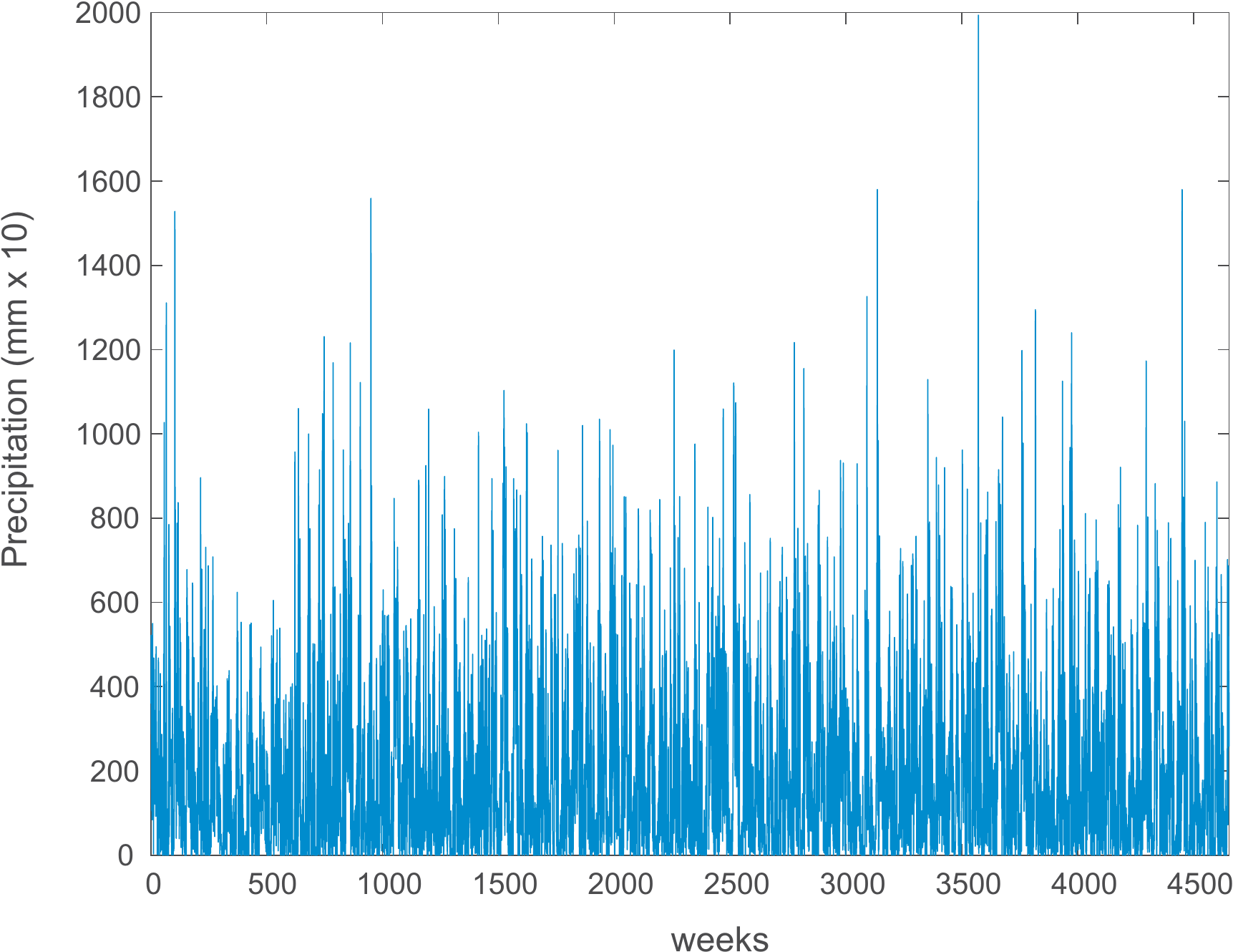}}
\end{center}
\caption{\label{Precipitation_series} Weekly averaged precipitation time series in A Coru\~na observatory, Galicia, Spain.}
\end{figure}

\begin{figure}[!ht]
\begin{center}
\subfigure[~sample/week, $s=5$]{\includegraphics[draft=false, angle=0,width=5cm]{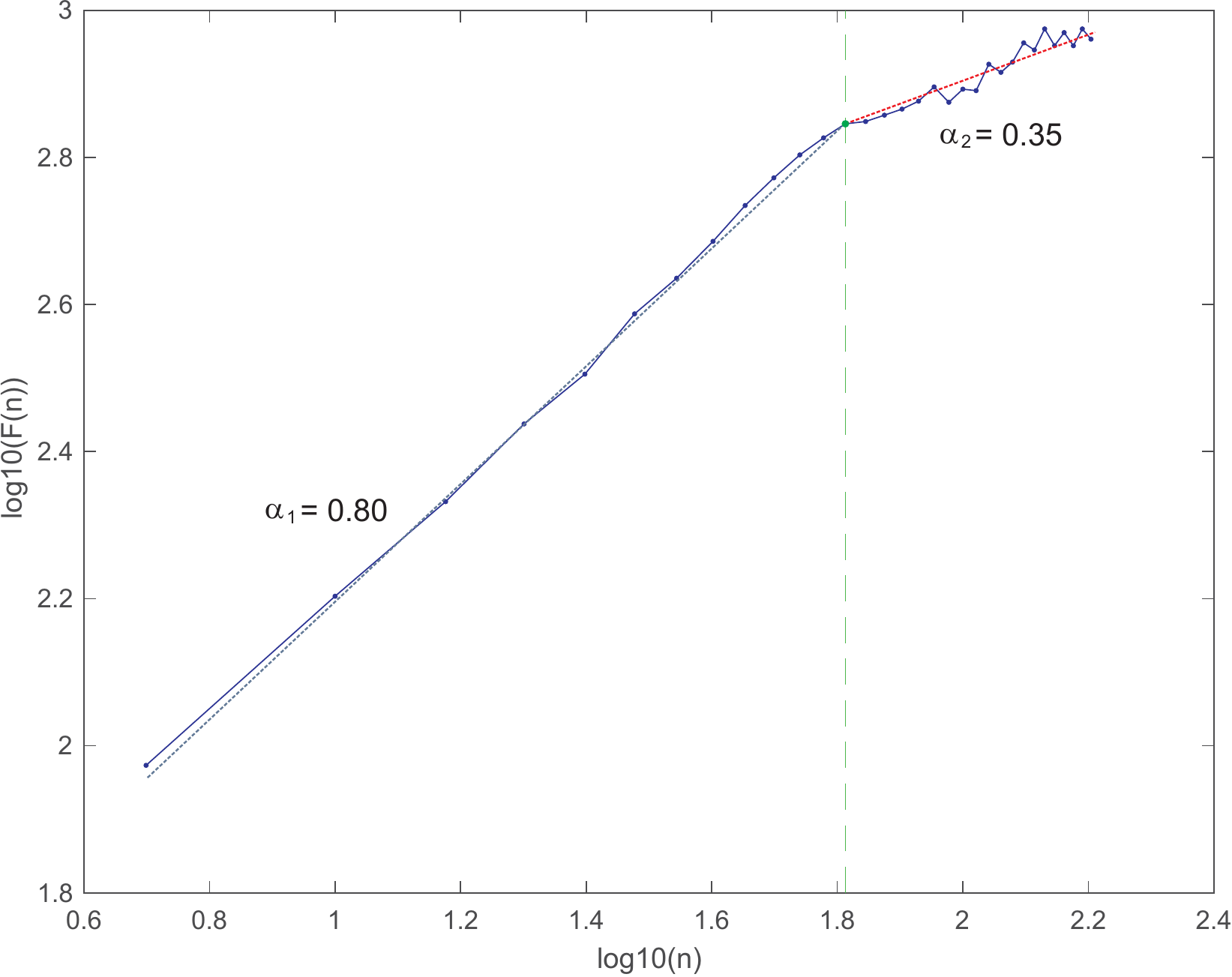}}
\subfigure[~sample/2 weeks, $s=5$]{\includegraphics[draft=false, angle=0,width=5cm]{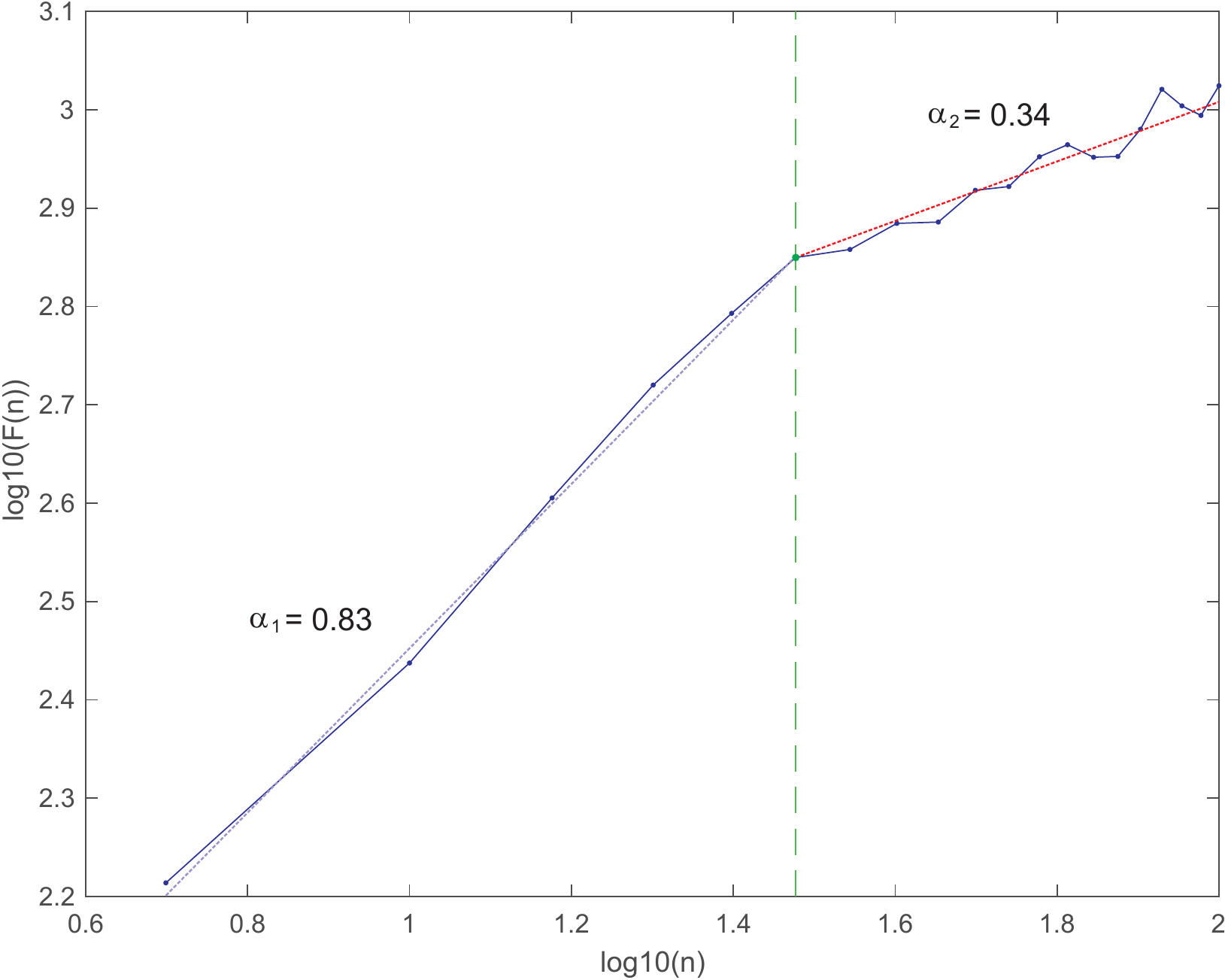}}
\subfigure[~sample/month, $s=5$]{\includegraphics[draft=false, angle=0,width=5cm]{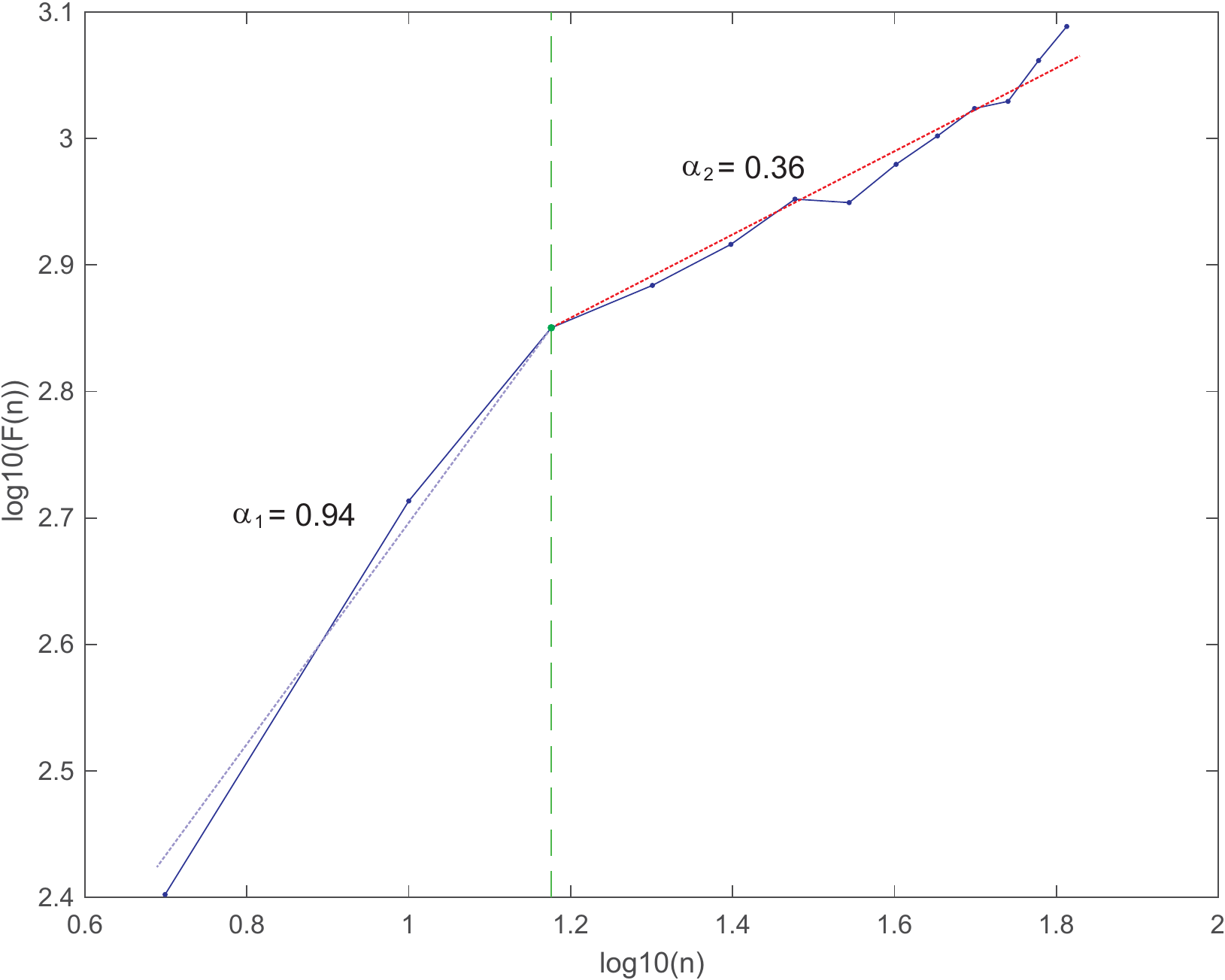}}
\end{center}
\caption{\label{Precipitation_Coruna} DFA (Double log plots of $F (s)$ vs $s$) for the precipitation time series at A Coru\~na observatory, Spain, considering different time resolutions and $s=5$; (a) (a sample each) 1 week; (b) 2 weeks; (c) 1 month.}
\end{figure}

Finally, we show an example of DFA application over precipitation data from the A Coru\~na observatory, Galicia, Spain. A very long measurement range of continuous data, from 1931 to 2019 is available. Figure \ref{Precipitation_series} shows the precipitation time series considered. The result of the DFA application (Double log plots of $F (s)$ vs $s$) to this precipitation series can be observed in Figure \ref{Precipitation_Coruna}. As in the previous cases, a two-ranges structure is shown, with $\alpha_1$ around 0.8/0.9, $\alpha_2$ around 0.35 and characteristic time around 15 months, time scale invariant, as obtained in all the atmospheric-related processes previously analyzed with the DFA algorithm. 

\subsection{Solar irradiation prediction with persistence-based models}

Solar energy is a clean, extremely abundant and renewable energy \cite{Ghimirea19}, with very low environmental impact. Solar is expected to be the future main renewable resource in the world, and its development is currently very strong in countries such as Australia, Middle-east, Southern Europe and Northern Africa, where the solar resource can be better exploited all year around \cite{Kalogirou14}. Maybe the main issue of solar energy (common to the majority of renewable energy sources) is that it is an energy resource intrinsically stochastic, so significant variations in solar energy production occur due to the presence of clouds, atmospheric dust or particles \cite{sovacool2009intermittency,zhou2018roles}. Thus, solar power prediction is extremely important to ensure a correct integration of this renewable source in the electric system.

\begin{figure}[!ht]
\begin{center}
\subfigure[~Reykjavik]{\includegraphics[draft=false, angle=0,width=6cm]{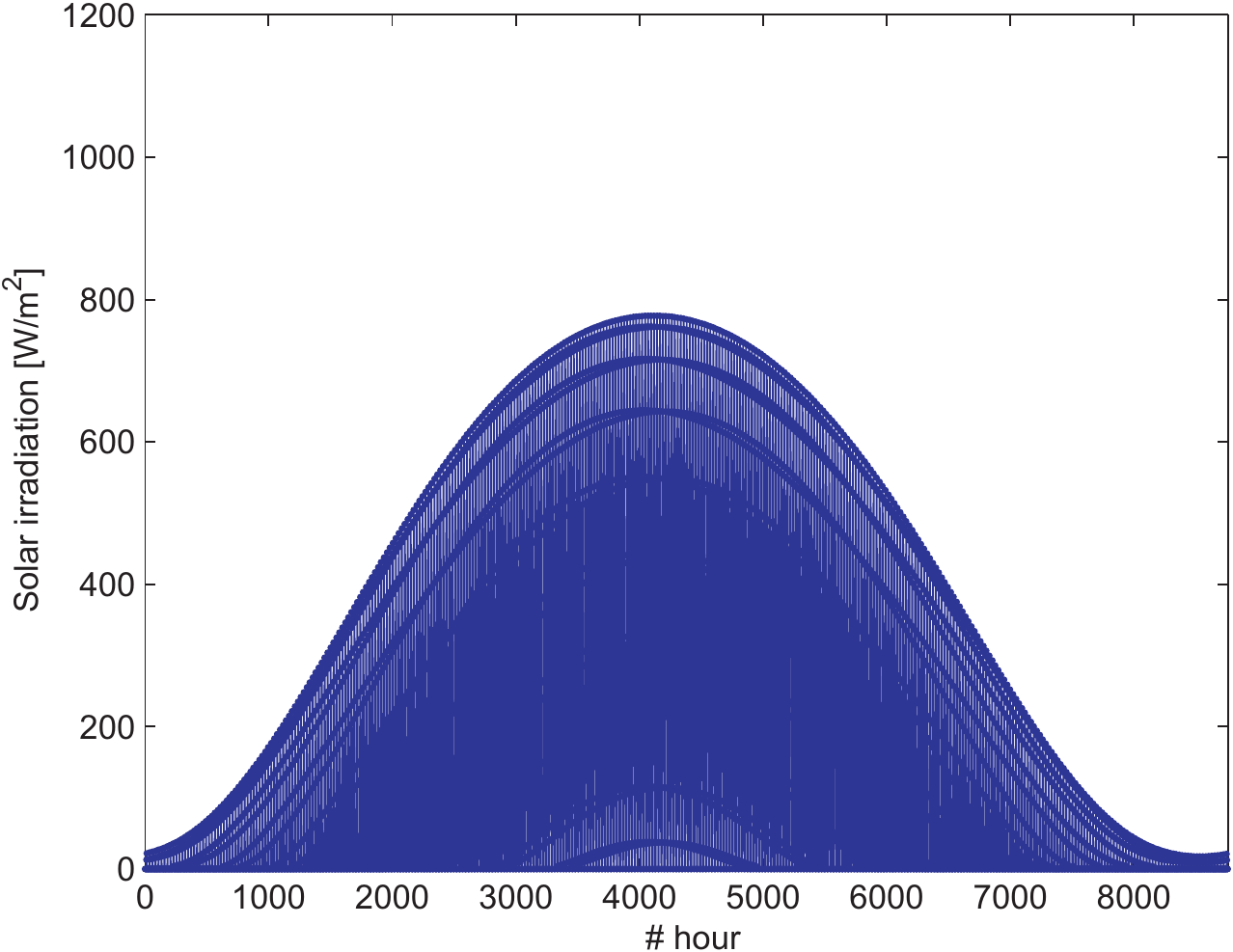}}
\subfigure[~Toledo]{\includegraphics[draft=false, angle=0,width=6cm]{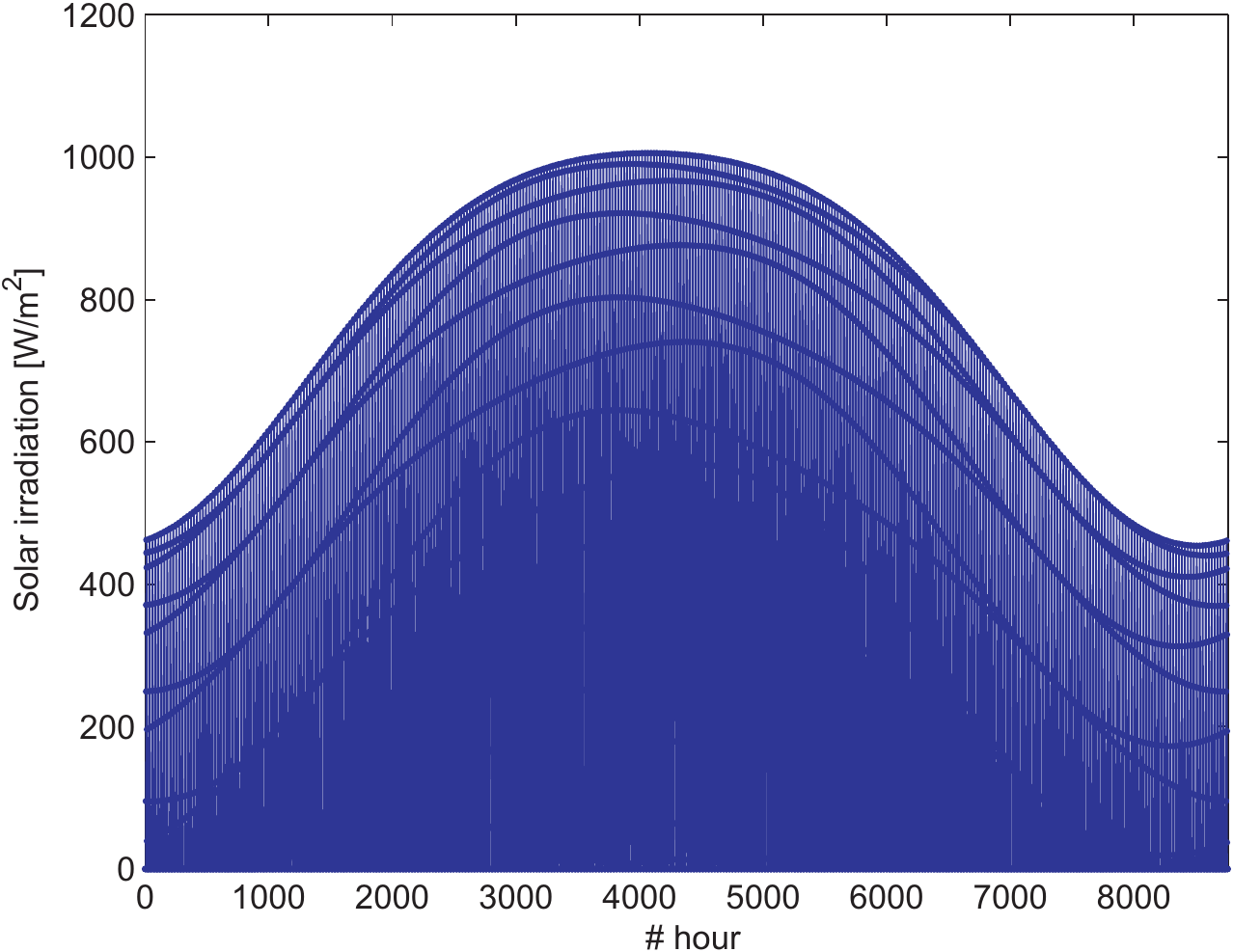}}\\
\subfigure[~Sydney]{\includegraphics[draft=false, angle=0,width=6cm]{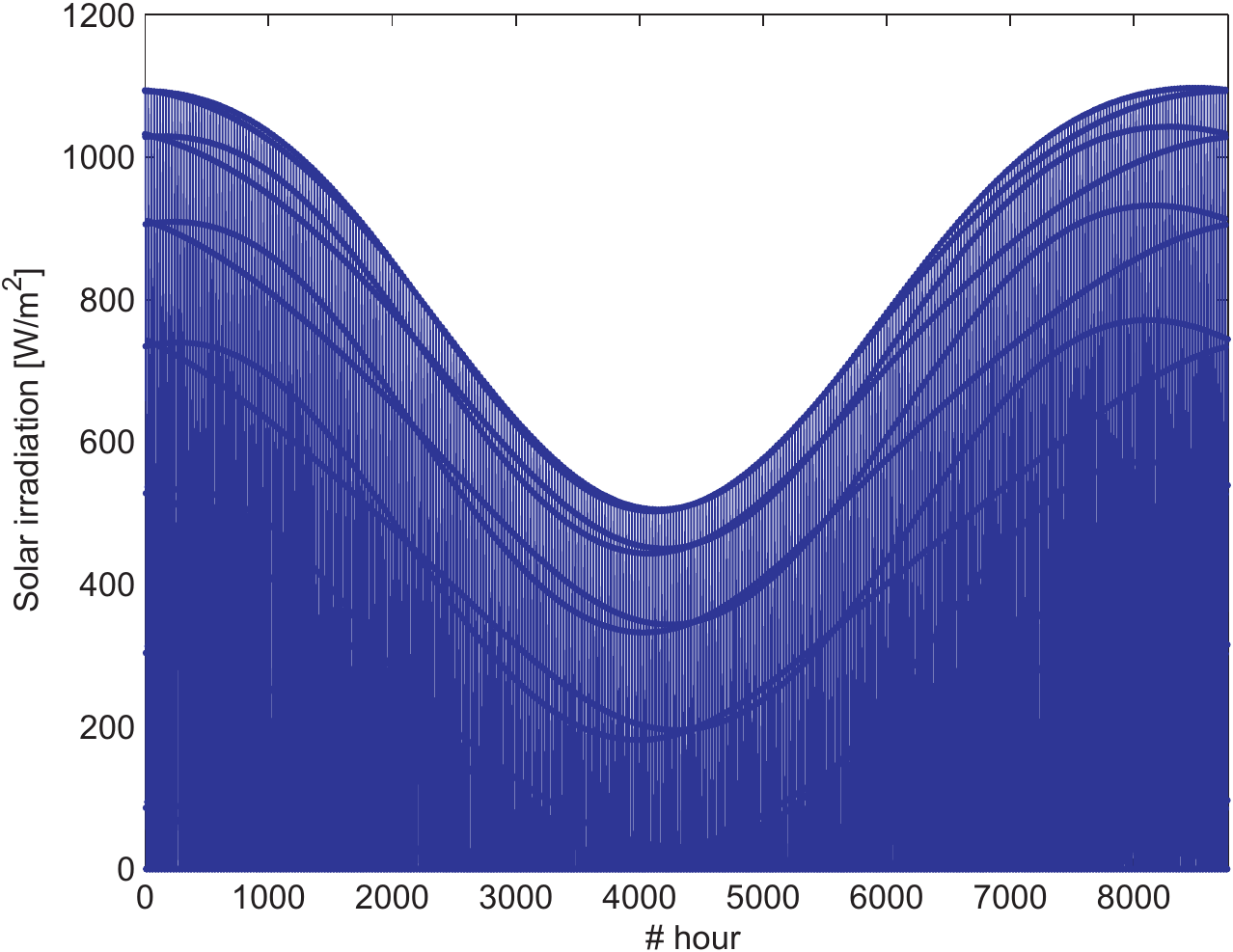}}
\subfigure[~Gulf of Guinea]{\includegraphics[draft=false, angle=0,width=6cm]{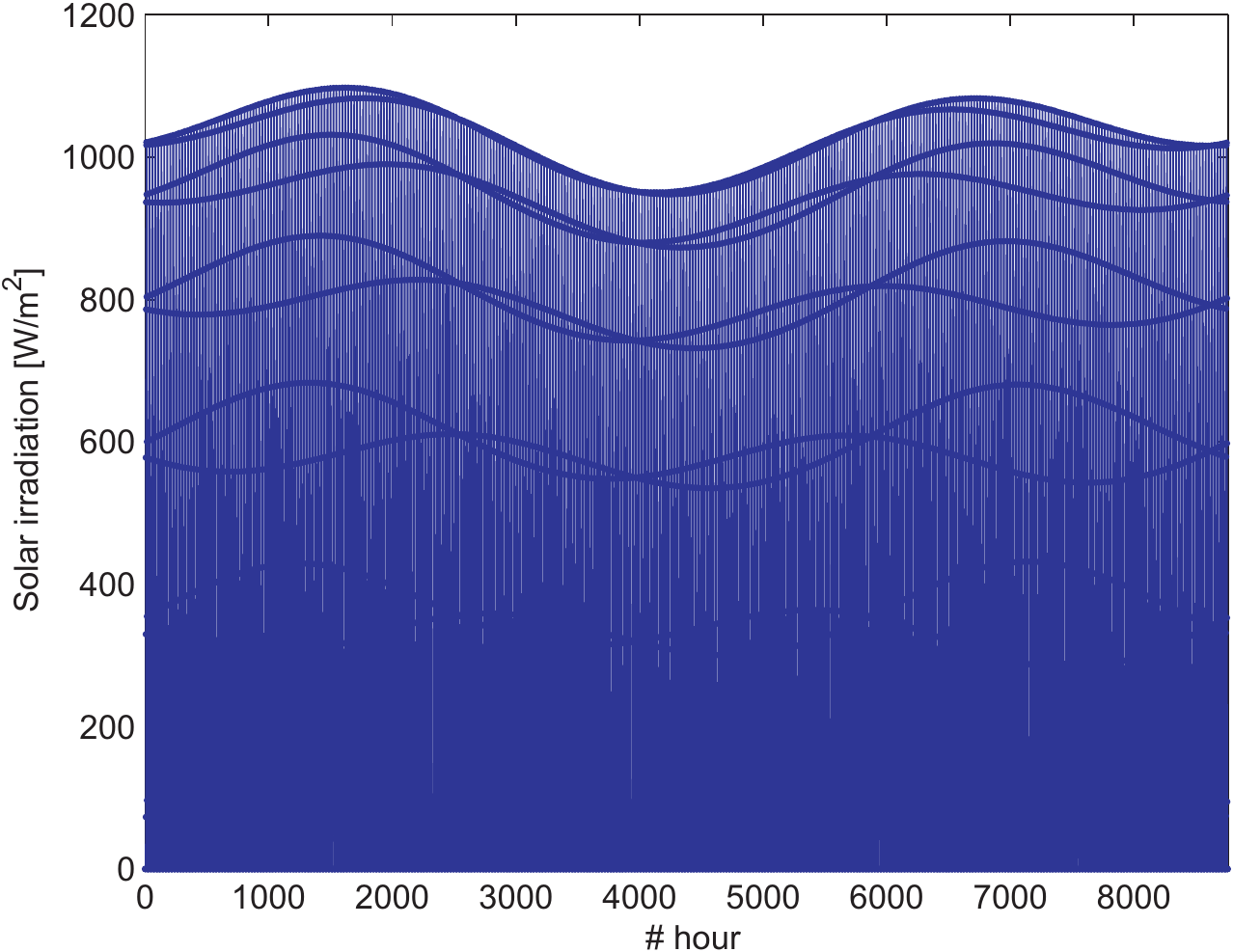}}\\
\subfigure[~Toledo (winter)]{\includegraphics[draft=false, angle=0,width=6.2cm]{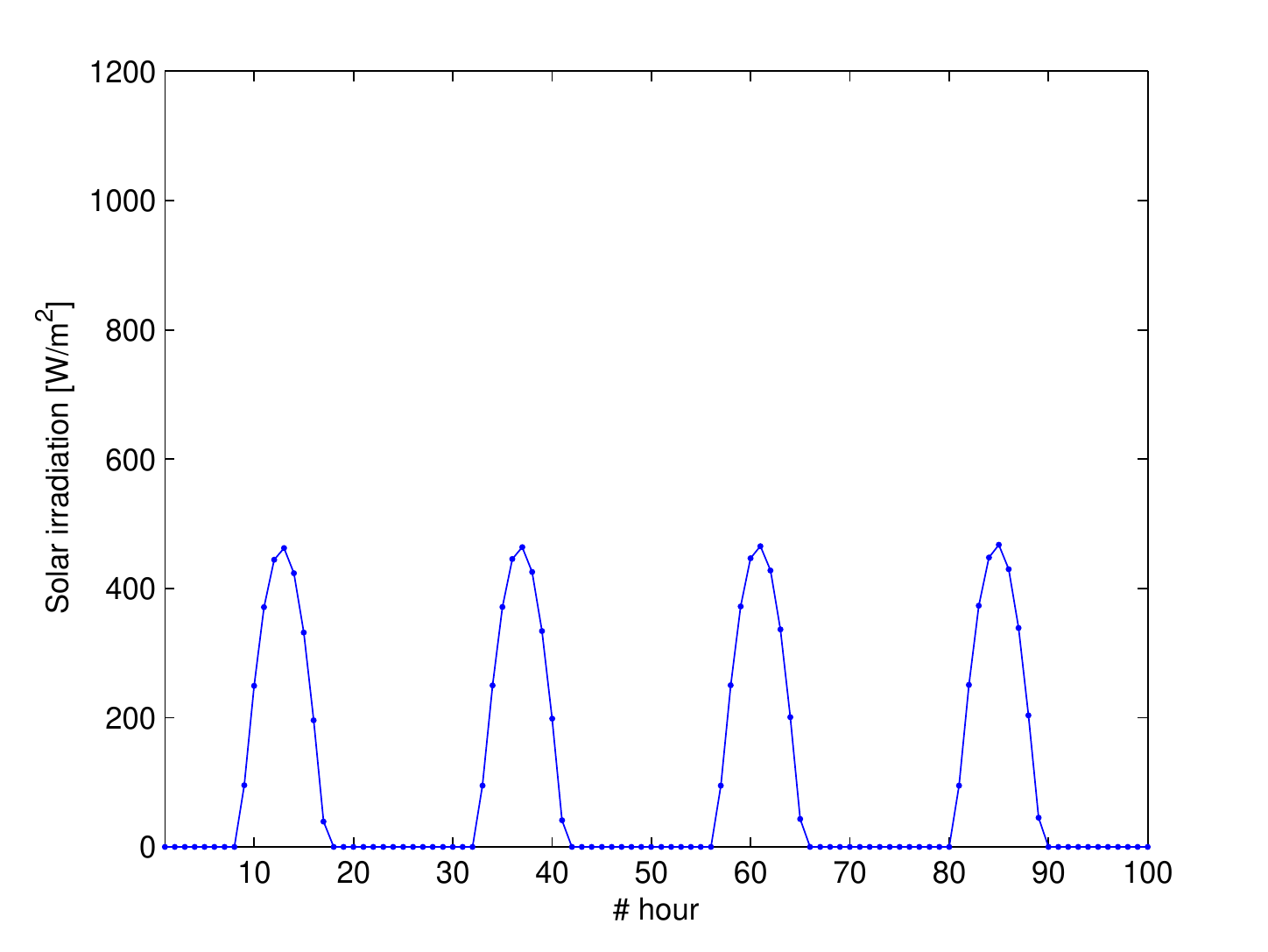}}
\subfigure[~Toledo (summer)]{\includegraphics[draft=false, angle=0,width=6.2cm]{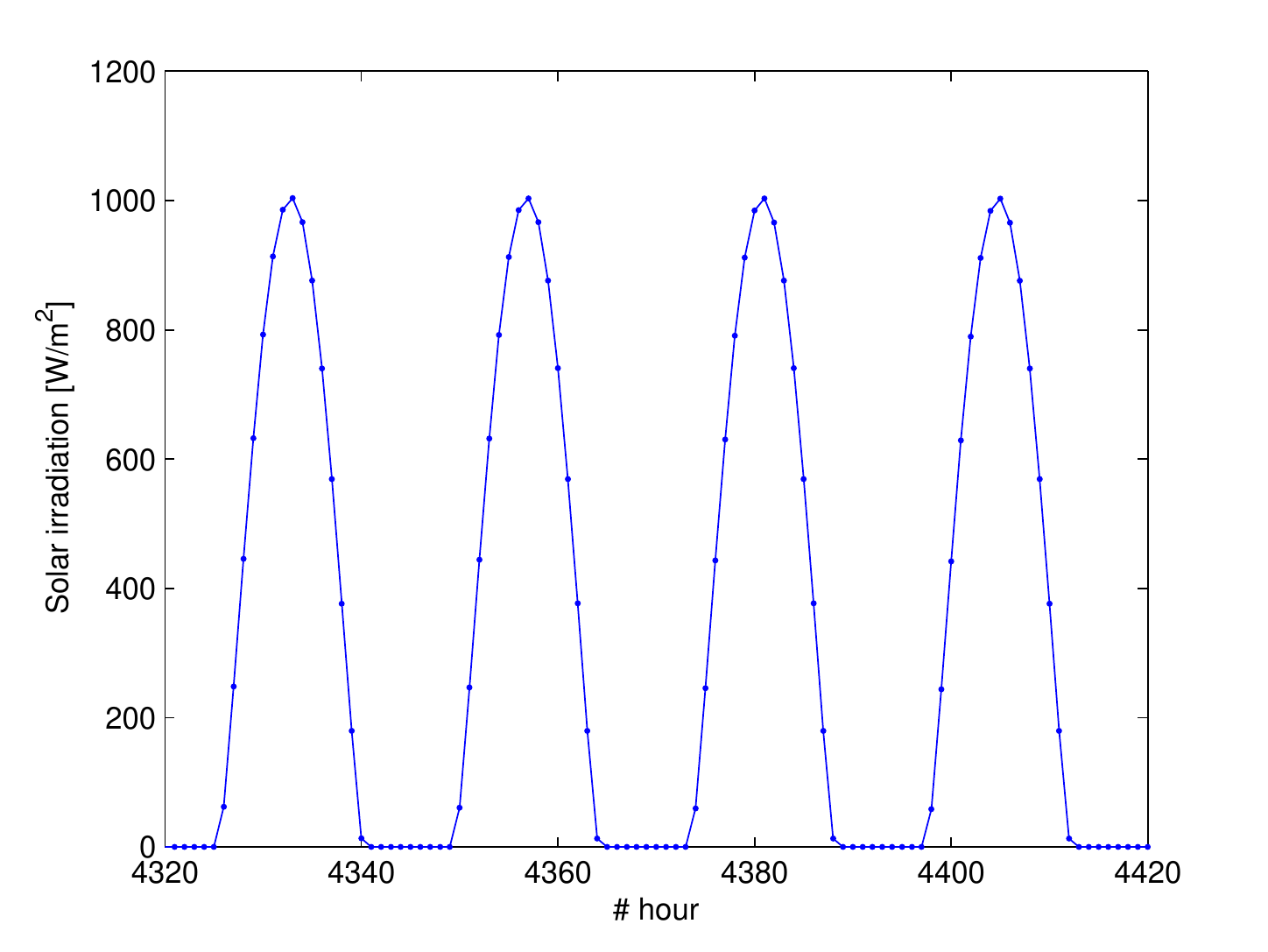}}\\
\end{center}
\caption{\label{CS_cities} Solar irradiation by the CS model (hourly basis), from January to December, for different sites in the Earth and zoom of the times series at Toledo (Spain); (a) Reykjavik; (b) Toledo (Spain); (c) Sydney; (d) Gulf of Guinea; (e) Zoom at Toledo (winter); (f) Zoom at Toledo (summer).}
\end{figure}

The basis to estimate solar radiation at any given location is to apply the classical astronomical equations \cite{Iqbal12}, also known as {\em Clear Sky} (CS) model \cite{Bird81}. CS model is based on the fact that, depending on the actual point of the Earth under study, the amount of solar radiation received (solar irradiation) without considering atmosphere, just depends on astronomical parameters. Of course, the effect of the atmospheric processes associated with solar irradiation at a given point is extremely important, and there are well-known different processes that modifies the solar irradiation, such as: molecular (or Rayleigh) scattering by the permanent gases, aerosol (or Mie) scattering due to particles and the abortion by different gases and finally the clouds presence \cite{Zenkai2008}. However, note that the CS can be seen as a kind of special {\em persistence} (astronomical) of solar energy, in such a way that we count on a maximum possible solar irradiation at a given point of the Earth, at any time of the year. Astronomical persistence given by the CS model takes into account the location of the point under study and the moment of the prediction (when in the year). Note that in this case of solar prediction, the classical persistence would take into account of atmospheric processes which affect the total radiation reaching to the ground. As an example of astronomical persistence of solar energy, Figure \ref{CS_cities} shows the CS (hourly) at four different study sites in the Earth (Reykjavik (Iceland), Toledo (Spain), Sidney (Australia) and (0Lat,0Long) point at Guinea Gulf. A detailed (zoomed) picture of the CS at Toledo (Spain), remaking the differences in solar irradiation in winter and summer, is given in Figures \ref{CS_cities} (e) and (f). As can be seen, the CS solar irradiance amount and its variability over the year depends on the location of the site under study on the Earth.

The use of a CS model is considered by the majority of prediction algorithms when dealing with solar irradiation prediction problems. This way, the persistence of the solar irradiation due to astronomical parameters is taken into account, and the algorithms can be fully focused on modeling the atmospheric effects to improve the prediction. There are different approaches to do this, such as the application of NWP models. These systems are able to model the dynamics of the atmosphere, as well as the physical processes involved, by employing a set of equations based on physical laws of motion and thermodynamics. For example, in \cite{Perez2013} the performance of different NWP models to forecast solar irradiance in the US, Canada and Europe is presented. More recently, \cite{Huang2017} carried out a study about the effects of the simulated regional weather variability at a relatively fine spatial resolution on the forecasting accuracy of solar irradiance. This was done using two NWP models: the Conformal Cubic Atmospheric Model (CCAM) and the Global Forecasting System (GFS). Additionally, NWP models have demonstrated their ability to provide useful solar radiation climate data sets, like the work of \cite{Perdigao2017} where a 60 years (1950-2010) climatology of incident shortwave downward solar radiation at the surface over the Iberian Peninsula is obtained from simulations performed using the WRF mesoscale model. ML approaches have become a reliable alternative/complement to NWP in solar energy prediction problems. ML approaches such as neural networks \cite{Mellit08,Amrouche2014}, Support Vector regression approaches \cite{Belaid16,Cornejo19}, Random Forests \cite{Ibrahim2017}, hybrid approaches \cite{Salcedo2018} or deep learning algorithms \cite{Ghimirec2019,Liu2019}, among others, have been successfully applied to model the atmospheric part (variables) related to solar irradiation prediction in different specific applications. A state-of-the-art review of ML models widely used in solar energy prediction problems can be found in recent works \cite{Voyant2017,Mosavi2019}. Both NWP and ML approaches are able to obtain highly accurate predictions at any time-horizon prediction levels, however, parsimonious models based on persistence are also able to obtain excellent results own their own, with an extremely low computational complexity \cite{Pedro12,Antonanzas16,Voyant18,FLIESS2018519}.

In this case study we will show an example of this, in a solar irradiation prediction problem at Toledo (Spain), by only considering persistence-based models. Following the approaches proposed in \cite{Voyant18}, we will tackle a problem of daily average solar radiation prediction with models, exclusively based on CS and modifications of na\"ive persistence. We will show how these models are able to obtain very reasonable prediction errors with extremely simple and low-complexity models.

\subsubsection{Solar irradiation data}

\begin{figure}[!ht]
\begin{center}
\includegraphics[draft=false, angle=0,width=9cm]{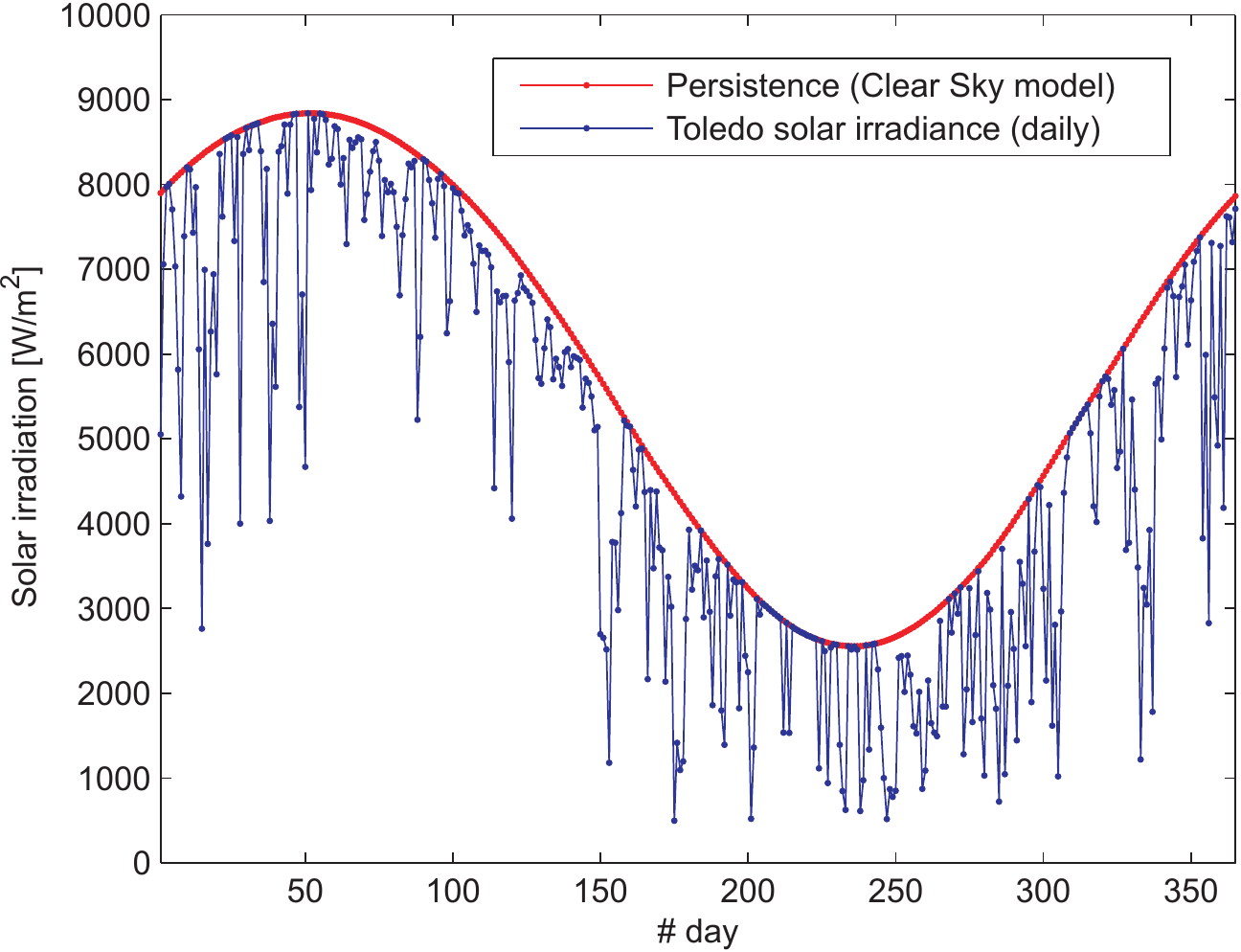}
\end{center}
\caption{\label{Medidas_y_CS_Toledo} Clear Sky model (daily aggregated) and measurements of solar irradiance at Toledo, Spain.}
\end{figure}

We consider solar irradiation at the radiometric station of Toledo, Spain ($39^\circ$ 53'N, $4^\circ$ 02'W). This meteorological station is part of the Spanish radiometric observing network, managed by the Meteorological State Agency of Spain (AEMET) and located in the South Plateau of the Iberian Peninsula, around 75 km south of Madrid (the capital city of Spain) at an altitude of 515m. We consider one year of daily average global solar irradiation data (from the 1st of May 2013 to the 30th of April 2014) measured at Toledo's station with a Kipp \& Zonen CM-11 pyranometer. All radiation measurements gathered by the AEMET are managed under a quality management system certified to ISO 9001:2008, which guarantees their accuracy and their compliance with the World Meteorological Organization (WMO) standards. Figure \ref{Medidas_y_CS_Toledo} shows the global solar irradiation measurements and the CS model for the period under study.

\subsubsection{Parsimonious persistence-based models for solar irradiation prediction}

The persistence-based models considered in this case are based on the application of the CS model $CS[n]$, to model astronomical persistence of solar irradiation, and na\"ive persistence to model the atmospheric processes, such as clouds, dust, particles, etc. We consider different time-steps $k$ to model this persistence $NV_k[n]$ as follows:

\begin{equation}
    NV_k[n]=\frac{1}{k}\sum_{j=1}^k I[n-j].
\end{equation}
where $I[n-j]$ stands for values of the solar irradiation time series, at times $n-j$. Note that with $CS[n]$ and $NV_k[n]$ terms, and following \cite{Voyant18}, we are able to obtain a reasonable prediction of solar irradiation ($\hat I$) at time $n$, based on $CS[n]$ and previous average values of irradiation given by $NV_k[n]$, depending on the number of time steps $k$ considered. We have used four different prediction expressions, to combine both types of persistence (CS-NV1 (Equation (\ref{Naive1})), CS-NV2 (Equation (\ref{Naive2})), CS-NV3 (Equation (\ref{Naive3})), CS-NV4 (Equation (\ref{Naive4})):
\begin{equation}\label{Naive1}
\hat I[n]= \frac{1}{k}\sum_{j=1}^k\frac{CS[n] + NV_k[n]}{2}
\end{equation}
\begin{equation}\label{Naive2}
\hat I[n]=\frac{1}{k}\sum_{j=1}^k\alpha CS[n] + \beta NV_k[n]
\end{equation}
\begin{equation}\label{Naive3}
\hat I[n]= \frac{1}{k}\sum_{j=1}^k \sqrt{CS[n] \cdot NV_k[n]}
\end{equation}
\begin{equation}\label{Naive4}
\hat I[n]=\frac{1}{k}\sum_{j=1}^k\sqrt{\alpha CS[n] \cdot \beta NV_k[n]}
\end{equation}

\subsubsection{Results}

\begin{figure}[!ht]
\begin{center}
\subfigure[~$k=1$]{\includegraphics[draft=false, angle=0,width=7cm]{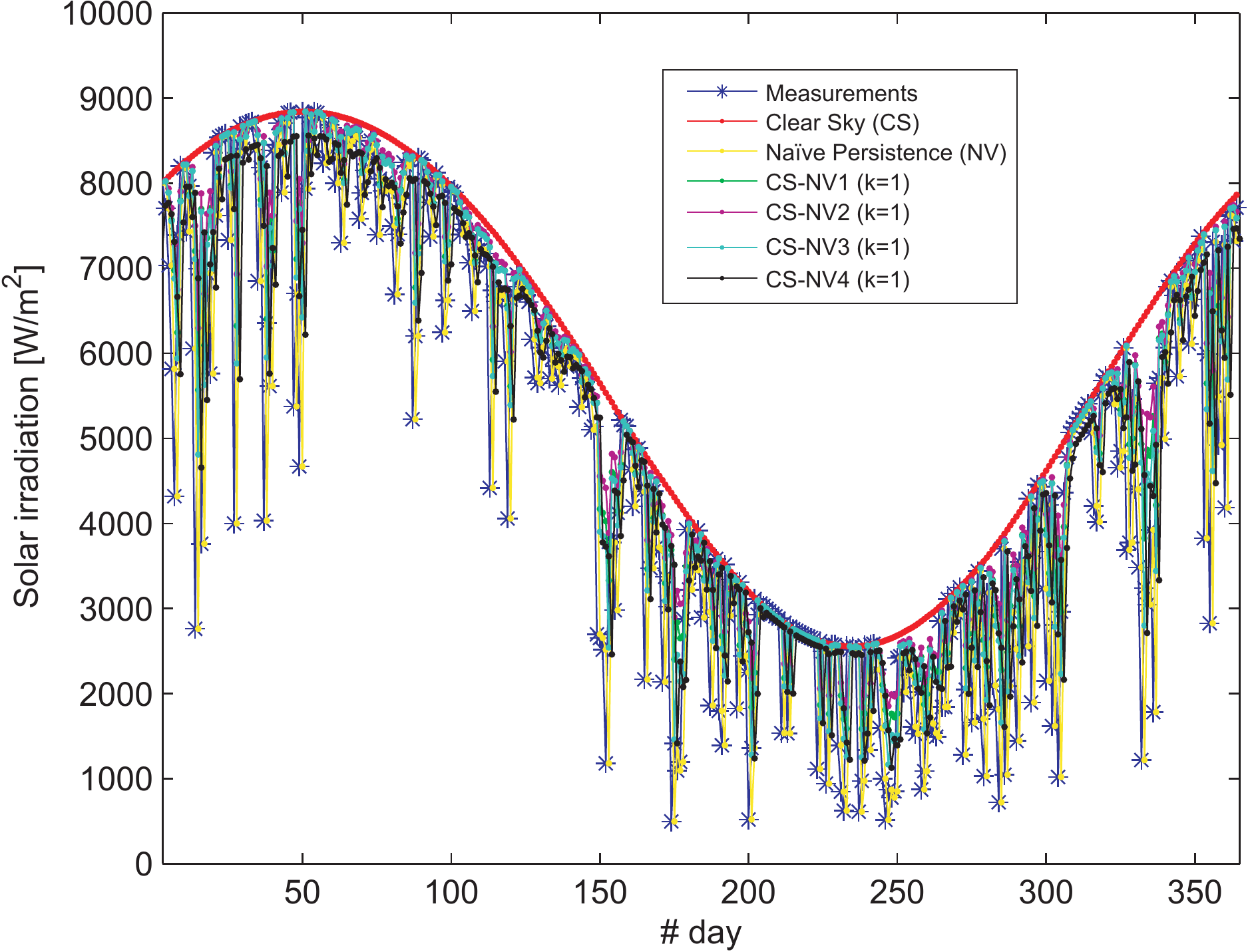}}
\subfigure[~$k=1$ zoom]{\includegraphics[draft=false, angle=0,width=7cm]{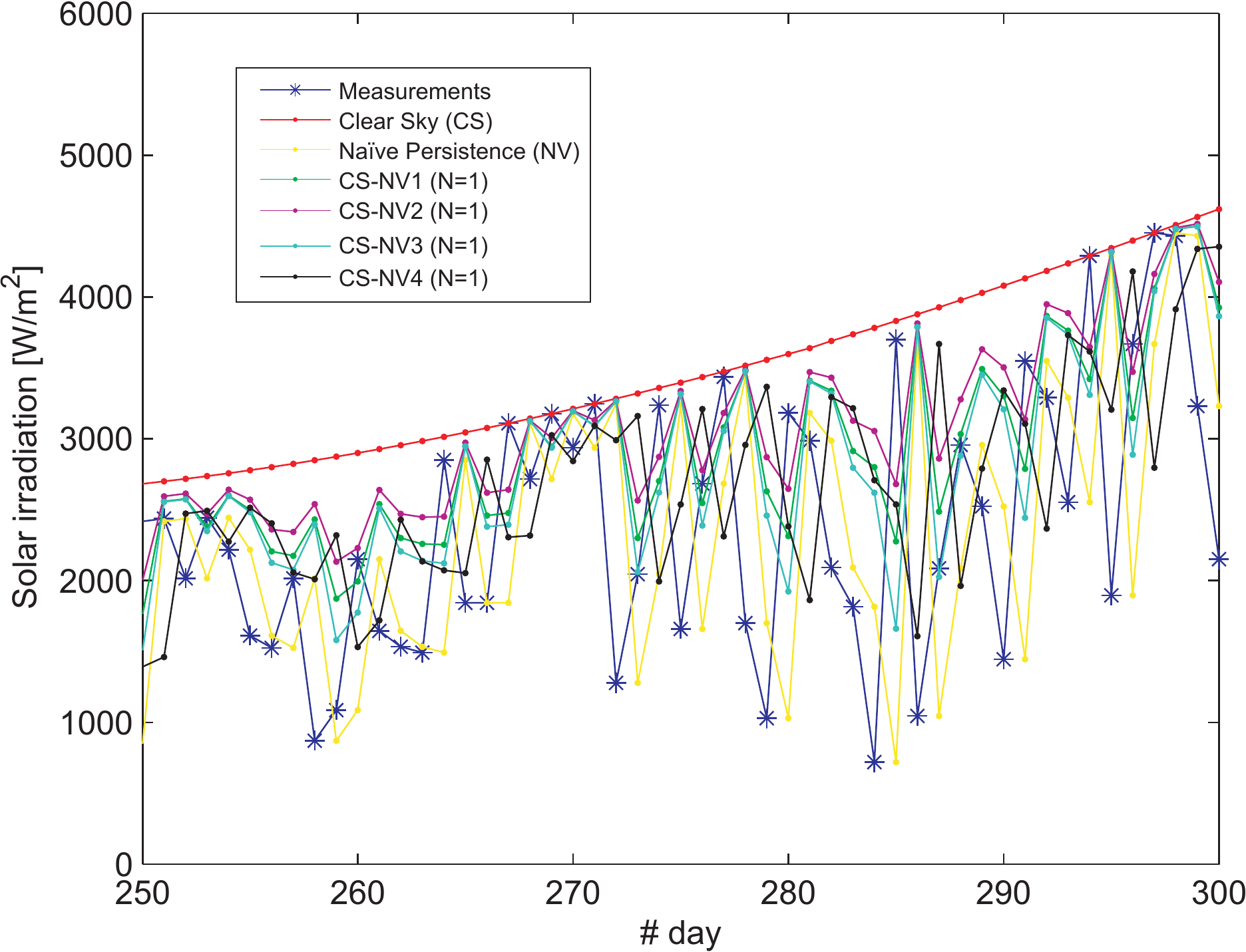}}\\
\subfigure[~$k=2$]{\includegraphics[draft=false, angle=0,width=7cm]{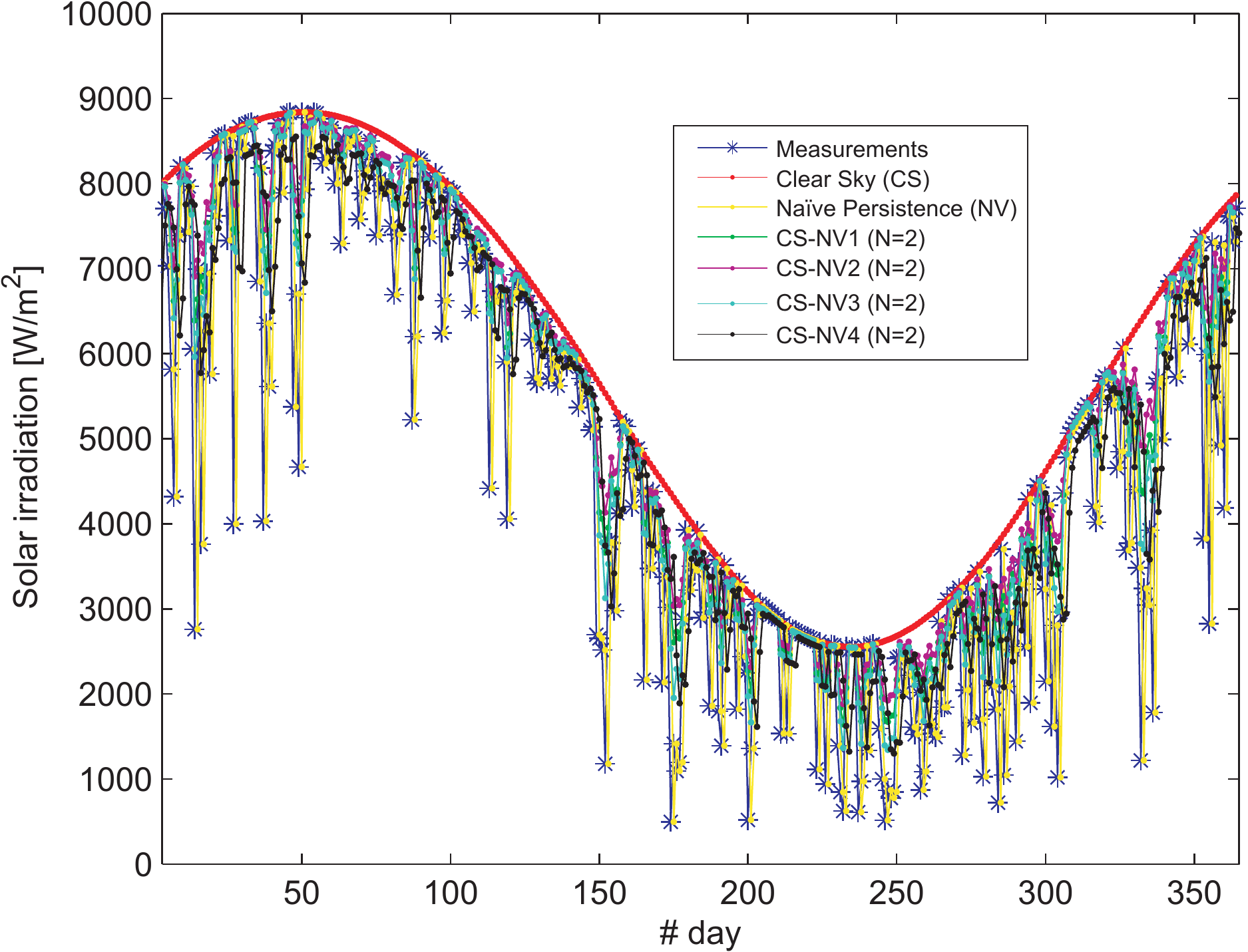}}
\subfigure[~$k=2$ zoom]{\includegraphics[draft=false, angle=0,width=7cm]{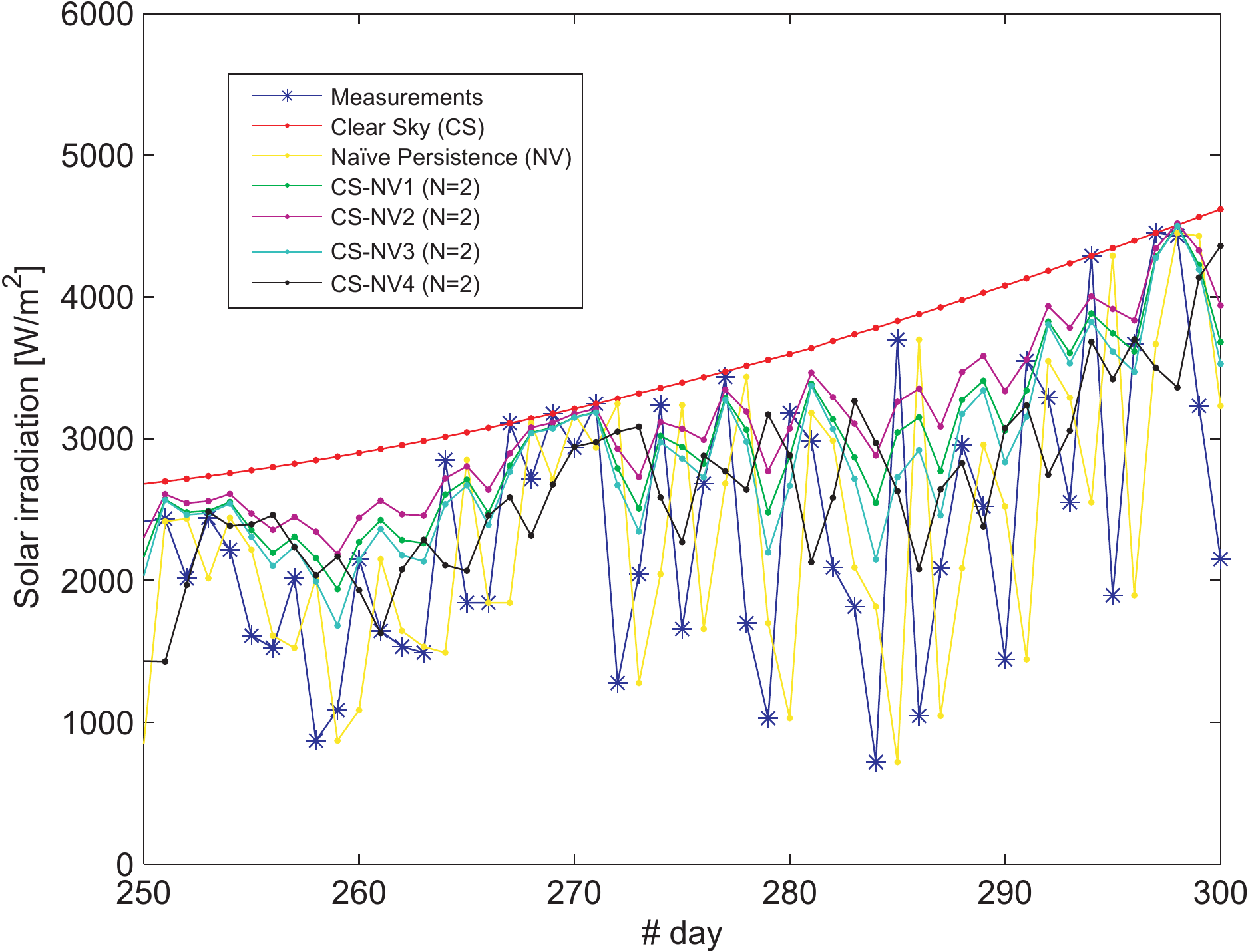}}\\
\subfigure[~$k=3$]{\includegraphics[draft=false, angle=0,width=7cm]{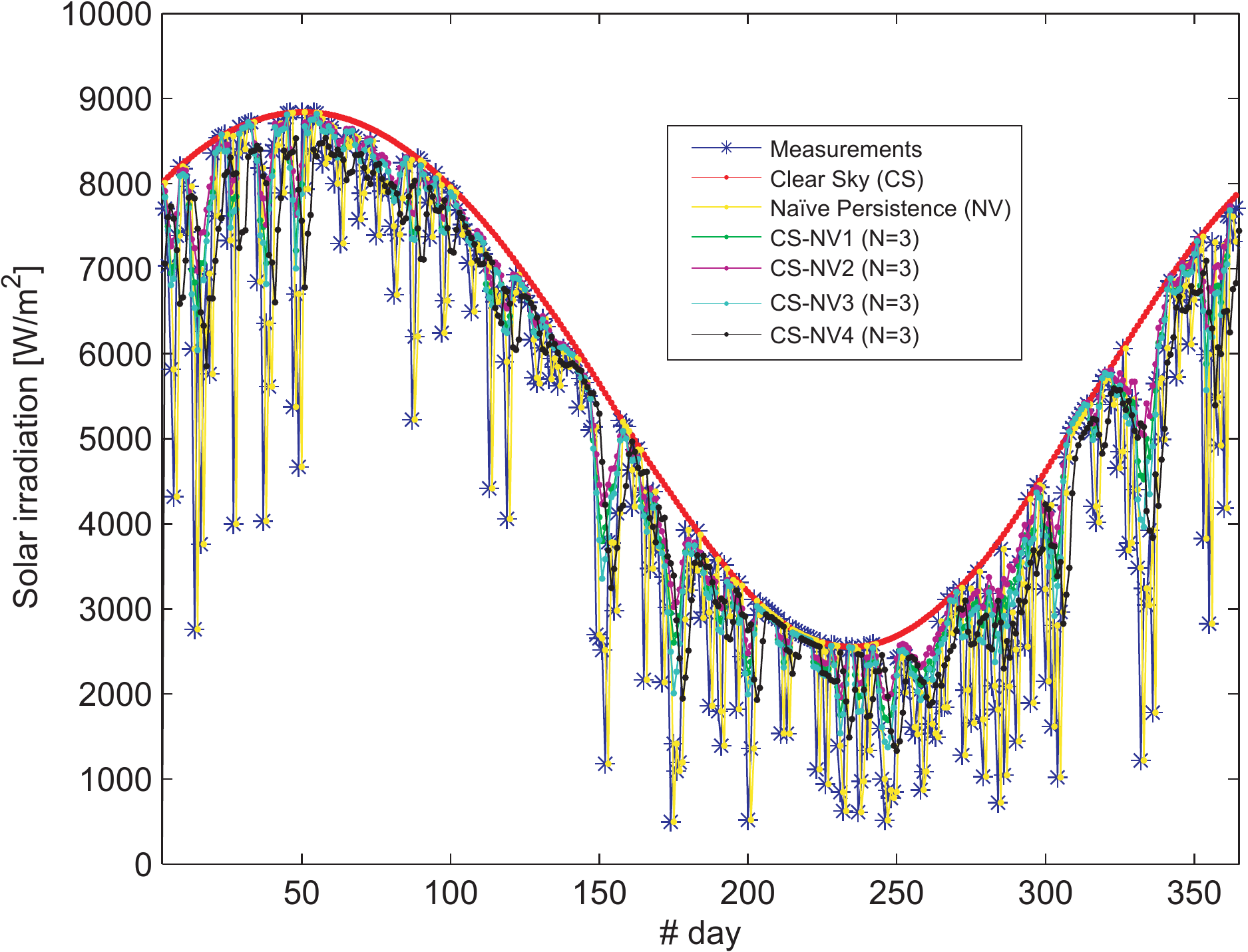}}
\subfigure[~$k=3$ zoom]{\includegraphics[draft=false, angle=0,width=7cm]{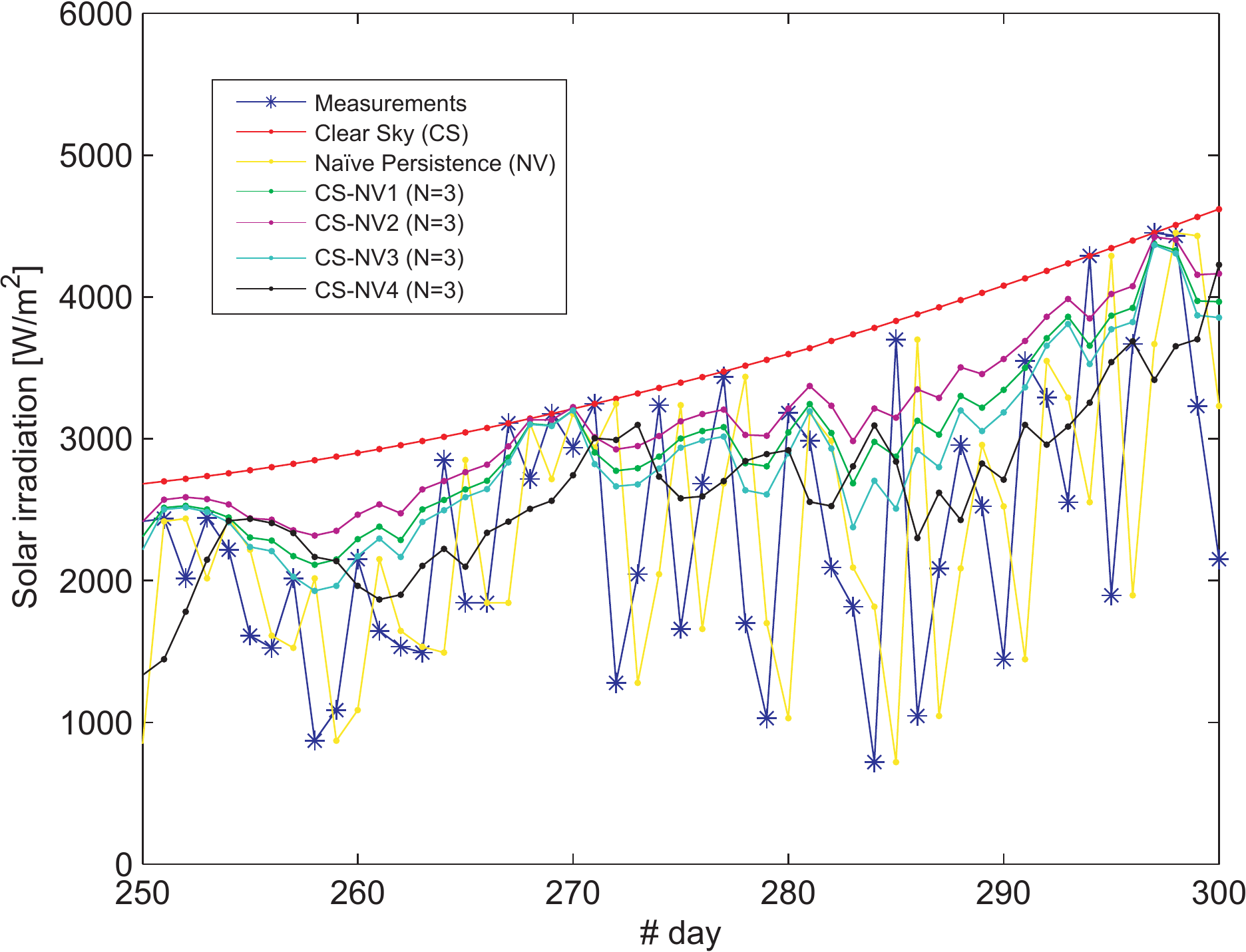}}\\
\end{center}
\caption{\label{Prediccion_6M} Prediction performance of the different persistence-based models in solar irradiation data at Toledo radiometric station.}
\end{figure}

Table \ref{solar_MAE} shows the results obtained with the persistence-based models CS-NV in the problem of daily solar irradiation prediction at Toledo (Spain), for different values of previous time-steps considered ($k$). As can be seen, the merging of na\"ive persistence and CS models to obtain persistence-based solar irradiation prediction models outperforms the na\"ive persistence and CS models on their own, in all cases tested. The differences are very significant, with a best MAE of 733 $W/m^2$ as best prediction result in the analyzed data, given by the CS-NV4 model with $k=2$. The detailed performance of these parsimonious models can be better analyzed in Figure \ref{Prediccion_6M}, where the performance of each model and a detailed zoom are shown. This figure is very illustrative of how each model works. In general all persistence-based models over-estimate the amount of global solar irradiation received, as expected, since they do not incorporate any atmospheric physics variable. Of course, numerical models and ML method with atmospheric variables are able to improve these predictions, but note that persistence-based methods are computationally much less costly. In addition, the results of these models can be incorporated as input variables in ML techniques, including valuable information which can be exploited by ML techniques to elaborate an accurate solar irradiation prediction.

\begin{table}[!h]
\small{}
\begin{center}
\caption{\label{solar_MAE} Results obtained with persistence-based models in the problem of daily solar irradiation prediction at Toledo, Spain. The best values for $\alpha$ and $\beta$ in CS-NV2 and CS-NV4 are 0.75 and 1.25, respectively.} \vspace{0.3cm}
\begin{tabular}{ccccc}
\hline
\vspace{0.1cm}
Method             & MAE ($W/m^2$) \\

\hline
\hline

Clear Sky model (CS)  \cite{Bird81}     &  950.86\\
Na\"ive Persistence (NV)                  &  837.39\\
\hline
$k=1$&\\
\hline
CS-NV1 & 749.40\\
CS-NV2&770.04\\
CS-NV3 & 757.30\\
CS-NV4&740.71\\
\hline
$k=2$&\\
\hline
CS-NV1 & 747.01\\
CS-NV2&769.79\\
CS-NV3& 743.87\\
CS-NV4&{\bf 733.61}\\
\hline
$k=3$&\\
\hline
CS-NV1& 752.89\\
CS-NV2&779.39\\
CS-NV3 & 747.34\\
CS-NV4&742.01\\
\hline
\end{tabular}
\end{center}
\end{table}

\subsection{Estimation of soil moisture persistence}

Soil Moisture (SM) memory or persistence can be defined as a measure of the time that a moisture anomaly is detectable, and during which it can influence the atmosphere. Characterization of this persistence has important implications in ecology \cite{Manzoni2012}, water management \cite{Rosenzweig2002}, and climate modeling \cite{Delworth89,Koster2001,Seneviratne2010}. The estimation of soil moisture persistence from Earth observation-based products is a promising area of intense current research (e.g. \cite{McColl2017,Ghannam2016,Shellito2018,Piles2021}). The motivation behind these works is that learning soil moisture persistence parameters from observations may reduce the sources of uncertainties in Earth system and climate models, where complex dynamics and memory effects are currently poorly represented, or even not considered at all. For instance, the impact of deeper soil layers on the upper layers in times of water stress \cite{Entin2000}, or the sign and strength of soil moisture precipitation feed-backs \cite{Koster2003} are relevant processes that could be better resolved with a proper soil moisture memory characterization. Further, due to its memory characteristics, soil moisture is an important contributor to climate persistence on land, impacting the development and duration of natural hazards such as droughts, floods, and heatwaves \cite{nicolai2016long}. As such, characterizing land-atmosphere interactions with observations is a key priority for improving the predictability of global water and climate processes \cite{Santanello2018}.

As introduced in Section \ref{ssec:efolding}, the autocorrelation function has been widely adopted to estimate short-term persistence in time series of atmospheric, land, or oceanic Earth system variables, for diverse application settings. Under certain assumptions, soil water dynamics can be approximated as a first order Markov process and hence the e-folding time scale can be a quite efficient measure of soil moisture memory \cite{Delworth88,Delworth89}. In this case study we first introduce the physical models describing soil water dynamics and discuss under which assumptions soil moisture time series can be reasonably represented as a red noise process. Subsequently, we show an example of soil moisture persistence estimation using ground-based measurements from the REMEDHUS network (Salamanca, Spain), as well as Earth Observation data over Europe. Following the approaches proposed in \cite{Piles2021}, we will tackle the problem of robust autocorrelation estimation for non-uniform time series and will then provide spatial-temporal descriptions of soil moisture persistence using e-folding times. 

\subsubsection{Modeling surface soil moisture}

Let us start with the physics of soil water dynamics and storage, which can be described with the surface water balance equation:
\begin{equation}
    \Delta z \frac{d\theta}{dt} = P(t) - (ET(t,\theta)+D(t,\theta)+Q(t,\theta)) = P(t) - L(t,\theta) 
\end{equation}
where $\theta$ is volumetric soil moisture (-), $t$ is time ($d$), $\Delta z$ is the depth of the soil volume below the surface ($mm$), $P(t)$ is the precipitation rate ($mm{\cdot}d^{-1}$), and $L(t,\theta)$ are the total water losses from the volume due to evapotranspiration $ET(t,\theta)$, drainage $D(t,\theta)$, and runoff $Q(t,\theta)$ ($mm{\cdot}d^{-1}$). While precipitation is an exogenous forcing, often modelled stochastically, the fluxes combined into the loss function are quasi-deterministic of the land surface \cite{Delworth88}. 

\begin{figure}
\begin{center}
\includegraphics[width=14cm]{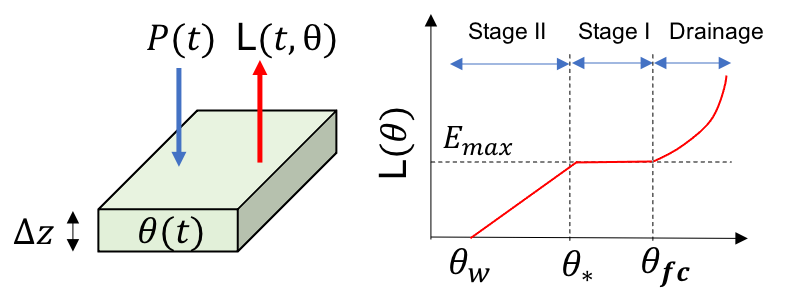}
\caption{Left: illustration of the surface water balance equation for a soil layer of depth $\Delta_z$, where the storage of water is represented by its volumetric soil water $\theta(t)$ and its dynamics are governed by precipitation inputs ($P(t)$) and drainage, runoff and evapotranspiration losses ($L(t,\theta)$). Right: Loss function under mean climate state, characterized by three regimes, 1) drainage or percolation dominant regime for soil wetness above the field-capacity $\theta_{fc}$, 2) Stage I evapotranspiration regime for intermediate soil wetness levels ($\theta_*\leq\theta\leq\theta_{fc}$), where there is not sensitivity to soil moisture, and 3) Stage II evapotranspitation or water-limited regime for soil moisture levels below the critical value $\theta_*$.}\label{fig:water_balance}
\end{center}
\end{figure}

The soil water loss function under mean climate state $L=L(\theta)$ is typically characterized by three regimes (Figure \ref{fig:water_balance}) that can be approximated by a three-stage piece-wise function \cite{Laio2001,Feng2017}. During and after precipitation, wet soils lose water rapidly due to drainage and runoff, following a power law until volumetric water content is lowered to the field capacity $\theta_{fc}$. Below this point, drainage and runoff become negligible and $L(\theta)$ is dominated by evapotranspiration, which is considered initially to be in ``stage I'' or energy-limited regime. During this period, $ET$ is invariant with respect to volumetric water content and is determined by the mean climatological potential evapotranspiration rate $E_{max}$ ($mm\cdot d^{-1}$), which depends on available energy and atmospheric evaporative demand \cite{Rodriguez-Iturbe2005}. The soil may dry further until it reaches the so-called critical soil moisture $\theta_*$, which is the transition point between water- and energy-limited hydrologic regimes. Below $\theta_*$,  evapotranspiration is limited by soil water availability and $L(\theta)$ decreases monotonically to zero as soil moisture decreases to the wilting point $\theta_w$. This phase is termed ``stage II'' or water-limited regime. The relative dominance of the water loss fluxes depends on how much time has passed since rainfall \cite{Laio2001}. Runoff, when present, ceases only minutes after rainfall, and drainage occurs on a timescale of hours. Assuming that 1) $P(t)$ can be represented as a white noise process; 2) at time scales on the order of days, water losses can then be assumed to be dominated by stage II or water limited conditions (see Figure \ref{fig:water_balance}), and the land water anomaly budget can be rewritten as:
\begin{equation}
    \frac{d\theta}{dt} = -\frac{\theta(t)-\theta_w}{\tau} + \epsilon(t),
\end{equation}
where $\epsilon$ is an independent and identically distributed random variable with a mean of zero and parameter $\tau=\frac{\Delta z(\theta_*-\theta_w)}{E_{max}}$. Noting that this expression is the definition of a red noise process, it follows that its e-folding time scale is $\tau$ and that it can reasonably represent the persistence of the time series \cite{Delworth88,Delworth89}. This approximation has been widely adopted to approximate SM dynamics in a variety of applications (e.g. \cite{Entin2000,Ghannam2016,Koster2001}). 

\subsubsection{Data and e-folding time estimation}

We considered hourly SM measurements acquired at 22 stations from the REMEDHUS SM measurement network \cite{Sanchez2012}, located in the central part of the Duero basin (41.1-41.5$^\circ$ N, 5.1-5.7$^\circ$ W), Spain. It covers a semi-arid continental-Mediterranean agricultural region of $35 \times 35$ km, and has extensively been used for calibration and validation of SM products (e.g. \cite{Sanchez2012,Polcher2016,Angel2019}). Each station within the network is equipped with capacitance probes (Hydra Probes of Stevens Water Monitoring System, Inc.) installed horizontally at a depth of 5 cm, with a reported accuracy of 0.003 $m^3m^{-3}$. We also considered satellite SM estimates over Europe (34-72$^\circ$N, 10$^\circ$W-36$^\circ$E) from the Eumetsat's C-band Advanced Scatterometer (ASCAT). The ASCAT product from ESA CCI Active v.3.2 product was chosen \cite{Dorigo2017}, with has a daily temporal resolution and a spatial resolution of 0.25$^\circ$. The study period covers six years, starting in June 2010. 

Soil texture from the Harmonized World Soil Database \cite{FAO2009} over Europe as well as ground-based measurements of soil textures at the 22 REMEDHUS stations have been used to investigate possible relations between the remotely-sensed soil moisture memory estimates (expressed with e-folding times) and soil sand and clay fractions. 

To cope with the non-uniform sampling of the observations, we estimated the autocorrelation function of each SM time series with the Lomb-Scargle periodogram and then applied an spatial-temporal denoising of the autocorrelation tensor before obtaining the e-folding time scale, following \cite{Piles2021}. Prior to analysis, we detrended the time series to avoid examining artificial memory associated with the climatological seasonal cycle of soil moisture. 

\begin{figure}[t!]
 \centering
     \includegraphics[trim=1cm 0.5cm 0.7cm 0.5cm,width=0.32\textwidth]{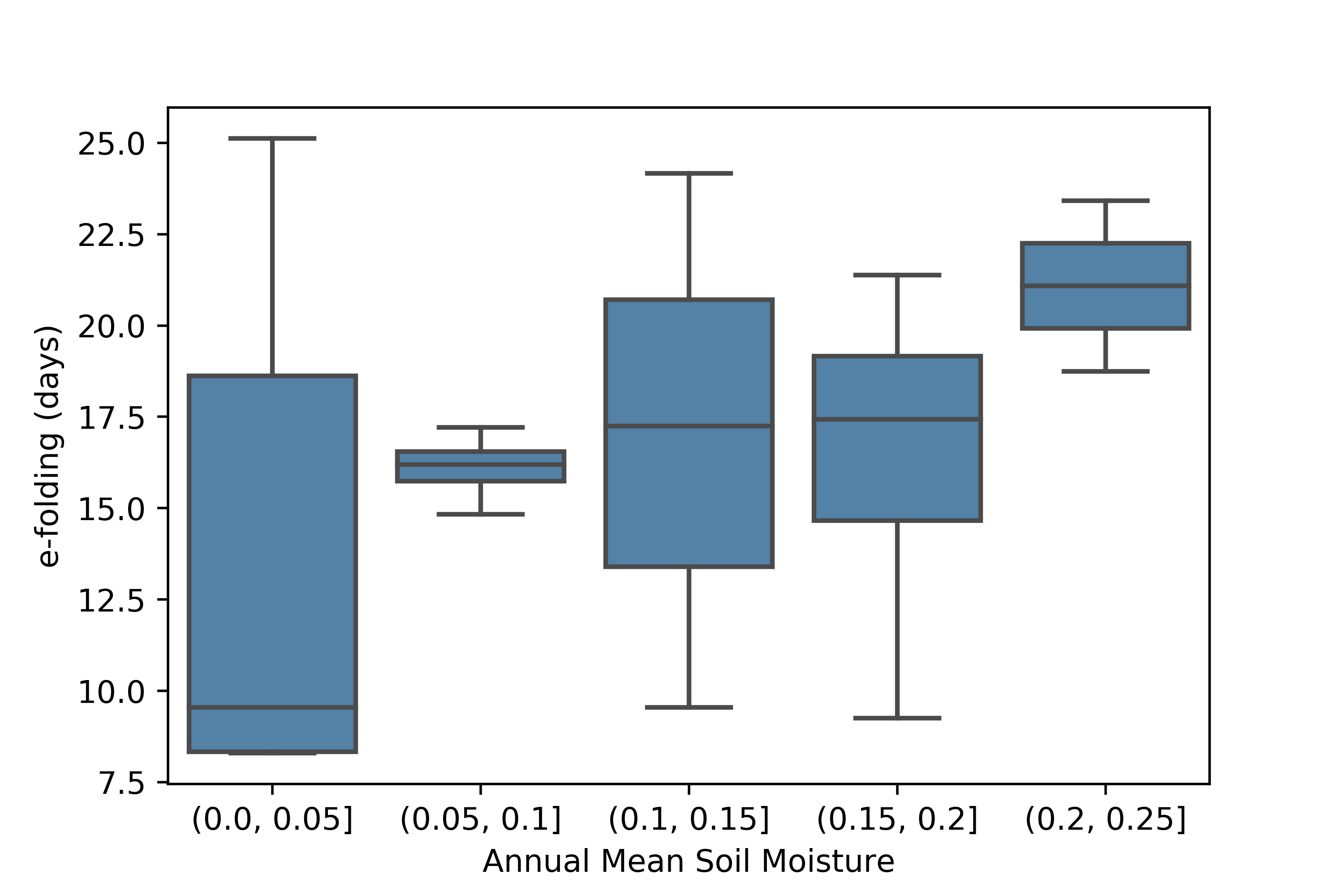}
    \includegraphics[trim=1cm 0.5cm 0.7cm 0.5cm, width=0.32\textwidth]{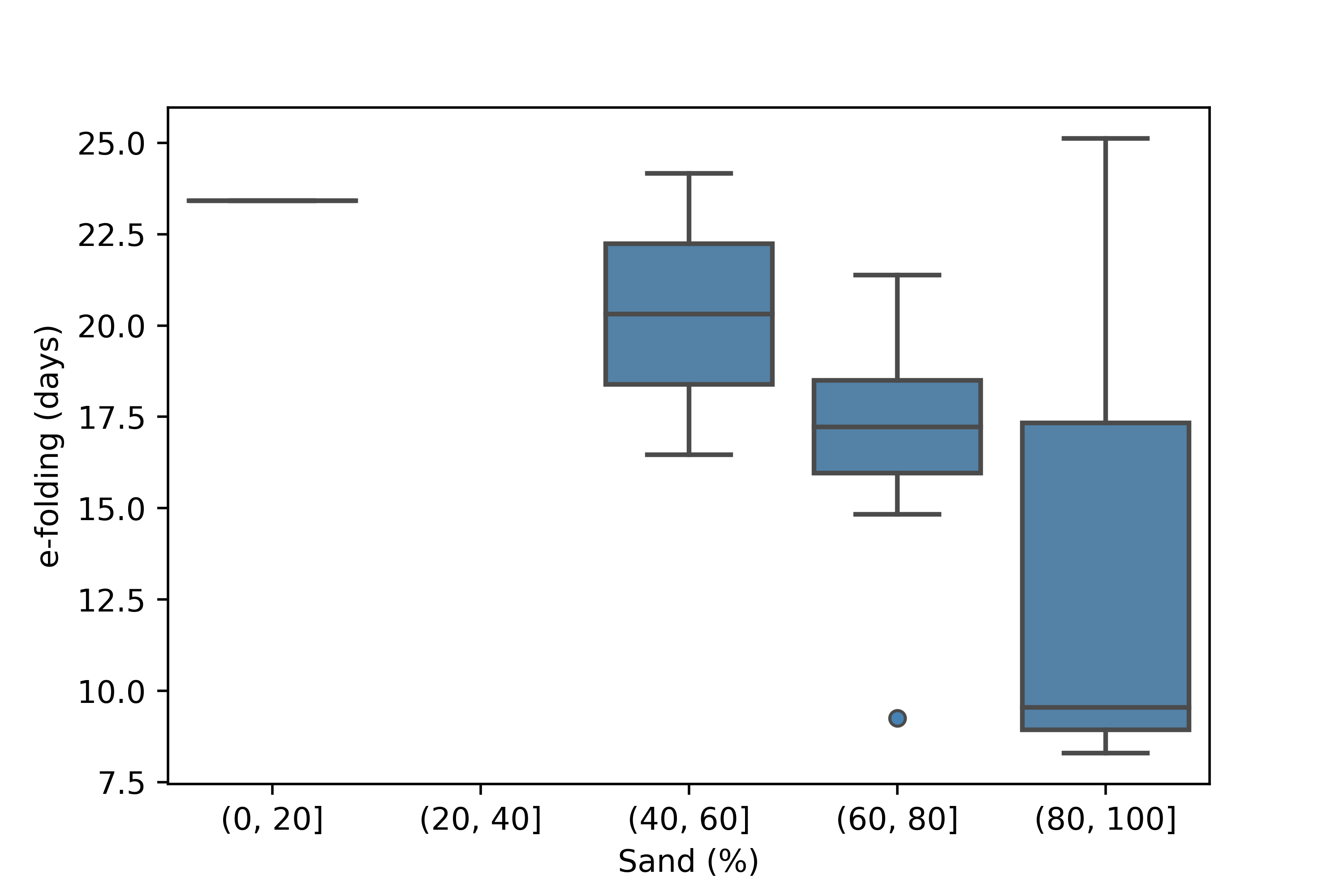}
    \includegraphics[trim=0.7cm 0.5cm 1cm 0.5cm, width=0.32\textwidth]{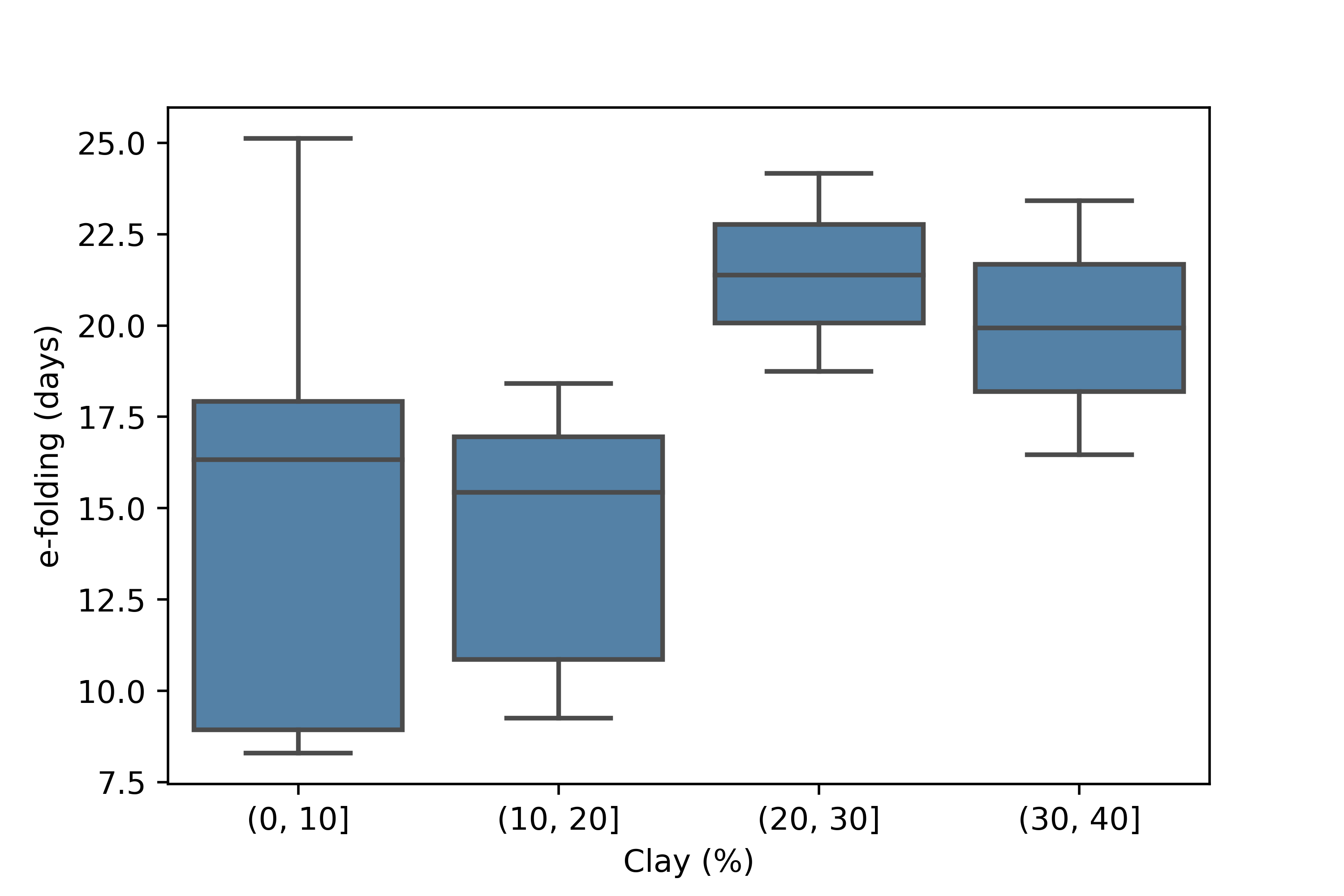}
 \caption{Boxplots of e-folding time (days) as a function of annual mean soil moisture (m$^3\cdot$m$^{-3}$) (left), the percentage of soil sand fraction (middle) and of soil clay fraction (right) for the REMEDHUS soil moisture measuring stations.}
 \label{fig:insitu_efolding}
\end{figure}

\subsubsection{Experimental results}

The SM e-folding times obtained at the 22 REMEDHUS stations are shown in Figure \ref{fig:insitu_efolding} as a function of the SM regime, characterized by the average SM content for the study period, and as a function of sand and clay fractions. Our results seem to indicate a proportional relationship of SM persistence with average SM content: the wetter the soil average condition, the higher the e-folding parameter, i.e. the longer it takes for the soil to dry out by a factor of $1/e$. Also, soil moisture memory at REMEDHUS stations appears to be highest for soils with low sand content and decreasing with sand content. This is in agreement with the general knowledge that coarse textured soils (sands and loamy sands) have faster percolation - due to large pores with limited ability to retain water. Also, a significant change is observed when the percentage of clay exceeds 20\%, so that the higher the percentage of clay, the greater the e-folding. 

\begin{figure}[H]
 \centering
    \includegraphics[width=0.6\textwidth]{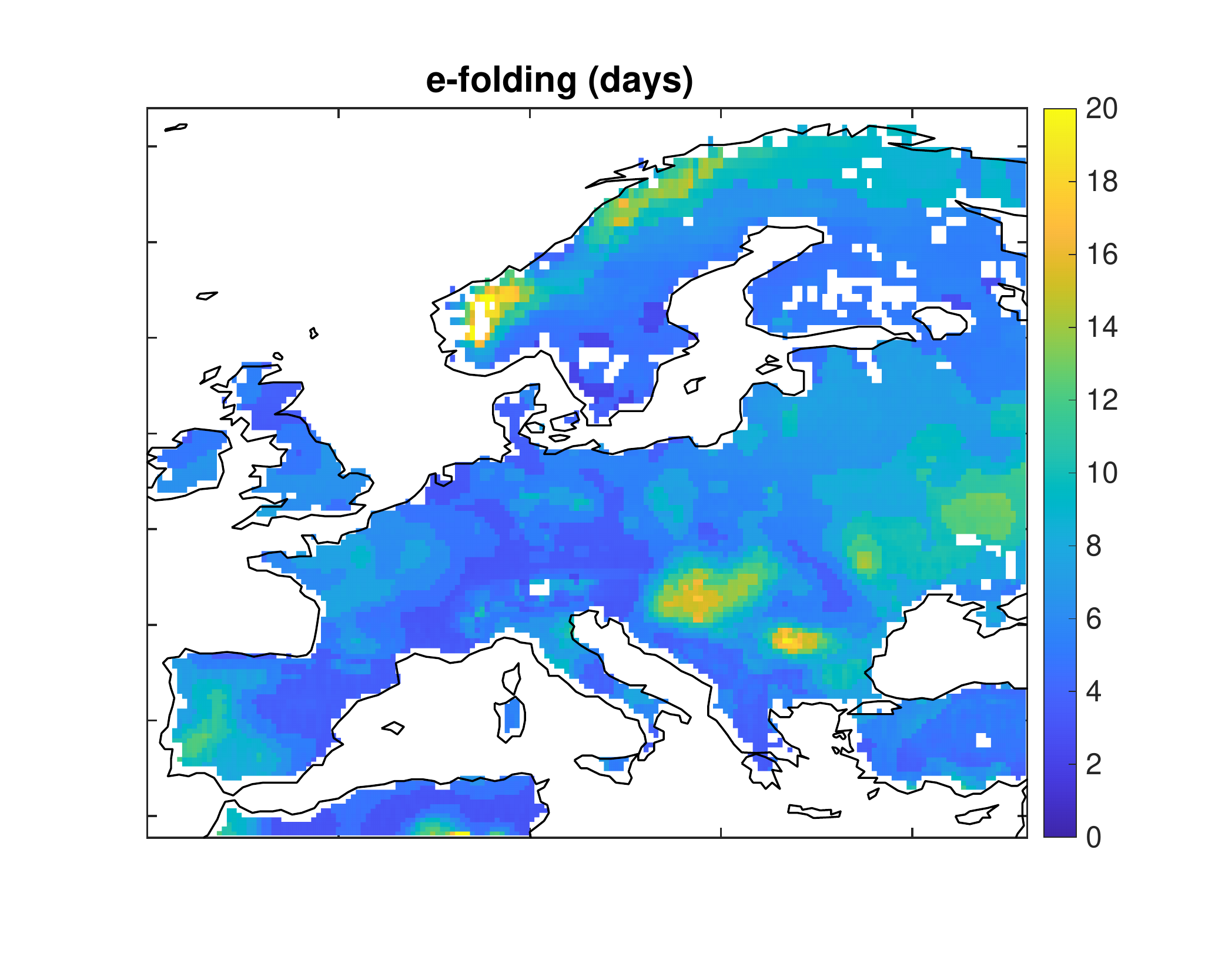}
    \vspace{-1cm}
 \caption{Map of e-folding time (days) as a measure of soil moisture persistence based on autocorrelation of ASCAT soil moisture time series for the period 2010-2016, following the approach in \cite{Piles2021}.}
 \label{fig:efolding_map}
\end{figure}

\begin{figure}[t!]
\begin{center}
\includegraphics[width=0.32\textwidth]{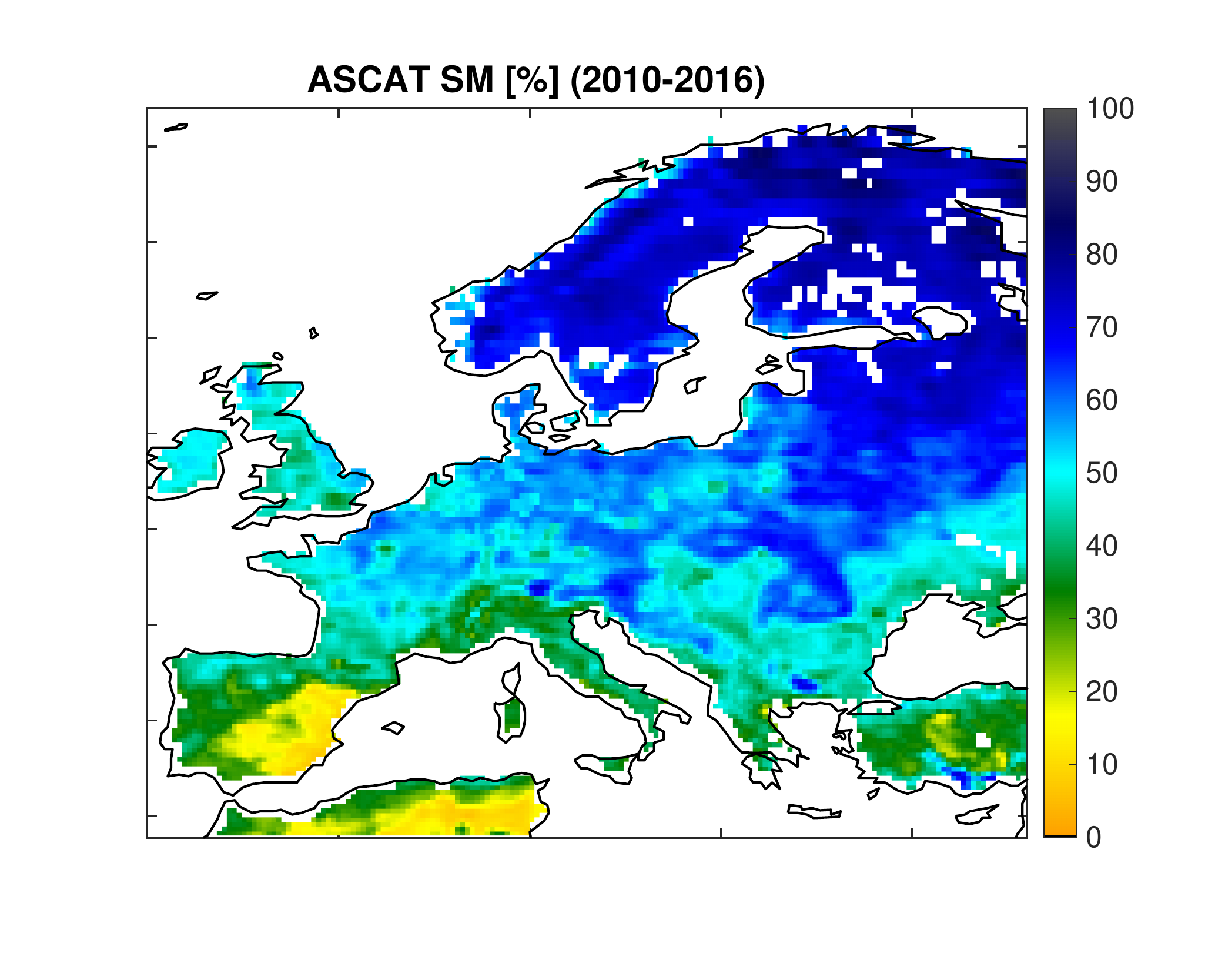}
\includegraphics[width=0.32\textwidth]{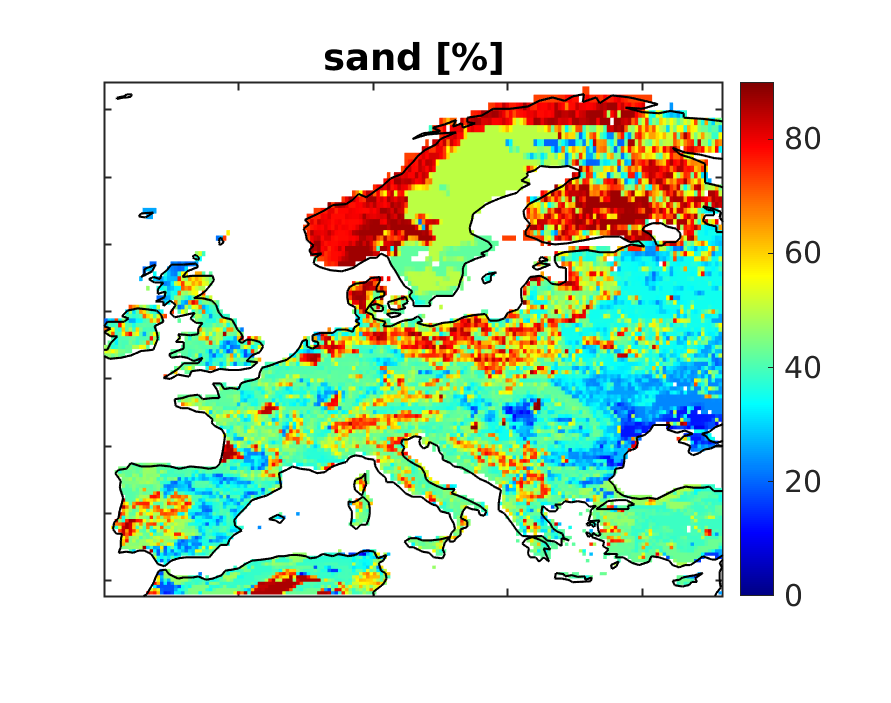}
\includegraphics[width=0.32\textwidth]{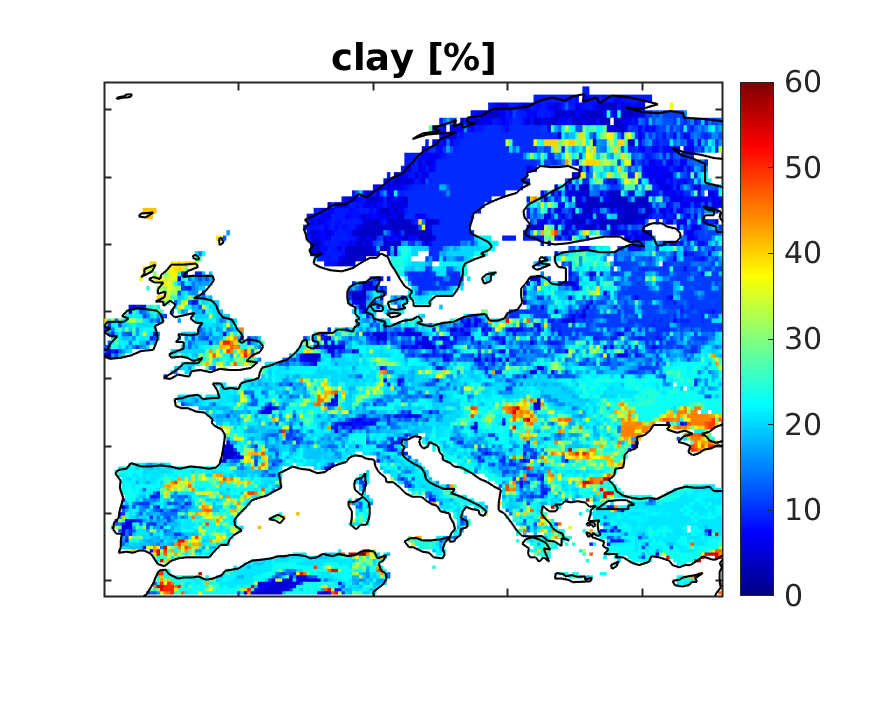}\\
\includegraphics[width=0.32\textwidth]{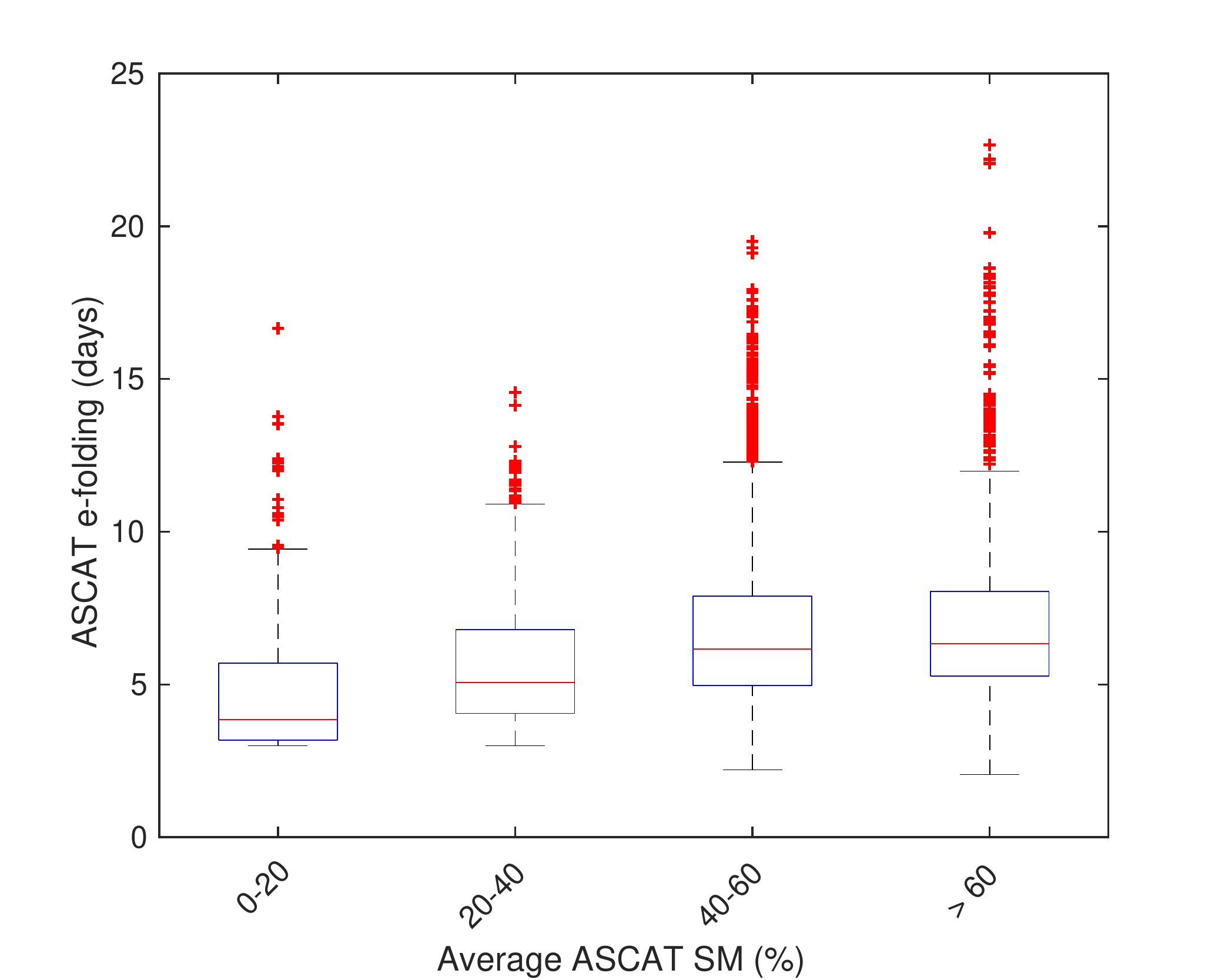}  
\includegraphics[width=0.32\textwidth]{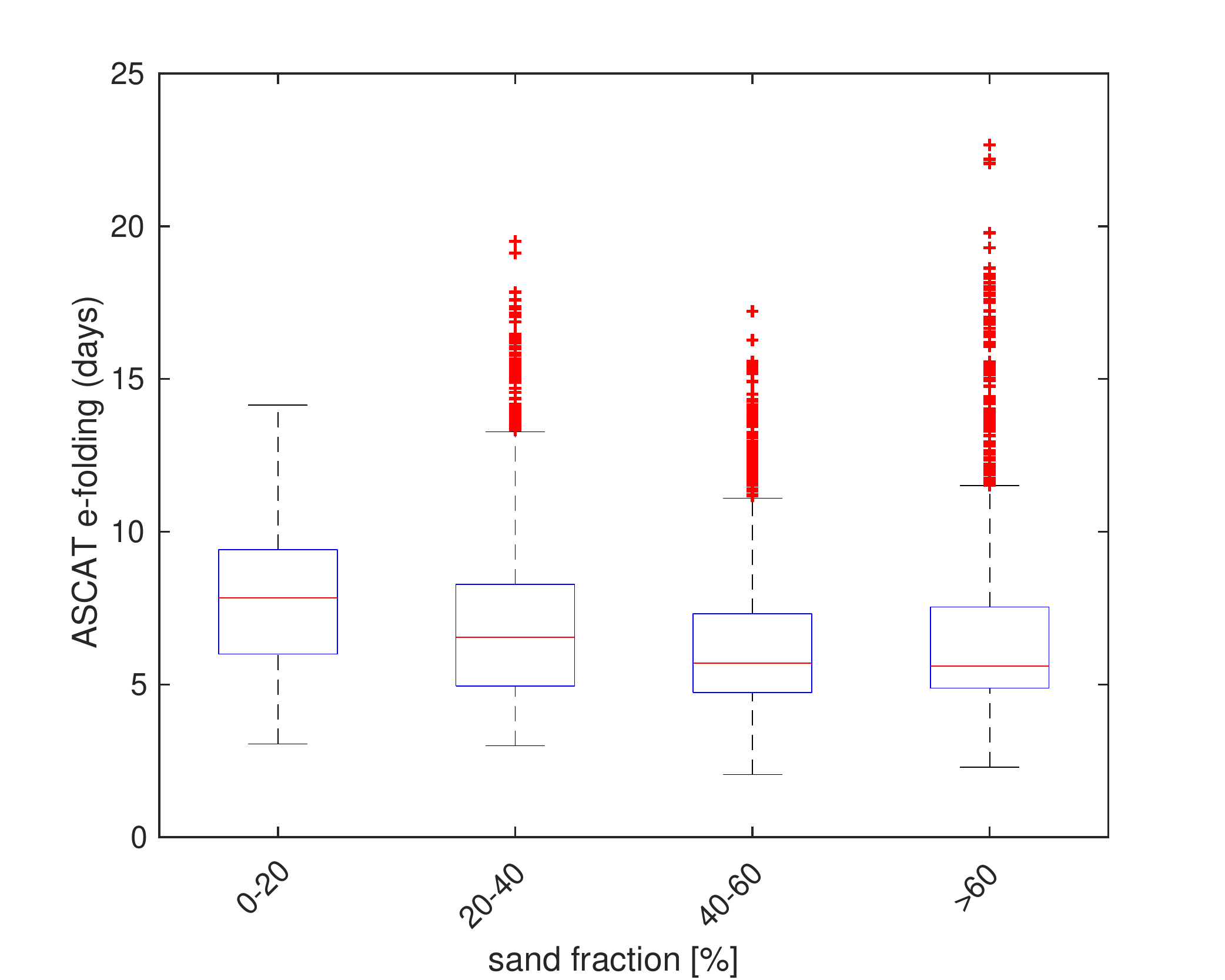}
\includegraphics[width=0.32\textwidth]{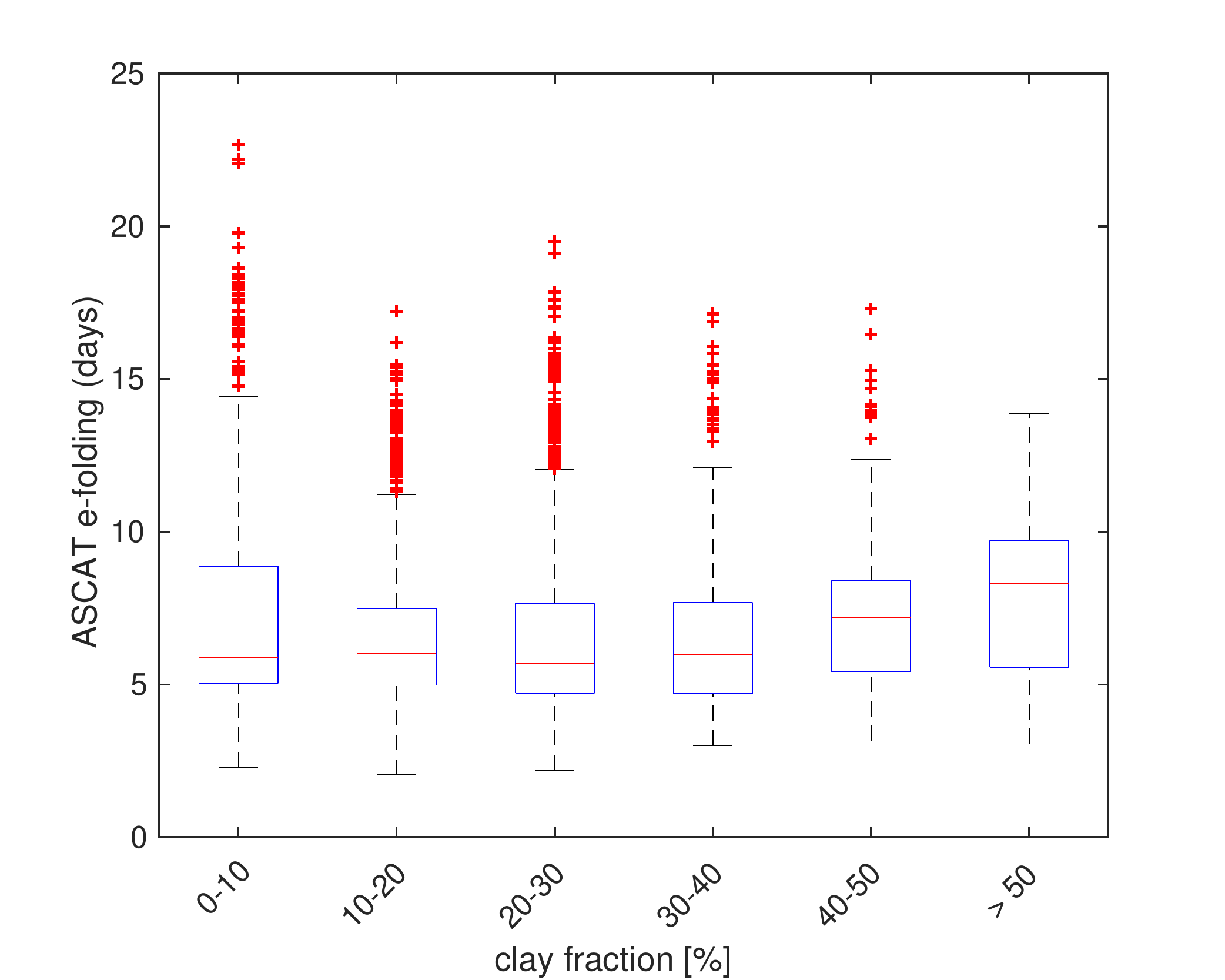}
\end{center}
    \vspace{-1cm}
\caption{Distribution of surface soil moisture, sand and clay soil fractions in Europe and their relations with surface soil moisture persistence. Top: Maps of average soil moisture from ASCAT for the period 2010-2016 (left), soil sand fraction (middle) and soil clay fraction (right). Bottom: relation between e-folding time scale (Fig.\ref{fig:efolding_map}), and long-term mean soil moisture (left column), soil sand fraction (middle column), and soil clay fraction (right column).}
\label{fig:explain_efolding}
\end{figure}

A map of e-folding time over Europe representing surface SM memory is shown in Figure \ref{fig:efolding_map}. SM e-folding times range from 2 and 20 days, with a median value of 5.92. The spatial patterns reveal areas of higher persistence such as the south-west of the Iberian Peninsula, the inner regions of the Balkan Peninsula, Eastern Europe and Norway. The obtained persistence patterns are further investigated in Figure \ref{fig:explain_efolding}, where the obtained e-folding times are compared to the distributions of average soil moisture, and soil texture classes. Satellite observations also show that lower soil moisture persistence is obtained for sandier soils and in regions that are more arid. This confirms the results obtained at REMEDHUS stations hold for a greater variety of soil and climate conditions. Also, our results agree well with previous studies that estimated soil moisture persistence by fitting an exponential model to drydown events \cite{McColl2017b,Ruscica2020}. 

Our results indicate that sandier soils are likely to exhibit lower memory, which is consistent with the fact that soil suction --and, therefore, the soil's ability to retain water in the presence of $ET$-- decreases with increasing sand content. Also, the observed decrease of soil moisture persistence with increasing aridity can be explained by the fact that the atmospheric demand for water is more intense in arid climates. For both relations, however, there is substantial unexplained variance that is likely due to vegetation effects and that should be subject of further research. This case study illustrates how the temporal autocorrelation of time series and their e-folding time can be used to characterize the persistence of variables such as the surface soil moisture under the assumption that they can be modeled as a stochastic, red noise process.

\subsection{Exploiting directional wind persistence in wind farms design}

Automatic turbine layout and planning in wind farms is an important topic in wind energy production, in order to obtain more efficient and profitable wind farms designs. Usually, the problem of turbine layout in wind farms has been tackled as a constrained optimization problem, in which meta-heuristic algorithms (mainly evolutionary computation techniques \cite{del2019bio}) have been proposed. The seminal paper in the use of evolutionary computation techniques in turbine layout problems was \cite{Mosetti94}, where a first genetic algorithm was proposed. The model used consisted of modeling the wind farm as a square divided into cells, in which turbines can be located. A wake model was proposed, and several experiments considering different average wind speed and direction were presented. The same binary model was used  some years later in \cite{Grady05}, in which an improved genetic approach was proposed to improve the results in \cite{Mosetti94}. Another improvement with the same model has been proposed in \cite{Emami10}. That paper proposed a modification of the objective function of the problem, to take into account deployment cost and efficiency of the turbines. The authors showed that this modification leaded to better layout designs than in previous approaches. Other works dealing with evolutionary computation in turbine layout problems are \cite{Riquelme07}, where a variable-length genetic algorithm with novel procedures of crossover was proposed, \cite{gonzalez2010optimization}, where a hybrid evolutionary algorithm was applied to wind farm designs, or the works \cite{Sisbot10} based on a multi-objective evolutionary algorithm and \cite{Saavedra11}, where evolutionary algorithms were hybridized with local heuristics to improve their performance in a turbine layout problem. Very recent works have also applied evolutionary computation techniques, such as in \cite{Wilson2018}, where evolutionary algorithms proposed for a wind turbine layout problem were described, or \cite{Long2020}, where a data-driven based evolutionary algorithm for wind turbine layout problems was described.

In all these and other previous approaches, it has been pointed out that the main issue to obtain a good quality turbine layout is that a realistic {\em wake model} must be considered \cite{Shakoor16}. Note that the computational complexity to solve a wind farm layout problem optimization primarily comes from the repetitive calculation of the wind speed loss known as {\em wake effect} caused by the presence of wind turbines upwind \cite{Tao20}. The vast majority of the literature on wind turbines layout optimization use the 1D Jensen's model \cite{jensen1983note,Shakoor16} assuming that the wake expands linearly behind a wind turbine, and that the far wake region has a cone form. Note that, in general, the evaluation of a wind farm layout considering a wind speed direction has a complexity of $O(n^2)$, where $n$ is the number of turbines considered. In other words, the computation of the influence matrix of a wind farm layout results in an evaluation time that is quadratic with the number of
turbines to be evaluated. In addition, this is only for a single wind speed direction, and normally Monte Carlo simulations are considered in wind layout problems, in which thousand of hours (with different wind direction and speeds) are simulated in order to obtain a final evaluation of the layout quality for a wind farm \cite{marmidis2008optimal}. Thus, in the last years, different approaches dealing with methods to reduce the computational complexity of the layout evaluation have been proposed. In \cite{Wagner2013} a fast local search algorithm is proposed for wind turbines layout optimization, in which the priority of the approach is to reduce the computational cost of the layout evaluation, by applying a reduce number of changes from a feasible solution with a local search approach. Other approaches bet for randomized algorithms in order to speed up the evaluation process, such as in \cite{feng2015solving}, where a random search algorithm with a strategy for reducing the evaluation cost of the layouts is proposed. In \cite{zergane2018optimization} also a random search approach is proposed with the help of a pseudo-random number generation, in order to reduce the computational complexity of the wind turbine layout problem. Finally, other works propose the use of data-driven approaches within meta-heuristics in order to reduce the computational complexity of the problem, such as \cite{Long2020}.

In this case study we show that the wind directional persistence can be exploited to reduce the computational complexity of a layout quality estimation in wind farms. First we will describe data from three real wind measurement towers in wind farms located in Spain. Long time series data of over 10 years are available for this study. We show that there is a marked directional persistence of wind speed in all the cases, and how this point can be used to reduce the computation time in the calculation of the quality of a given layout.

\subsubsection{Real wind speed data considered and directional wind persistence}

\begin{figure}[!ht]
\begin{center}
\includegraphics[draft=false, angle=0,width=8cm]{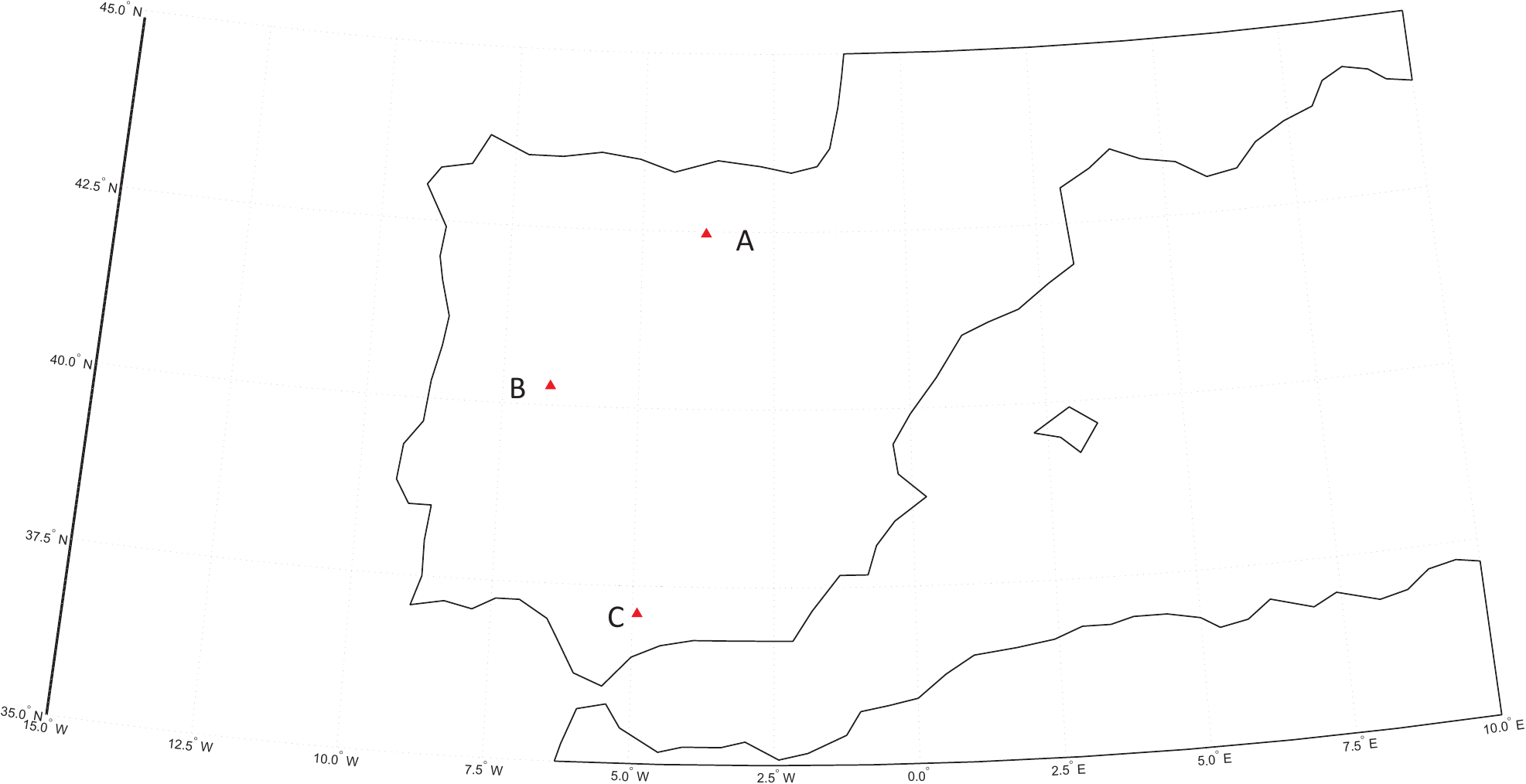}
\end{center}
\caption{\label{fig:maps_WF} Location of the three wind speed measurement towers considered in Spain.}
\end{figure}

Three wind measurement towers in wind farms has been considered (Figure \ref{fig:maps_WF} labeled ``A'', ``B'' and ``C''). Note that the wind farms selected cover different parts of Spain, north, center and south, characterized by different wind regimes. Long time series of hourly wind speed have been measured in each wind farm: in ``A'', data ranges from 11/01/2002 to 29/10/2012, in ``B'' data ranges from 23/11/2000 to 17/02/2013, and finally in ``C'', the data available ranges between 02/03/2002 and 30/06/2013 (over 10 years of hourly data in all cases).

\begin{figure}[htpb]
    \centering
      \subfigure[]{
        \label{Mon_mod}
        \includegraphics[width=0.32\textwidth]{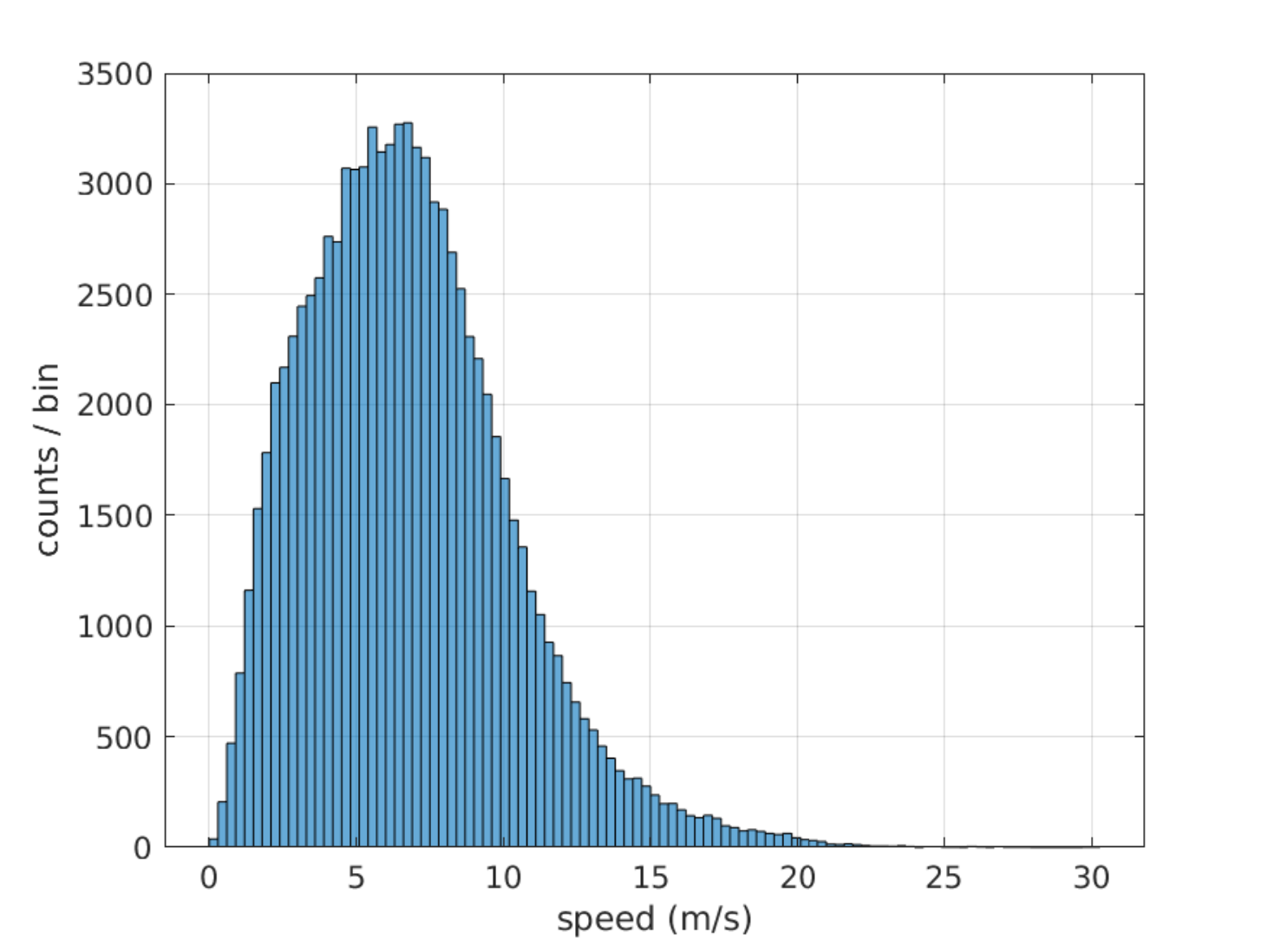}}
      \subfigure[]{
        \label{Pena_mod}
        \includegraphics[width=0.32\textwidth]{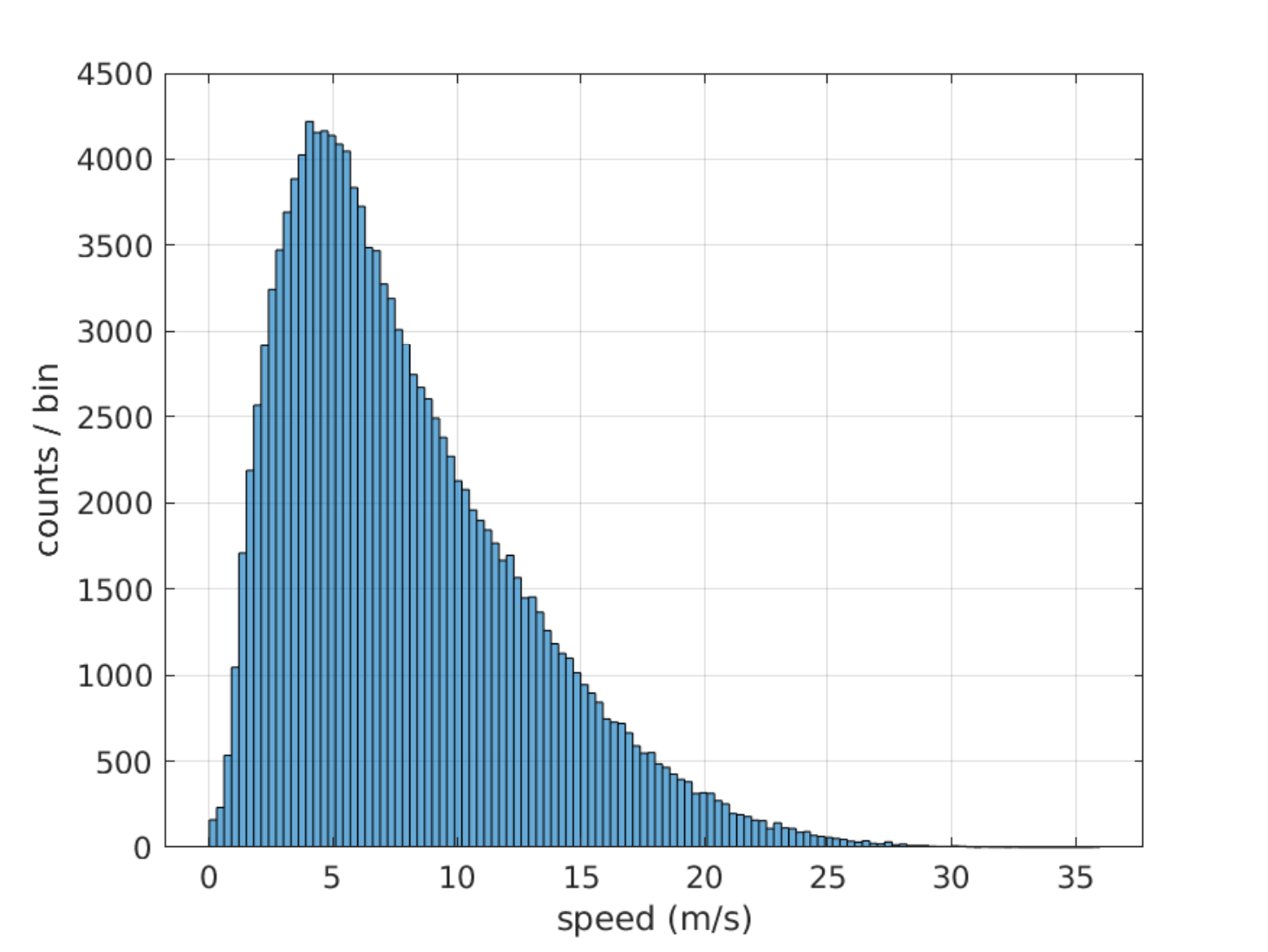}}
      \subfigure[]{
        \label{Jara_mod}
        \includegraphics[width=0.32\textwidth]{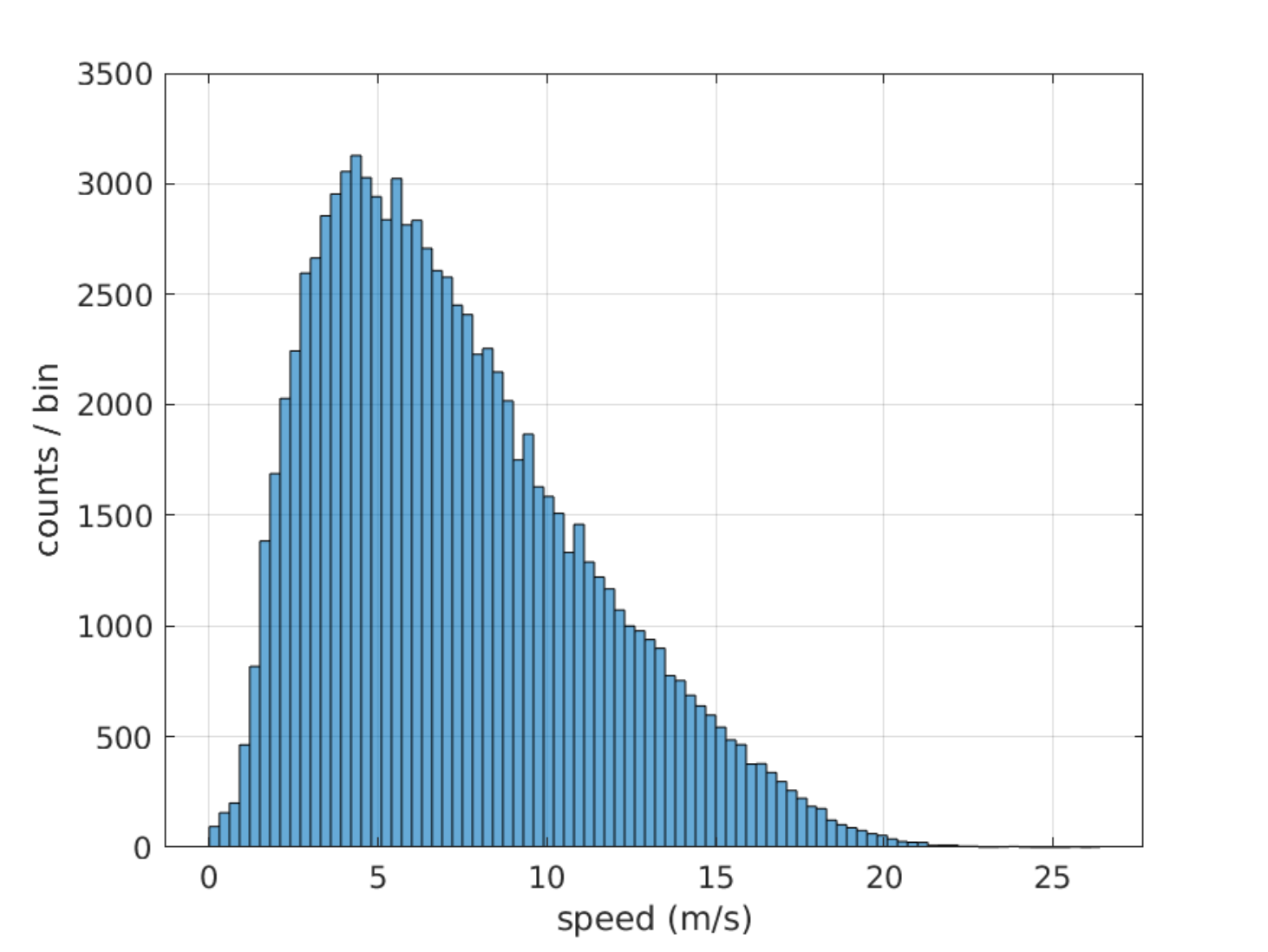}}
    \caption{Histograms of wind speed module for the three wind farms: A, B and C. Observe that this modulus can be modelled by a Weibull distribution.}
    \label{fig:histMod}
\end{figure}

Available data consist of measurements of wind speed module $v$ and wind direction $\theta_v$ with hourly temporal resolution. Wind speed module has a resolution of $0.01$ m/s and wind direction a resolution of $0.1$ degree. Figure \ref{fig:histMod} show the histograms of wind speed module for the three wind farms considered. Table \ref{tab:W_wind-param} shows the parameters of the Weibull distributions for the three wind farms.

\begin{table}[htpb]
    \centering
    \begin{tabular}{ccc}
        \hline
        Wind Farm & scale ($\lambda$) & shape ($k$) \\ \hline\hline
         A & $7.6553$ & $2.0494$ \\ \hline
         B & $8.8942$ & $1.7132$ \\ \hline
         C & $8.2415$ & $1.9473$ \\ \hline
    \end{tabular}
    \caption{Weibull parameters ($\lambda$ and $k$) estimation for the three wind farms considered.}
    \label{tab:W_wind-param}
\end{table}

\begin{figure}[htpb]
    \centering
      \subfigure[]{
        \label{Mon_ang}
        \includegraphics[width=0.32\textwidth]{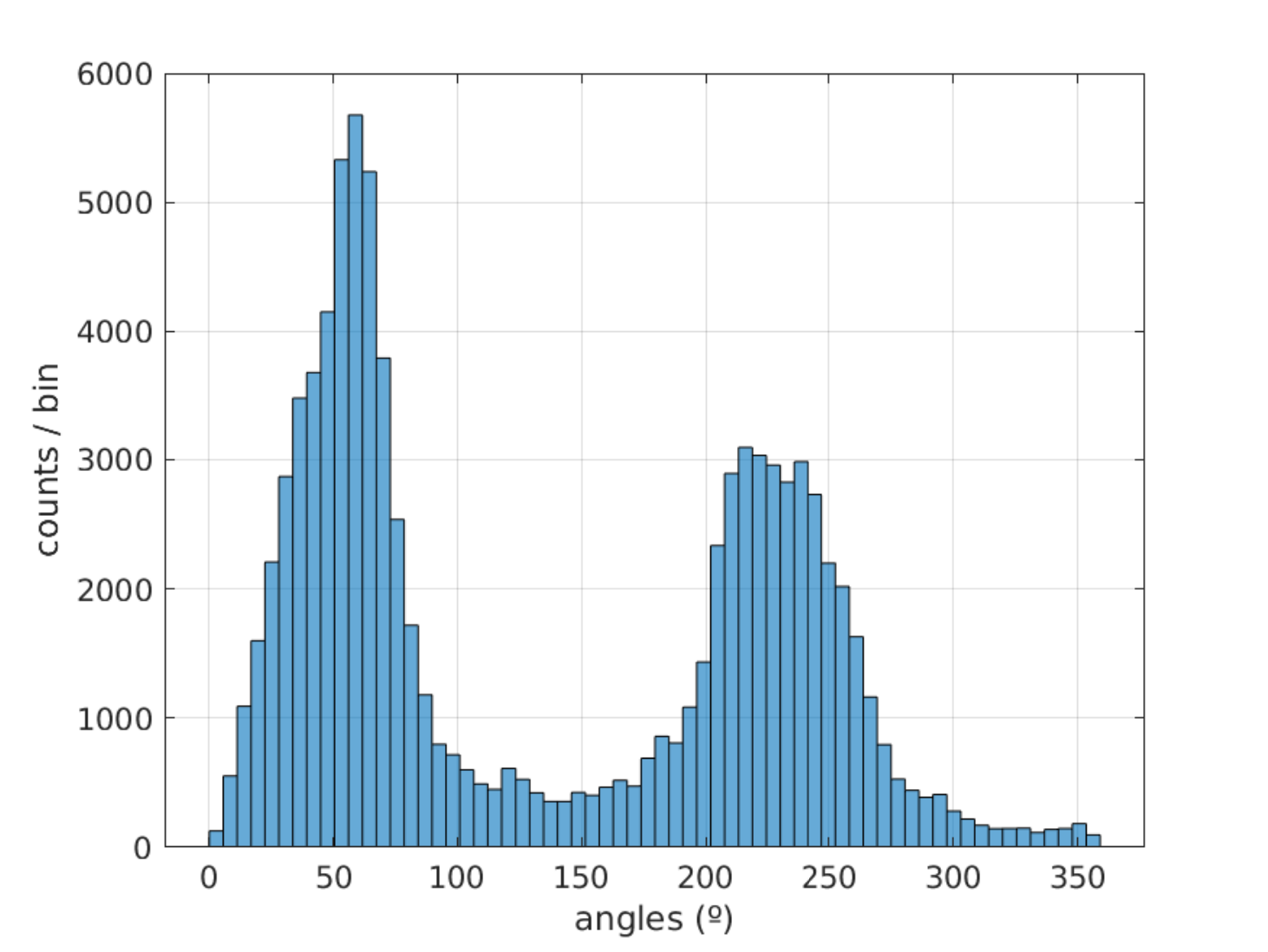}}
      \subfigure[]{
        \label{Pena_ang}
        \includegraphics[width=0.32\textwidth]{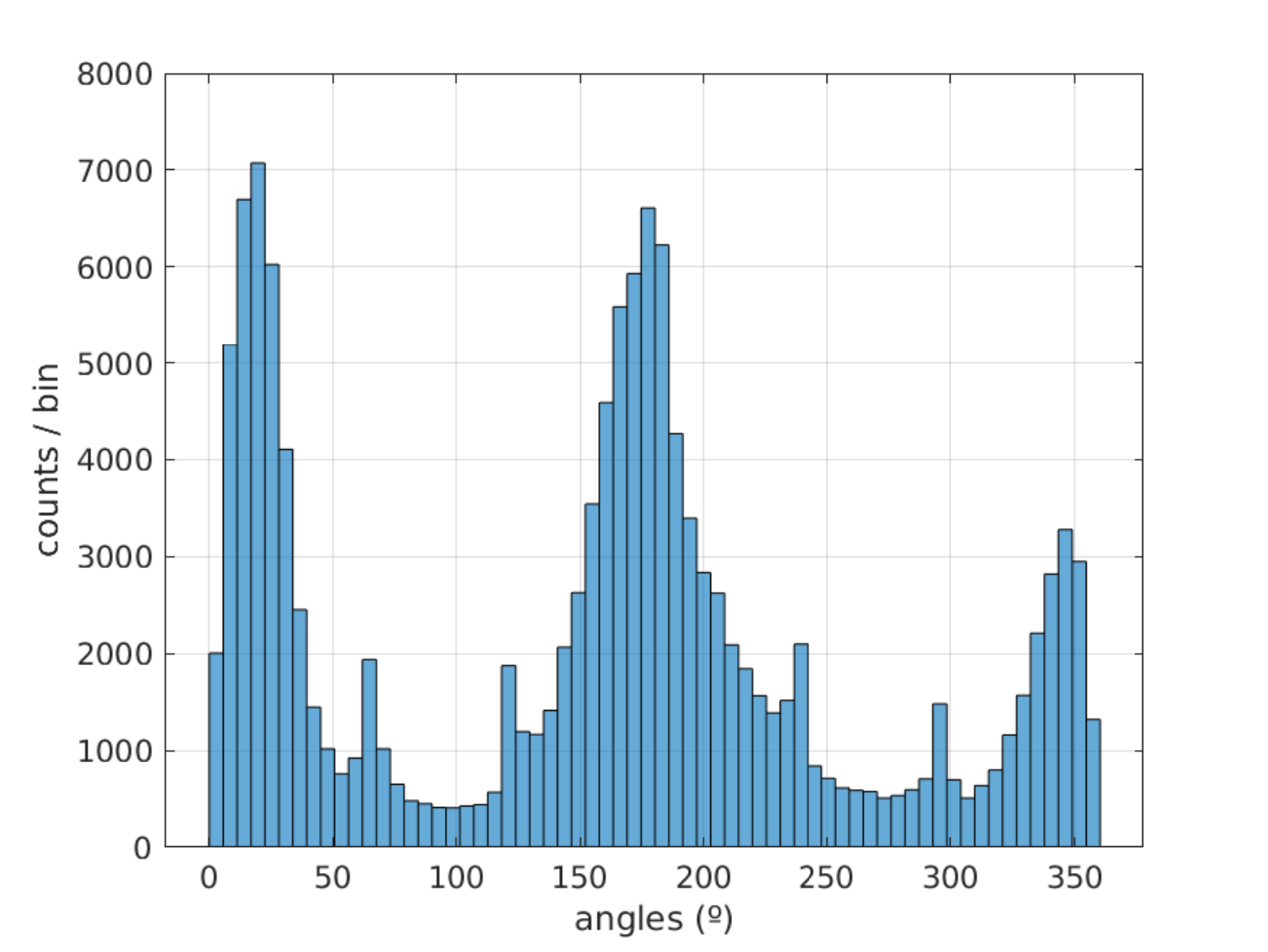}}
      \subfigure[]{
        \label{Jara_ang}
        \includegraphics[width=0.32\textwidth]{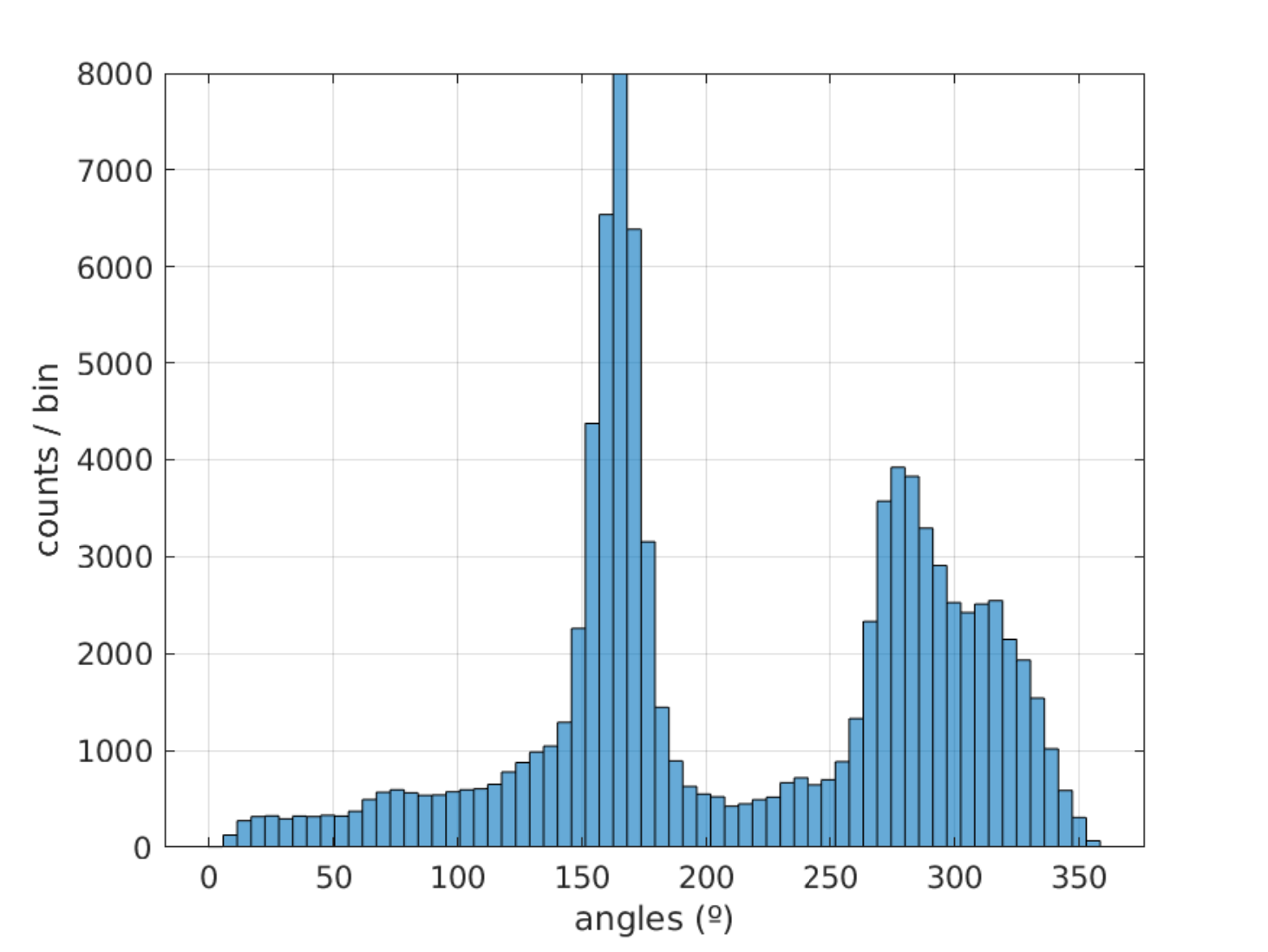}}
        \subfigure[]{
        \label{Mon_wr}
        \includegraphics[width=0.32\textwidth]{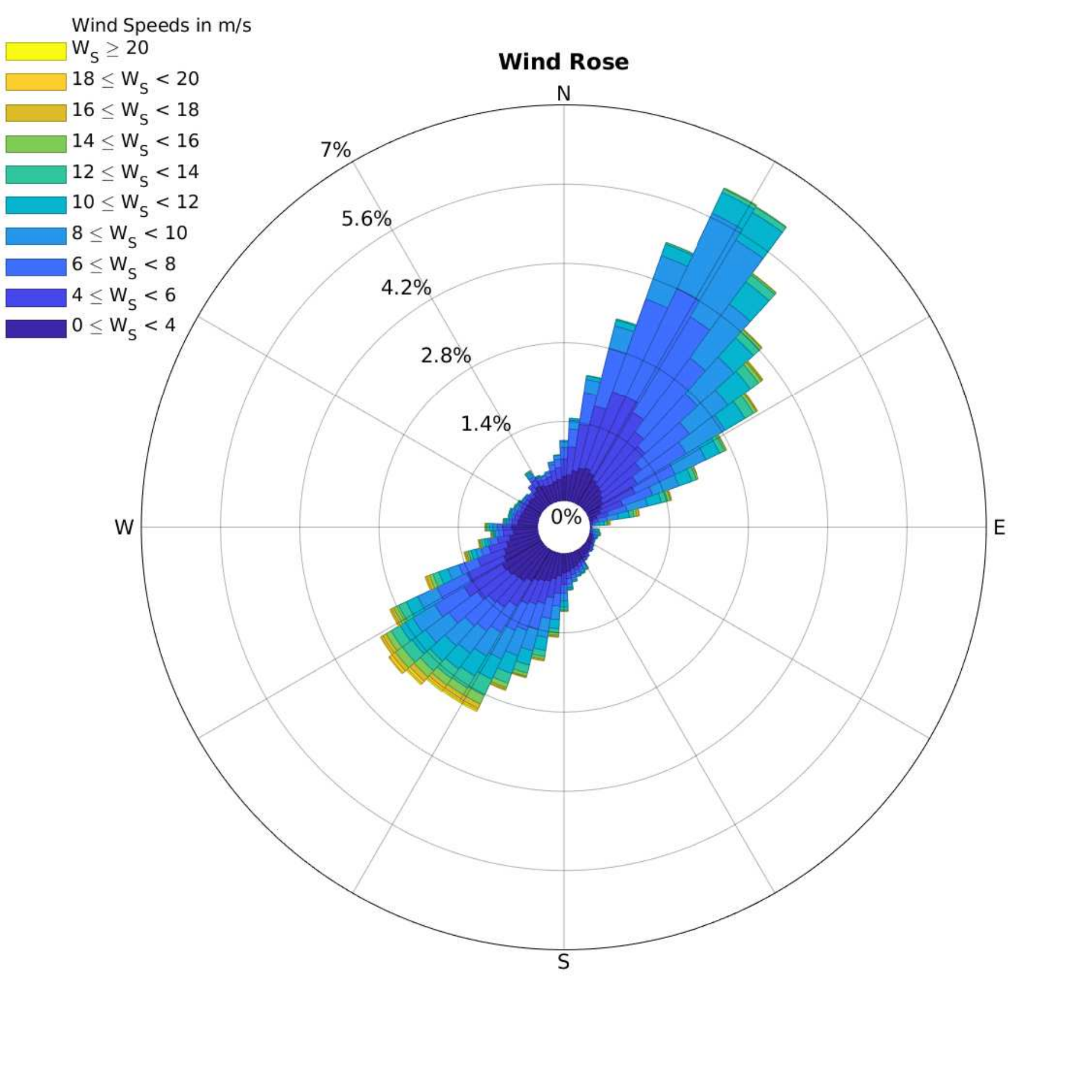}}
      \subfigure[]{
        \label{Pena_wr}
        \includegraphics[width=0.32\textwidth]{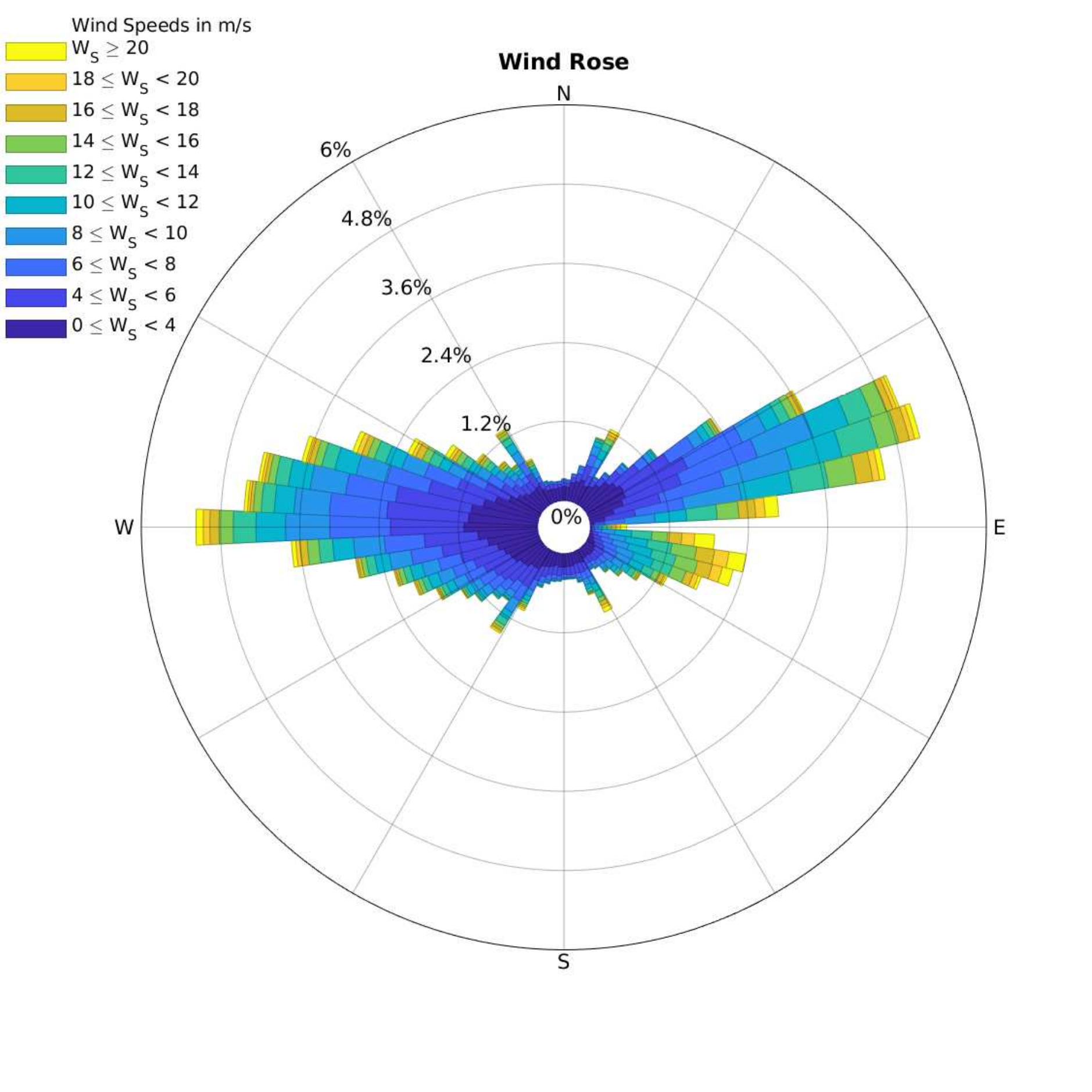}}
      \subfigure[]{
        \label{Jara_wr}
        \includegraphics[width=0.32\textwidth]{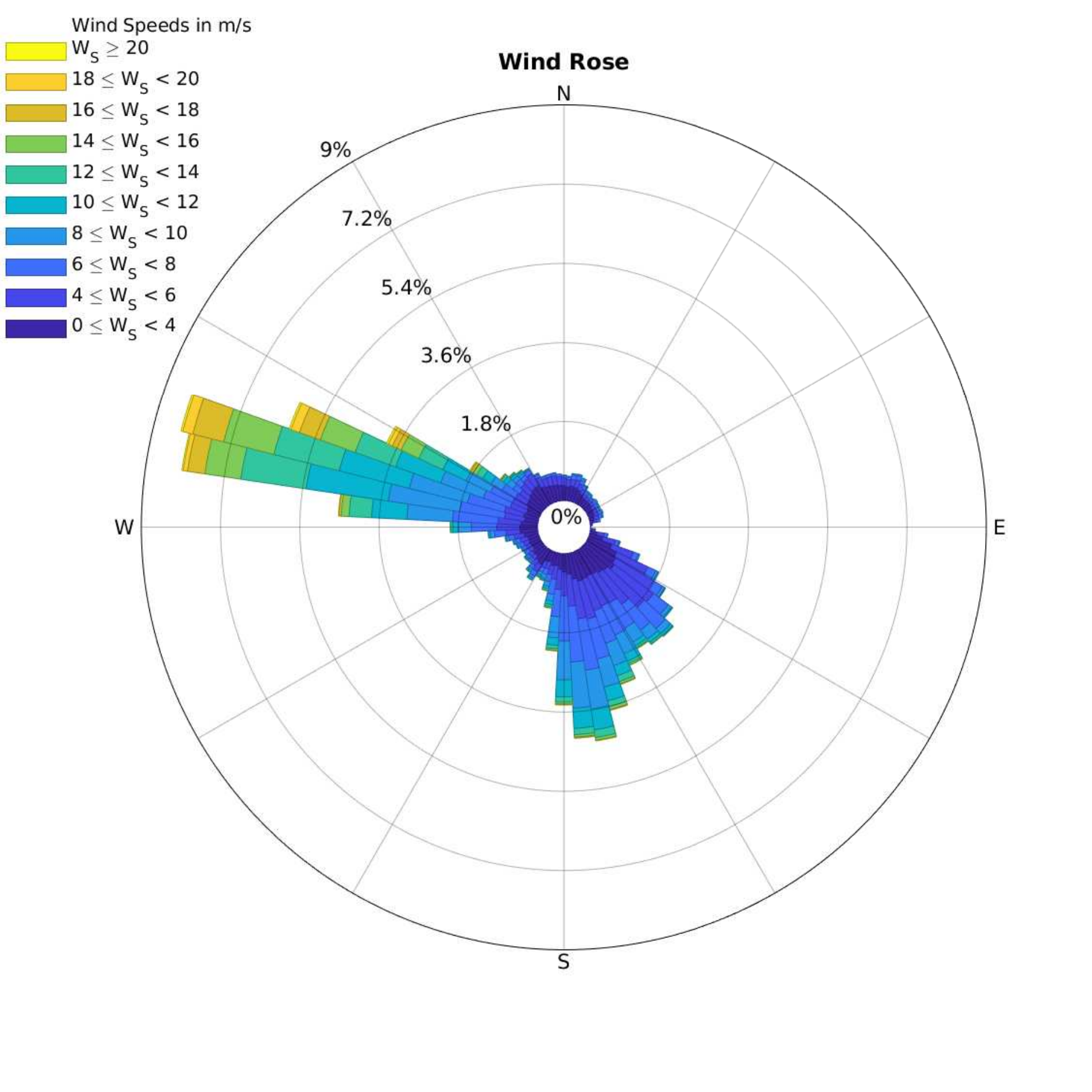}}
    \caption{Wind direction histograms and wind roses for the three wind farms considered; (a) Wind direction histogram for wind farm A; (b) Wind direction histogram for wind farm B; (c) Wind direction histogram for wind farm C; (d) Wind rose for wind farm A; (e) Wind rose for wind farm B; (f) Wind rose for wind farm C.}
    \label{fig:histangwr}
\end{figure}

In our case, is it indeed more interesting the estimation of the direction histograms and wind roses, since they suggest a clear directional persistence of the wind in the three wind farms considered. Figure \ref{fig:histangwr} shows the wind direction histograms and wind roses for the three wind farms considered in this case study. As can be seen, all the wind farms present directional persistence of wind, with some directions presenting much more frequency than others. In fact, note that there are some directions in which the frequency of wind blowing is nearly negligible. We will exploit this fact of directional persistence to improve the computation time for evaluating a given turbine layout. 

\subsubsection{Computational time reduction by exploiting wind directional persistence}

\begin{figure}[htpb]
    \centering
      \subfigure[]{
        \includegraphics[width=0.25\textwidth]{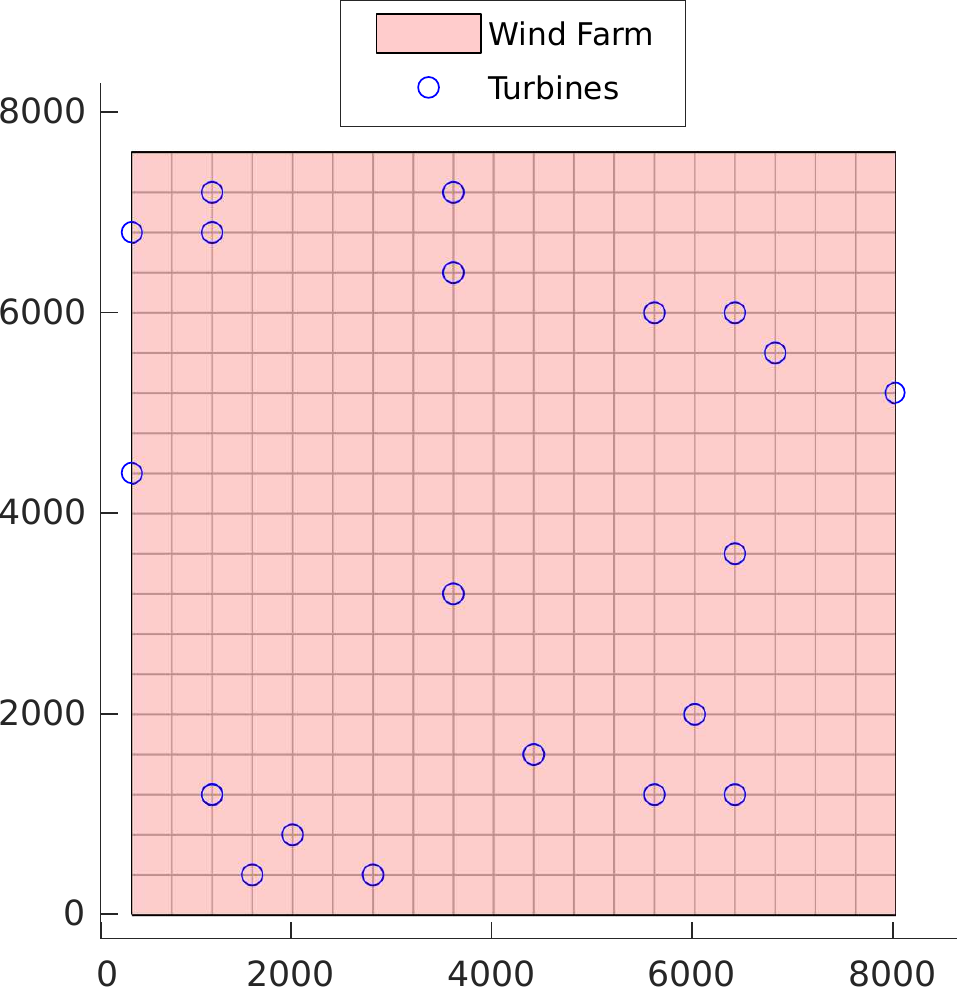}}
      \subfigure[]{
         \includegraphics[width=0.3\textwidth]{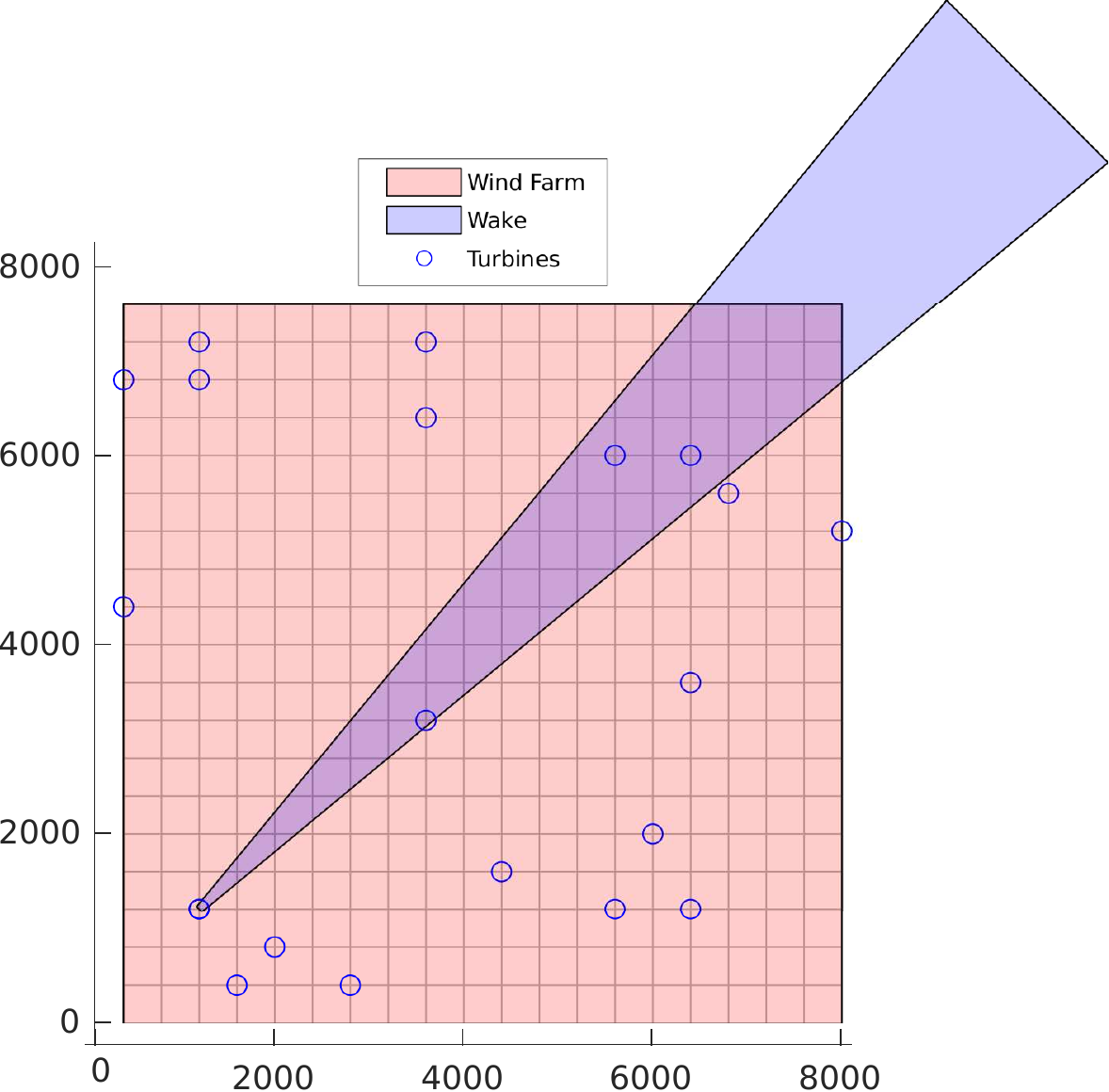}}\\
             \subfigure[]{
          \includegraphics[width=0.28\textwidth]{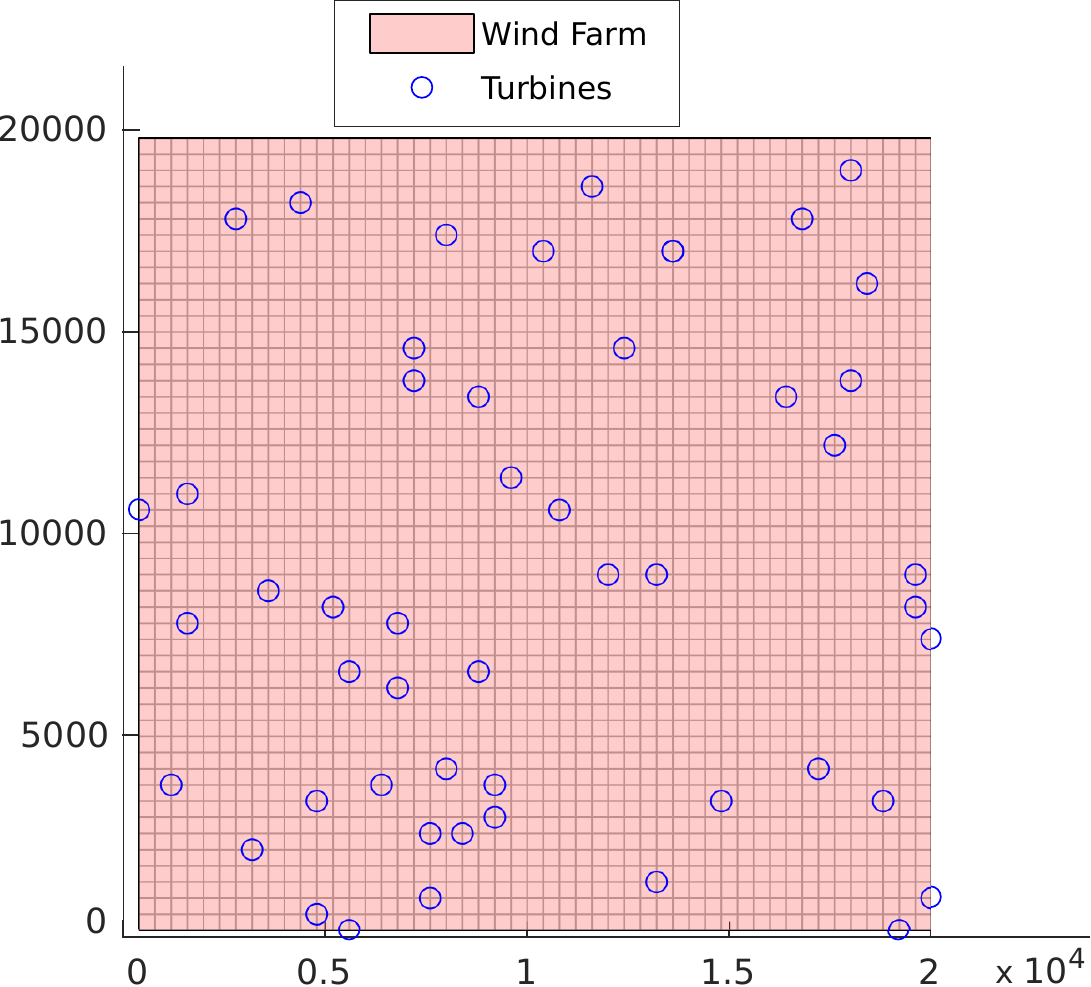}}
      \subfigure[]{
          \includegraphics[width=0.3\textwidth]{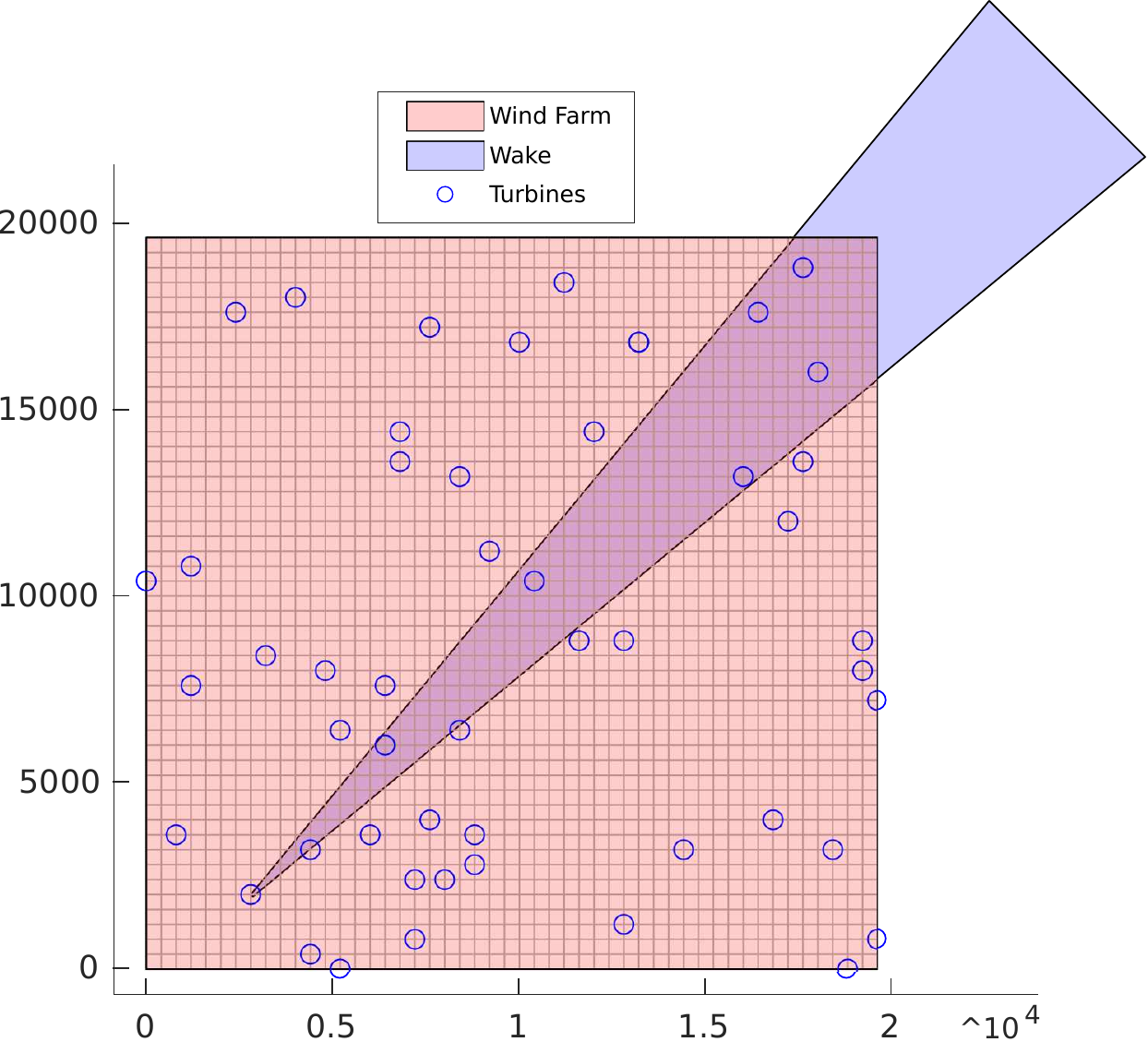}}
    \caption{Examples of two layouts in a wind farm; (a) $20 \times 20$ layout; (b) Wake effect for a turbine in the $20 \times 20$ layout; (c) $50 \times 50$ layout; (d) Wake effect for a turbine in the $20 \times 20$ layout.}
    \label{fig:layout}
\end{figure}

\begin{figure}
    \centering
    \includegraphics[width=0.4\textwidth]{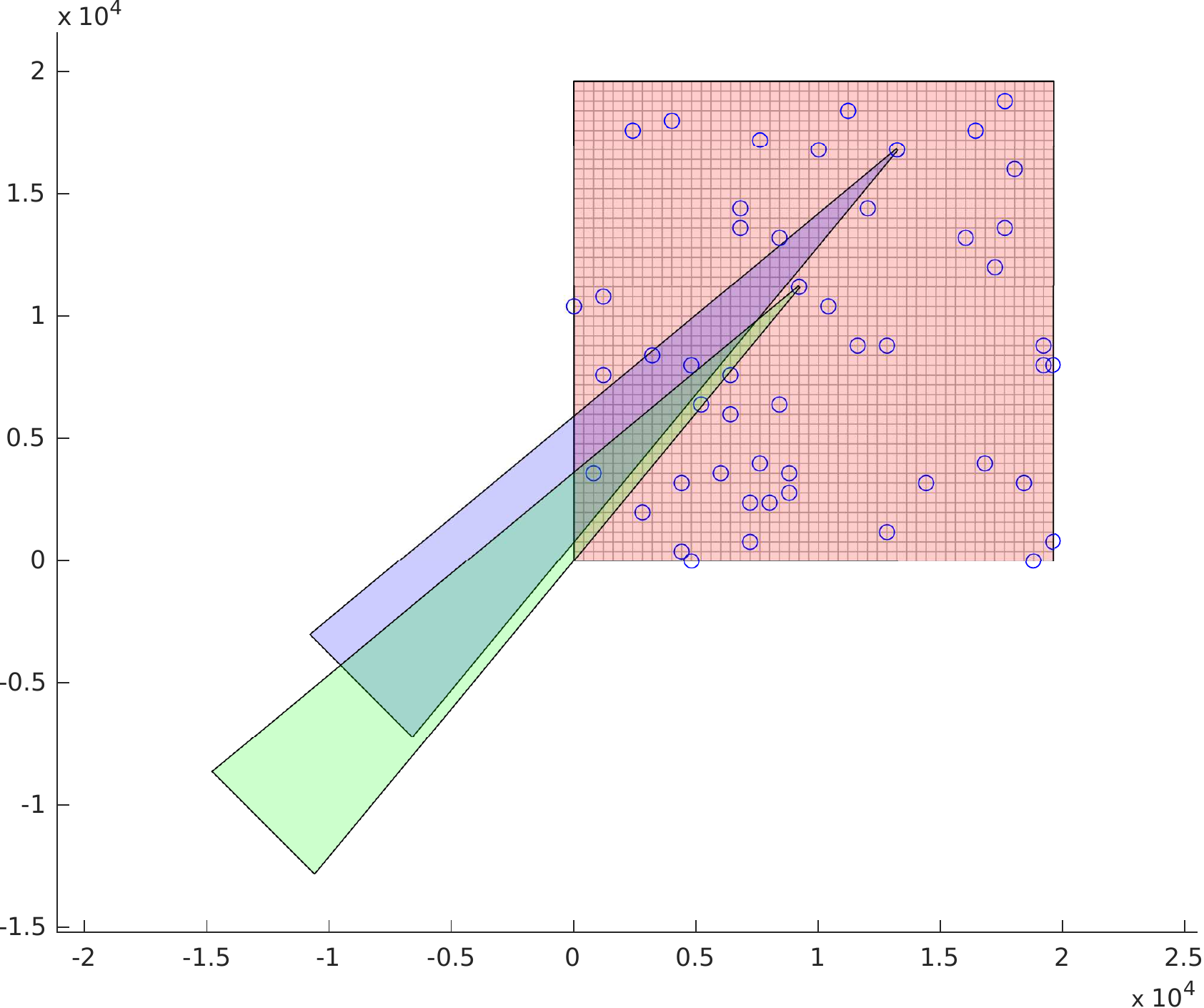}
    \caption{Wakes overlapping example.}
    \label{fig:WakeOverlap}
\end{figure}

Different models for wind farms design can be found in the literature, but the most commonly used consist of a square grid, considering that a turbine can be located in a cell of the grid. This model is very convenient to encode solutions to the layout problem, which can be easily managed by meta-heuristic algorithms: note that a square binary matrix can be considered, so a $1$ in a cell stands for a wind turbine to be located, and a $0$ means no turbine in that cell. Figure \ref{fig:layout} shows two examples of possible layouts. Figure \ref{fig:layout} (a) shows a layout in a wind farm of $20 \times 20$ cells (each cell is considered to have a length of $5D$ to avoid  constraints related to proximal location of turbines). In this case $20$ turbines will be considered. In turn, Figure \ref{fig:layout} (c) shows a wind $50 \times 50$ cells wind farm, with $50$ turbines to be located. Figure \ref{fig:layout} (b) and Figure \ref{fig:layout} (d) show two possible wake effects in each wind farm, for a given turbine. In both Figures \ref{fig:layout} (b) and \ref{fig:layout} (d) a wind with South-East direction is considered.
Note that in both cases a given turbine may affect by its wake to many other turbines in the wind farm. Note that wakes from two different turbines may overlap, adding their wake effect to other turbines, as shown in Figure \ref{fig:WakeOverlap}.

Note that the original wind direction data is obtained with a resolution of $0.1$ degree. Thus, given a layout to be studied (calculating the power in a given time period), it is necessary to calculate the final wind speed in each point of the wind farm, including the wake effect of all turbines. Since the wake effect mainly depends on the wind direction, there will be a different wake for each wind direction considered in the wind farm. A first option to reduce the computational cost of this calculation is to discretize the wind direction to a number of directions, typically $8$ or $16$ wind direction will provide a good quality approximation at this stage. From this point, we will show how the directional persistence of the wind can be exploited to further reduce the computation time of the wake effect.

\begin{figure}[htpb]
    \centering
    \includegraphics[width=0.5\textwidth]{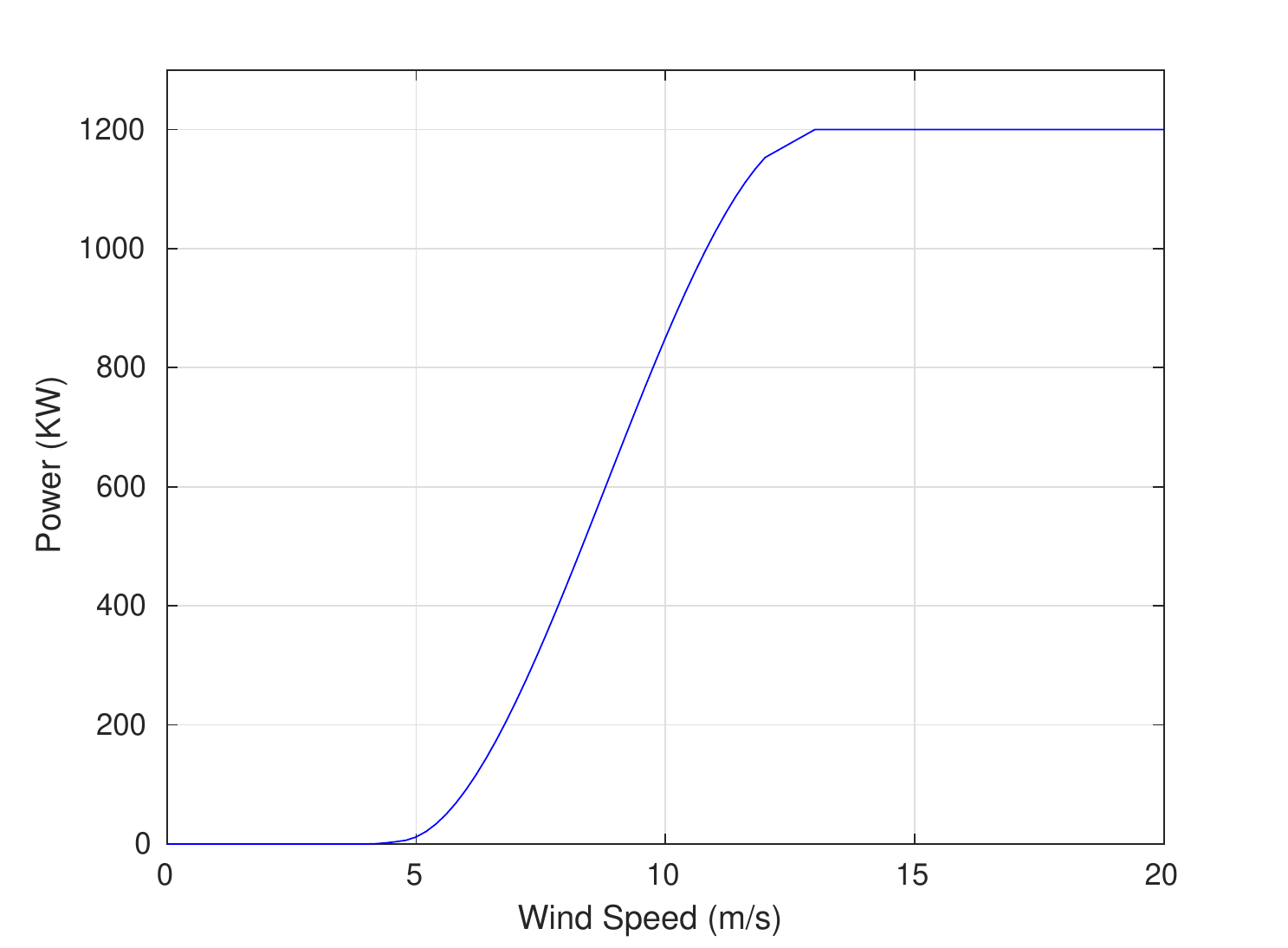}
    \caption{Wind turbine power curve considered.}
    \label{fig:powerCurver}
\end{figure}

Let us suppose, without loss generality, that we deal with a wind rose discretized to $8$ directions, that is, we consider that the wind direction can be in $8$ different states, $s_1, \ldots, s_8$. However, note that if there is wind directional persistence, some of these states will be very few times visited, i.e. in this case, the persistence consists of the prevalence of some of the states over the rest. This can be easily seen in the wind rose of Figure \ref{fig:histangwr}. For example in Figure \ref{fig:histangwr} (d), wind rose of wind farm A, note that the system is mainly in states $s_1, s_2, s_5, s_6$, whereas states $s_3,s_4,s_7, s_8$ are much less frequent. We propose to exploit this fact to reduce the computation time in the wake effect calculation for a given turbine layout in a wind farm. Note that given a turbine layout, the wind speed in the wind farm is modified due to the presence of turbines (because their wake effect), and therefore the power produced in the turbines, which is usually calculated by means of a curve power (non-linear relationship between wind speed and power generation in a wind turbine or wind farm). Figure \ref{fig:powerCurver} shows an example of power curve which will be used in the experiments of this case study. 

\subsubsection{Experimental results}
In a turbine layout problem, optimization algorithms look for the best location of turbines in wind farms which maximize the total amount of power generated by the facility. As previously mentioned, this process requires, in the majority of cases, the execution of a huge number of iterations to find the best places for each farm, which severely affect the computational cost of turbine layout processes.

We present here an experiment to show that the computation time of a wind turbine layout power outcome can be reduced if the wind farm considered presents directional persistence of the wind. We have evaluated the time reduction improvements of considering wind directional persistence for the three wind farm whose wind roses were shown in Figure~\ref{fig:histangwr}. First of all, we calculate the computation time of the layout shown in Figures \ref{fig:layout} (a) and (c), with $20 \times 20$ and $50 \times 50$, respectively. We use a mid-range microprocessor Intel(R) Core(TM) i7-8700 CPU for all the computations carried out. In this first experiment, the computation of the wind speed in the wind farm is carried out with the original wind data (module and direction) for the three wind farms considered. Table~\ref{tab:times20x20} (Experiment 1) shows the computation and the total amount of generated power in the whole period available for three wind farms A, B, C and the grid $20\times20$ layout, (Figure~\ref{fig:layout} (a)). Table~\ref{tab:times50x50} (Experiment 1) shows the same results for the $50\times50$ layout, (Figure~\ref{fig:layout} (c)). As can be seen in both tables, the computation time of the layouts is extremely high (about $200$ seconds for the case of the $20\times20$ layout and over $2000$ seconds for the $50\times50$ layout. Recall, that an optimization process would need a high number of optimization steps like this one, and this computation time makes any optimization process unfeasible. 

Can we somehow reduce the computation time for this problem? The answer is yes, if we assume some reasonable cost on accuracy. To this end, we will consider discretization of the wind data, and also we will exploit the directional persistence of the wind in the wind farms considered. First, Experiment $2$ in Tables~\ref{tab:times20x20} and \ref{tab:times50x50} shows the results obtained in the calculation of the total power in each wind farm, when the number of directions has been discretized to $8$, i.e. we consider here $8$ possible states of the system, in terms of the wind speed direction. As can be seen, in this case the computation time is already reduced a lot (without even consider persistence), to $4$ seconds in the case of $20\times20$ wind farm and around $8$ for the case of the $50\times50$ wind farm. If we consider as well a further discretization, not only in direction, but also in wind speed module, we have Experiment 3, as can be seen in both tables, with results around $0.5$ seconds of computation time for the $20\times20$ layout and around $1.5$ second for the $50\times50$ layout. Note that with discretization, the calculation of the produced power produces a small error, which seems in both cases assumable, considering the large improvement in computation time.

At this point, we can obtain further computation time improvements in the calculation of the power for a given layout by considering the wind direction persistence. The idea is to discard those directions which contribute the less to the generation of wind power in the wind farm. For this, we calculate the following parameter:
\begin{equation}
    Pw_i = \bar d(s_i) \cdot \bar p_i(s_i) \cdot n_i,\quad 1 \le i \le 8
\end{equation}
where $Pw_i$ stands for the mean power contribution by system state (each of the directions), $d(s_i)$ is the mean duration of the wind speed in a given state (direction), $\bar p_i$ is the mean power generated in a given wind direction, and finally, $n_i$ stands for the frequency of each state $s_i$. 

We can use parameter $Pw_i$ to decide whether taking into account (or not) a given state $s_i$ (wind speed direction), in the calculation of the power generated in the wind farm. For this, we define a threshold $\gamma$ in percentage, and we establish that we consider all the systems states such that $Pw_i>\gamma$, discarding the rest of states. In Tables~\ref{tab:times20x20} and \ref{tab:times50x50} we have set $\gamma=1\%$, which means a power loss due to persistence of less than $1\%$, with an improvement of the computation time. This can be seen in Experiments 4 (direction discretization and directional persistence) and 5 (direction discretization, module discretization and directional persistence). Note that the best result is obtained in Experiment 5, with a final computation time of around $400$ ms in the $20 \times 20$ layout and around 1 second for the case of $20 \times 20$, with a very reasonable error in power production compared to the results of Experiment 1.

\begin{table}[!ht]
    \centering
    \begin{tabular}{|c|c|c|c|c|}
        \hline
         Experiment & Wind farm &  Number of samples & Total Power & Computation Time (s) \\ \hline \hline
          \multirow{3}{*}{1}& A & $91920$ & $5.6886\times10^8$ & $216.00$ \\ \cline{2-5}
          & B & $136023$ & $1.2111\times10^9$ & $262.13$ \\\cline{2-5}
          & C & $92713$ & $6.9802\times10^8$ & $29.74$ \\ 
          \hline
          \hline
          \multirow{3}{*}{2}& A & $91920$ & $4.9465\times10^8$ & $4.15$ \\ \cline{2-5}
          & B & $136023$ & $1.2044\times10^9$ & $4.51$ \\  \cline{2-5}
          & C & $92713$ & $6.6962\times10^8$ & $3.70$ \\
          \hline
          \hline
          \multirow{3}{*}{3}& A & $91920$ & $4.9586\times10^8$ & $0.50$ \\ \cline{2-5}
          & B & $136023$ & $1.2039\times10^9$ & $0.63$ \\\cline{2-5}
          & C & $92713$ & $6.7063\times10^8$ & $0.46$ \\
          \hline
          \hline
          \multirow{3}{*}{4}& A & $91920$ & $4.8201\times10^8$ & $3.82$ \\   \cline{2-5}
          & B & $136023$ & $1.1950\times10^9$ & $3.90$ \\  \cline{2-5}
          & C & $92713$ & $6.4326\times10^8$ & $2.81$ \\
          \hline
          \hline
          \multirow{3}{*}{5}& A & $91920$ & $4.8325\times10^8$ & $0.38$ \\   \cline{2-5} 
          & B & $136023$ & $1.1944\times10^9$ & $0.49$ \\  \cline{2-5}
          & C & $92713$ & $6.4411\times10^8$ & $0.37$ \\
          \hline
    \end{tabular}
    \caption{Computation time for the grid $20 \times 20$}
    \label{tab:times20x20}
\end{table}

\begin{table}[!ht]
    \centering
    \begin{tabular}{|c|c|c|c|c|}
        \hline
         Experiment & Wind farm &  Number of samples & Total Power & Computation Time (s) \\ \hline \hline
          \multirow{3}{*}{1}& A & $91920$ & $1.4173\times10^9$ & $2711.57$ \\ \cline{2-5}
          & B & $136023$ & $3.0275\times10^9$ & $4015.32$ \\\cline{2-5}
          & C & $92713$ & $1.7410\times10^9$ & $2761.52$ \\ 
          \hline
          \hline
          \multirow{3}{*}{2}& A & $91920$ & $1.5279\times10^9$ & $6.89$ \\ \cline{2-5}
          & B & $136023$ & $3.0109\times10^9$ & $10.05$ \\  \cline{2-5}
          & C & $92713$ & $1.8298\times10^9$ & $7.07$ \\
          \hline
          \hline
          \multirow{3}{*}{3}& A & $91920$ & $1.5288\times10^9$ & $1.22$ \\ \cline{2-5}
          & B & $136023$ & $3.0100\times10^9$ & $1.63$ \\\cline{2-5}
          & C & $92713$ & $1.8326\times10^9$ & $1.22$ \\
          \hline
          \hline
          \multirow{3}{*}{4}& A & $91920$ & $1.4967\times10^9$ & $6.53$ \\   \cline{2-5}
          & B & $136023$ & $2.9646\times10^9$ & $9.29$ \\  \cline{2-5}
          & C & $92713$ & $1.7581\times10^9$ & $6.06$ \\
          \hline
          \hline
          \multirow{3}{*}{5}& A & $91920$ & $1.4977\times10^9$ & $1.08$ \\   \cline{2-5}
          & B & $136023$ & $2.9633\times10^9$ & $1.47$ \\  \cline{2-5}
          & C & $92713$ & $1.7604\times10^9$ & $0.97$ \\
          \hline
    \end{tabular}
    \caption{Computation times for the grid $50 \times 50$}
    \label{tab:times50x50}
\end{table}

\subsection{Persistence in Machine Learning over data streams}
In this final case study, we show how the lack of persistence can yield a notable degradation of performance in ML predictive modeling over data streams (see Section \ref{Pers_ML}). To this end, we will first discuss a synthetically generated data stream, composed by two substreams, in which the values of the target variable $x[n,\boldsymbol{\theta}]$ to be predicted based on the observed values of $\boldsymbol{\theta}$ are given by:
\begin{equation}\label{eq:concepts}
x[n,\boldsymbol{\theta}] = \left\lbrace
\begin{array}{ll}
     \sqrt{\theta_1^2 + (\theta_2 \theta_3 - 1 / (\theta_2\theta_4))^2} + X_n & \mbox{if $1\leq n \leq 15000-\Delta]$ (concept A),} \\
     \arctan{((\theta_2\theta_3 - (1/(\theta_2 \theta_4)))/\theta_1)} + X_n) & \mbox{if $15000+\Delta\leq n\leq 30000]$ (concept B),} \\
\end{array}\right.
\end{equation}
where $X_n \sim N(0, \sigma)$ is a realization of a Gaussian random variable with standard deviation $\sigma$, and the entries $\{\theta_1,\theta_2,\theta_3,\theta_4\}$ of the variable vector $\boldsymbol{\theta}$ fulfill $0\leq\theta_1\leq 100$, $40\pi\leq\theta_2\leq 560\pi$, $0\leq\theta_3\leq 1$ and $1\leq\theta_4\leq 11$. These two concepts correspond to the so-called Friedman 2 (concept A) and Friedman 3 (concept B) benchmark problems proposed in \cite{Friedman:91} and often used thereafter to build synthetic datasets for regression problems. A transition between both concepts is set at time $n=15000-\Delta$, and modeled by a sigmoid function $f[n] = 1/(1+\exp(-4(n-15000)/\Delta))$ that establishes the probability that the new sample for $n>15000-\Delta$ belongs to concept B. Parameter $width$ establishes the duration of this transition, which is set to $200$ samples in the results shown in what follows. We also set $\sigma=1$.
 
A predictive model $M_{\bm{\vartheta}}$, therefore, aims to characterize the relation between $\boldsymbol{\theta}$ and $x[n,\boldsymbol{\theta}]$ by incrementally learning from the observed values over the stream. Therefore, we assume that once the model has produced a prediction $\widehat{x}[n,\boldsymbol{\theta}]$ for a given time, the true value $x[n,\boldsymbol{\theta}]$ is made available for incrementally updating its knowledge. This common practice in data stream mining is referred to as \emph{test-then-train}, and is a valid assumption for studies focused on error performance statistics disregarding any eventual data supervision latency. It should be also clear that only some particular modeling flavors $M_{\bm{\vartheta}}$ allow for incremental updates based on new supervised data. 

We first proceed by analyzing the performance shown by an incrementally learnable model that does not incorporate any countermeasures for accommodating the transition between concepts (concept drift). Specifically, we use a windowed version of the CART algorithm \cite{Loh:11} to grow a tree regression model, wherein the quality of a split in the tree is measured by the mean squared error between the output of the split and that of the stream measured over the instances of the window. It is important to note that the capability of this model to forget and adapt its knowledge to that in the stream is given by the size $Wtree$ of the window of instances from where it is learned. This size also establishes the amount of supervised information stored for learning the model. Therefore, there is a trade-off between performance and storage requirements as per the selected window size. 

Alongside this model we consider a specific flavor of this same learning model suited for drifting data streams, called Fast Incremental Model Tree with Drift Detection (FIMT-DD) \cite{Ikonomovska:11}. This model features several changes with respect to the naive tree induction algorithm: 1) the use of the Hoeffding bound to decide whether the best splitting variable can be decided on a given set of instances; 2) the adoption of the Page-Hinckley change detection test to examine whether error increases at each node of the tree are a symptomatic sign of a concept drift; and 3) a tree regrowing mechanism that operates at every node of the tree that has been positively marked by the aforementioned test. By virtue of all these algorithmic ingredients, FIMT-DD is expected to perform better than the windowed version of the CART regression tree, mainly due to explicit methods to detect and adapt possible changes over the stream. Furthermore, FIMT-DD does not require any window of past samples to be stored, thereby making it compliant with stream settings.
\begin{figure}[!t]
\centering
\includegraphics[width=0.9\textwidth]{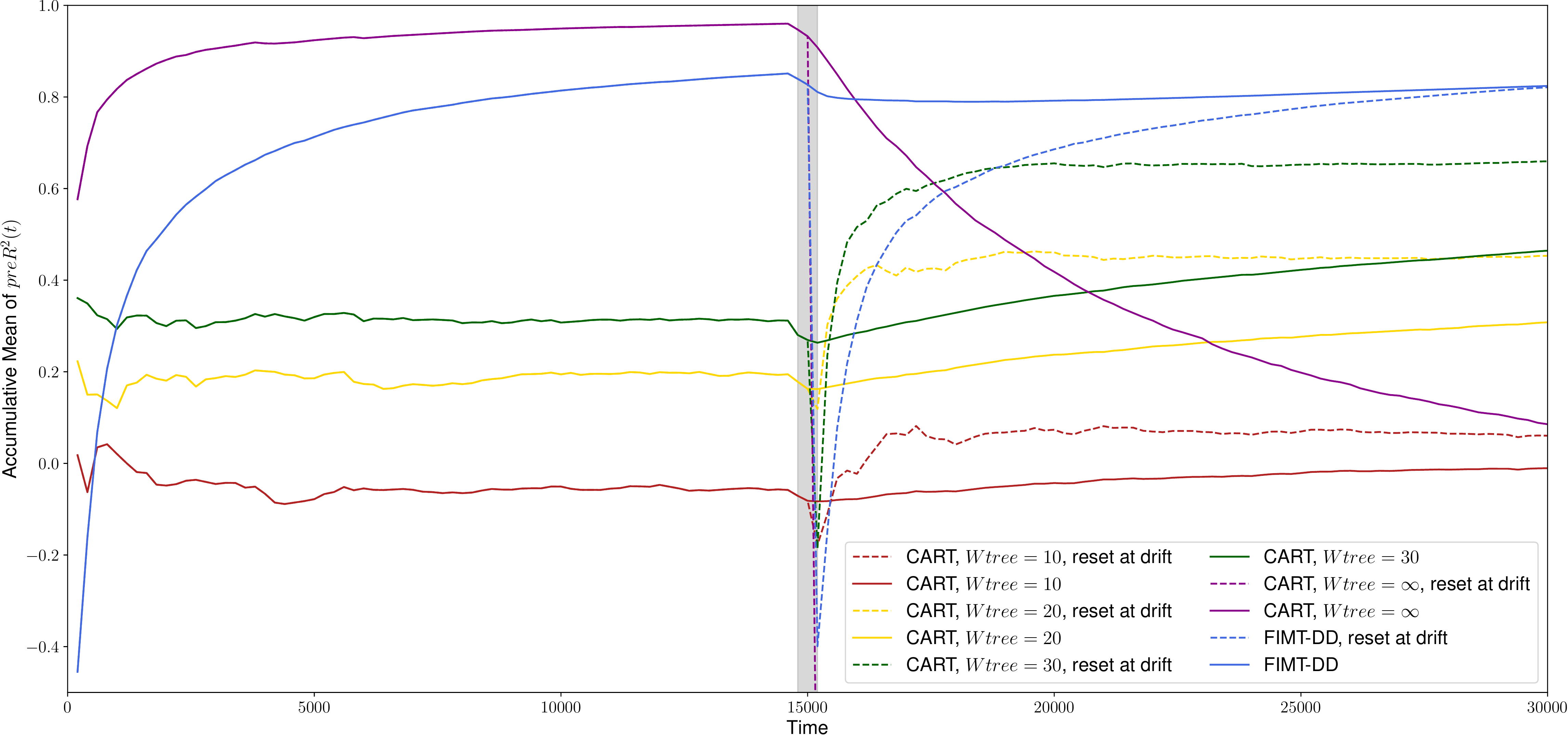}
\caption{Cumulative mean of different learning models, $preR^2$, when processing a data stream composed by two different concepts. A gradual drift occurs within $14800<t<15200$, changing from concept A to concept B as per \eqref{eq:concepts}. The lack of persistence in the characterized relationship between the observed variables $\boldsymbol{\theta}$ and the temporal variable $x(t,\boldsymbol{\theta})$ gives rise to a decay in the precision of the model, which is counteracted under different strategies. On one hand, CART models are fed only with the last $Wtree$ instances, which provides an effective forgetting mechanism at a penalty in performance when no drift occurs. This penalty becomes lower the higher the window is made, but makes the model less compliant with space constraints typically imposed in stream settings, and asymptotically makes the overall model not adapt suitably to concept changes ($Wtree=\infty$). On the other hand, FIMT-DD not only achieves better performance when the pattern between observed variables and the temporal target is stable, but also reacts smoothly against the drift without storing past samples.}
\label{drift_synthetic}
\end{figure}

Figure \ref{drift_synthetic} summarizes the results obtained for this first experimentation. Specifically, the plot depicts the accumulated mean of the prequential coefficient of determination obtained by different models along time, given by:
\begin{equation}\label{preqr2}
preR^2[n] = 1 - \frac{\sum_{i=n-200}^{n} (x[i,\boldsymbol{\theta}] - \widehat{x}[i,\boldsymbol{\theta}])^2}{\sum_{i=n-200}^{n} \left(x[i,\boldsymbol{\theta}] - \frac{1}{200}\sum_{j=n-200}^n x[j,\boldsymbol{\theta}]\right)^2},
\end{equation}
i.e. as the coefficient of determination between the prediction $\widehat{x}[n,\boldsymbol{\theta}]$ and the true value $\widehat{x}[n,\boldsymbol{\theta}]$ of the temporal variable to be predicted. The time at which the transition from concept A to concept B occurs is marked with a dashed vertical line. For any model, a second dashed curve is included in the plot, which corresponds to the accumulated mean reinitialized to $0$ at the time the drift is held. This permits to visually inspect the reaction of the model when facing the drift between concepts, i.e. how the precision of the model recovers \emph{after} the drift disregarding its behavior \emph{before} the drift. The models in this benchmark include windowed CART regression trees with different $Wtree$ values, namely, $10$, $20$, $30$ and $\infty$ (i.e. infinite window size). Despite not suitable for stream mining due to its unbounded storage requirements, this last case provides empirical evidence of the need for adapting the models to the lack of persistence in the characterized relationship between observed and target variables. The benchmark also considers the aforementioned FIMT-DD approach, which incorporates more sophisticated drift adaptation mechanisms than the forgetting capability implicitly provided by windowing over time.

This being said, we pause at Figure \ref{drift_synthetic} to discuss on several interesting results. To begin with, the evolution of the CART regression tree with $Wtree=\infty$ conforms to intuition: inside concept A, the model attains the best performance, but once drift occurs, its performance degrades catastrophically. CART-based approaches with finite window sizes behave as expected, showing a better predictive performance as the size of the window of past instances from where the model is learned is higher. This, however, produces a model that may not be entirely compliant with the storage requirements imposed by the most stringent stream learning scenarios. On the other hand, we note that the lack of persistence in the modeled pattern emerges from a decay in the performance of the model when the drift occurs. Focusing on the dashed curves, we see that the speed at which the model evolves and characterizes the new concept depends on its effectiveness to forget the previous concept and capture the new one. While in CART-based approaches this effectiveness depend on the size of the window (at the cost of a worse performance once the concept has stabilized), FIMT-DD exhibits a superior performance and a good reaction to the drift, which comes at no cost in terms of data storage.

Departing from the insights gained with this first experiment, we now focus our attention on a second set of experiments dealing with the real-world \texttt{weather} dataset, widely used in the literature related to concept drift \cite{Elwell:11,Ditzler:12}. In detail, \texttt{weather} is part of the data published by the National Oceanic and Atmospheric Administration (NOAA). The dataset consists of 18,159 daily observations, each comprising 8 weather-related measurements (e.g. temperature, dew point, sea level pressure, visibility or average wind speed) collected at the Offutt Air Force Base located in Bellevue, Nebraska. Based on these observed variables, the goal is to predict whether it will rain ($x[t,\boldsymbol{\theta}]=1$) or not ($x[t,\boldsymbol{\theta}]=0$), hence we deal with a binary classification problem. The dataset spans over more than 50 years, featuring interesting sources of drift during this period such as short-term seasonal changes and long-term climate trends. Unfortunately, there is no ground truth in regards to the presence of verified drifts, so experiments performed over this dataset aim at showing the predictive accuracy of different learning models over the stream, each endowed with different mechanisms of detection and/or adaptation to such eventual drifting periods. 

\begin{figure}[!h]
\centering
\includegraphics[draft=false, angle=0,width=0.7\textwidth]{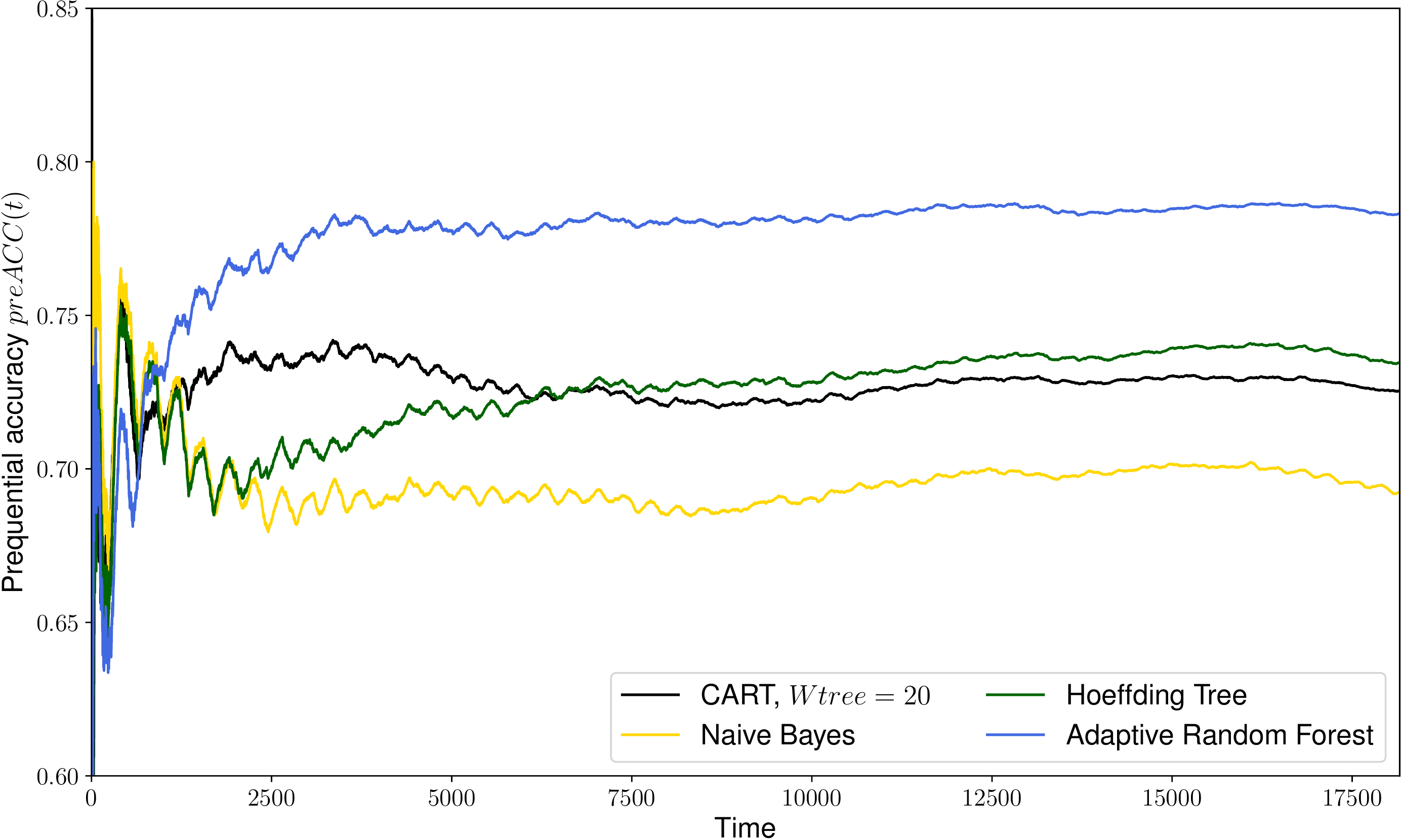}
\caption{Prequential accuracy $preACC(t)$ achieved by several modeling choices over the real-world \texttt{weather} dataset. The accuracy degradation of models without drift adaptation in the first part of the stream (Gaussian Naive Bayes and Hoeffding Tree) expose severe drifting events in this part of the stream. The windowed CART and the Adaptive Random Forest, on the other hand, exhibit better accuracy levels in this stage of the stream. As more stream samples are fed to the models, the Hoeffding Tree and Adaptive Random Forest outperform the rest of schemes in the plot, as a result of their more efficient learning algorithms, that permit to retain knowledge without resorting to restrictive forgetting methods.}
\label{drift_real}
\end{figure}

This is indeed the purpose of the plot in Figure \ref{drift_real}, where we depict the prequential accuracy $preACC[n]$ of different incremental classifiers when progressively fed with the instances of the aforementioned dataset. Similarly to the prequential coefficient of determination in Expression \eqref{preqr2}, $preACC[n]$ computes incrementally the average accuracy of the model along time, yet on a per data instance basis:
\begin{equation}\label{preqACC}
    preACC[n] = preACC[n-1]+\frac{ACC[n]+preACC[n-1]}{n-1},\quad n > 1
\end{equation}
where $preACC(0)=0$, and $ACC(t)=0$ if $\hat{x}(t,\boldsymbol{\theta})=x(t,\boldsymbol{\theta})$ (i.e. good prediction) and $0$ otherwise. We consider 4 different classification models: i) a Naive Gaussian Bayes classifier whose statistics are updated incrementally based on arriving samples; ii) a windowed CART model similar to those utilized in the first set of experiments, configured with $Wtree=20$ samples; iii) a Hoeffding Tree Classifier (also known as Very Fast Decision Trees, VFDT \cite{Hulten:01}), which induces classification trees from data streams by relying on the Hoeffding bound to grow or shrink their nodes, and which does not deal with possible drifts anyhow; and iv) an Adaptive Random Forest \cite{Gomes:17}, which is an ensemble of tree classifiers that accommodate drifts by adopting a twofold strategy: an online bootstrap aggregating sampling procedure with random feature selection to induce diversity in the ensemble, and the use drift detectors for every tree learner in the ensemble, which trigger selective tree resets when a drift is detected.

The worst performance is shown by the Gaussian Naive Bayes classifier, whose lack of adaptation mechanisms to learn from drifting data hinders further its reduced modeling capability with respect to the rest of models in the benchmark. Next, the windowed CART model attains good performance scores in the early stages of the data stream, but progressively degrades as a result of the low number of samples from where it is learned. The Hoeffding Tree, however, performs comparatively worse in the first part of the stream. This last observation, along with the results for the Gaussian Naive Bayes counterpart, suggests that the \texttt{weather} dataset contains severe concept changes, as can be inferred from the degradation shown by classifiers without adaptation like these ones. As the stream advances, the performance of the Hoeffding Tree improves to eventually surpass that of the windowed CART approach, as a result of its better capability to incrementally learn from the stream. Finally, the Adaptive Random Forest dominates the benchmark, as it is the only model in the benchmark that incorporates efficient incremental learning, a superior modeling capability and explicit mechanisms for drift adaptation. All these algorithmic ingredients yield an incremental learning model that adapts better to drifting data streams.

\section{Conclusions, perspective and outlook}\label{Conclusions}

Persistence is an elusive concept, but can be broadly defined as the average time taken for a certain variable of a complex system to change from one state to a different one. Persistence is an important characteristic of any system, which describes part of its statistical physics behaviour. The study of complex systems' persistence involves different definitions and uses different techniques, depending on whether we consider short-term or long-term persistence. Short-term persistence has been described by considering definitions focused on measuring the average time spent at a given state, or the probability of changing to another one. 
Markov chain methods or auto-regressive models have been applied, among others, to describe the characteristics of systems with short-term persistence. 
On the other hand, long-term persistence is usually defined based on the autocorrelation function of the time series describing the complex system at hand. A long-range persistent process exhibits a power-law scaling of its autocorrelation function. Long-term persistence can also be defined in terms of the power spectral density of the time series (which forms a Fourier par with its autocorrelation function). A process can be defined as long-range persistent if its power spectral density scales asymptotically as a power law for frequencies close to the origin. Some methods used to describe the characteristics of long-term persistent systems have been described in this article, such as Hurst R/S analysis, DFA or methods based on Wavelets.

Concepts related to persistence in arbitrary systems appear in a wide spectrum and diversity of the scientific literature. 
In this review we 
provided a structured literature review in this work, highlighting the most important fields where persistence of complex systems has been 
studied. Without doubt, the persistence of systems in Earth and atmospheric sciences has been profusely studied over the last fifty years. Short-term and long-term persistence of atmospheric processes such as precipitation, global temperature, droughts and heat waves, soil moisture, climate models performance, sea level, air pollution and different hydrologic-related processes have been reviewed. In connection with atmospheric and climate, the persistence of renewable energy resources has been studied in different works, which have been discussed in this paper. Wind speed and solar persistence are the two renewable resources which have attracted the most attention in the last years. Note that solar resource is in general much more persistent than wind, as discussed in the literature review carried out. We have also revised some works dealing with the persistence of systems in Earth science, such as geophysics and seismology. 
We have also dealt with the review of persistence in complex networks, where a number works dealing with network dynamics and persistence, or link persistence, among others, have been discussed. The research on persistence in Economics has been intense in the last years. We have reviewed here a number of works dealing with persistence in inflation, assets prices, stock market indices or exchange rates. We have next revised different works on persistence of non-equilibrium thermodynamics systems, most of them characterized by probabilities of system's state change which decay as power law of the form $P_o(t) \sim t^{-\theta}$. We have also revised different techniques to construct synthetic time series with long-term persistence properties, some of them specific for hydrology applications. The last large block of works discussed in this paper is related to persistence in optimization and planning. In this type of problems, persistence is referred to the optimal solution found, which should maintain a certain degree of performance versus changes in the original definition of the problem. We have closed this review of existing works about persistence in complex systems by reviewing some previous studies in biomedical applications and sport science.

Specific case studies and applications have been described in the last part of this paper: The study of time scales in DFA over different time series, how persistence-based methods can be accurate in solar radiation prediction problems, an study of soil moisture persistence and its analysis, how to exploit directional persistence of wind to reduce the computational cost in turbine layout problems for wind farms design and the effect of persistence in ML algorithms when dealing with data streams, are the problems exhaustively discussed in this work, to show the importance of studying persistence of real complex systems.

The study of complex systems' persistence is expected to continue and be a major topic of research interest in the near future. We anticipate a cross-fertilization of the field, and suggest here some prospective uses and applications of persistence analysis for complex systems:
\begin{itemize}
\item {\em Merging with ML algorithms.} ML is currently the paradigm of computational techniques, applied to the study of almost every known system or open problem. The development of ML is being huge in the last years, specially with the raise of deep learning as the last and maybe most powerful tool for knowledge discovering and analysis from data. However, in the study of different highly complex real systems, sometimes the solutions obtained by ML techniques (exclusively working on data) do not respect the most elementary laws of physics, like mass or energy conservation \cite{reichstein2019deep}. Some researchers have pointed out the necessity of mixing ML algorithms with numerical models dealing with the physical equations of the phenomenon under study, so the physical meaning of ML solution is ensured \cite{karpatne2017theory}. The inclusion of the system's persistence as one of the rules that ML algorithms must keep within their outcome could be basic in order to obtain solutions with physical meaning out of ML approaches.  

\item {\em Models of memory and persistence with a biological inspiration.} Quantifying persistence in nonlinear systems requires nonlinear models of memory and adaptation, and artificial neural networks excel into that. Neural nets have found inspiration in biological processes in the brain system or parts thereof, such as the visual, auditory or olfactory brain. Such bio-inspiration has resulted in a great success in applied fields of engineering, such as computer vision, speech or natural language processing. There is no reason why this should not continue in the future. In order to maintain a persistent activity, it has been  argued that a positive feedback is strictly needed \cite{hebb1949organization}. There are also some evidences that neuronal responses can instead be maintained by a purely feedforward mechanism that maintains short-term memory in which both feedforward and feedback processes interact to govern network behavior. Such mechanisms are currently accounted explicitly in recurrent neural networks, but one could gain in computational efficiency by accounting for {\em memory without feedbacks} \cite{goldman2009memory}. The framework of {\em attractor neural networks} \cite{barak2007persistent} could be extended to represent time-dependent stimuli in artificial neural nets. Actually, persistent activity states (attractors) in biological neurons is believed to be the basis of the working memory \cite{barak2014working}. The existence of persistent states was originally suggested in \cite{little1974existence}, and gave rise to exciting developments in multiscale neural networks \cite{Werbos90:BTW}. Yet, also relevant computational principles of memory are concerned here \cite{chaudhuri2016computational}. After all, the ability to store and use information is essential for a variety of adaptive behaviors, including learning, generalization, prediction and inference, which are familiar concepts in ML (and deep learning). 

\item {\em Extreme events.} The persistence of extreme events has not been deeply studied in the literature so far, but in some specific cases related to meteorology and climate, such as droughts or heat waves. Even in these cases, the concept is still far from being fully understood and modelled, and it seems that persistence of these events is fully dependent on each specific case. It is necessary a more intense research to characterize extreme events persistence, which is difficult, due to the lack of observational data because the low rate of occurrence of these events.

\item {\em Persistence and physical drivers.} In many works, persistence of a complex system is described from a statistical physics point of view, i.e. they try to characterize the statistics of the time series representing an event, without a clear connection with the physical processes (in a broad sense) driving the system. 
We have touched upon these 
topics in Section \ref{Pers_ML}, where the relationship of persistence in ML with the exogenous variables describing the system has been discussed. In a more general framework, to study the relationship between the persistence of a complex system, and the physical variables of the phenomenon (exogenous variables) would be extremely interesting and also challenging. This point clearly connects with topics such as causal inference, 
observed physical processes, and persistence. 

\item {\em Persistence of new complex systems.} Persistence analysis has proven to be a very useful way of characterizing part of the statistical properties of complex systems. The techniques revised in this work, and other alternatives that may arise in the next years, could be applied to any new complex system of interest, in any research area. Some specific new problems in which persistence analysis is extremely interesting arise, for example, in mass media research (why a given news item remains at newscasts during more time than others? how is this process of news items importance decaying over time? and in turn, how does the audience interest in a news item decay over time?), social research and politics, sports science, further problems related to persistence in optimization of systems, and many others.
\end{itemize}
The study of persistence in complex systems has reached a point of maturity in theoretical developments and exciting applications. We reviewed the literature, gave empirical evidence of performance in several case studies, and outlined a collective agenda for the future research and developments. We anticipate a wider adoption of the techniques by newcomers and experienced researchers willing to make sound advances in persistence characterization. 
And vice-versa, we expect to spark new ideas across fields. Our insights and prospects held over this field pave a promising path plenty of opportunities for persistence related studies, including a synergistic interplay with other domains of knowledge. 

\section*{Acknowledgments}
This research has been partially supported by the project PID2020-115454GB-C21 of the Spanish Ministry of Science and Innovation (MICINN). This research has also been partially supported by Comunidad de Madrid, PROMINT-CM project (grant ref: P2018/EMT-4366). J. Del Ser would like to thank the Basque Government for its funding support through the EMAITEK and ELKARTEK programs (3KIA project, KK-2020/00049), as well as the consolidated research group MATHMODE (ref. T1294-19). GCV work is supported by the European Research Council (ERC) under the ERC-CoG-2014 SEDAL Consolidator grant (grant agreement 647423) and the ERC Synergy Grant `Understanding and Modelling the Earth System with Machine Learning' (USMILE) under the Horizon 2020 research and innovation programme (Grant agreement No. 855187).


\end{document}